\pdfoutput=1
\documentclass{article}
\usepackage{lmodern}
\usepackage{graphicx}
\usepackage{amsmath}
\usepackage{amsthm}
\usepackage{amsfonts}
\usepackage{mathrsfs}
\usepackage{amssymb}
\usepackage{enumerate}
\usepackage{lscape}
\usepackage{longtable}
\usepackage{rotating}
\usepackage{multirow}
\usepackage{color}
\usepackage{tikz}
\usepackage{url}
\usepackage{subfigure}
\usepackage{rotating}
\usepackage{comment}
\usepackage[utf8]{inputenc}
\usepackage{acro}
\usepackage{color}
\usepackage{booktabs}
\usepackage{float}
\usepackage{placeins}
\newfloat{algorithm}{tbp}{loa}
\floatname{algorithm}{Algorithm}
\newenvironment{algorithmic}[1][0]{\begin{enumerate}}{\end{enumerate}}
\newcommand{\State}{\item}
\newcommand{\Require}{\item[\textbf{Require:}]}

\newcommand{\If}[1]{\item \textbf{if} #1 \textbf{then}\begin{list}{}{\leftmargin=1.5em\topsep=0pt\itemsep=0pt}}
\newcommand{\Else}{\end{list}\item \textbf{else}\begin{list}{}{\leftmargin=1.5em\topsep=0pt\itemsep=0pt}}
\newcommand{\EndIf}{\end{list}}
\newcommand{\While}[1]{\item \textbf{while} #1 \textbf{do}\begin{list}{}{\leftmargin=1.5em\topsep=0pt\itemsep=0pt}}
\newcommand{\EndWhile}{\end{list}}
\newcommand{\For}[1]{\item \textbf{for} #1 \textbf{do}\begin{list}{}{\leftmargin=1.5em\topsep=0pt\itemsep=0pt}}
\newcommand{\EndFor}{\end{list}}
\usepackage[hidelinks]{hyperref}
\hypersetup{pdftitle={Riemannian Real-Time Optimization for Robotics: Active Set SQP Parametric Pathfollowing},pdfauthor={Vyacheslav Kungurtsev},pdfsubject={Riemannian optimization, SQP, pathfollowing, robotics}}
\usepackage{geometry}  
\usepackage[most]{tcolorbox}
\usepackage{environ}
\usepackage{xcolor}

\newtheorem{theorem}{Theorem}[section]
\newtheorem{corollary}{Corollary}
\newtheorem{lemma}[theorem]{Lemma}
\newtheorem{proposition}[theorem]{Proposition}

\newtheorem{assumption}[theorem]{Assumption}
\newtheorem{definition}[theorem]{Definition}

\title{Riemannian Real-Time Optimization for Robotics\\
\large Active Set SQP Parametric Pathfollowing}
\author{Vyacheslav Kungurtsev\footnote{Oak Ridge National Laboratory \texttt{kungurtsevv@ornl.gov} \\ Notice: This document has been authored in part by UT-Battelle, LLC under Contract No. DE-AC05-00OR22725 with the U.S. Department of Energy. The United States Government retains and the publisher, by accepting the article for publication, acknowledges that the United States Government retains a non-exclusive, paid-up, irrevocable, world-wide license to publish or reproduce the published form of this manuscript, or allow others to do so, for United States Government purposes. The Department of Energy will provide public access to those results of federally sponsored research in accordance with the DOE Public Access Plan \url{http://energy.gov/downloads/doe-public-access-plan} }}
\date{November 2026}
\begin{document}
\maketitle
\begin{abstract}
Robot optimal control often combines Euclidean states with configurations on
Lie groups and other smooth manifolds. This paper develops a parametric
Riemannian Newton pathfollowing method based on homotopy quadratic programs and
parallel strict-working-set tasks. It specifies the event handling,
working-set exchanges, matrix factor updates, synchronization, and residual
checks used by the algorithm. Tracking guarantees are given under mildly degenerate regularity conditions and the Riemannian strong second-order
sufficient condition. Drone navigation and arm manipulation are formulated
as geometric optimal-control problems, including finite-force grasp
acquisition, transport and supported release. Numerical examples illustrate
drone navigation and arm pick-and-place in closed loop with MuJoCo
simulations, through KKT and feasibility residuals at every feedback sample
and along the QP homotopies. They include a perturbed drone plant and the
full-horizon arm problems, for which the controller does not yet run in
real time.
\end{abstract}
\section{Introduction}
\label{sec:introduction}
Robot kinematics naturally involve geometric configuration spaces, including rotation Lie groups and products of revolute-joint circles~\cite{selig2005geometric}. Riemannian optimization provides a framework for enforcing these geometric constraints while optimizing motion~\cite{sato2021riemannian}. Nonlinear model predictive control (NMPC) adds a second requirement: the optimal control problem changes with every measurement, and an improved solution must be available within one sampling interval. This paper develops a Riemannian predictor--corrector method for that setting. It coordinates strict working sets, inertia-aware linear algebra and active-set quadratic programming (QP) within one feedback update of NMPC. A strict working set is a set of active inequality constraints whose gradients, together with those of the equality constraints, are linearly independent, so that it determines one multiplier vector. Inertia-aware linear algebra means that each factorization of a Karush--Kuhn--Tucker (KKT) matrix also returns the numbers of its positive, negative and zero eigenvalues, from which the second-order conditions are checked.

The discretized optimal control problem (OCP) uses four principal variable groups and one explicit algebraic feature block. For a positive integer $M$, write $[M]=\{1,\ldots,M\}$. The integers $M^u$, $M^v$ and $M^q$ count the corresponding variable blocks; $d_u,d_v,d_x,d_q,d_z$ denote their aggregate ambient dimensions:
\begin{enumerate}
    \item Geometric control variables:
    \[
    u_k \in \mathbb{U}:=\prod_{m\in [M^u]}\mathbb{U}_{m}\subseteq \prod_{m\in [M^u]}\mathcal{M}^u_{m} \subset \mathbb{R}^{d_u}
    \]
    These controls lie on a product manifold whose factors represent the kinematic geometry. The set $\mathbb{U}_{m}$ specifies physical actuator limits on coordinate block $m$, such as a prescribed spherical cap or an intersection of hemispheres.
    \item Euclidean control variables:
    \[
    v_k\in \mathbb{V}:=\prod_{m\in [M^v]}\mathbb{V}_{m}\subset \prod_{m\in [M^v]}\mathbb{R}^{d^v_{m}} = \mathbb{R}^{d_v}
    \]
    These include torque, thrust, and other actuator amplitudes that enter the manifold dynamics.

    \item Euclidean state variables:
    \[
    x_k \in \mathbb{R}^{d_x}
    \]
    These include position and velocity in space.

    \item Geometric state configuration variables:
    \[
    q_k \in \mathbb{Q}:= \prod_{m\in [M^q]} \mathbb{Q}_{m}\subseteq \prod_{m\in [M^q]}\mathcal{M}^q_{m} \subset \prod_{m\in [M^q]}\mathbb{R}^{d^q_{m}}=\mathbb{R}^{d_q}
    \]
    These variables describe the robotic configuration.

     \item Algebraic feature variables $z_k\in \mathbb{R}^{d_z}$ describe the Euclidean forward kinematics of the geometric states,
     \[
     z_k = h^z(x_k,q_k)
     \]
     These include, for example, the end-effector position of a robotic arm composed of revolute joints.
\end{enumerate}
The feature variables are redundant, since they are functions of the states. They are kept as explicit variables because collision and tracking constraints are simple functions of $z_k$, whereas they are complicated functions of $q_k$.

    The Euclidean state variables evolve according to
    \[
    x_{k+1}=f^x(x_k,q_k,u_k,v_k)
    \]
    which may rotate a body-frame input into world coordinates before integration.
As a simple illustration, consider an ideal velocity-controlled model of drone dynamics with geometric control
$u(t)\in SO(3)$ and prescribed speed $\nu\geq0$,
\[
 \dot x(t)=\nu u(t)e_3,\qquad e_3=(0,0,1)^\top,
 \qquad
 x(t+h)=x(t)+\int_t^{t+h}\nu u(s)e_3\,ds.
\]
Here $x\in\mathbb R^3$ is position and $h>0$ is the integration interval, so $\|\dot x\|=\nu$. (The drone model of Section~\ref{sec:drone-example} is different: there the attitude is a state and the control is a direction on $S^2$.)

    The geometric state dynamics are described by
    \[
    q_{k+1}=f^q(q_k,x_k,u_k,v_k)
    \]

A retraction can be used in the numerical integrator $f^q$. A retraction maps
a tangent increment to a manifold point and has differential equal to the
identity at the zero increment. For a revolute
joint let $\theta$ be an unwrapped local angle and $q=[\theta]\in S^1$.
An ideal angular-velocity actuator has
\[
 \dot\theta(t)=v(t),\qquad
 q_{k+1}=R_{q_k}(h v_k)=[\theta_k+h v_k],\qquad
 \theta_\ell\leq\theta_k\leq\theta_u
\]
on a declared chart with $\theta_u-\theta_\ell<2\pi$. The maps $f^x$, $f^q$ and $h^z$ may depend on the stage; they are then written $f^x_k$, $f^q_k$ and $h^z_k$.

Let $\hat t$ denote a feedback-sample index. The horizon length
$H_{\hat t}$ counts \emph{control intervals}. There are $H_{\hat t}+1$
state/feature nodes, including the measured initial node and terminal node. Set
\[
 \mathcal N_{\hat t}=\{\hat t,\ldots,\hat t+H_{\hat t}\},\qquad
 \mathcal I_{\hat t}=\{\hat t,\ldots,\hat t+H_{\hat t}-1\}.
\]
There are exactly $H_{\hat t}+1$ state/feature nodes and $H_{\hat t}$ controls
and shooting intervals. For a task ending at absolute temporal knot $N$,
$H_{\hat t}=N-\hat t$ (a \emph{shrinking horizon}), while a fixed-horizon controller instead keeps $H$ fixed.
In relative indices the states are $0,\ldots,H$ and the controls are
$0,\ldots,H-1$; Sections~\ref{sec:geometric-case-studies}--\ref{sec:pathfollowing} use relative indices, so that node $0$ is always the measured node. An OCP with controls requires $H\geq1$; $H=0$ is terminal-only.

The horizon length counts intervals, not time. The physical length of an
interval is specified separately by the integrator, and it may itself be an
optimization variable, as in the drone problem of
Section~\ref{sec:drone-example}.

The resulting full-space nonlinear program (NLP) is
\begin{equation}\label{eq:riemnlp}
\begin{aligned}
 \min_{\substack{(x_k,q_k,z_k)_{k\in\mathcal N_{\hat t}}\\
                  (u_k,v_k)_{k\in\mathcal I_{\hat t}}}}
 &\sum_{k\in\mathcal I_{\hat t}}
 J(u_k,x_k,v_k,q_k,z_k,k)
 +V(x_{\hat t+H_{\hat t}},q_{\hat t+H_{\hat t}},z_{\hat t+H_{\hat t}}),\\
 \text{subject to}\quad
 &x_{k+1}=f^x_k(x_k,q_k,u_k,v_k),
              &&k\in\mathcal I_{\hat t},\\
 &q_{k+1}=f^q_k(q_k,x_k,u_k,v_k),
              &&k\in\mathcal I_{\hat t},\\
 &z_k-h^z_k(x_k,q_k)=0,
              &&k\in\mathcal N_{\hat t},\\
 &c_k(x_k,z_k)\geq0,
              &&k\in\mathcal N_{\hat t}\setminus\{\hat t\},\\
 &(u_k,v_k)\in\mathbb U\times\mathbb V,
              &&k\in\mathcal I_{\hat t},\\
 &q_k\in\mathbb Q,
              &&k\in\mathcal N_{\hat t}\setminus\{\hat t\},\\
 &(x_{\hat t},q_{\hat t})=(\hat x(\hat t),\hat q(\hat t)).
\end{aligned}
\end{equation}
This is a multiple-shooting transcription: the states at all nodes are variables, and the dynamics are imposed as equality constraints between consecutive nodes.
The inequalities $c_k(x_k,z_k)\geq0$ and the configuration bounds impose
the specified plant requirements at future prediction nodes. The measured
initial state is a problem parameter: the last line of \eqref{eq:riemnlp}
is not a constraint passed to the solver. The measured values are
substituted for $(x_{\hat t},q_{\hat t})$ before optimization, so the
initial node carries no variables. Operational state, contact and phase
constraints start at the next predicted node, because the optimizer cannot
change the measured state, which may itself violate them. The first shooting
equation therefore computes a correction from the actual measured state,
including any mismatch between the model and the measured state. Control
bounds apply to the first command as well as to subsequent commands. Whether
the measured state and the command about to be applied are physically
acceptable is tested by the supervising controller; these admission checks
are separate from the constraints imposed on the prediction.

In NMPC, the first computed control pair is applied, the next state is measured, and
$\hat t$ is incremented. The measured state $(\hat x(\hat t),\hat q(\hat t))$ is thus the parameter of a family of problems \eqref{eq:riemnlp}, and the controller follows the solution of this family as the parameter changes.

Related work includes online Riemannian optimization for convex objectives~\cite{wang2023online} and sequential quadratic programming (SQP) on Riemannian manifolds~\cite{obara2022sequential}. The former treats problems with no additional operational constraints. SQP is central to real-time MPC. Here these components are coordinated for optimal-control models that couple geometric and Euclidean variables and whose operational constraints change their active set from one sample to the next.

Many robotic factors have small intrinsic dimension: $S^2$ has dimension two,
$SO(3)$ three and $SE(3)$ six. Explicit matrix representations and inexpensive
component retractions make geometric feasibility practical. The problems at consecutive feedback samples differ only slightly and often have the same or nearly the same active set. This motivates the classical active-set approach to real-time SQP-based parametric optimization, in which the active set of the previous solution is the starting guess for the next.

Constraint degeneracy is a practical concern in navigation and contact tasks, where several active constraints may have dependent gradients. Earlier predictor--corrector work~\cite{kungurtsev2017predictor} established tracking guarantees under the Mangasarian--Fromovitz constraint qualification (MFCQ), allowing dependencies excluded by the linear independence constraint qualification (LICQ); see also the NMPC stability analysis in~\cite{suwartadi2020fast}. That procedure estimates from residuals which constraints are strongly active (active with a nonzero multiplier) and which are weakly active (active with a zero multiplier), and it solves an auxiliary linear program (LP). Identification from residuals compares each constraint value with a threshold derived from the KKT residual. Such a threshold, like the regularization parameter, is valid only up to an error-bound constant that depends on the problem and is rarely known. The auxiliary LP is a second kind of subproblem, which adds a computational layer to an embedded implementation otherwise built on QPs.

The method developed here handles mild degeneracy through strict working representations and coordinated QP paths. It combines linear algebra, active-set QP and pathfollowing SQP, together with a globalized method for initialization and recovery. Under MFCQ the multiplier set at a solution can contain more than one element; each strict working set selects one of them, and the parallel tasks keep several such selections available until the nonlinear acceptance test.

\section{Background}
\label{sec:background}
\subsection{Real-time correction and active-set linear algebra}
In NMPC the problem \eqref{eq:riemnlp} must be solved again at every
sample. Real-time iteration (RTI) does not iterate to convergence at each
sample; it normally budgets one SQP correction per measurement
\cite{diehl2001real,diehl2005real}, warm started from the shifted previous
solution, so that the iterates track the moving solution. The advanced-step
approach~\cite{zavala2009advanced} instead prepares an NLP solution using a
predicted next state and corrects it when the measurement arrives.

In this paper, at sampling index $\hat t$ the data, horizon dimension and
constraint catalogue are fixed during one numerical correction. The
previous trajectory is shifted or restricted before that correction
starts. Several corrector rounds at the same sample are possible. Each
round has finitely many operational working-set tasks, which may execute
asynchronously. A candidate is selected after every admitted task has
completed or has been explicitly terminated.

\paragraph{Inertia-revealing factorizations.}
Optimization linear algebra commonly uses dense direct factorizations for small matrices, sparse direct factorizations for structured problems, and iterative methods when direct factorization is too costly. Optimal control problems have substantial sparsity; the exposition here concentrates on dense factorizations to make the working-set operations explicit.

A classical dense technique for symmetric indefinite Karush--Kuhn--Tucker (KKT) systems is the
$LBL^\top$ factorization~\cite{gill2021numerical}. 
It provides both a representation for solving the system and its inertia: the numbers of positive, negative and zero eigenvalues. 

Let $H\in\mathbb{R}^{n\times n}$ represent a Hessian and let $A\in \mathbb{R}^{m\times n}$ be the active-constraint Jacobian; the same letters are used for the working KKT systems of Section~\ref{sec:pathfollowing}. A symmetric indefinite KKT matrix admits the triangular and block-diagonal factorization
\begin{equation}\label{eq:lblt}
 K=\begin{pmatrix}H&A^\top\\A&0\end{pmatrix},\qquad
 \Pi^\top K\Pi=L\operatorname{diag}(B_1,\ldots,B_{m_B})L^\top,
 \qquad L_{ii}=1,\quad L_{ij}=0\ (i<j),
\end{equation}
where $\Pi$ is a permutation matrix, $L$ is unit lower triangular, and each $B_j$ is a real symmetric $1\times 1$ or $2\times 2$ block. The block-diagonal factor is written $B=\operatorname{diag}(B_1,\ldots,B_{m_B})$, which gives the factorization its name. The eigenvalue signs of each block determine its contribution to the inertia.

Write $\operatorname{In}(K)=(n_+,n_-,n_0)$ for the positive, negative and zero eigenvalue counts of a symmetric matrix $K$. When $A$ has full row rank and the reduced Hessian $Z^\top H Z$ is positive definite, the KKT inertia is $(n,m,0)$: $n$ positive and $m$ negative eigenvalues, with no zero eigenvalues. Here the columns of $Z$ form a basis of $\ker A$. These conditions ensure a well-defined reduced Newton solve, and comparing the computed inertia with $(n,m,0)$ is how the algorithm tests them. More generally, if $\operatorname{rank}A=\rho$, then
\[
 \operatorname{In}\begin{pmatrix}H&A^\top\\A&0\end{pmatrix}
 =(\rho,\rho,m-\rho)+\operatorname{In}(Z^\top HZ),
\]
where the columns of $Z\in\mathbb R^{n\times(n-\rho)}$ are an orthonormal
basis of $\ker A$.  The addition is componentwise.  In particular, dependent
rows produce $m-\rho$ dual zero directions even when the reduced Hessian is
positive definite: the KKT matrix is singular whenever the rows of $A$ are
linearly dependent, whatever the curvature. This is the reason each task of
the method keeps a working set whose rows are independent.

The practical computation of an $LBL^{\top}$ factorization requires symmetric
pivoting. Unlike Gaussian elimination for nonsymmetric matrices, arbitrary row
interchanges would destroy symmetry. Instead, a symmetric permutation matrix
$\Pi$ is chosen so that the factorization is performed on
\[
\Pi^\top K\Pi = LBL^\top,
\]
where $\Pi$ records a sequence of symmetric row--column interchanges.
Classical algorithms, such as the Bunch--Parlett and Bunch--Kaufman methods,
select pivots to balance numerical stability with preservation of the symmetric
structure. Depending on the size of the available pivot, the next block of $B$
is chosen as either a $1\times1$ pivot or a $2\times2$ pivot. The latter is
required when no sufficiently stable diagonal pivot is available and preserves
backward stability while avoiding excessive element growth.

\paragraph{$LBL^\top$ factorization.}


At the $k$th step, let $S$ denote the current symmetric Schur complement, with entries $s_{ij}$.  A
symmetric pivoting strategy chooses either a scalar pivot or a $2\times2$ pivot
from $S$.  Write
\[
\lambda := \max_{i>1}|s_{i1}|,\qquad
r\in\arg\max_{i>1}|s_{i1}|,
\]
and, if $\lambda>0$,
\[
\sigma_r := \max_{j\neq r}|s_{jr}|.
\]
For the Bunch--Kaufman choice, with
\[
\alpha := \frac{1+\sqrt{17}}{8},
\]
the pivot block $B_j$ is selected by the first of the following tests that holds, in the order listed:
\[
\begin{array}{ll}
|s_{11}|\geq \alpha\lambda
&\Longrightarrow\quad
B_j = [s_{11}],\\[0.6em]
|s_{11}|\sigma_r\geq \alpha\lambda^2
&\Longrightarrow\quad
B_j = [s_{11}],\\[0.6em]
|s_{rr}|\geq \alpha\sigma_r
&\Longrightarrow\quad
\text{interchange }1\leftrightarrow r,\qquad
B_j=[s_{rr}],\\[0.6em]
\text{otherwise}
&\Longrightarrow\quad
\text{interchange }2\leftrightarrow r,\qquad
B_j =
\begin{pmatrix}
s_{11} & s_{r1}\\
s_{r1} & s_{rr}
\end{pmatrix}.
\end{array}
\]

Once $B_j$ is chosen, the elimination step is a block congruence.  For a
$1\times1$ pivot,
\[
S =
\begin{pmatrix}
\beta & a^\top\\
a & C
\end{pmatrix},
\qquad
\beta=s_{11},
\]
we write
\[
S =
\begin{pmatrix}
1 & 0\\
\ell & I
\end{pmatrix}
\begin{pmatrix}
\beta & 0\\
0 & C-aa^\top/\beta
\end{pmatrix}
\begin{pmatrix}
1 & \ell^\top\\
0 & I
\end{pmatrix},
\qquad
\ell = a/\beta .
\]
Thus the next diagonal block is
\[
B_j=[\beta],
\]
and the next Schur complement is
\[
S^+ = C-\frac{aa^\top}{\beta}.
\]

For a $2\times2$ pivot, with $G$ denoting the block below the pivot,
\[
S =
\begin{pmatrix}
E & G^\top\\
G & C
\end{pmatrix},
\qquad
E =
\begin{pmatrix}
s_{11} & s_{12}\\
s_{12} & s_{22}
\end{pmatrix},
\]
the corresponding elimination step is
\[
S =
\begin{pmatrix}
I_2 & 0\\
M & I
\end{pmatrix}
\begin{pmatrix}
E & 0\\
0 & C-GE^{-1}G^\top
\end{pmatrix}
\begin{pmatrix}
I_2 & M^\top\\
0 & I
\end{pmatrix},
\qquad
M = GE^{-1}.
\]
Hence
\[
B_j=E,\qquad
S^+=C-GE^{-1}G^\top.
\]

Each elimination step is a congruence, so by Sylvester's law of inertia
the inertia of the factored matrix is accumulated blockwise:
\[
\operatorname{In}(K)
=
\sum_j \operatorname{In}(B_j).
\]
For a scalar block,
\[
\operatorname{In}([\beta]) =
\begin{cases}
(1,0,0), & \beta>0,\\
(0,1,0), & \beta<0,\\
(0,0,1), & \beta=0,
\end{cases}
\]
while for
\[
E=
\begin{pmatrix}
a & b\\
b & d
\end{pmatrix},
\]
one may compute
\[
\det(E)=ad-b^2,\qquad \operatorname{tr}(E)=a+d.
\]
In particular,
\[
\det(E)<0 \quad\Longrightarrow\quad \operatorname{In}(E)=(1,1,0),
\]
while for $\det(E)>0$ both eigenvalues have the sign of the trace, so that
$\operatorname{In}(E)=(2,0,0)$ if $\operatorname{tr}(E)>0$ and
$\operatorname{In}(E)=(0,2,0)$ if $\operatorname{tr}(E)<0$. Finally,
\[
\det(E)=0
\]
flags a singular pivot block.

\paragraph{Warm-started active-set quadratic programming.}
Active-set QP methods traverse a sequence of independent working sets selected from the active constraints~\cite{wong2011active}. They are well suited to the related QPs that arise in feedback control, where nearby problems often share an optimal active set. qpOASES was developed specifically for warm and hot starts in optimal control~\cite{ferreau2014qpoases}. Interior-point methods~\cite{gill2024projected} and shifted primal--dual QP methods such as ODYN~\cite{rojas2026odyn} provide related approaches to warm-started optimization, while active-set solvers remain widely used~\cite{gill2014user}. Robotics provides a useful setting for studying their behavior under changing geometric and operational constraints.

The strict-working-set framework for degenerate SQP~\cite{wright2002modifying} explains fast local convergence under mild degeneracy. It was motivated by the observed effectiveness of active-set QP methods within SQP even when LICQ fails. A working set retains an independent subset of the active constraints. On such a subset the KKT matrix is nonsingular, so the multiplier associated with the working set is unique even when the set of all multipliers is not. Its careful maintenance gives well-defined linear systems and supports fast Newton-type local convergence.

\subsection{Riemannian geometry and robotics}

\paragraph{Riemannian optimization.}

Riemannian optimization concerns numerical minimization over a Riemannian manifold. We recall the definitions needed for the SQP construction, following \cite{absil2008optimization} and \cite{boumal2023introduction}. Two facts are used throughout. The gradient and the Hessian of a smooth function at a point are a vector and a symmetric linear map on the tangent space at that point, so a second-order model is an ordinary quadratic function of a tangent vector. A retraction maps a tangent vector to a nearby point of the manifold. An SQP method can therefore pose its QP in tangent coordinates and take its step with the retraction.

Holonomic constraints restrict the configuration itself, for example by a
fixed link length, and many holonomic robot configuration spaces can be
described locally by regular level sets. Nonholonomic constraints, such as
rolling without slipping, instead restrict velocities and generally do not
integrate to configuration equalities; they must be retained in the
kinematic/dynamic equations \cite[Section~2.4]{lynch2017modern}. Intrinsic formulations offer two practical benefits. Retractions preserve the geometric feasibility needed for meaningful actuator commands, and tangent models expose the geometry used by Newton-type methods. Their convergence constants and attraction neighborhoods depend on the chosen metric, representation and derivative bounds~\cite{da2025thorough}. 

Throughout, $\mathcal M$ is a smooth manifold~\cite{lee2003smooth}. Its tangent space at $x$ is $T_x\mathcal M$, and its tangent bundle is $\bigcup_{x\in\mathcal M}T_x\mathcal M$. 

Riemannian algorithms use the gradient and Hessian defined by the manifold metric and connection. In Euclidean space, the components of $\nabla f$ are directional derivatives along the standard basis vectors $e_i$. A manifold has no such global basis, so the gradient is defined through the metric: it is the tangent vector whose inner product with any direction equals the directional derivative in that direction. 

For a scalar-valued function $f$ on a Riemannian manifold and a point $x\in\mathcal{M}$, the gradient is defined as follows.

\begin{definition}\label{def:riemgrad}
    The \emph{Riemannian gradient} $\operatorname{grad} f(x)$ is the unique element of $T_x\mathcal{M}$ satisfying
    \[
    \langle \operatorname{grad} f(x),\xi\rangle:= g(\operatorname{grad} f(x),\xi) := \sum g_{ij}\left( \operatorname{grad} f(x)\right)^i \xi^j = D f(x)[\xi],\,\,\, \forall \xi\in T_x\mathcal{M}
    \]
    where $g$ is the metric associated with $\mathcal{M}$.
\end{definition}

The differential $Df(x)$ is a covector in $T_x^*\mathcal M$, whereas
$\operatorname{grad}f(x)=(Df(x))^\sharp$ is a tangent vector. Here $\sharp$ is the \emph{musical isomorphism} mapping a covector to its corresponding tangent vector, also called index raising.

The negative Riemannian gradient defines the steepest descent direction. For $\operatorname{grad} f(x)\neq 0$,
\[
-\frac{ \operatorname{grad} f(x)}{\| \operatorname{grad} f(x)\|} = \arg\min\limits_{\xi\in T_x\mathcal{M},\|\xi\|=1} D f(x)[\xi]
\]

The Riemannian Hessian is defined using an affine connection.
\begin{definition}
    Let $\mathcal{V}(\mathcal{M})$ be the set of smooth vector fields on $\mathcal{M}$. An \emph{affine} connection $\nabla$ on a manifold $\mathcal{M}$ is a mapping:
    \[
    \nabla : \mathcal{V}(\mathcal{M})\times \mathcal{V}(\mathcal{M}) \to \mathcal{V}(\mathcal{M})
    \]
    denoted $\nabla_v w$ for $v,w\in \mathcal{V}(\mathcal{M})$, satisfying
\[
\nabla_{f w+gy} v = f\nabla_{w} v+g\nabla_y
v,\,\,\,\nabla_w(av+by) = a\nabla_w v+b\nabla_w y,\,\,\,\nabla_w(f v) = (wf)v+f\nabla_w v
\]
for scalar fields $f,g$, scalars $a,b$ and vector fields $v,w,y\in \mathcal{V}(\mathcal{M})$.
\end{definition}
A Riemannian metric determines the following connection.
\begin{definition}
    The \emph{Levi-Civita connection} for a Riemannian manifold $\mathcal{M}$ is the unique affine connection that is symmetric and compatible with the Riemannian metric,
    \[
    (\nabla_w v-\nabla_v w) f= [w,v]f := w(vf)-v(wf),\,\,\,\,y\langle w,v\rangle=\langle \nabla_y w,v\rangle+\langle w,\nabla_y v\rangle 
    \]
    This connection satisfies the Koszul formula
    \[
    2\langle \nabla_y w,v\rangle = y\langle w,v\rangle+w\langle v,y\rangle - v\langle y,w\rangle-\langle y,[w,v]\rangle+\langle w,[v,y]\rangle+\langle v,[y,w]\rangle
    \]
\end{definition}
The Riemannian Hessian is then defined by covariant differentiation of the gradient.
\begin{definition}
    The \emph{Riemannian Hessian} of $f$ at a point $x\in\mathcal{M}$ is the linear mapping 
    \[
    \operatorname{Hess} f(x):T_x\mathcal{M}\to T_x\mathcal{M}
    \]
    defined by
    \[
    \operatorname{Hess} f(x)[\xi]=\nabla_\xi \operatorname{grad} f
    \]
    for all $\xi\in T_x\mathcal{M}$, where $\nabla$ is the Riemannian (Levi-Civita) connection on $\mathcal{M}$. 
\end{definition}
The Koszul formula gives the identity
    \[
    \langle\operatorname{Hess} f[\xi],\eta\rangle = \xi(\eta f)-(\nabla_\xi \eta) f
    \]
    for tangent vectors extended locally to smooth vector fields $\xi,\eta$. The value at the point is independent of the chosen extensions.

Suppose the manifold is represented by regular ambient equalities, $\mathcal M=\{x:h(x)=0\}$ with $Dh$ of full row rank. At a stationary point the Riemannian Hessian then agrees, on tangent directions, with the Hessian of the ambient Lagrangian $f+\nu^\top h$. The term with the multiplier $\nu$ of the defining equations is the contribution of the curvature of the manifold in its normal directions. More generally, the exponential pullback satisfies~\cite[Proposition 5.5.4]{absil2008optimization}
\[
\operatorname{Hess} f(x) = \operatorname{Hess} (f\circ \operatorname{Exp}_x)(0_x)
\]
where the right-hand side is the Euclidean coordinate Hessian on $T_x\mathcal{M}$, endowed with the inner product induced by the Riemannian metric on $\mathcal{M}$ at $x$. 

The SQP step in the Riemannian real-time iteration is computed in tangent coordinates. If a component $r_m$ of the decision trajectory lies on a manifold $\mathcal M_m$, then the QP variable is the tangent increment
\[
    \Delta r_m \in T_{\bar r_m}\mathcal M_m.
\]
The local QP is therefore posed over the product tangent space of the current trajectory. After the QP is solved, the update is performed by a retraction,
\[
    \bar r_m^+
    =
    R_{\bar r_m}(\Delta r_m^\star).
\]
For Euclidean components this reduces to ordinary addition, while for geometric components such as $SO(3)$, $SE(3)$, $S^2$, or $T^n$, the retraction enforces exact feasibility with respect to the underlying configuration manifold.

Retractions are particularly useful in robotics: they preserve valid rotations, unit quaternions and circular joint configurations throughout the numerical updates. RTI on manifolds retains the real-time structure of classical NMPC while maintaining this geometric consistency.

\paragraph{Robotic kinematics.}

Lie groups describe the rotational and rigid-body motions arising in robotic kinematics. A revolute joint, for example, permits rotation about a fixed axis. 

Formally, a Lie group $(\mathcal{M},\circ)$ is a smooth manifold $\mathcal{M}$ with a group law $\circ$ compatible with the manifold structure. Thus for $x,y\in\mathcal{M}$, $(x,y)\mapsto x\circ y$ and $x\mapsto x^{-1}$ are smooth maps. Let $\iota\in\mathcal{M}$ be the identity element. The tangent space $T_{\iota}\mathcal{M}$ is the Lie algebra $\mathbb{M}$. For each $x\in\mathcal{M}$, left translation, right translation and conjugation, together with their derivatives, are given by
\[
\begin{array}{l}
\mathsf L_x(y) := x\circ y,\,\,\mathsf R_x(y):=y\circ x,\,\, \mathsf C_x(y):=x\circ y \circ x^{-1},\\
J\mathsf L_x(y):T_y\mathcal{M}\to T_{x\circ y}\mathcal{M},\,J\mathsf R_x(y):T_y\mathcal{M}\to T_{y\circ x}\mathcal{M},\,\,\,J\mathsf C_x(y):T_y\mathcal{M}\to T_{x\circ y \circ x^{-1}}\mathcal{M}
\end{array}
\]
The sans-serif letters distinguish these translations from the retraction $R_x$ and the Lagrangian $L$ used in the rest of the paper. Left translation is what identifies every tangent space of a matrix Lie group with its Lie algebra, $T_Y\mathcal M=Y\mathbb M$, as in the formulas of Section~\ref{subsec:background-configuration-spaces}.

Variational symplectic integrators preserve geometric structure in the discretized dynamics; see, for example, \cite{chyba2009role}. Automatic differentiation on manifolds is also available in numerical libraries~\cite{lidec2024end}, building on the foundations of Lie-group integration~\cite{iserles2000lie,celledoni2000approximating}. Applications to path planning are developed in~\cite{de2024geometry,teng2025optimization,kobilarov2011discrete}.

Geometric control constraints require particular care in numerical MPC. Model mismatch and finite optimization tolerances produce residuals in predicted dynamics, whereas a physical command still needs a valid geometric representation. The manifold update preserves that representation during the correction. 

The convergence theory of Riemannian optimization provides manifold analogues of many standard numerical methods. Their construction uses a \emph{retraction}, which maps a tangent increment to a manifold point with the correct first derivative, and often parallel transport between tangent spaces. The definitions and computational cost of these operations depend on the manifold. This cost is part of the real-time iteration budget. 

The robotics applications considered here use a small collection of standard kinematic configuration manifolds.

\subsection{Kinematic configuration spaces}
\label{subsec:background-configuration-spaces}
We now describe standard configuration manifolds in robotic kinematics, the motions they represent, and computational formulas for retractions and their domains. The exposition follows~\cite{lynch2017modern,selig2005geometric,absil2008optimization}. For each manifold we give the three items an implementation needs: the tangent space, a retraction, and the metric with its connection, from which the Hessian is computed. The case studies of Section~\ref{sec:geometric-case-studies} use the torus $T^n$, the sphere $S^2$, the rotation group $SO(3)$ and the rigid-motion group $SE(3)$. The symplectic and unitary groups are included for completeness; they arise in Hamiltonian and quantum control rather than in rigid-body kinematics.

\textbf{Remark.} Rigid chains can also have nonholonomic constraints. Their additional velocity restrictions require a separate treatment beyond the configuration geometries considered here.
\subsubsection{Spheres and Tori \texorpdfstring{$T^n$}{Tn}}
The torus is a product of circles,
\[
T^n=S^1\times\dots \times S^1
\]
and describes revolute joints and multi-joint manipulators. 

These spaces can be represented simply by angles $\theta_i\in \mathbb{R}/ (2\pi \mathbb{Z})$. The tangent space is then naturally defined by $T_{\theta}T^n=\mathbb{R}^n$.

The exponential retraction is componentwise addition modulo $2\pi$,
\[
R_{\theta}(\psi) = (\theta+\psi)\operatorname{mod}(2\pi)
\]
A joint limit is a simple bound on the unwrapped angle: a connected interval $\theta_{i\ell}\leq\theta_i\leq\theta_{iu}$ of width less than $2\pi$. An endpoint is active when it is attained. When an active endpoint is included in the operational working set, the corresponding tangent coordinate is fixed at zero, exactly as for a bound in Euclidean space. 

With the flat metric $g_\theta(\xi,\eta)=\xi^\top\eta$ on tangent vectors, the Levi-Civita connection satisfies $\Gamma^k_{ij}=0$ and the Riemannian Hessian is just the ordinary Hessian in angle coordinates.

For $n>1$, the sphere geometry is given by
\[
\begin{aligned}
 S^n&=\{x\in\mathbb R^{n+1}:\|x\|=1\},&
 T_xS^n&=\{v\in\mathbb R^{n+1}:x^\top v=0\},\\
 \nabla^{S^n}_{\eta}\xi&=P_xD\xi(x)[\eta],&&\eta\in T_xS^n,
\end{aligned}
\]
where $\xi$ is a smooth local tangent vector field and $P_x=I-xx^\top$. 

A retraction is given by normalization,
\[
R_x(v) = \frac{x+v}{\|x+v\|_2}
\]

For the sphere $S^2$, spherical coordinates $\theta,\phi$ of a point $x=(x_1,x_2,x_3)$ can be chosen as
\[
x_1=\cos\theta,\quad x_2=\sin\theta\cos\phi,\quad x_3=\sin\theta\sin\phi,\quad 0<\theta<\pi.
\]
In these coordinates the metric $ds^2=d\theta^2+\sin^2\theta d\phi^2$ has Christoffel symbols:
\[
\Gamma^{\theta}_{\phi\phi}=-\sin\theta\cos\theta,\,\,\Gamma^{\phi}_{\theta\phi} = \Gamma^{\phi}_{\phi\theta}=\cot\theta,
\]
and zero otherwise. 

A simple constraint on $S^n$ has the form $a^\top x\geq b$ and describes a
spherical cap. Its boundary $a^\top x=b$ is nondegenerate when $|b|<\|a\|$; on
$S^2$ it is a small circle. The boundaries of several caps may intersect in the
empty set, a point, two points, or a higher-dimensional set, depending on the
rank of the normals and on the offsets. When several such constraints are
active, only those in the selected working set, whose normals are linearly
independent, remove directions from the tangent basis;
Section~\ref{subsec:riemannian-bounds} gives the normal and retraction
equations.

The unit quaternions $S^3$ arise in computer vision and quadcopter control. Their sphere geometry is represented analogously to $S^2$.

\subsubsection{The Orthogonal Groups \texorpdfstring{$O(n)$ and $SO(n)$}{O(n) and SO(n)}}
The orthogonal group
\[
O(n):=\{Y\in \mathbb{R}^{n\times n}: Y^\top Y=I_n\}
\]
consists of the rotations and reflections of $\mathbb R^n$. For two dimensions,
\[
O(2)=\left\{
 \begin{pmatrix}\cos\theta&-\varepsilon\sin\theta\\
 \sin\theta&\varepsilon\cos\theta\end{pmatrix}:
 \theta\in\mathbb R/(2\pi\mathbb Z),\ \varepsilon\in\{-1,1\}\right\}
\]
contains anticlockwise rotation by $\theta$ about the origin and reflection across a line making angle $\theta/2$ with the first coordinate axis. Thus it is in bijection with $[0,2\pi)\times\{-1,1\}$. This group provides the foundation for the rotational configurations used below. 

The tangent space of this manifold at a member $Y\in O(n)$ can be written as:
\[
T_YO(n) :=\left\{ Y \Omega,\,\,\Omega\in\mathbb{R}^{n\times n},\, \Omega^\top=-\Omega\right\}
\]
One computational choice is the \emph{Cayley retraction}, evaluated by solving a linear system. For $V=Y\Omega \in T_YO(n)$,
\[
R_Y(V) = Y(I-\frac{1}{2} \Omega)^{-1}(I+\frac{1}{2} \Omega)
\]
Alternatively, let $A=\widetilde Q\widetilde R$ be a QR factorization and
write $\mathcal Q_{QR}(A)=\widetilde Q$. The corresponding retraction is
\[
R_Y(V=Y\Omega) = \mathcal{Q}_{QR}(Y+V)=\mathcal{Q}_{QR}(Y(I+\Omega))
\]
Use the QR convention with positive diagonal entries of $\widetilde R$ near
the base point, so that $R_Y(0)=Y$ and the $SO(n)$ component is preserved.
Householder or Givens implementations may be used; the choice concerns
the cost and structure of the linear algebra for this retraction. 

The special orthogonal group is the determinant-one submanifold of $O(n)$, 
\[
SO(n):=\left\{Y\in O(n),\,\det(Y)=1\right\}
\]
This selects the proper rotations used to represent rigid-body orientation. 

For general $n$, the retractions above apply to $O(n)$. In the three-dimensional case $SO(3)$, an explicit exponential is also available. The skew-symmetric matrices defining $T_{Y}SO(3):=\{Y\Omega\}$ are 
\[
 \widehat\omega:=\begin{pmatrix}
 0&-\omega_3&\omega_2\\\omega_3&0&-\omega_1\\-\omega_2&\omega_1&0
 \end{pmatrix},\qquad
 \mathfrak{so}(3)=\{\widehat\omega:\omega\in\mathbb R^3\}.
 \]
This is the Lie algebra of $SO(3)$. In general the tangent space is $T_Y SO(n)=Y\mathfrak{so}(n)$, where $\mathfrak{so}(n)$ is the space of skew-symmetric $n\times n$ matrices.

Let $\Omega=\widehat\omega$ and $\theta=\left\|\begin{pmatrix} \omega_1 \\ \omega_2 \\ \omega_3 \end{pmatrix}\right\|$. The exponential retraction follows the geodesic with the prescribed initial tangent (Rodrigues' formula):
\[
R_{Y}(Y\Omega) = Y\left(I+\frac{\sin\theta}{\theta}\Omega+\frac{1-\cos\theta}{\theta^2}\Omega^2\right)
\]

The natural metric for $\mathfrak{so}(n)$ is the Frobenius-induced bi-invariant metric,
\[
\langle \Omega_1,\Omega_2\rangle =-\frac{1}{2}\operatorname{tr}\left(\Omega_1\Omega_2\right)
\]
For $\mathfrak{so}(3)$, let $\Omega_1,\Omega_2,\Omega_3$ be a basis with structure constants $\rho_{ijk}$, that is, $[\Omega_i,\Omega_j]=\sum_k\rho_{ijk} \Omega_k$. The Christoffel symbols and Levi-Civita connection are then
\[
\Gamma^k_{ij} = \frac{1}{2} \rho_{ijk},\,\,\,\nabla_{\Omega_i} \Omega_j = \frac12[\Omega_i,\Omega_j]=\frac12\sum_k\rho_{ijk}\Omega_k
\]

\subsubsection{The Real Symplectic Group \texorpdfstring{$Sp(2n)$}{Sp(2n)}}
The real symplectic group preserves the nondegenerate skew-symmetric bilinear form $a^\top J b$ on $\mathbb R^{2n}$, where
\[
J:=\begin{pmatrix} 0 & I \\ -I & 0 \end{pmatrix}
\]
and is defined by
\[
Sp(2n):=\left\{S\in\mathbb{R}^{2n\times 2n},\,S^\top J S = J\right\}
\]
The same construction applies after a change of coordinates $J'=P^\top J P$ with $P$ invertible. 

The tangent space at $Y\in Sp(2n)$ is equal to
\[
T_YSp(2n)=\left\{\Xi=Y\Omega,\,\Omega^\top J+J\Omega=0\right\}
\]
In two dimensions the Lie algebra is
\[
\mathfrak{sp}(2) = \left\{\begin{pmatrix} \alpha & \beta \\ \gamma & -\alpha\end{pmatrix}: \,\alpha,\beta,\gamma\in\mathbb{R}\right\}
\]
with $T_YSp(2)=Y\mathfrak{sp}(2)$.

The Cayley transform gives the computable retraction
\[
R_Y(\Xi)=R_Y(Y\Omega) = Y\left(I-\frac{1}{2}\Omega\right)^{-1}\left(I+\frac{1}{2}\Omega\right)
\]
on the neighborhood where $I-\Omega/2$ is invertible. Its dense solve has
cubic complexity in the component size; the full horizon KKT cost also
depends on the shooting structure and working set.

The noncompact real symplectic group admits a left-invariant Riemannian metric obtained by extending a positive-definite inner product on its Lie algebra. Its connection follows from the Koszul formula~\cite{selig2005geometric}.

\subsubsection{The Unitary Groups \texorpdfstring{$U(n)$ and $SU(n)$}{U(n) and SU(n)}}
The unitary and special unitary groups are the complex analogues of the orthogonal and special orthogonal groups:
\[
U(n):=\{Y\in\mathbb{C}^{n\times n},\,Y^* Y = I\},\,\,\,SU(n):=\{Y\in U(n),\,\det(Y)=1\}
\]
Here
\[
 T_YU(n)=\{Y\Omega:\Omega^*=-\Omega\},\qquad
 T_YSU(n)=\{Y\Omega:\Omega^*=-\Omega,\ \operatorname{tr}\Omega=0\}.
\]
The exponential $Y\exp(\Omega)$ preserves both requirements.  A generic
unitary Cayley or QR update may leave $SU(n)$. For a candidate unitary
update $U$, a local determinant correction is $U\exp(-\log(\det U)/n)$,
using a continuous logarithm branch near the identity.  That scalar factor is understood to multiply $I_n$.
Here $Y^*$ denotes the conjugate transpose. 

The special unitary group $SU(n)$ has a natural bi-invariant metric $\langle \Omega_1,\Omega_2\rangle=-\operatorname{tr}(\Omega_1\Omega_2)$. Multiplying the conventional Hermitian (generalized) Gell-Mann matrices by $i$ gives an anti-Hermitian trace-free basis, commonly used in quantum control and computing. 

\subsubsection{Rigid Body Motion \texorpdfstring{$SE(n)$}{SE(n)}}
Rigid body motions are defined by a combined rotation and translation:
\[
SE(n):=\left\{\begin{pmatrix} Y & z \\ 0 & 1\end{pmatrix} ,\,Y\in SO(n),\,z\in\mathbb{R}^n\right\}
\]
In three dimensions, this has the corresponding Lie algebra,
\[
\mathfrak{se}(3)=\left\{\begin{pmatrix}
    \widehat{\omega} & v \\ 0 & 0
\end{pmatrix} :\omega,v\in\mathbb{R}^3\right\}
\]
and in general the tangent space is $T_XSE(n) = X\mathfrak{se}(n)$. An element $(\omega,v)$ of $\mathfrak{se}(3)$ is called a \emph{twist}: $\omega$ is an angular velocity and $v$ a linear velocity.

The Lie bracket is
\[
[(\omega_1,v_1),(\omega_2,v_2)] = (\omega_1\times\omega_2,\omega_1\times v_2-\omega_2\times v_1)
\]
Using the same metrics as for $\mathfrak{so}(n)$ and $\mathbb{R}^n$ gives a left-invariant metric on $SE(n)$; it is not bi-invariant, and $SE(n)$ admits no positive-definite bi-invariant Riemannian metric. Its connection is
\[
\nabla_{(\omega,v)}(\eta,w) = \left(\frac{1}{2}\omega\times \eta,\omega\times w\right)
\]
for left-invariant vector fields.
\subsubsection{The Screw Subgroup \texorpdfstring{$H_p(n)$}{Hp(n)}}
The subgroup $H_p(n)\subset SE(n)$ represents a screw operation. 

For a concrete example, consider $H_2(3)\subset SE(3)$, the screw subgroup with translation $p=2$ per full revolution about the first coordinate axis. 

An element of $ H_2(3)$ is parameterized by a single angle $\theta$, and its Lie algebra has the generator $\xi$:
\[
H_2(3):=\left\{\begin{pmatrix} 1 & 0 & 0 & \theta/\pi \\
0 & \cos\theta & -\sin\theta & 0 \\
0 & \sin\theta & \cos\theta & 0 \\
0 & 0 & 0 & 1\end{pmatrix}:\theta\in\mathbb R\right\},\,\,\xi =\begin{pmatrix} 0 & 0 & 0 & 1/\pi \\ 0 & 0 & -1 & 0 \\ 0 & 1 & 0 & 0 \\ 0 & 0 & 0 & 0 \end{pmatrix}
\]
The angle $\theta\in\mathbb R$ is unwrapped: a $2\pi$ turn also translates
the body by $p$, so it is not identified with zero when $p\ne0$.  The retraction
adds $\psi$ to this unwrapped angle.  The screw pitch per radian used in
\cite{lynch2017modern} is $p/(2\pi)$ under this convention. 

\subsection{Intrinsic optimization on robotic configuration manifolds}
\label{subsec:why-riemannian-robotics}

Many robotic configuration spaces admit Euclidean embeddings and can therefore be
written as nonlinear equality-constrained sets.  Examples include
\[
 S^n=\{q\in\mathbb R^{n+1}:q^\top q=1\},\qquad
 SO(n)=\{R\in\mathbb R^{n\times n}:R^\top R=I,\ \det R=1\},
\]
and
\[
 Sp(2n)=\{S\in\mathbb R^{2n\times 2n}:S^\top JS=J\}.
\]
A conventional nonlinear-programming method can therefore be applied in the
ambient vector space. The embedding introduces variables, constraints and
directions associated with the representation. The intrinsic Riemannian
formulation eliminates these representation directions before forming the SQP
subproblem, retaining the physical degrees of freedom.

The distinction enters the second-order model. Suppose that
\[
 \mathcal M=\{r\in\mathbb R^N:h(r)=0\}
\]
is a regular embedded manifold.  An equality-constrained Euclidean formulation
uses the KKT matrix
\[
 K_{\rm emb}=
 \begin{pmatrix}
  \nabla^2_{rr}\mathcal L(r,\nu)&Dh(r)^\top\\
  Dh(r)&0
 \end{pmatrix}.
\]
The physically relevant second-order operator is instead the restriction of the
Lagrangian Hessian to
\[
 T_r\mathcal M=\ker Dh(r).
\]
If \(Z_r\) spans this tangent space, the relevant reduced Hessian is
\[
 Z_r^\top\nabla^2_{rr}\mathcal L(r,\nu)Z_r.
\]
At a constrained stationary point, with the induced Euclidean metric and
its corresponding geometric multipliers, this is the Riemannian Hessian
in the basis $Z_r$. At a nonstationary point the appropriate geometric
multiplier is determined by
\[
 \nu_{\mathrm{geom}}=-(Dh\,Dh^\top)^{-1}Dh\,\nabla\bar f(r),\qquad
 \operatorname{grad}_{\mathcal M}f
 =\nabla\bar f(r)+Dh(r)^\top\nu_{\mathrm{geom}},
\]
where $\bar f$ is an ambient extension. A second-derivative SQP model includes
the corresponding weighted constraint-Hessian term
$\sum_i(\nu_{\mathrm{geom}})_i\nabla^2h_i(r)$, which accounts for the curvature
of the manifold. The Riemannian formulation
performs this geometric reduction before the active-set and inertia calculations.
The ambient KKT system carries tangent coordinates, normal coordinates and
multipliers for the geometric identities, so its conditioning can depend
strongly on a redundant or poorly scaled representation.

This distinction is central to Newton-type methods.  A representation symmetry may
create a continuum of equivalent ambient points even when the physical solution is
isolated.  The standard example is a homogeneous eigenvector objective: scaling the
vector does not change the physical eigendirection, so the ambient problem has
nonisolated minimizers and a singular or misleading Newton model.  Restriction to
the sphere, or passage to the corresponding quotient manifold, removes the symmetry
and its continuous scaling direction. Isolation additionally requires,
for example, a simple eigenvalue; a repeated eigenvalue can retain a genuine
family of solutions. In robotics, normalization and orthogonality remove known
representation-normal directions~\cite{absil2008optimization}.

\paragraph{Exact geometric feasibility.}
A conventional constrained method satisfies nonlinear manifold equations to a
prescribed numerical tolerance. Intermediate matrices may therefore lose
orthogonality, quaternions may lose unit norm, and periodic joint coordinates
may cross an inappropriate chart branch. Valid representations are needed for
evaluating robotic states and controls. The intrinsic method computes
\[
 \eta_k\in T_{r_k}\mathcal M
\]
and updates by a retraction,
\[
 r_{k+1}=R_{r_k}(\eta_k).
\]
In exact arithmetic every defined retraction lies on the underlying manifold.
In floating-point arithmetic, the acceptance tests of
Section~\ref{sec:pathfollowing} monitor numerical geometric residuals.
A difference between the measured state and the model prediction is not a
geometric error: it enters the next feedback problem through the measured
initial node.
Geometric consistency is important when evaluating kinematics,
dynamics, and actuator commands during real-time iteration.

\paragraph{Spheres and unit quaternions.}
For \(S^n\), the ambient equation \(q^\top q=1\) introduces a normal direction
proportional to \(q\), whereas the physical tangent space is
\[
 T_qS^n=\{\eta:q^\top\eta=0\}.
\]
The intrinsic Newton equation uses only these \(n\) physical coordinates.  This
applies to directional controls on \(S^2\) and to unit quaternions on \(S^3\). The antipodal symmetry $q\sim-q$ of unit quaternions also requires care: $q$ and $-q$ represent the same rotation, so a large change in quaternion coordinates can represent no change in physical orientation.

\paragraph{Tori and revolute joints.}
For
\[
 q\in\mathbb T^m=(S^1)^m,
\]
a Euclidean angle vector ignores periodicity.  Two vectors differing by \(2\pi k\)
represent the same physical configuration, although an ordinary Euclidean norm may
treat them as far apart.  The torus formulation instead uses tangent increments in
\(\mathbb R^m\), wrapped differences, and the retraction
\[
 R_q(\eta)=(q+\eta)\bmod 2\pi.
\]
Consequently the local model depends on physical angular displacement rather than
the arbitrary location of a chart cut.

\paragraph{\(O(n)\) and \(SO(n)\).}
The direct representation of \(O(n)\) uses \(n^2\) matrix entries together with
\(n(n+1)/2\) orthogonality equations.  The physical dimension of \(SO(n)\) is only
\(n(n-1)/2\), and its tangent space is
\[
 T_RSO(n)=\{R\Omega:\Omega^\top=-\Omega\}.
\]
Thus an \(SO(3)\) update uses three tangent coordinates, eliminating the nine-entry ambient
representation and its orthogonality multipliers. Exponential, Cayley and QR-based
retractions preserve orthogonality in exact arithmetic.  The intrinsic Hessian consequently acts only on
rotational generator directions, while normal directions associated with
\(R^\top R=I\) never enter the QP factorization.

\paragraph{\(SE(n)\).}
For a rigid-body pose
\[
 X=\begin{pmatrix}R&t\\0&1\end{pmatrix}\in SE(n),
\]
the intrinsic tangent variable is a twist
\[
 \xi=\begin{pmatrix}\Omega&v\\0&0\end{pmatrix},
 \qquad \Omega^\top=-\Omega.
\]
The tangent coordinates coincide with the physical rotational and translational
degrees of freedom.  An ambient matrix formulation would additionally require
orthogonality, determinant selection, and fixed affine-row equations.  Eliminating
these identities before SQP gives a smaller intrinsic system. It removes the
conditioning effects of those representation constraints. 
The exact condition-number improvement depends on the chosen metric,
scaling and remaining physical constraints.

\paragraph{\(U(n)\), \(SU(n)\), and \(Sp(2n)\).}
For \(U(n)\) and \(SU(n)\), the intrinsic generators are skew-Hermitian, with an
additional trace-zero condition for \(SU(n)\).  For the symplectic group,
\[
 T_SSp(2n)=\{S\Omega:\Omega^\top J+J\Omega=0\}.
\]
In both cases the ambient equality representation introduces normal matrix
directions that violate the defining geometry.  The intrinsic formulation removes
them and leaves only admissible generator coordinates.  This is especially useful
for \(Sp(2n)\), whose noncompactness and metric choice can make an ambient
representation strongly scaling dependent.

The Riemannian formulation eliminates known geometric constraints before the
optimization step. The active-set layer then handles dynamics, collision avoidance,
contact conditions and actuator limits whose activity changes along the robotic
trajectory.

\section{Geometric Optimal-Control Case Studies}
\label{sec:geometric-case-studies}
We illustrate the intrinsic variables of Section~\ref{sec:introduction} with two robotic
optimal-control problems. The drone problem couples attitude,
angular-direction control and free arrival time. The arm problem couples
revolute-joint geometry, rigid-body motion and frictional wrench generation.

\subsection{Drone Navigation}
\label{sec:drone-example}

We consider a drone with position $x^p(t)\in\mathbb R^3$, translational
velocity $x^v(t)\in\mathbb R^3$, and attitude $q^R(t)\in SO(3)$.
The columns of $q^R$ define the drone's body axes expressed in world coordinates.  Thus the
state blocks are
\[
 x(t)=\begin{pmatrix}x^p(t)\\x^v(t)\end{pmatrix}\in\mathbb R^6,
 \qquad q(t)=q^R(t)\in SO(3).
\]
The superscript $R$ identifies the rotational state; $\dot q^R$ denotes its
time derivative.

The geometric control $u(t)\in S^2$ specifies the \emph{body-frame}
direction of angular velocity, whose magnitude $\omega>0$ is fixed.
The Euclidean scalar control $v(t)$ is the total thrust magnitude in newtons,
with $0\le v(t)\le v_{\max}$.  We use an \emph{ideal kinematic actuator
model}: an inner controller realizes the commanded angular velocity
$\omega u$ and thrust without actuator lag.  More detailed geometric models can include angular-momentum dynamics,
rotor mixing, aerodynamic drag and motor saturation.  

\begin{figure}[tbp]
\centering
\begin{tikzpicture}[>=stealth, every node/.style={font=\small},
 arrow/.style={->,thick}, box/.style={draw,rounded corners,align=center,
 inner sep=5pt}]
 \coordinate (O) at (-3.2,-1.1);
 \coordinate (C) at (0,0);
 \draw[thick] (-1.15,0)--(1.15,0);
 \draw[thick] (0,-0.65)--(0,0.65);
 \foreach \p in {(-1.15,0),(1.15,0),(0,-0.65),(0,0.65)}
   {\draw[thick] \p circle (0.15);}
 \draw[fill=gray!15,thick] (C) circle (0.26);
 \node[above left] at (-0.3,0.4) {$q^R(t)\in SO(3)$};
 \draw[arrow] (O)--(-2.3,-1.1) node[right] {$e_1$};
 \draw[arrow] (O)--(-3.2,-0.25) node[above] {$e_3$};
 \node[below] at (O) {world};
 \draw[arrow,blue] (O)--(-0.25,-0.1);
 \node[blue,above] at (-1.9,-0.75) {$x^p(t)$};
 \draw[arrow,purple] (0.5,-0.22)--(2,-0.65);
 \node[purple,below] at (1.5,-0.7) {$x^v(t)$};
 \draw[arrow,red] (0,0.25)--(0,1.7);
 \node[red,above,align=center] at (0,1.7)
   {$v(t)q^R(t)e_3$\\thrust force};
 \draw[arrow] (-0.1,-0.28)--(-0.1,-1.65);
 \node[below] at (-0.1,-1.65) {$-mg e_3$};
 \draw[thick] (3.5,0.25) circle (0.63);
 \draw[dashed] (2.87,0.25) arc (180:360:0.63 and 0.19);
 \draw[dashed] (4.13,0.25) arc (0:180:0.63 and 0.19);
 \draw[arrow,orange] (3.5,0.25)--(3.92,0.66);
 \node[above right,orange] at (3.82,0.75) {$u(t)\in S^2$};
 \node[below,align=center] at (3.5,-0.5)
   {body angular direction\\angular speed $\omega$};
 \node[box] at (0.7,-2.75) {Euclidean state $x=(x^p,x^v)\in\mathbb R^6$;
 geometric state $q=q^R\in SO(3)$\\
 $\dot x^p=x^v,\quad \dot x^v=-ge_3+(v/m)q^Re_3,\quad
 \dot q^R=q^R\widehat{\omega u}$\\
 $0\le v\le v_{\max},\quad u\in S^2,\quad c(x^p)\ge0$};
\end{tikzpicture}
\caption{Drone states and controls under the ideal kinematic actuator model.
The hat in the attitude equation is the skew map
\eqref{eq:example-hat}, and the two force arrows distinguish thrust from gravity.}
\label{fig:drone_variables}
\end{figure}
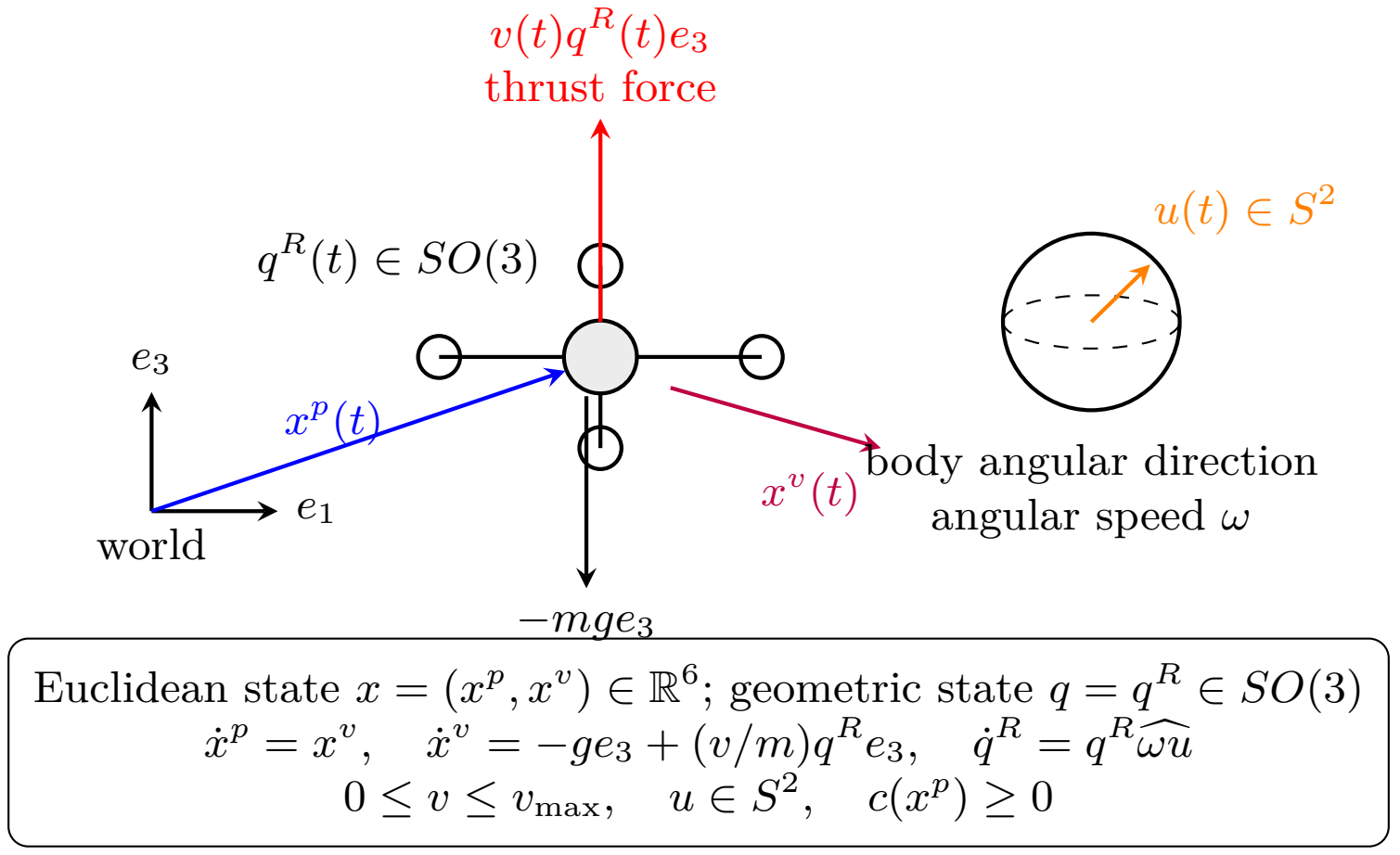

Let $e_3=(0,0,1)^\top$ point upward in the world frame, let $m>0$ be
the vehicle mass, and let $g>0$ be gravitational acceleration.  For
$a=(a_1,a_2,a_3)^\top\in\mathbb R^3$, define
\begin{equation}\label{eq:example-hat}
 \widehat a:=\begin{pmatrix}
 0&-a_3&a_2\\ a_3&0&-a_1\\-a_2&a_1&0
 \end{pmatrix},\qquad \widehat a\,b=a\times b.
\end{equation}
This is the skew-matrix operation denoted $[a]$ in
\cite[Section~3.2.2]{lynch2017modern}.  The body-angular-velocity convention gives
$(q^R)^\top\dot q^R=\widehat{\omega u}$, so the coupled dynamics are
\begin{equation}\label{eq:ex1dyn}
 \begin{aligned}
 \dot x^p(t)&=x^v(t),\\
 \dot x^v(t)&=-g e_3+\frac{v(t)}m q^R(t)e_3,\\
 \dot q^R(t)&=q^R(t)\widehat{\omega u(t)}.
 \end{aligned}
\end{equation}
The attitude equation is a $3\times3$ matrix equality. The vector
$q^Re_3$ gives the unit thrust direction in world coordinates, and
$vq^Re_3$ is the corresponding thrust \emph{force}.  Figure~\ref{fig:drone_variables}
uses the same notation and includes gravity.

\paragraph{Free final time and collision constraints.}
For the initial plan, choose the arrival time $T>0$ as a decision variable
and minimize
\begin{equation}\label{eq:ex1cost}
 J(v,T)=\int_0^T \rho_v |v(t)|^2\,dt+V(T),
 \qquad V(T)=\rho_T T,\qquad \rho_v,\rho_T>0.
\end{equation}
The weights reconcile the units of the two terms.  The boundary conditions are
\begin{equation}\label{eq:ex1initfinal}
 x^p(0)=a,\qquad x^v(0)=0,\qquad q^R(0)=q^R_0,
 \qquad x^p(T)=b,\qquad x^v(T)=0.
\end{equation}
Terminal attitude is free.  No zero terminal angular velocity is imposed,
because that would contradict the fixed positive angular speed.

For a stationary obstacle $\mathcal O_j$, let $d_{\mathcal O_j}$ denote
signed distance, positive outside.  A known enclosing radius $a_d>0$
contains the drone and its swept rotor disks in every attitude.  Then
\begin{equation}\label{eq:ex1ineqcon}
 c_j(x^p):=d_{\mathcal O_j}(x^p)-a_d-d_{\rm safe}\ge0,
 \qquad j=1,\ldots,N_{\rm obs},
\end{equation}
is a conservative finite-body collision-avoidance constraint.  Ground is included by
$e_3^\top x^p-a_d-d_{\rm ground}\ge0$ when required.  Smooth local
clearance representations and between-node verification are discussed
below for the arm.  The kinematic feature block is empty, so $d_z=0$; position and velocity
are differential states in $x$, while $z$ is absent.

\paragraph{Time reparameterization and shrinking horizon.}
Let $\hat t=0,1,\ldots,N-1$ be the \emph{integer feedback-sample index},
let $t_{\hat t}$ be the current physical time, and set
\[
 H_{\hat t}=N-\hat t,\qquad H_0=N,\qquad
 \mathcal I_{\hat t}=\{0,\ldots,H_{\hat t}-1\},\qquad
 \mathcal N_{\hat t}=\{0,\ldots,H_{\hat t}\}.
\]
At sample $\hat t$, the remaining physical duration $T_{\hat t}$ is an
additional shared Euclidean decision variable.  With normalized time
$\sigma\in[0,1]$,
\begin{equation}\label{eq:drone-time-scaling}
 t=t_{\hat t}+T_{\hat t}\sigma,\quad
 \frac{d x^p}{d\sigma}=T_{\hat t}x^v,\quad
 \frac{d x^v}{d\sigma}=T_{\hat t}\left(-ge_3+\frac vm q^Re_3\right),
 \quad
 \frac{d q^R}{d\sigma}=T_{\hat t}q^R\widehat{\omega u}.
\end{equation}
For $H_{\hat t}\geq1$, use $\sigma_k=k/H_{\hat t}$ and physical interval length
$h_{\hat t}=T_{\hat t}/H_{\hat t}$.  Impose
$0<h_{\min}\le h_{\hat t}\le h_{\max}$, with prescribed integration
and actuation limits.  This is a shrinking-horizon, variable-sampling
implementation. The number of remaining intervals decreases by one at every
sample, while their common length is reoptimized at every sample: the first
control is applied for the optimized $h_{\hat t}^\star$, and
$t_{\hat t+1}=t_{\hat t}+h_{\hat t}^\star$.
At $H_{\hat t}=0$ no interval remains, so no optimization problem for the
duration and controls is formed; only the terminal conditions are tested on
the measured state. 
The distinction between a path and its time scaling follows
\cite[Chapter~9]{lynch2017modern}; the shrinking feedback policy is
specified here for this example.

For held controls, \eqref{eq:ex1dyn} admits an exact interval map.  Define $\Omega_k=\widehat{\omega u_k}$ and
\[
 \Phi_1(h,\Omega)=\int_0^h e^{s\Omega}\,ds,\qquad
 \Phi_2(h,\Omega)=\int_0^h(h-s)e^{s\Omega}\,ds.
\]
These matrix functions can be evaluated by their analytic expressions or
convergent power series.  With $h=h_{\hat t}$, the discrete maps are
\begin{equation}\label{eq:drone-discrete-map}
 \begin{aligned}
 q^R_{k+1}&=q^R_k e^{h\Omega_k},\\
 x^v_{k+1}&=x^v_k-gh e_3+\frac{v_k}{m}q^R_k\Phi_1(h,\Omega_k)e_3,\\
 x^p_{k+1}&=x^p_k+h x^v_k-\frac12gh^2e_3
             +\frac{v_k}{m}q^R_k\Phi_2(h,\Omega_k)e_3.
 \end{aligned}
\end{equation}
The maps preserve $SO(3)$ and are smooth in the states, controls and
positive duration.  A consistently differentiated Lie-group integrator
is an alternative.

\begin{samepage}
Suppressing the prediction superscript $\hat t$ on individual nodes,
the full finite-dimensional problem is
\begin{equation}\label{eq:ex1nlp}
\begin{aligned}
 \min_{\substack{T_{\hat t},\ (u_k,v_k)_{k\in\mathcal I_{\hat t}}\\
                         (x_k,q^R_k)_{k\in\mathcal N_{\hat t}}}}
 &\quad h_{\hat t}\sum_{k\in\mathcal I_{\hat t}}\rho_v v_k^2
                      +\rho_T T_{\hat t}\\
 \textnormal{s.t.}\quad
 &(x_{k+1},q^R_{k+1})\text{ satisfy \eqref{eq:drone-discrete-map}},
                         &&k\in\mathcal I_{\hat t},\\
 &c(x^p_k)\ge0,&&k=1,\ldots,H_{\hat t},\\
 &q^R_k\in SO(3),&&k\in\mathcal N_{\hat t},\\
 &u_k\in S^2,\quad 0\le v_k\le v_{\max},&&k\in\mathcal I_{\hat t},\\
 &H_{\hat t}h_{\min}\le T_{\hat t}\le H_{\hat t}h_{\max},\\
 &x_0=\hat x(\hat t),\quad q^R_0=\hat q^R(\hat t),\\
 &x^p_{H_{\hat t}}=b,\quad x^v_{H_{\hat t}}=0.
\end{aligned}
\end{equation}
\end{samepage}
The terminal attitude $q^R_{H_{\hat t}}$ remains a free decision variable.  The first measured state has
$\hat x(0)=(a,0)$ and $\hat q^R(0)=q^R_0$.
Elapsed physical time adds only a constant to the linear arrival-time
penalty, so penalizing the remaining $T_{\hat t}$ gives the same
optimizer at this sample.  The duration $T_{\hat t}$ is a single variable shared by all intervals. It
therefore appears in every shooting equation, and its first and second
derivatives enter the gradient, Jacobian and Hessian assembled for the SQP
step; \eqref{eq:ex1nlp} is the free-final-time extension of the
fixed-duration structure in Section~\ref{sec:introduction}.  Each control pair is applied
only once before the state is measured again.  The measured node is substituted as data before the KKT system is formed;
state-clearance constraints begin at the next predicted node. At the
terminal node the position is fixed by the terminal equality, so
$c(x^p_{H_{\hat t}})=c(b)$ is a constant: it is verified once on the data
and is not carried as a constraint row. After the
first control has been applied and the next state measured, the remaining
part of the previous solution (its tail) is shifted by one interval,
$H_{\hat t}$ decreases by one, and the remaining
duration is reoptimized on the new fixed layout. 

A sphere retraction is used for SQP updates of $u$, and an $SO(3)$
retraction for updates of $q^R$.  The identities
$(q^R)^\top q^R=I_3$ and $\det(q^R)=1$ are enforced by the manifold
representation, with no additional ambient NLP equations.

\subsection{Arm Navigation and Grasping}
\label{sec:arm-example}

This example uses the UR5 six-revolute-joint kinematics of
\cite[Example~4.5 and Figure~4.6]{lynch2017modern}, together with a
specified parallel-jaw gripper and a rigid cuboidal payload.  The mechanics below follow the book's
product-of-exponentials, rigid-body dynamics and frictional-contact
conventions.  Link inertias, actuator limits, and friction coefficients
are supplied calibration data.

We formulate three optimal-control problems, solved successively in
closed-loop MPC. The arm first approaches and grasps the cuboidal object,
then carries it to the target, and finally places and releases it before
retreating. The task is split because the set of constraints changes with
the contact state: before the grasp the box is a fixed obstacle, during
transport it is rigidly attached to the hand, and after release it is a
separate body resting on a support. Within each problem the constraint
list is fixed, so each problem is a smooth NLP to which the method of
Section~\ref{sec:pathfollowing} applies; the switches are made by the
supervisor on measured tests.

The basic arm variables are the joint configuration $q\in T^6=(S^1)^6$,
joint rate $x^v\in\mathbb R^6$, jaw aperture $\gamma\in\mathbb R$, joint
torque $v^q\in\mathbb R^6$, and aperture-rate command
$v^\gamma\in\mathbb R$. A freely moving payload adds a pose
$X_b\in SE(3)$ and body twist $\mathcal V_b\in\mathbb R^6$.
Table~\ref{tab:arm-layout-variables} summarizes the decision blocks.
During a rigid grasp the object pose is determined by the arm, while
during supported release it becomes an independent geometric state.
The aperture and its rate are substituted as calibrated constants while
the force servo holds the grasp. Contact-force allocations are auxiliary primal variables for wrench
balance: the optimizer chooses how the load is distributed over the contact
points so that the balance holds, but these forces are not commanded. The
physical controls remain the joint torques and aperture command.

\begin{table}[H]
\centering\small
\setlength{\tabcolsep}{4pt}
\begin{tabular}{p{.17\linewidth}p{.24\linewidth}p{.26\linewidth}p{.21\linewidth}}
\toprule
Problem/layout & Geometric state & Euclidean state & Controls\\
\midrule
$\mathsf A$: acquire & $q\in T^6$ & $(x^v,\gamma)\in\mathbb R^7$ & $(v^q,v^\gamma)\in\mathbb R^7$\\[3pt]
$\mathsf B$: carry & $q\in T^6$ & $x^v\in\mathbb R^6$ & $v^q\in\mathbb R^6$\\[3pt]
$\mathsf C$: held placement & $q\in T^6$ & $x^v\in\mathbb R^6$ & $v^q\in\mathbb R^6$\\[3pt]
$\mathsf C$: release/retreat & $(q,X_b)\in T^6\times SE(3)$ & $(x^v,\gamma,\mathcal V_b)\in\mathbb R^{13}$ & $(v^q,v^\gamma)\in\mathbb R^7$\\
\bottomrule
\end{tabular}
\caption{State and control variables of the three arm OCPs. The two
listed layouts of $\mathsf C$ belong to the same placement--release
problem. A measured layout is fixed during each solve.}
\label{tab:arm-layout-variables}
\end{table}

\begin{figure}[H]
\centering
\resizebox{\linewidth}{!}{%
\begin{tikzpicture}[>=stealth,line cap=round,line join=round,
 every node/.style={font=\small}]
 \begin{scope}
  \node[font=\bfseries] at (1.95,2.75) {(a) Arm, hand and payload};
  \draw[thick] (-.5,-.27)--(4.9,-.27);
  \fill[gray!30] (-.35,-.25) rectangle (.35,.02);
  \draw[thick] (-.35,-.25) rectangle (.35,.02);
  \draw[line width=11pt,gray!65] (0,.08)--(0,.75)--(1.1,1.9)--(2.35,1.65)--(2.8,1.3)--(3.2,1.3)--(3.52,1.3);
  \draw[line width=6pt,blue!35] (0,.08)--(0,.75)--(1.1,1.9)--(2.35,1.65)--(2.8,1.3)--(3.2,1.3)--(3.52,1.3);
  \foreach \x/\y in {0/.08,0/.75,1.1/1.9,2.35/1.65,2.8/1.3,3.2/1.3}
   {\draw[fill=white,line width=1pt] (\x,\y) circle (.12);}
  \node[right] at (.15,.1) {$\theta_1$};
  \node[left] at (-.12,.75) {$\theta_2$};
  \node[above] at (1.1,2.04) {$\theta_3$};
  \node[above] at (2.35,1.81) {$\theta_4$};
  \node[below] at (2.8,1.13) {$\theta_5$};
  \node[below] at (3.25,1.13) {$\theta_6$};
  \draw[fill=blue!18,thick] (3.48,.99) rectangle (3.66,1.61);
  \draw[line width=3pt,blue!65!black] (3.62,1.56)--(4.12,1.56)--(4.12,1.43);
  \draw[line width=3pt,blue!65!black] (3.62,1.04)--(4.12,1.04)--(4.12,1.17);
  \draw[fill=orange!22,thick] (3.98,1.17) rectangle (4.58,1.43);
  \draw[->,thick] (3.55,1.3)--(3.55,1.99) node[above] {$e^e_3$};
  \draw[->,thick] (3.55,1.3)--(4.25,1.3) node[below right] {$e^e_1$};
  \node[above] at (4.35,2.18) {hand frame $\{e\}$};
  \node[right] at (4.61,1.55) {$X_b$};
  \draw[->,thick] (-1.02,-.15)--(-1.02,.48) node[above] {$e_3$};
  \draw[->,thick] (-1.02,-.15)--(-.42,-.15) node[below] {$e_1$};
  \node[below,align=center] at (2.0,-.43)
   {$q\in T^6$, $x^v\in\mathbb R^6$, torques $v^q$\\aperture $\gamma$ and rate $v^\gamma$};
 \end{scope}
 \begin{scope}[xshift=7.0cm]
  \node[font=\bfseries] at (1.7,2.75) {(b) Finite-pad wrench geometry};
  \draw[fill=orange!12,thick] (.8,.4) rectangle (2.5,1.7);
  \draw[fill=orange!18,thick] (.8,1.7)--(.98,1.86)--(2.68,1.86)--(2.5,1.7)--cycle;
  \draw[fill=orange!25,thick] (2.5,.4)--(2.68,.56)--(2.68,1.86)--(2.5,1.7)--cycle;
  \draw[fill=blue!20,thick] (.35,.22) rectangle (.66,1.99);
  \draw[fill=blue!20,thick] (2.79,.22) rectangle (3.10,1.99);
  \foreach \x/\y in {.8/.65,.8/1.45,.98/.81,.98/1.61,2.5/.65,2.5/1.45,2.68/.81,2.68/1.61}
   {\filldraw[red!80!black] (\x,\y) circle (.035);}
  \fill[green!13] (.8,.65)--(1.45,.39)--(1.45,.91)--cycle;
  \draw[->,green!50!black,thick] (.8,.65)--(1.45,.39);
  \draw[->,green!50!black,thick] (.8,.65)--(1.45,.91);
  \fill[green!13] (2.5,1.45)--(1.85,1.19)--(1.85,1.71)--cycle;
  \draw[->,green!50!black,thick] (2.5,1.45)--(1.85,1.19);
  \draw[->,green!50!black,thick] (2.5,1.45)--(1.85,1.71);
  \draw[->,orange!90!black,line width=1.3pt] (.8,.65)--(1.41,.76);
  \node[orange!90!black] at (1.35,.18) {$f_j$};
  \draw[->,blue!75!black] (2.5,1.45)--(1.98,1.45);
  \fill (1.65,1.05) circle (.035);
  \node[below right] at (1.65,1.05) {$\{b\}$};
  \draw[->,purple,dashed,thick] (1.65,1.05)--(.8,.65) node[midway,above left] {$\rho_j$};
  \node[left,align=center] at (.3,1.1) {pad\\4 points};
  \node[right,align=center] at (3.15,1.1) {pad\\4 points};
  \node[above,green!40!black] at (1.7,2.0) {inscribed friction-cone rays};
  \node[below,align=center] at (1.7,-.10)
   {$\mathcal W_j=\operatorname{col}(\rho_j\times f_j,f_j)$\\
    $\operatorname{rank}G_{\rm ray}=6$, $G_{\rm ray}\alpha=0$, $\alpha>0$};
 \end{scope}
\end{tikzpicture}}
\caption{Geometric interpretation of the arm problem. (a) A schematic
projection shows the six revolute joints, gripper and rigid payload;
the precise kinematics are specified by Table~\ref{tab:robotarm} and
\eqref{eq:arm-poe}. (b) Each of the two pads has four spatially
separated reactions. Their moment arms generate torque as well as
force. Green arrows show cone rays, blue arrows the inward normal, and
the orange arrow an admissible contact force. Full wrench rank and a
strictly positive balanced ray combination jointly express
three-dimensional force closure in the object-frame grasp model.}
\label{fig:arm-robot-and-force-closure}
\end{figure}
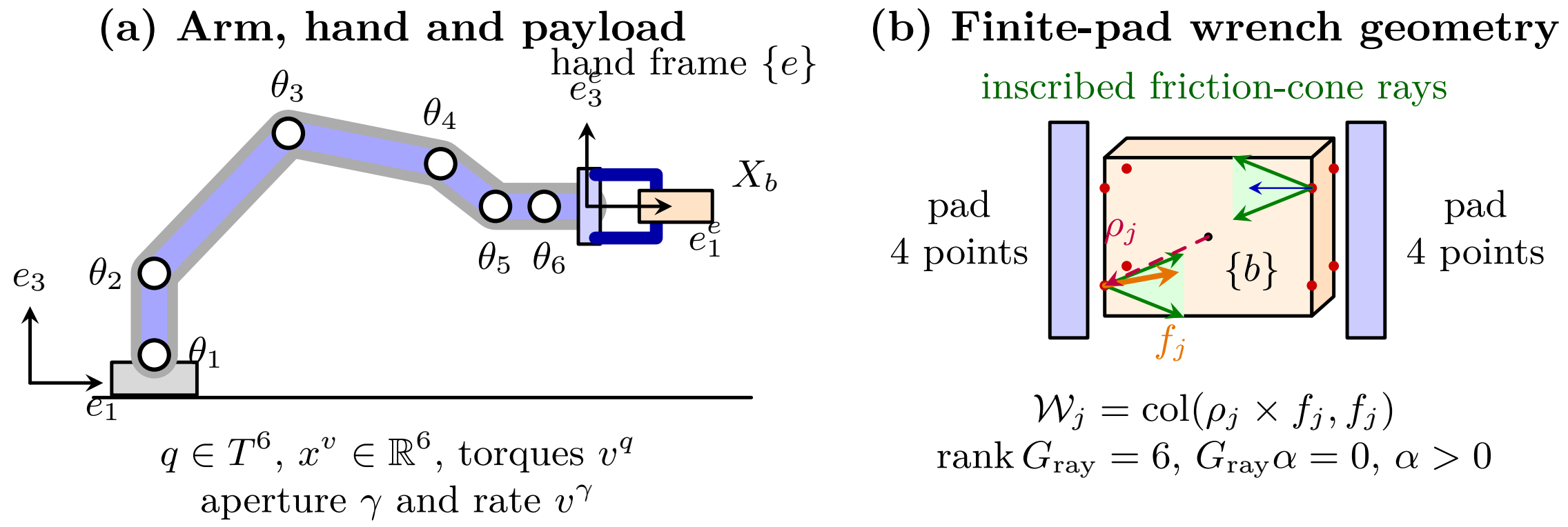

\subsubsection{Joint geometry, controls and forward kinematics}
\paragraph{Axis points and twists.}
All home-configuration axis data are expressed in the fixed base frame
$\{0\}$.  The vector $p_i\in\mathbb R^3$ is a \emph{point on joint
axis $i$}, and $\omega_i$ is its unit direction.  The associated
zero-pitch revolute screw twist and its matrix are
\begin{equation}\label{eq:arm-screw-definition}
 \mathcal S_i=\begin{pmatrix}\omega_i\\-\omega_i\times p_i\end{pmatrix},
 \qquad
 [\mathcal S_i]=\begin{pmatrix}
 \widehat{\omega_i}&-\omega_i\times p_i\\0&0
 \end{pmatrix}.
\end{equation}
Changing $p_i$ by a multiple of $\omega_i$ leaves the same axis and
same twist.  Table~\ref{tab:robotarm} gives representative axis points
whose derived linear twist components reproduce the cited UR5 table.
In particular, the linear component of the first twist is zero. The link
dimensions $L_1,L_2,W_1,W_2,H_1,H_2$ are the textbook's symbols; they are
used only in Table~\ref{tab:robotarm} and in \eqref{eq:arm-poe} and are
unrelated to the horizon $H$ and the other uses of these letters.

\begin{table}[tbp]
\centering
\small
\setlength{\tabcolsep}{4pt}
\begin{tabular}{c c c c}
\toprule
$i$ & $\omega_i$ & axis point $p_i$ & $-\omega_i\times p_i$\\
\midrule
1 & $(0,0,1)$ & $(0,0,0)$ & $(0,0,0)$\\
2 & $(0,1,0)$ & $(0,0,H_1)$ & $(-H_1,0,0)$\\
3 & $(0,1,0)$ & $(L_1,0,H_1)$ & $(-H_1,0,L_1)$\\
4 & $(0,1,0)$ & $(L_1+L_2,0,H_1)$ & $(-H_1,0,L_1+L_2)$\\
5 & $(0,0,-1)$ & $(L_1+L_2,W_1,H_1)$ & $(-W_1,L_1+L_2,0)$\\
6 & $(0,1,0)$ & $(L_1+L_2,W_1,H_1-H_2)$ & $(H_2-H_1,0,L_1+L_2)$\\
\bottomrule
\end{tabular}
\caption{UR5 home axes: points $p_i$ and the corresponding screw twists
are distinct quantities.  In the textbook example,
$(L_1,L_2,W_1,W_2,H_1,H_2)=(425,392,109,82,89,95)$ mm
\cite[Example~4.5]{lynch2017modern}.}
\label{tab:robotarm}
\end{table}

Let $M_6$ be the textbook home flange pose and let $X_{6e}\in SE(3)$ be
the known fixed flange-to-gripper transform.  Then
\begin{equation}\label{eq:arm-poe}
 \begin{aligned}
 M_6&=\begin{pmatrix}
 -1&0&0&L_1+L_2\\0&0&1&W_1+W_2\\
 0&1&0&H_1-H_2\\0&0&0&1
 \end{pmatrix},\qquad M_e=M_6X_{6e},\\
 F_e(q)&=\left(\prod_{i=1}^6 e^{[\mathcal S_i]\theta_i}\right)M_e
       =\begin{pmatrix}R_e(q)&p_e(q)\\0&1\end{pmatrix},\\
 F_i(q)&=\left(\prod_{j=1}^i e^{[\mathcal S_j]\theta_j}\right)M_i.
 \end{aligned}
\end{equation}
The $M_i$ are the known home poses of the selected link geometry frames;
products are ordered in increasing joint index from left to right.
The gripper frame $\{e\}$ is attached to the palm and remains fixed as
the jaws translate.
All lengths used in numerical evaluation are converted to meters.

\paragraph{States, controls, and simple bounds.}
The geometric state is
\[
 q=(\theta_1,\ldots,\theta_6)\in T^6.
\]
Write $x^v\in\mathbb R^6$ for joint angular rates, so that
$\dot\theta=x^v$ in an unwrapped joint chart, and define
\begin{equation}\label{eq:arm-state-control}
 x=\begin{pmatrix}x^v\\\gamma\end{pmatrix}\in\mathbb R^7,
 \qquad
 v=\begin{pmatrix}v^q\\v^\gamma\end{pmatrix}\in\mathbb R^7.
\end{equation}
Here $\gamma$ is the full inner-jaw separation in meters,
$v^q$ is the six-vector of joint torques, and $v^\gamma$ is the aperture
rate in meters per second.  There is no geometric control: the
$u$-product is empty.  Joint torques are assumed to be realized by an
inner servo.  Aperture motion is an ideal quasistatic velocity servo
when the hand is free; after measured contact, a force-regulating inner servo
maintains the prescribed grip. The normal-force setting is an inner-servo parameter, distinct from the
aperture-rate command. The active servo regulates either velocity or force
in the constrained closing direction.
This distinction follows the motion/force-control discussion in
\cite[Sections~11.3--11.6]{lynch2017modern}.

The admissible variables satisfy
\begin{equation}\label{eq:arm-simple-bounds}
 \begin{gathered}
 \theta_i^\ell\le\theta_i\le\theta_i^u,\qquad
 \underline x_i^v\le x_i^v\le\overline x_i^v,\qquad
 \underline v_i^q\le v_i^q\le\overline v_i^q,\quad i=1,\ldots,6,\\
 0\le\gamma\le\gamma_{\max},\qquad
 -\overline v^\gamma\le v^\gamma\le\overline v^\gamma.
 \end{gathered}
\end{equation}
Joint limits are evaluated on the chosen physical unwrapped joint
intervals.  Locally these intervals are charts of the circular joint
factors.  If a motor can turn through more than one revolution and successive turns
must be distinguished, the angle is not a point of the circle; the
unwrapped encoder value, which is its lift to the real line, is then
retained as the variable.  The free aperture evolves according to 
\begin{equation}\label{eq:arm-aperture}
 \dot\gamma=v^\gamma,\qquad
 \gamma_{k+1}=\gamma_k+h v^\gamma_k.
\end{equation}
The optimization problem enforces its travel and rate limits through
\eqref{eq:arm-simple-bounds}.

The discrete index $k$ labels task nodes, with fixed physical
sampling interval $h>0$ for this arm example.

\subsubsection{Rigid payload, finite pads and force closure}
Let the exactly known box in its body frame $\{b\}$ be
\begin{equation}\label{eq:arm-box-geometry}
 \mathcal B_0=\{\xi\in\mathbb R^3:|\xi_i|\le a_i,\ i=1,2,3\},
 \quad a_i>0,\qquad
 \mathcal B(X_b)=p_b+R_b\mathcal B_0,
 \quad X_b=\begin{pmatrix}R_b&p_b\\0&1\end{pmatrix}.
\end{equation}
We take the center of mass to be the box center and the body axes to
be face-aligned.  Mass $m_b>0$ and inertia $I_b\succ0$ about this center
are known additional data.
For a homogeneous box, for example,
\[
 I_b=\frac{m_b}{3}\operatorname{diag}
 (a_2^2+a_3^2,\ a_1^2+a_3^2,\ a_1^2+a_2^2).
\]
A measured noncentral center of mass can instead be handled by the
usual rigid change of spatial inertia and contact moment arms
\cite[Section~8.2.3]{lynch2017modern}.

The hand has two opposing flat contact pads on the faces $\xi_1=\pm a_1$.
Each pad is represented by four frictional point reactions at known
interior face locations
\begin{equation}\label{eq:arm-contact-geometry}
 \rho_j=(\sigma a_1,\epsilon d_2,\delta d_3)^\top,
 \quad n_j=-\sigma e_1,\quad
 (\sigma,\epsilon,\delta)\in\{-1,1\}^3,
 \quad 0<d_2<a_2,\quad0<d_3<a_3.
\end{equation}
Thus $N_c=8$, and $n_j$ is the unit normal pointing \emph{into} the
object.  These distributed reactions describe the two finite pads. Their spatial
extent supplies resistance to moments about the line between the pads,
which two point contacts alone cannot provide in three dimensions
\cite[Section~12.2.3.1]{lynch2017modern}.

Let $X_g={}^bT_e$ denote the fixed pose of the gripper frame expressed
in box coordinates.  It maps gripper-coordinate points into box
coordinates.  Let $\rho_j^e(\gamma)$ be the corresponding pad points in
the gripper frame.  The calibrated grasp geometry satisfies
\begin{equation}\label{eq:arm-contact-compatibility}
 X_g\begin{pmatrix}\rho_j^e(\gamma_g)\\1\end{pmatrix}
   =\begin{pmatrix}\rho_j\\1\end{pmatrix},\qquad
 \gamma_g=2a_1,\qquad
 0<\gamma_g<\gamma_{\rm open}\le\gamma_{\max}.
\end{equation}
Writing $X_g=(R_g,p_g;0,1)$, one explicit parallel-jaw contact-point map is
\[
 \rho_j^e(\gamma)=R_g^\top\!\left(
 \rho_j-p_g+\frac{\sigma(j)}2(\gamma-\gamma_g)e_1\right).
\]
The equality $\gamma_g=2a_1$ uses aligned parallel jaws and an aperture
measured between the contact faces; pad thickness is already included
in the calibration.  Consequently zero aperture is \emph{not} the
grasp state.  Pose compatibility is $F_e(q)=X_bX_g$.  In world
coordinates the pad points are
$c_j^0(q,\gamma)=p_e(q)+R_e(q)\rho_j^e(\gamma)$, and sticking contact
requires
\[
 c_j^0(q,\gamma_g)=p_b+R_b\rho_j,\qquad
 \dot c_j^0=\dot p_b+\omega_b^0\times(R_b\rho_j),
\]
where $\widehat{\omega_b^0}=\dot R_bR_b^\top$.
The latter is a no-relative-slip condition, not a new independent
state equation when rigid attachment is imposed.

All contact forces $f_j$, normals $n_j$, moment arms $\rho_j$, twists
and wrenches below are expressed in the \emph{same body frame} $\{b\}$.
With known coefficient $\mu_f>0$, Coulomb admissibility is
\begin{equation}\label{eq:arm-friction-cone}
 f_j\in\mathcal K_j:=\left\{f\in\mathbb R^3:
 n_j^\top f\ge0,\quad
 \|f-(n_j^\top f)n_j\|_2\le\mu_f n_j^\top f\right\}.
\end{equation}
Following \cite[Sections~3.4 and 12.2.1]{lynch2017modern}, the wrench
ordering is \emph{moment first, force second}:
\begin{equation}\label{eq:arm-grasp-map}
 \mathcal W_j=\begin{pmatrix}\rho_j\times f_j\\f_j\end{pmatrix},
 \qquad
 \mathcal W_{\rm pads}=Gf,\qquad
 G=\begin{pmatrix}
 \widehat{\rho_1}&\cdots&\widehat{\rho_{N_c}}\\
 I_3&\cdots&I_3
 \end{pmatrix},\quad
 f=\operatorname{col}(f_1,\ldots,f_{N_c}).
\end{equation}
A world-frame force is transformed by $R_b^\top$, with its moment
taken about the stated body-frame origin.  

For a conservative finite test, one can approximate each circular cone by the
\emph{inscribed} cone generated by
\[
 d_{j\ell}=n_j+\mu_f(t_{j1}\cos\phi_\ell+t_{j2}\sin\phi_\ell),
 \qquad \phi_\ell=2\pi\ell/L,\quad L\ge3,
\]
where $(t_{j1},t_{j2})$ is an orthonormal tangent basis of face $j$.
The finite ray-wrench matrix is
$G_{\rm ray}=[\operatorname{col}(\rho_j\times d_{j\ell},d_{j\ell})]_{j,\ell}$.
For this finite-ray contact model, force closure is characterized by
$\operatorname{rank}G_{\rm ray}=6$ and $\varepsilon_*>0$, where
\begin{equation}\label{eq:arm-force-closure}
 \begin{aligned}
 \varepsilon_*:=\max_{\alpha,\varepsilon}\;&\varepsilon\\
 \textnormal{s.t.}\;&G_{\rm ray}\alpha=0,\quad
              \mathbf1^\top\alpha=1,\quad
              \alpha\ge\varepsilon\mathbf1.
 \end{aligned}
\end{equation}
This is a normalized version of the strictly-positive null-combination
test in \cite[Equation~(12.14) and Section~12.2.3]{lynch2017modern}.
Problem \eqref{eq:arm-force-closure} is a linear program, and its optimal
value $\varepsilon_*$ is the largest margin by which a normalized null
combination can be kept positive. Together with the rank condition,
$\varepsilon_*>0$ certifies positive spanning of wrench space by the finite rays,
and hence force closure for the original larger circular cones.
Strict positivity of the null combination is essential to this
frictional force-closure certificate.

\paragraph{Rank as a Riemannian constraint.}
The matrices of full row rank form an open set, so ``rank six'' cannot be
imposed as an equation on the entries of the wrench matrix. The thin QR
factorization gives smooth coordinates on that set, so the rank condition
becomes membership of the two factors in two manifolds, together with an
ordinary matching equation.

Let $m_r=8L$. The wrench scaling is
$S_w=\operatorname{diag}(\ell_b^{-1}I_3,I_3)$, with characteristic
length $\ell_b>0$.  Write $W(\vartheta)=S_wG_{\rm ray}(\vartheta)$,
where $\vartheta$ denotes the smooth contact-geometry data; during
rigid transport it is fixed at the acquired geometry.  The rank
condition is enforced through the matrix manifold
\begin{equation}\label{eq:arm-rank-manifold}
 \begin{aligned}
 \mathcal R_6&=\{W\in\mathbb R^{6\times m_r}:\operatorname{rank}W=6\},\\
 W(\vartheta)^\top&=U_gT_g,\qquad
 U_g\in\operatorname{St}(m_r,6),\quad
 T_g\in\mathcal T_{++}(6).
 \end{aligned}
\end{equation}
Here $\operatorname{St}(m_r,6)=\{U:U^\top U=I_6\}$, and
$\mathcal T_{++}(6)$ is the manifold of upper triangular matrices
with positive diagonal.  The positive-diagonal thin QR factors are
unique and give a smooth factor representation of the full-row-rank stratum.
They serve as auxiliary geometric coordinates and become dependent
functions of the physical contact geometry in the reduced graph. The physical state and control blocks remain those of the arm model.
The dimensions of the two factor manifolds add to $6m_r$, the dimension of
the open set of full-row-rank matrices.  Thus rank is represented by a
Riemannian constraint, in the same way that the robot attitudes are
represented on $SO(3)$, while
\begin{equation}\label{eq:arm-rank-matching}
 K_g(\vartheta,U_g,T_g):=W(\vartheta)^\top-U_gT_g=0
\end{equation}
ties the factors to the physical contact geometry.  Intrinsic matrix
manifold constructions follow \cite{absil2008optimization,boumal2023introduction}.

For the symmetric finite pads in \eqref{eq:arm-contact-geometry},
the uniform coefficients $\alpha_*=(8L)^{-1}\mathbf1$ balance the
contact rays. Proposition~\ref{prop:arm-analytic-force-closure} gives the
full rank calculation and its explicit perturbation margin. That
calculation clarifies why the spatial extent of the pads supplies the
moment directions needed for three-dimensional force closure.

\paragraph{Normalized closure geometry.}
In the rank representation, the normalized closure equality is
\begin{equation}\label{eq:arm-geometric-closure}
 H_g(U_g,\alpha)=
 \begin{pmatrix}U_g^\top\alpha\\\mathbf1^\top\alpha-1\end{pmatrix}=0,
\end{equation}
with $\alpha>0$.
Since $T_g$ is invertible, this is equivalent to the null combination
in \eqref{eq:arm-force-closure}.

Let $B_g=[U_g^\top;\mathbf1^\top]$, a $7\times m_r$ matrix. The certificate
$\alpha$ is not unique: \eqref{eq:arm-geometric-closure} is a system of
seven equations in $m_r$ unknowns. Free directions of this kind would be
directions of zero curvature in the NLP, so they are removed. The
optimization eliminates these
artificial certificate fibers by selecting seven columns $J_g$ with
invertible $B_{g,J_g}$; $J_g^c$ denotes the complementary columns.
Their reference coefficients $\bar\alpha_{J_g^c}$ are fixed, giving
\begin{equation}\label{eq:arm-closure-section}
 \alpha_{J_g}=B_{g,J_g}^{-1}
 \left[\begin{pmatrix}0_6\\1\end{pmatrix}
       -B_{g,J_g^c}\bar\alpha_{J_g^c}\right].
\end{equation}
This is a smooth positive section near the certified point, that is, a
smooth choice of one positive certificate for each nearby contact geometry.
It is evaluated by a small solve and removes $m_r-7$ purely auxiliary
directions.
In the nominal rigid mode $\alpha_*$ and the factors are fixed
geometric data and are substituted into the NLP as constants.

The tangent geometry is
$U_g^\top\dot U_g+\dot U_g^\top U_g=0$ with upper triangular
$\dot T_g$.  For example, a Stiefel polar retraction is
$R_{U_g}(\dot U_g)=(U_g+\dot U_g)
 (I_6+\dot U_g^\top\dot U_g)^{-1/2}$; the positive diagonal
of $T_g$ has local logarithmic coordinates.  

We use the graph-reduced representation: the physical contact geometry is
independent, while its QR factors and normalized certificate are reconstructed.
In other words, the auxiliary variables $(U_g,T_g,\alpha)$ are treated as
functions of the physical contact geometry $\vartheta$ and substituted,
instead of being kept as independent variables with matching constraints.
From here on the selected column set $J_g$ is written $J$.
The following choices specify the metric, connection and second-order model.
Let $\mathcal P$ be the product manifold of independent physical geometry
variables. Locally, its factors are Euclidean coordinates $a\in\mathbb R^{d_a}$,
joint angles $\theta\in T^{d_\theta}$, rotations $R_\ell\in SO(3)$ and,
when independently variable, unit contact directions $n_j\in S^2$.
Factors absent from a particular phase are omitted. In particular, an
independent payload pose is represented by $(R_b,p_b)\in SO(3)\times\mathbb R^3$;
rigidly attached poses are reconstructed from the retained joint variables.
The operational contact equalities remain OCP constraints on this product.

Write a physical tangent vector as
$\xi=(\xi_a,\xi_\theta,\{R_\ell\Omega_\ell\},\{\xi_{n_j}\})$,
where $\Omega_\ell^\top=-\Omega_\ell$ and $n_j^\top\xi_{n_j}=0$.
Choose constant matrices $M_a\succ0$, $M_\theta\succ0$ and constants
$c_\ell>0$, $d_j>0$, fixed throughout the local homotopy. The physical
metric is
\begin{equation}\label{eq:grasp-physical-metric}
\begin{split}
 g^{\mathcal P}_\vartheta(\xi,\eta)
 ={}&\xi_a^\top M_a\eta_a
       +\xi_\theta^\top M_\theta\eta_\theta\\
 &+\sum_\ell\frac{c_\ell}{2}
       \operatorname{tr}(\Omega_\ell^\top\Lambda_\ell)
       +\sum_jd_j\xi_{n_j}^\top\eta_{n_j},
\end{split}
\end{equation}
where the rotational block of $\eta$ is $R_\ell\Lambda_\ell$.
The weights specify the scaling of quantities with different physical
units. The rotational weights are scalar, so these rotational metrics
are bi-invariant.

For a smooth tangent field $Y$, the Levi--Civita connection is the
product connection
\begin{equation}\label{eq:grasp-physical-connection}
\begin{aligned}
 (\nabla^{\mathcal P}_\xi Y)_a&=DY_a[\xi],&
 (\nabla^{\mathcal P}_\xi Y)_\theta&=DY_\theta[\xi],\\
 (\nabla^{\mathcal P}_\xi Y)_{R_\ell}
 &=R_\ell\!\left(D\Lambda_\ell[\xi]
                 +\tfrac12[\Omega_\ell,\Lambda_\ell]\right),\\
 (\nabla^{\mathcal P}_\xi Y)_{n_j}
 &=(I_3-n_jn_j^\top)DY_{n_j}[\xi],
 \qquad [A,B]:=AB-BA.
\end{aligned}
\end{equation}
Here $Y_{R_\ell}=R_\ell\Lambda_\ell$, so $D\Lambda_\ell[\xi]$
differentiates the body-coordinate field, including its dependence on all
physical variables. The angle formula is evaluated in a local unwrapped
chart. Equivalently, the rotational connection is
$R_\ell\operatorname{skew}(R_\ell^\top DY_{R_\ell}[\xi])$.

A concrete physical retraction is the product exponential
\begin{equation}\label{eq:grasp-physical-retraction}
\begin{aligned}
 R^{\mathcal P}_\vartheta(\xi)
  &=\left(a+\xi_a,\ \theta+\xi_\theta\pmod{2\pi},\
        \{R_\ell e^{\Omega_\ell}\},\
        \{\mathcal E_{n_j}(\xi_{n_j})\}\right),\\
 \mathcal E_n(v)&=\cos(\|v\|_2)n+
             \frac{\sin(\|v\|_2)}{\|v\|_2}v,
 \qquad \frac{\sin 0}{0}:=1.
\end{aligned}
\end{equation}
Thus $\nabla^{\mathcal P}_{\dot c}\dot c(0)=0$ for
$c(t)=R^{\mathcal P}_\vartheta(t\xi)$.
For a pose this gives $(R_b e^{\Omega_b},p_b+\xi_{p_b})$.
It is the exponential of the declared product metric; an $SE(3)$ group
exponential is a different map. Three-dimensional rotation exponentials
are evaluated by Rodrigues' formula.

Let $\mathcal U_J\subset\mathcal P$ be the open domain where
$W(\vartheta)$ has full row rank, $B_{g,J}$ is invertible, and the
certificate in \eqref{eq:arm-closure-section} is strictly positive.
On a smooth contact branch, assume $W$ is $C^3$.
For fixed complementary coefficients $\bar\alpha_{J^c}$, define
\begin{equation}\label{eq:grasp-graph-map}
 \Psi_J(\vartheta)
  =\bigl(\vartheta,U_g(\vartheta),T_g(\vartheta),
                      \alpha^{(J)}(\vartheta)\bigr),\qquad
 \mathcal G_J=\Psi_J(\mathcal U_J).
\end{equation}
The positive-diagonal thin QR factorization and the section solve make
$\Psi_J$ smooth and single valued. Its first component is the identity,
so it is an embedding even when the physical wrench map has redundant
parameter directions. With $\pi(\vartheta,U,T,\alpha)=\vartheta$,
the graph metric is defined explicitly by
\begin{equation}\label{eq:grasp-graph-metric}
 g^{\mathcal G_J}_z(\widehat\xi,\widehat\eta)
   =g^{\mathcal P}_{\pi(z)}(D\pi\widehat\xi,D\pi\widehat\eta),
 \quad
 z=\Psi_J(\vartheta),\quad
 \widehat\xi=D\Psi_J\xi,\quad\widehat\eta=D\Psi_J\eta.
\end{equation}
Consequently $\Psi_J$ is an isometry. The auxiliary factors and
certificate inherit their metric through their dependence on the physical
variables: lengths on the graph are measured through the physical variables
alone. Adding separate squared norms of the changes in $U_g$, $T_g$ or
$\alpha$ would define a different metric, and therefore a different
gradient and Hessian.

For physical vector fields $X,Y$, the graph connection, gradient and
Hessian are
\begin{equation}\label{eq:grasp-graph-connection}
\begin{aligned}
 \nabla^{\mathcal G_J}_{(\Psi_J)_*X}(\Psi_J)_*Y
   &=(\Psi_J)_*(\nabla^{\mathcal P}_X Y),\\
 \operatorname{grad}_{\mathcal G_J}\ell
   &=(\Psi_J)_*\operatorname{grad}_{\mathcal P}(\ell\circ\Psi_J),\\
 \operatorname{Hess}_{\mathcal G_J}\ell[D\Psi_J\xi]
   &=D\Psi_J\,
       \operatorname{Hess}_{\mathcal P}(\ell\circ\Psi_J)[\xi].
\end{aligned}
\end{equation}
Here $\ell$ may be the scalar Lagrangian with all multipliers held fixed.
These formulas specify the intrinsic geometry associated with the
physical metric on the graph.

For an explicit matrix formula, let $x=(x^1,\ldots,x^d)$ be physical
local coordinates, let $G_{ij}=g^{\mathcal P}(\partial_i,\partial_j)$,
and write the ambient components of the reconstruction as $z(x)$.
Its Christoffel symbols and the Hessian bilinear-form matrix are
\begin{equation}\label{eq:grasp-graph-hessian}
\begin{aligned}
 \Gamma^k_{ij}
   &=\frac12\sum_lG^{kl}
       (\partial_iG_{jl}+\partial_jG_{il}-\partial_lG_{ij}),\\
 Q_{ij}
   &=D^2\widetilde\ell(z)[z_i,z_j]
       +D\widetilde\ell(z)
         \left[z_{ij}-\sum_k\Gamma^k_{ij}z_k\right],\\
 z_i&=\partial_i z,\qquad z_{ij}=\partial_i\partial_jz.
\end{aligned}
\end{equation}
Here $G^{kl}$ are the entries of $G^{-1}$ and
$\widetilde\ell$ is any smooth ambient extension of $\ell$ near the
graph point. The resulting $Q$ is independent of that extension.
The coordinate gradient and Hessian operator are respectively
$G^{-1}D(\ell\circ\Psi_J)$ and $G^{-1}Q$.
Equation~\eqref{eq:grasp-graph-hessian} includes both the reconstruction
second derivatives and the physical connection correction. If the
Lagrangian is evaluated directly in physical variables, the same matrix
is $Q_{ij}=\partial_{ij}(\ell\circ\Psi_J)
-\sum_k\Gamma^k_{ij}\partial_k(\ell\circ\Psi_J)$.

\paragraph{Differentiating the reconstruction.}
The needed derivatives follow from the factorization and section
equations by triangular and $7\times7$ solves. Along a smooth physical
curve set $A=W^\top=U_gT_g$ and use dots for ordinary curve derivatives.
Define $\operatorname{low}(D)$ as the strictly lower-triangular part.
Then
\begin{equation}\label{eq:grasp-qr-first}
\begin{aligned}
 Z&=\dot A T_g^{-1},\qquad D=U_g^\top Z,\qquad
 \Omega=\operatorname{low}(D)-\operatorname{low}(D)^\top,\\
 \dot U_g&=(I-U_gU_g^\top)Z+U_g\Omega,\qquad
 \dot T_g=(D-\Omega)T_g.
\end{aligned}
\end{equation}
Right multiplication by $T_g^{-1}$ denotes a triangular solve.
These equations imply
$U_g^\top\dot U_g+\dot U_g^\top U_g=0$ and upper-triangular
$\dot T_g$. For the second derivative, put
\begin{equation}\label{eq:grasp-qr-second}
\begin{aligned}
 E&=(\ddot A-2\dot U_g\dot T_g)T_g^{-1},\qquad
 K=U_g^\top E,\qquad C_U=\dot U_g^\top\dot U_g,\\
 L_{ii}&=K_{ii}+(C_U)_{ii},\\
 L_{ij}&=K_{ij}+K_{ji}+2(C_U)_{ij}\quad(i<j),
 \qquad L_{ij}=0\quad(i>j),\\
 \ddot U_g&=E-U_gL,\qquad \ddot T_g=LT_g.
\end{aligned}
\end{equation}
Indeed, differentiating $A=U_gT_g$ twice gives the last line, while
differentiating $U_g^\top U_g=I$ twice gives
$U_g^\top\ddot U_g+\ddot U_g^\top U_g=-2C_U$ and determines $L$.
The derivatives of $A=W(\vartheta)^\top$ include all physical
first and second derivatives of the contact geometry.
Mixed coordinate derivatives follow by polarization of these
directional second derivatives in a physical coordinate chart.

For $B_g=[U_g^\top;\mathbf1^\top]$ and
$b_0=\operatorname{col}(0_6,1)$, differentiate
$B_g\alpha^{(J)}=b_0$ with the complementary coefficients fixed:
\begin{equation}\label{eq:grasp-certificate-derivatives}
\begin{aligned}
 B_{g,J}\dot\alpha_J&=-\dot B_g\alpha,
       &\dot\alpha_{J^c}&=0,\\
 B_{g,J}\ddot\alpha_J&=-\ddot B_g\alpha-2\dot B_g\dot\alpha,
       &\ddot\alpha_{J^c}&=0.
\end{aligned}
\end{equation}
The lower row of $\dot B_g$ and $\ddot B_g$ is zero. Thus all
reconstruction terms in \eqref{eq:grasp-graph-hessian} are defined
by explicit linear equations.

\paragraph{Graph retraction and changes of section.}
The update used for the reduced Newton step is
\begin{equation}\label{eq:grasp-graph-retraction}
 R^{\mathcal G_J}_{\Psi_J(\vartheta)}(D\Psi_J\xi)
       :=\Psi_J\!\left(R^{\mathcal P}_\vartheta(\xi)\right),
 \qquad R^{\mathcal P}_\vartheta(\xi)\in\mathcal U_J.
\end{equation}
It preserves physical matching and the selected positive certificate
exactly in exact arithmetic. By the connection identity in
\eqref{eq:grasp-graph-connection}, its curves have zero covariant
acceleration at zero. Therefore it is a second-order retraction and
\begin{equation}\label{eq:grasp-second-order-pullback}
 D^2(\ell\circ R^{\mathcal G_J}_z)(0)
           [\widehat\xi,\widehat\eta]
  =g^{\mathcal G_J}_z
       (\widehat\xi,\operatorname{Hess}_{\mathcal G_J}\ell[
                                                     \widehat\eta]).
\end{equation}
This is the Hessian bilinear form used in the Newton model
\cite[Section~5.10]{boumal2023introduction}.
The reconstructed ambient acceleration in \eqref{eq:grasp-qr-second}
can be nonzero despite the vanishing intrinsic acceleration.

Explicit validity margins for a section are
\begin{equation}\label{eq:grasp-section-margins}
 \rho_W=\|T_g^{-1}\|_F^{-1},\qquad
 \rho_J=\|B_{g,J}^{-1}\|_F^{-1},\qquad
 \rho_\alpha=\min_i\alpha_i^{(J)}.
\end{equation}
The selected neighborhood has these three margins bounded below by
declared positive constants. A trial point must retain the corresponding
positive margins; inverse norms can be evaluated by solves.
On an overlap $\mathcal U_J\cap\mathcal U_K$, the representation
transition is
\begin{equation}\label{eq:grasp-section-transition}
 \mathcal T_{KJ}=\Psi_K\circ\pi:\mathcal G_J\longrightarrow\mathcal G_K,
 \qquad
 D\mathcal T_{KJ}(D\Psi_J\xi)=D\Psi_K\xi.
\end{equation}
The domain of this map is the graph over the physical overlap.
It is an isometry for \eqref{eq:grasp-graph-metric}.
The new complementary coefficients can be chosen as the current
$\alpha_{K^c}$, preserving the current certificate at the change.
Different sections represent the same physical geometry; their
certificate derivatives can differ. The physical objective and
constraints are independent of this artificial certificate choice,
and their transformed derivatives obey
\eqref{eq:grasp-graph-connection}--\eqref{eq:grasp-graph-hessian}.
In a rigid nominal mode these geometric quantities are fixed data;
their reconstruction adds no decision variables.

\paragraph{Rigid transport and invariance of closure.}
Let $X_b=(R_b,p_b;0,1)$ map object-frame coordinates into world
coordinates. With the moment-first wrench ordering of
\eqref{eq:arm-grasp-map} and world moments taken about the world origin,
the wrench transformation is
\begin{equation}\label{eq:grasp-wrench-transform}
 \begin{pmatrix}\tau^0\\f^0\end{pmatrix}
 =\mathcal A(X_b)\begin{pmatrix}\tau^b\\f^b\end{pmatrix},
 \qquad
 \mathcal A(X_b)=
 \begin{pmatrix}R_b&\widehat p_bR_b\\0&R_b\end{pmatrix}
 =\operatorname{Ad}_{X_b}^{-\top}.
\end{equation}
For a rigidly maintained grasp, the ray-wrench matrix $G_{\rm ray}^b$
is constant in the object frame. Consequently
\begin{equation}\label{eq:grasp-rigid-closure-invariance}
\begin{aligned}
 G_{\rm ray}^0&=\mathcal A(X_b)G_{\rm ray}^b,\qquad
 W^0=S_w\mathcal A(X_b)S_w^{-1}W^b,\\
 \operatorname{rank}W^0&=\operatorname{rank}W^b=6,\qquad
 W^0\alpha=0\quad\text{whenever }W^b\alpha=0.
\end{aligned}
\end{equation}
Since $\det\mathcal A(X_b)=1$, the transformation is invertible and
preserves the rank. The null space of $W^b$ is contained in that of $W^0$,
so the same normalized positive
certificate, with $\mathbf1^\top\alpha=1$ and $\alpha>0$, remains valid
throughout rigid transport and the held part of placement.
If world-coordinate moments are instead taken about the moving object
origin, $\mathcal A(X_b)$ is replaced by
$\operatorname{diag}(R_b,R_b)$, with the same rank and certificate
invariance. The following finite-force constraints specify the
loads that the grasp can actually support.

\subsubsection{Finite grip strength and coupled payload mechanics}
Finite grip strength restricts the realizable loads within the
geometric contact cones \cite[Section~12.2.3.1]{lynch2017modern}.
Let $\overline N_j$ and $\overline N_\sigma$ be the available local
and total jaw normal-force limits.  Loaded transport uses
preload reserves $\underline N_j>0$ and a conservative coefficient
$0<\mu_{\rm safe}<\mu_f$.  Let $\mathcal K_j^{\rm safe}$ be
\eqref{eq:arm-friction-cone} with $\mu_f$ replaced by $\mu_{\rm safe}$.
Define the realizable, reserved wrench set
\begin{equation}\label{eq:arm-finite-wrench-set}
 \mathbb W_{\rm grip}:=\left\{Gf:
 \begin{array}{l}
 f_j\in\mathcal K_j^{\rm safe},\quad\underline N_j\le n_j^\top f_j\le\overline N_j,\\
 \displaystyle\sum_{j:\,\sigma(j)=\sigma}n_j^\top f_j
                    \le\overline N_\sigma,\quad \sigma\in\{-1,1\},\\
 A_{\rm hand}f\le b_{\rm hand}
 \end{array}\right\}.
\end{equation}
The supplied linear limits $A_{\rm hand}f\le b_{\rm hand}$ describe
any additional pressure or actuator coupling of the hand.  For the ideal pad-reaction model these may be empty; a physical hand
uses the coupling limits of its actual actuator and pressure model.  At stationary acquisition/unloading the preload lower
limits are ramped to or from zero, and the corresponding feasible
contact/support forces are checked separately.

A sufficient-normal-strength certificate for a declared operating envelope is
\begin{equation}\label{eq:arm-load-certificate}
 \mathbb W_{\rm req}\subset\operatorname{int}\mathbb W_{\rm grip}.
\end{equation}
Here $\mathbb W_{\rm req}$ includes gravity, inertial loads and the specified
bounded disturbance wrenches. In the three MPC problems below, this inclusion
is verified on a finite prescribed polytope of load scenarios, or the
instantaneous wrench balance and finite-force constraints are retained in
the transport NLP. The force-closure condition and these finite-force constraints jointly
specify the required grasp: the former determines the wrench directions,
while the latter impose preload and load capacity.

During unsupported transport the box pose is
$X_b(q)=F_e(q)X_g^{-1}$.  Define the body Jacobian by
$[J_e^b(q)x^v]=F_e(q)^{-1}\dot F_e(q)$ and set
\begin{equation}\label{eq:arm-payload-twist}
 J_b(q)=\operatorname{Ad}_{X_g}J_e^b(q),\qquad
 \mathcal V_b=J_b(q)x^v,\qquad
 \mathcal G_b=\operatorname{diag}(I_b,m_bI_3).
\end{equation}
Here $[\mathcal V_b]=X_b^{-1}\dot X_b$, with angular components first,
and for $X=(R,p;0,1)$ the convention is
$\operatorname{Ad}_X=\left(\begin{smallmatrix}R&0\\\widehat p R&R\end{smallmatrix}\right)$.
The required contact wrench follows
\cite[Equations~(8.40) and (12.31)]{lynch2017modern}:
\begin{equation}\label{eq:arm-payload-wrench}
 \begin{aligned}
 \mathcal W_{\rm req}
 &=\mathcal G_b\dot{\mathcal V}_b
     -\operatorname{ad}_{\mathcal V_b}^{\top}\mathcal G_b\mathcal V_b
     -\mathcal W_g,\\
 \mathcal W_g&=\begin{pmatrix}0\\-m_bgR_b^\top e_3\end{pmatrix},\qquad
 \dot{\mathcal V}_b=J_b\dot x^v+\dot J_b x^v,\qquad
 Gf=\mathcal W_{\rm req}.
 \end{aligned}
\end{equation}
For $\mathcal V=\operatorname{col}(\mathcal V^\omega,\mathcal V^p)$,
$\operatorname{ad}_{\mathcal V}=
\left(\begin{smallmatrix}\widehat{\mathcal V^\omega}&0\\
\widehat{\mathcal V^p}&\widehat{\mathcal V^\omega}\end{smallmatrix}\right)$;
these twist components are not the control vector $v$.
All terms in \eqref{eq:arm-payload-wrench} use body coordinates about
the center of mass.  In particular gravity is included once, with
its world direction rotated into that frame.

Let $M_0(q,\gamma)$, $C_0(q,x^v;\gamma)x^v$ and $g_0(q,\gamma)$
be the calibrated unloaded-arm quantities, including the hand.
For fixed aperture, payload kinetic and potential energy add
\begin{equation}\label{eq:arm-payload-mass}
 M_b(q)=J_b(q)^\top\mathcal G_bJ_b(q),\qquad
 U_b(q)=m_bg e_3^\top p_b(q),\qquad
 g_b(q)=\nabla_\theta U_b(q).
\end{equation}
The velocity term $C_b(q,x^v)x^v$ follows from the Lagrange equations
for $M_b$ and is retained in the loaded model. In joint coordinates,
\[
 (C_b)_{ij}=\frac12\sum_{\ell=1}^6
 \left(\partial_\ell(M_b)_{ij}+\partial_j(M_b)_{i\ell}
                          -\partial_i(M_b)_{j\ell}\right)x^v_\ell.
\]
Hence unsupported transport obeys
\begin{equation}\label{eq:arm-loaded-dynamics}
 (M_0+M_b)\dot x^v+(C_0+C_b)x^v+g_0+g_b=v^q,
 \qquad \dot\theta=x^v.
\end{equation}
Equivalently, the payload contribution is $J_b^\top\mathcal W_{\rm req}$
in unloaded inverse dynamics.  These are equivalent ways of accounting
for the payload: the loaded matrices already contain its wrench.
Joint impacts and actuator electrical dynamics are outside the rigid-body,
ideal-servo prediction model. Aperture motion is quasistatic; the pad-force
controller performs preload and unloading at the measured interfaces below.
The moving jaw masses are neglected in this reference model; the unloaded
arm inertia includes the fixed hand casing and is independent of aperture.
Kinematics and collision geometry retain their full aperture dependence.
The arm prediction always retains its acceleration equation, including at a
nominal stopping or holding configuration. A small measured joint velocity is retained as the initial condition and
damped by the controller; zero-velocity equalities apply at the designated
future stopping nodes.

\subsubsection{Three measured manipulation problems}
\label{subsec:arm-three-mpc}
In what follows the index $\mathsf P\in\{\mathsf A,\mathsf B,\mathsf C\}$
denotes the problem; it selects the dynamics, features, costs and
constraint sets that carry it as a subscript or superscript.
The manipulation task is decomposed into three consecutive MPC problems:
$\mathsf A$ approaches and establishes a certified grasp, $\mathsf B$
transports the rigidly grasped box to a pre-placement pose above the
destination support, and $\mathsf C$ sets the box down, transfers its
load to that support, unloads, opens and retreats while penalizing
object disturbance. These are three top-level OCPs; placement and
release are prescribed, measured internal contact layouts of the third.
Measured contact and wrench tests determine each transition between
problems. A predicted terminal pose is therefore followed by the measured
acceptance tests before the supervisor advances.
For each problem $\mathsf P\in\{\mathsf A,\mathsf B,\mathsf C\}$, fix a current
terminal sample $N_{\mathsf P}$ and write
\begin{equation}\label{eq:arm-event-indices}
 H_{\hat t}^{\mathsf P}=N_{\mathsf P}-\hat t,\qquad
 \mathcal I=\{0,\ldots,H-1\},\quad
 \mathcal N=\{0,\ldots,H\},\quad
 \mathcal N^+=\{1,\ldots,H\},\quad H=H_{\hat t}^{\mathsf P}.
\end{equation}
Thus $H$ counts control intervals and the remaining duration is $hH$.
The control NLP is formed for $H\geq1$; at $H=0$ the supervisor checks its
measured success gate before terminating or extending that problem.
If the measured success test remains unsatisfied at $N_{\mathsf P}$, the
horizon is extended or the problem is replanned. Minimum
closing/opening durations obey $hH\geq
|\gamma_{\rm final}-\hat\gamma|/\overline v^\gamma$.
The nominal node zero equals the measured state and is substituted before
forming tangent systems. State bounds, geometric contact conditions and
operational inequalities are imposed on $\mathcal N^+$; actuator bounds apply
on $\mathcal I$. The supervisor performs sensor and physical consistency checks at the
measured node, which is fixed data for the optimizer.

Let $X_b^{(0)}$ be the measured initially supported box pose and
$X_b^{(d)}$ its desired upright placement pose.  Define a collision-free
pre-placement pose $X_b^{\rm pre}$ with the same intended upright attitude
and $p_b^{\rm pre}=p_b^{(d)}+h_{\rm pre}e_3$, $h_{\rm pre}>0$.
The supports have known
horizontal planes $e_3^\top p=h_s$ and convex contact polygons. The box's center
of mass must project strictly inside the intended support polygon, with a
specified reserve. The nominal box pose used in each prediction is
\begin{equation}\label{eq:arm-box-modes}
 X_{b,k}=\begin{cases}
 X_b^{(0)},&\mathsf P=\mathsf A,\\
 F_e(q_k)(X_g^{\mathsf B})^{-1},&\mathsf P=\mathsf B
                    \text{ or held placement in }\mathsf C,\\
 X_{b,k}\text{ is an independent predicted object state},&
                    \text{supported release/retreat in }\mathsf C.
 \end{cases}
\end{equation}
At the $\mathsf A$--$\mathsf B$ interface set
$X_g^{\mathsf B}=\hat X_b^{-1}F_e(\hat q)$ and
$\gamma_{\mathsf B}=\hat\gamma$ from the established measured grasp;
recompute its contact geometry and certificates before admission.
During transport the relative-transform error and relative velocity are
monitored against this calibration. A detected slip initiates recovery with the calibrated $X_g^{\mathsf B}$
retained until a new grasp is measured and admitted.
At the $\mathsf B$--$\mathsf C$ interface, the held placement
trajectory inherits the measured grasp and pre-placement state. The
independent object pose and twist are introduced at $\mathsf C$'s measured
supported-release interface, retaining the measured pose and velocity.

\paragraph{A common arm integrator and explicit features.}
A semi-implicit Euler step on each interval $k\in\mathcal I$ is
\begin{equation}\label{eq:arm-discrete-dynamics}
 \begin{aligned}
 a_k&=M_{\mathsf P}(q_k,\gamma_k)^{-1}
       \bigl(v^q_k-C_{\mathsf P}(q_k,x^v_k;\gamma_k)x^v_k
                       -g_{\mathsf P}(q_k,\gamma_k)+\tau^{\rm ext}_{\mathsf P,k}\bigr),\\
 x^v_{k+1}&=x^v_k+ha_k,\qquad q_{k+1}=R_{q_k}(h x^v_{k+1}),\\
 \gamma_{k+1}&=\gamma_k+h v^\gamma_k,\qquad
 R_q(\psi)=(q+\psi)\bmod2\pi.
 \end{aligned}
\end{equation}
Problem $\mathsf A$ and the detached part of $\mathsf C$ use the unloaded
$M_0,C_0,g_0$; Problem $\mathsf B$ and the held placement part of
$\mathsf C$ use $M_0+M_b,C_0+C_b,g_0+g_b$.
The aperture integration equation applies in $\mathsf A$ and after
unloading in $\mathsf C$.  In $\mathsf B$ and held placement the aperture
is a held force-servo parameter and is omitted
from the dynamic state and control blocks; its measured consistency is
checked by the supervisor.
The external torque is zero in free approach and unsupported transport.
During $\mathsf C$'s support transfer the loaded model instead has
$\tau^{\rm ext}_{\mathsf C,k}=J_b^\top\mathcal W_{S,k}$, so that
the external support wrench reduces the required hand load.  After
detachment it is the measured/predicted residual hand reaction specified
below. A different consistently differentiated integrator may replace
\eqref{eq:arm-discrete-dynamics}. The acceleration equation is kept even
when the arm is nominally at rest: a holding phase is not modeled by the
identity update $q_{k+1}=q_k$. The rate equation may equivalently be imposed
in inverse-dynamics form,
$M_{\mathsf P}(q_k,\gamma_k)(x^v_{k+1}-x^v_k)/h+C_{\mathsf P}x^v_k
+g_{\mathsf P}-\tau^{\rm ext}_{\mathsf P,k}-v^q_k=0$,
which defines the same discrete map and the same feasible set and avoids
differentiating $M_{\mathsf P}^{-1}$; the multipliers of these rows are then
expressed in the residual basis $M_{\mathsf P}(q_k,\gamma_k)/h$. The
implementation of Section~\ref{sec:results1019B} uses this form.
At every future node retain explicit kinematic features
\begin{equation}\label{eq:arm-explicit-z}
 z_k=h^z_{\mathsf P}(x_k,q_k,X_{b,k}),\qquad
 C^h_k=z_k-h^z_{\mathsf P}(x_k,q_k,X_{b,k})=0,\quad k\in\mathcal N^+.
\end{equation}
The stack contains gripper/link vertices, pad points and the eight box
vertices. Node-zero features are evaluated from the measurement. In
$\mathsf A$, $\mathsf B$, and held placement the box argument is the map
\eqref{eq:arm-box-modes}; in supported release/retreat it is the
independent object pose.
Rotation entries in $z$ are Euclidean features whose dependence on the
physical attitudes is imposed by $C^h$.

\subsubsection{Full-body collision constraints and reference costs}
Let $\mathcal K^{\rm rob}_\ell(q,\gamma)$ be the known convex components of the
links, palm and jaws, and let $\mathcal O_j$ be known static obstacles,
including noncontact portions of the supports.  Meshes may be split
into convex components.  A conservative sphere cover of a noncontact
link component is also admissible, with center
$p_{i,a}(q)$ and radius $a_{i,a}$; the cover must contain the actual
link volume.  The grasped box uses its full geometry in
\eqref{eq:arm-box-geometry}.

With body coordinates $y=R_b^\top(p-p_b)$ and
$\delta=|y|-a$, the exact point-to-box signed distance is
\begin{equation}\label{eq:arm-box-distance}
 d_{\mathcal B(X_b)}(p)
   =\|\max(\delta,0)\|_2+\min\{\max_i\delta_i,0\}.
\end{equation}
Thus a noncontact collision sphere satisfies, for example,
\[
 d_{\mathcal B(X_b(k))}(p_{i,a}(q_k))
       -a_{i,a}-d_{\rm box}\ge0.
\]
The corresponding sphere--obstacle and sphere--ground conditions are
\begin{equation}\label{eq:arm-sphere-clearance}
 d_{\mathcal O_j}(p_{i,a})-a_{i,a}-d_{\rm safe}\ge0,
 \qquad e_3^\top p_{i,a}-a_{i,a}-d_{\rm ground}\ge0,
\end{equation}
for pairs not designated as contact pairs.  For nonadjacent
self-collision spheres a smooth test is
\[
 \|p_{i,a}-p_{j,b}\|_2^2
       -(a_{i,a}+a_{j,b}+d_{\rm self})^2\ge0.
\]
The allowed-pair catalogue excludes physically connected and
intentionally contacting pairs; all other pairs involving the same links
remain in the collision checks.

An exact box signed-distance function is not globally twice
continuously differentiable at closest-feature changes.  For the
smooth NLP in this paper we therefore use a fixed, certified
separating-direction catalogue on each planned phase/corridor.
For convex components $\mathcal A_k,\mathcal B_k$, let
$n_{AB,k}$ be a prescribed unit vector pointing from $\mathcal A$
to $\mathcal B$.  For their vertices $a_i(z_k),b_j(z_k)$ enforce
\begin{equation}\label{eq:arm-smooth-separation}
 n_{AB,k}^\top\bigl(b_j(z_k)-a_i(z_k)\bigr)
       -d_{AB,k}\ge0\qquad\text{for every vertex pair }(i,j).
\end{equation}
For a sphere, replace its extremal projection by the center projection
plus/minus its radius.  These finitely many smooth inequalities imply
separation of the whole convex bodies; the chosen direction may be
conservative even though the box geometry is exact.  The normals are
geometric planning data.  Direction or closest-feature changes use a
catalogue update and a new local linearization.  This
is a specific smooth-optimization construction using the geometric
collision viewpoint of \cite[Section~10.2.2]{lynch2017modern}.

The pair catalogue includes every link and gripper component against
the environment, all nonadjacent self-collision pairs, all
\emph{noncontact} robot--box pairs, and the box against the environment
in each of the three problems, using its corresponding box-pose model.  Only the designated
pad--box pairs and box--support pairs are allowed zero clearance.
A nonpenetration inequality is retained only while it is independent of the
chosen mode equations. A gap fixed to zero by a rigid grasp or support mode
is removed from that mode's NLP and checked independently against measurements
and the physical geometry.  The pad components use their actual flat contact geometry, avoiding
spurious penetration from enclosing-sphere approximations. Arm--box and
payload--environment checks remain active throughout the grasp.

For a horizontal support of height $h_s$, let
$\xi_\ell\in\{(\pm a_1,\pm a_2,\pm a_3)\}$ be the eight box vertices.
A scheduled lift-clearance profile $\ell_k\ge0$ is imposed by
\begin{equation}\label{eq:arm-payload-clearance}
 e_3^\top\bigl(p_b(k)+R_b(k)\xi_\ell\bigr)-h_s
                   \ge\ell_k,\qquad\ell=1,\ldots,8.
\end{equation}
In Problem $\mathsf B$ use positive clearance at future lift/transport nodes and
at its pre-placement terminal target, with a plateau at $h_{\rm lift}$.
In $\mathsf C$'s held placement segment lower this profile to zero at the
supported placement target.
An endpoint gap implied by an imposed placement or support equality is
omitted from that layout. For different support
heights the corresponding local surface/corridor constraints are
used.  This enforces clearance at the \emph{lowest point of the full box}.
Payload separation from all other obstacles is imposed independently.

Let $c_{{\rm env},k}(z_k)\ge0$ collect this fixed phase-specific
catalogue. The following interval test complements the nodal checks with
continuous-time clearance verification.
\paragraph{A computable continuous-clearance certificate.}
For each noncontact pair, extend its fixed separating normal over a physical
control interval $I=[a,b]$ and let
$g_{ij}(t)=n^\top(b_j(t)-a_i(t))-d_{AB}$ be every required vertex-pair margin.
A validated continuous integrator (or an explicitly declared interpolation
model) encloses the positions and their velocities on $I$. Let
$L_{ij,I}\geq\sup_{t\in I}|\dot g_{ij}(t)|$ be a validated derivative
bound and $[g_{ij}(t_c)]$ an interval value at $t_c=(a+b)/2$. The test
\begin{equation}\label{eq:continuous-clearance-certificate}
 \underline{g_{ij}(t_c)}-L_{ij,I}(b-a)/2\geq0
 \qquad\text{for every required pair and vertex pair}
\end{equation}
proves separation throughout $I$.  A certified upper value below zero at any
point proves violation. An unresolved interval is bisected and tested
again, with the remaining unresolved intervals retained.  On a compact interval with a strictly positive true
clearance margin and convergent derivative enclosures, subdivision eventually
certifies the interval.  At grazing contact, designated contacts are checked by their explicit
nonpenetration and sticking identities. An unclassified noncontact interval
retains status \texttt{UNRESOLVED}, and control application is withheld;
finite termination of subdivision requires the strictly positive margin
stated above.
For a known plant-state error tube of radius $e_I$, subtract the validated
margin Lipschitz bound $L_{r,ij,I}e_I$ in
\eqref{eq:continuous-clearance-certificate}.  The error tube extends the nominal-model certificate to the declared
plant uncertainty. In its absence, the certificate applies to the nominal
continuous model.
The continuous-clearance procedure therefore has an explicit sufficient
test, refinement rule, and unresolved status.

\paragraph{Smooth references and a common tracking cost.}
References are supplied separately for each problem. Problem $\mathsf A$
ends at the calibrated pad-contact geometry; $\mathsf B$ follows a lift
and transfer corridor ending above the destination support; $\mathsf C$
lowers to placement, transfers support, opens at the measured supported
state and then retreats to a collision-free home.
For moving unwrapped-joint waypoints, a smooth reference interpolant is
\begin{equation}\label{eq:arm-reference-interpolation}
 \theta^{\rm ref}(t_a+\xi\Delta t)
 =\theta_a+(10\xi^3-15\xi^4+6\xi^5)(\theta_b-\theta_a),
 \qquad 0\leq\xi\leq1.
\end{equation}
The discretization-consistent feedforward quantities are
\begin{equation}\label{eq:arm-reference-feedforward}
 \begin{aligned}
 x^{v,\rm ref}_{k+1}&=(\theta^{\rm ref}_{k+1}-\theta^{\rm ref}_k)/h,\\
 v^{q,\rm ref}_k&=M_{\mathsf P}
       (x^{v,\rm ref}_{k+1}-x^{v,\rm ref}_k)/h
           +C_{\mathsf P} x^{v,\rm ref}_k+g_{\mathsf P}-\tau^{\rm ext,ref}_{\mathsf P,k},\\
 v^{\gamma,\rm ref}_k&=(\gamma^{\rm ref}_{k+1}-\gamma^{\rm ref}_k)/h.
 \end{aligned}
\end{equation}
All quantities in the second equation are evaluated at the reference node.
When a future reference node is required to have zero velocity, the
reference repeats the same position at the preceding node (an
equal-position predecessor), so that the finite-difference velocity in
\eqref{eq:arm-reference-feedforward} vanishes there. These nominal conditions apply to the reference trajectory; the initial
velocity remains the measured value. Admissible references satisfy
the actual torque, rate and full-body interval-clearance limits.
With positive control weights and positive semidefinite tracking weights,
define
\begin{equation}\label{eq:arm-tracking-stage-cost}
 \begin{aligned}
 J_k^{\mathsf P}={}&\|p_e(z_k)-p_{e,k}^{\rm ref}\|_{Q_p}^2
 +\|\operatorname{vec}(R_e(z_k)-R_{e,k}^{\rm ref})\|_{Q_R}^2
 +\|x_k^v-x_k^{v,\rm ref}\|_{Q_v}^2\\
 &+\rho_{\gamma,x}|\gamma_k-\gamma_k^{\rm ref}|^2
 +\|v_k^q-v_k^{q,\rm ref}\|_{R_q}^2
 +\rho_{\gamma,v}|v_k^\gamma-v_k^{\gamma,\rm ref}|^2.
 \end{aligned}
\end{equation}
Here $\|a\|_Q^2=a^\top Qa$. The chordal attitude cost is defined even for
infeasible feature trials $z\ne h^z$. At a feasible rotation,
$\frac12\|R-R^{\rm ref}\|_F^2=2(1-\cos\theta)$ for relative angle $\theta$.

\subsubsection{\texorpdfstring{Problem $\mathsf A$: approach and force-capable acquisition}{Problem A: approach and force-capable acquisition}}
\label{subsec:arm-mpc-a}
Let $\mathcal F_\delta$ be the force allocations in
\eqref{eq:arm-finite-wrench-set}, with a prescribed strictly positive
reserve $\delta$ in each finite-force and hand-coupling limit.
Circular friction constraints are implemented smoothly, for positive
$N_j=n_j^\top f_j$, as
$\mu_{\rm safe}^2N_j^2-\|f_j-N_jn_j\|^2\geq\delta_f^2$.
All reserves have declared force scales. A finite family of load scenarios
$w^1,\ldots,w^{N_w}$ is chosen with
$\mathbb W_{\rm req}\subseteq\operatorname{conv}\{w^1,\ldots,w^{N_w}\}$;
the scenario polytope includes the chosen disturbance reserve.

For a candidate terminal hand, let
$\rho_j^\mathsf A=R_b^{(0)\top}(c_j^0(q_H,\gamma_H)-p_b^{(0)})$ and construct
$G_\mathsf A$ and $G_{{\rm ray},\mathsf A}$ from
\eqref{eq:arm-grasp-map} and the fixed target cone rays using these moment
arms. They describe the candidate contact geometry, with terminal force
allocations used to evaluate its load capacity. These terminal forces are
variables of the optimization only: they test whether the candidate
geometry could carry each load scenario, and no force is commanded during
the approach. Physical forces are applied after contact.
A strictly positive normalized null witness $\alpha^{\rm fc}$ is fixed
from \eqref{eq:arm-force-closure} for the target geometry. The wrench
scaling is $S_w=\operatorname{diag}(\ell_b^{-1}I_3,I_3)$.
A precise terminal acquisition cost is
\begin{equation}\label{eq:arm-acquisition-terminal-cost}
 \begin{aligned}
 \Phi_\mathsf A={}&\sum_{j=1}^{N_c}
       \|c_j^0(q_H,\gamma_H)-(p_b^{(0)}+R_b^{(0)}\rho_j)\|_{Q_c}^2
       +\|x_H^v\|_{Q_a}^2+\rho_a|\gamma_H-\gamma_g|^2\\
 &+\rho_{\rm fc}\|S_wG_{{\rm ray},\mathsf A}\alpha^{\rm fc}\|^2
       +\rho_{\rm cap}\sum_{\ell=1}^{N_w}
                \|S_w(G_\mathsf A f^\ell-w^\ell)\|^2,
 \qquad f^\ell\in\mathcal F_\delta.
 \end{aligned}
\end{equation}
The resulting approach OCP is
\begin{equation}\label{eq:arm-ocp-a}
\begin{aligned}
 \min_{\substack{(q_k,x_k^v,\gamma_k,z_k)_{k\in\mathcal N^+}\\
                 (v_k^q,v_k^\gamma)_{k\in\mathcal I},\ (f^\ell)_{\ell=1}^{N_w}}}
 &h\sum_{k\in\mathcal I}J_k^\mathsf A+\Phi_\mathsf A,\\
 \text{subject to}\quad
 &(q_0,x_0^v,\gamma_0)=(\hat q,\hat x^v,\hat\gamma)
                       \quad\text{as substituted data},\\
 &\eqref{eq:arm-discrete-dynamics}\text{ with }(M_{\mathsf P},C_{\mathsf P},g_{\mathsf P})
       =(M_0,C_0,g_0), &&k\in\mathcal I,\\
 &C_k^h=0,\quad c_{{\rm env},k}(z_k)\geq0,
       \quad\text{state bounds in \eqref{eq:arm-simple-bounds}},
                                      &&k\in\mathcal N^+,\\
 &(v_k^q,v_k^\gamma)\in\mathbb V_{\mathsf A},
                                      &&k\in\mathcal I,\\
 &f^\ell\in\mathcal F_\delta,&&\ell=1,\ldots,N_w.
\end{aligned}
\end{equation}
Here $\mathbb V_{\mathsf A}$ consists of the torque and aperture-rate
bounds in \eqref{eq:arm-simple-bounds}, and $C_k^h$ is
\eqref{eq:arm-explicit-z}. The candidate rank factors are the dependent
geometric coordinates of \eqref{eq:arm-rank-manifold}--\eqref{eq:arm-rank-matching};
if retained explicitly, their variables and matching equations are included in
\eqref{eq:arm-ocp-a}. The terminal variables $f^\ell$ enter the
capacity term of $\Phi_\mathsf A$.
The approach retains pad--box nonpenetration inequalities and uses a
terminal acquisition cost before measured sticking is established. The null-witness term promotes
the certified target geometry, and the scenario terms promote finite-force
capacity. The measured test below establishes actual force closure through explicit
rank, positive-combination and finite-capacity checks.  The
terminal candidate matrix is represented on
\eqref{eq:arm-rank-manifold} with physical matching
\eqref{eq:arm-rank-matching}, or by its unique dependent factors.
The positive null equality \eqref{eq:arm-geometric-closure} becomes
active at the handoff, once the contact geometry has been established;
free approach uses the acquisition cost.

After pad contact, the low-level force servo establishes preload while the
arm damps its measured rate. At the initially supported box, a continuous
load-transfer parameter $\chi\in[0,1]$ gives the quasistatic target
\begin{equation}\label{eq:arm-support-balance}
 \mathcal W_{\rm pads}=-\chi\mathcal W_g,\qquad
 \mathcal W_{\rm support}=-(1-\chi)\mathcal W_g,\qquad
 v^q=g_0+\chi g_b.
\end{equation}
This interface servo uses measured wrench feedback while node zero
remains substituted measurement data. For a nonzero measured motion,
include the inertial wrench from \eqref{eq:arm-payload-wrench} in the
force-balance check. The supervisor advances to $\mathsf B$ only after
all of the following hold for a prescribed dwell $T_g>0$: both opposing pads are detected within
pose tolerance $\varepsilon_g$; relative pad/object speed is at most
$\varepsilon_v$; the measured geometry has rank six and
$\varepsilon_*\geq\varepsilon_{\rm fc}$ with
$0<\varepsilon_{\rm fc}<(8L)^{-1}$ in
\eqref{eq:arm-force-closure}; every load scenario admits
$Gf^\ell=w^\ell$ with $f^\ell\in\mathcal F_\delta$; measured preload is
within its limits; and full hand load acceptance has been verified at
$\chi=1$ without object slip. Failure retains or replans $\mathsf A$, or invokes the offline backup,
with the measured pose and velocity preserved.

\subsubsection{\texorpdfstring{Problem $\mathsf B$: hard force closure, lift and transport}{Problem B: hard force closure, lift and transport}}
\label{subsec:arm-mpc-b}
Once the measured acquisition test passes, the rigid sticking condition is
\begin{equation}\label{eq:arm-event-compatibility}
 F_e(q_k)=X_{b,k}X_g^{\mathsf B},\qquad
 \gamma_k=\gamma_{\mathsf B},\qquad k\in\mathcal N^+.
\end{equation}
When both poses are retained, compatibility is imposed by six local pose
equations; the intrinsic alternative substitutes $X_{b,k}$ through
\eqref{eq:arm-box-modes}. The aperture and its zero rate are
eliminated as known held quantities. Thus the dynamic arm state in
$\mathsf B$ is $(q,x^v)$, while aperture measurements are checked as
force-servo data. A discrepancy outside its certified calibration envelope
withholds transport and reinitializes the contact model. At the initial node this relation
is verified and calibrated from the measured established grasp; no
additional contact equation is imposed on the measured parameters.
Future predictions retain the exact rigid relation.
Hard force closure consists of the rigid relation together with
\eqref{eq:arm-rank-manifold}--\eqref{eq:arm-geometric-closure} on
the acquired contact geometry.  The acquired rank and positive certificate are maintained by the fixed
object-frame grasp geometry. Its equivalent intrinsic reduction substitutes
the fixed object-frame factors and certificate and retains the
independent rigid-grasp equations.  The resulting equality system contains the independent rigid-grasp rows.
In other words, once the grasp is established the contact geometry no
longer changes in the frame of the box. The rank factors and the
certificate computed at acquisition therefore remain valid and enter
Problem $\mathsf B$ as constants, and only the rigid attachment of the box
to the hand is imposed by equations.

For every transport interval retain contact-force allocations
$f_k\in\mathcal F_\delta$ and impose
\begin{equation}\label{eq:arm-transport-force-balance}
 Gf_k=\mathcal W_{{\rm req},k},\qquad
 \mathcal W_{{\rm req},k}=\mathcal G_b(J_ba_k+\dot J_bx_k^v)
 -\operatorname{ad}_{J_bx_k^v}^\top\mathcal G_bJ_bx_k^v-\mathcal W_g.
\end{equation}
The acceleration $a_k$ is that of the loaded model
\eqref{eq:arm-discrete-dynamics}. Thus pad forces certify the required
payload wrench without adding the payload inertia twice to the arm dynamics.
If the envelope inclusion \eqref{eq:arm-load-certificate} has been certified
for the full allowed state/control domain, these allocations may instead
be eliminated using that certificate.
The transport tracking cost has terminal term
\begin{equation}\label{eq:arm-transport-terminal-cost}
 \Phi_\mathsf B=
 \|\operatorname{vee}\log((X_b^{\rm pre})^{-1}X_{b,H})\|_{Q_b}^2
       +\|J_b(q_H)x_H^v\|_{Q_{b,v}}^2.
\end{equation}
The loaded transport OCP is therefore
\begin{equation}\label{eq:arm-ocp-b}
\begin{aligned}
 \min_{\substack{(q_k,x_k^v,z_k)_{k\in\mathcal N^+}\\
                 (v_k^q,f_k)_{k\in\mathcal I}}}
 &h\sum_{k\in\mathcal I}J_k^\mathsf B+\Phi_\mathsf B,\\
 \text{subject to}\quad
 &(q_0,x_0^v)=(\hat q,\hat x^v)\quad\text{as substituted data},\\
 &\eqref{eq:arm-discrete-dynamics}\text{ with }
 (M_{\mathsf P},C_{\mathsf P},g_{\mathsf P})
 =(M_0+M_b,C_0+C_b,g_0+g_b),&&k\in\mathcal I,\\
 &X_{b,k}=F_e(q_k)(X_g^{\mathsf B})^{-1},\quad C_k^h=0,
                                               &&k\in\mathcal N^+,\\
 &c_{{\rm env},k}(z_k)\geq0,
       \quad\eqref{eq:arm-payload-clearance},
       \quad\text{joint/rate bounds},&&k\in\mathcal N^+,\\
 &\underline v^q\leq v_k^q\leq\overline v^q,\quad
 f_k\in\mathcal F_\delta,\quad Gf_k=\mathcal W_{{\rm req},k},
                                               &&k\in\mathcal I.
\end{aligned}
\end{equation}
The aperture is the measured parameter $\gamma_{\mathsf B}$ and is
absent from the decision vector. Hard geometric closure is retained as
the acquired rank and positive-certificate data in
\eqref{eq:arm-rank-manifold}--\eqref{eq:arm-geometric-closure}.
If both object and gripper poses are retained instead of substituting
$X_{b,k}$, the equivalent pose equations are
\eqref{eq:arm-event-compatibility}.
Logarithms are evaluated in a fixed normal neighborhood with declared
angular/length scaling. Future lift/transfer nodes satisfy the full-body
clearance catalogue and \eqref{eq:arm-payload-clearance}.  An optional
hard terminal target is $X_{b,H}=X_b^{\rm pre}$, using independent
local pose equations and the endpoint rank check below.  The transition
to $\mathsf C$ requires the measured pre-placement pose and velocity
tolerances, retained grasp/force-capacity margins, positive support
clearance, and an admitted $\mathsf C$ placement trajectory.  Placement begins
in the next problem.

\subsubsection{\texorpdfstring{Problem $\mathsf C$: placement, support transfer, release and retreat}{Problem C: placement, support transfer, release and retreat}}
\label{subsec:arm-mpc-c}
The first internal layout retains the same rigid grasp, geometric rank
manifold, positive closure certificate, held aperture, loaded dynamics,
and actual force limits as $\mathsf B$.  It lowers from
$X_b^{\rm pre}$ to $X_b^{(d)}$ along the fixed placement corridor.
The held-layout cost uses $J_k^{\mathsf C}$ with
$J_{{\rm release},k}=0$ and terminal term
\begin{equation}\label{eq:arm-placement-terminal-cost}
 \Phi_{\mathsf C}^{\rm place}=
 \|\operatorname{vee}\log((X_b^{(d)})^{-1}X_{b,H})\|_{Q_b}^2
 +\|J_b(q_H)x_H^v\|_{Q_{b,v}}^2.
\end{equation}
If placement is a hard future pose equality, remove its implied
zero support gaps and verify its independent endpoint differential.
The transition from the placement cost $\Phi_{\mathsf C}^{\rm place}$
to the next internal layout requires measured support acquisition. Over a dwell $T_s>0$, verify the supported
pose tolerance, center-of-mass reserve in the support polygon,
$\|\hat{\mathcal V}_b\|\leq\varepsilon_s$, and measured force balance
\begin{equation}\label{eq:arm-support-switch-test}
 \|S_w(\hat{\mathcal W}_{\rm pads}+
       \hat{\mathcal W}_{\rm support}-\hat{\mathcal W}_{\rm req})\|
        \leq\varepsilon_W.
\end{equation}
Here $\hat{\mathcal W}_{\rm req}$ includes gravity and measured inertial
terms, as in \eqref{eq:arm-payload-wrench}.
Before support acquisition set $\mathcal W_{S,k}=0$.  In the admitted
held-support layout retain attachment and impose the independent
equations $C_S(X_b(q_k))=0$ from
\eqref{eq:arm-release-support-mode}, omitting their implied gaps.
The support-friction law \eqref{eq:arm-release-support-friction} is
evaluated at the attached pose and $\mathcal V_{b,k}=J_b(q_k)x_k^v$.
The supported wrench balance is
\begin{equation}\label{eq:arm-held-support-balance}
 Gf_k+\mathcal W_{S,k}=\mathcal W_{{\rm req},k}.
\end{equation}
Thus actual hand capacity is tested against
$\mathcal W_{{\rm req},k}-\mathcal W_{S,k}$, while the loaded
acceleration equation includes $+J_b^\top\mathcal W_{S,k}$ once.
While the rigid grasp remains engaged, the interface controller reverses
\eqref{eq:arm-support-balance} from $\chi=1$ to $\chi=0$; the arm torque compensates the support wrench using the loaded-support
dynamics during unloading. Release requires a feasible
support-wrench allocation with reserve and measured net hand load below
$\varepsilon_{\rm unload}$. Individual pad normal forces, support reactions and relative motion are
monitored throughout the transfer. These interval checks establish the
load-transfer history in addition to the final support-contact test.
During scheduled unloading, the loaded preload and friction reserves
decrease to zero with $\chi$, replacing the positive transport lower
bounds.  At nominal zero pad load eliminate
the unloaded pad-force variables and their cone/bound rows in the new
release layout; measured residual load is handled by the residual-wrench
model below.  The zero-force apex is thereby removed from the smooth friction-cone
constraint system of the release layout. At zero force all cone and bound
constraints of a contact are active at once and their gradients are
linearly dependent; eliminating the force variable avoids this degenerate
point.

After these placement and unloading tests, remove
\eqref{eq:arm-event-compatibility} and let the box have its own
$SE(3)$ pose $X_{b,k}$ and body twist $\mathcal V_{b,k}$, initialized
from the measured supported state. The rigid box is allowed small
predicted tangential motion: it may slide and turn about the vertical on
its support, while its height and its two tilt angles are fixed. This fixed
upright-support mode imposes three
independent future equations
\begin{equation}\label{eq:arm-release-support-mode}
 C_S(X_{b,k})=\begin{pmatrix}
 e_3^\top p_{b,k}-h_s-a_3\\e_3^\top R_{b,k}e_1\\e_3^\top R_{b,k}e_2
 \end{pmatrix}=0,\qquad k\in\mathcal N^+,
 \quad e_3^\top R_{b,k}e_3>0.
\end{equation}
The bottom-contact gaps implied by these equations are omitted from the
NLP. Positive normal support reactions and all other geometric clearances
are still checked. The support geometry and available force capacity must
admit the small measured initial motion under the finite-step model.

For bottom contact points $\rho_i^S=(\pm a_1,\pm a_2,-a_3)$, let
$N_{i,k}^S\in[0,\overline N_i^S]$ be reaction variables and define world
contact velocities and a smooth bounded friction model by
\begin{equation}\label{eq:arm-release-support-friction}
 \begin{aligned}
 u_{i,k}^{\parallel}&=(I-e_3e_3^\top)R_{b,k}
       (\mathcal V_{b,k}^p+\mathcal V_{b,k}^{\omega}\times\rho_i^S),\\
 f_{i,k}^S&=R_{b,k}^\top\left(N_{i,k}^Se_3-
       \mu_SN_{i,k}^S\frac{u_{i,k}^{\parallel}}
               {\sqrt{v_S^2+\|u_{i,k}^{\parallel}\|^2}}\right),
       \qquad v_S>0,\\
 \mathcal W_{S,k}&=\sum_i\operatorname{col}(\rho_i^S\times f_{i,k}^S,f_{i,k}^S).
 \end{aligned}
\end{equation}
This is an explicitly regularized sliding-friction prediction law; its
parameters must be identified for the support. The regularization speed
$v_S$ also has to be compatible with the explicit object step
\eqref{eq:arm-release-object-dynamics} below. For a tangential speed small
compared with $v_S$ the friction force is nearly linear in that speed, with
slope $\mu_SN^S_{i,k}/v_S$, and one explicit step multiplies the yaw rate
of the supported box by $1-\kappa_S$ with
$\kappa_S=h\mu_Sm_bg\,r_S^2/(I_{b,3}v_S)$. Here $r_S$ is the horizontal
distance of the support points from the center of mass and $I_{b,3}$ is the
moment of inertia about the vertical axis; for a homogeneous box supported
at its bottom corners $\kappa_S=3h\mu_Sg/v_S$. The step is stable only for
$\kappa_S<2$. A small speed such as $v_S=0.01$~m/s gives
$\kappa_S\approx35$ for $\mu_S=0.6$ and $h=20$~ms, so the predicted yaw
motion is amplified at every step and the KKT matrix becomes numerically
singular. The experiments of Section~\ref{sec:results1019B} therefore use a
speed proportional to the sampling interval,
$v_S=(37.5\,{\rm m/s^2})\,h$, for which $\kappa_S<1/2$. It satisfies the Coulomb
force bound and opposes tangential motion. Normal reactions are determined
with the rigid support equations and Newton--Euler dynamics. With a declared
residual-pad wrench prediction $\mathcal W_{P,k}=\eta_k\hat{\mathcal W}_{P,0}$,
$0\leq\eta_k\leq1$, use
\begin{equation}\label{eq:arm-release-object-dynamics}
 \begin{aligned}
 \mathcal V_{b,k+1}&=\mathcal V_{b,k}+h\mathcal G_b^{-1}
   \bigl(\operatorname{ad}_{\mathcal V_{b,k}}^\top\mathcal G_b\mathcal V_{b,k}
              +\mathcal W_g+\mathcal W_{S,k}+\mathcal W_{P,k}\bigr),\\
 X_{b,k+1}&=X_{b,k}\exp(h[\mathcal V_{b,k+1}]).
 \end{aligned}
\end{equation}
Here $\eta_k$ is a prescribed smooth unloading profile validated against
individual pad-force measurements; it is zero after the measured pad-force
release. For example, with carried measured unloading progress
$s_{\hat t}\in[0,1]$ and duration $T_u>0$, set
$s_k=\min\{1,s_{\hat t}+kh/T_u\}$ and
$\eta_k=1-10s_k^3+15s_k^4-6s_k^5$. While $\eta_k>0$, set
$v_k^\gamma=0$ and use the force servo. The first aperture command also
remains zero until measured pad normal forces are below $\varepsilon_P$;
thereafter the aperture velocity servo is enabled.
Residual-wrench uncertainty is included in the stated disturbance envelope.
The reaction on the unloaded arm is
$\tau^{\rm ext}_{\mathsf C,k}=-J_{eb,k}^\top\mathcal W_{P,k}$, where
$J_{eb,k}=\operatorname{Ad}_{X_{b,k}^{-1}F_e(q_k)}J_e^b(q_k)$.
The equal and opposite hand/object wrench therefore enters both models once.
The first update retains the measured arm and object velocities. If the measured support mode is inconsistent or the
reaction problem is infeasible, release is withheld and support recovery
is invoked.

Let $X_b^{\rm placed}$ be the measured placement pose and define
$e_{b,k}=\operatorname{vee}\log((X_b^{\rm placed})^{-1}X_{b,k})$.
In addition to $J_k^\mathsf C$, minimize
\begin{equation}\label{eq:arm-release-motion-cost}
 J_{{\rm release},k}=
       \|e_{b,k}\|_{Q_{b,p}}^2+\|\mathcal V_{b,k}\|_{Q_{b,v}}^2
       +\eta_k\|J_{eb,k}x_k^v\|_{Q_{\rm drag}}^2
       +\rho_o|\gamma_k-\gamma_{\rm open}|^2.
\end{equation}
The first two terms penalize residual object displacement and motion.
The third discourages hand motion while residual preload can transmit
it to the object. When the measured pad wrench is zero, the nominal detached object
is driven only by its support and the declared disturbance envelope.

An internal measured clearance flag $b_{\rm sep}\in\{0,1\}$ governs the
reference and catalogue for $\mathsf C$. While $b_{\rm sep}=0$, command
opening with $v_k^\gamma\geq0$ once the pressure-release test passes,
and impose only on future nodes
\begin{equation}\label{eq:arm-release-before-retreat}
 e_3^\top(p_e(q_k)-p_e(\hat q_{\mathsf C}))\leq\varepsilon_{\rm hold},
 \qquad
 \|p_e(q_k)-p_e(\hat q_{\mathsf C})\|^2\leq\delta_{\rm hold}^2.
\end{equation}
Here $\hat q_{\mathsf C}$ is the measured configuration at release entry;
these bounded hold regions permit damping small measured rates. The
reference holds position and opens the hand. The flag changes to
$b_{\rm sep}=1$ only when
both pads and every hand component have measured/certified clearance
at least $d_{\rm sep}>0$ from the object, individual pad normal forces are
below $\varepsilon_P$, and support-wrench feasibility and object-motion
tests remain satisfied for a dwell $T_{\rm sep}>0$.
The full-body separating inequalities \eqref{eq:arm-smooth-separation}
with margin $d_{\rm sep}$ implement this clearance test; aperture alone
is insufficient. Then remove the hold region, activate the upward/home
retreat reference, and retain positive hand--object separation on every
future retreat node and its certified intervals. A failed force or support
check returns to support recovery before further opening or upward motion.

The formal placement--release OCP uses a measured internal layout
$m$ fixed during each solve. Write $y_k^{\mathsf C,m}$ for the
corresponding state in Table~\ref{tab:arm-layout-variables}, and
$\varpi_k^{\mathsf C,m}$ for its physical force allocations. Then
\begin{equation}\label{eq:arm-ocp-c}
\begin{aligned}
 \min_{\substack{(y_k^{\mathsf C,m},z_k)_{k\in\mathcal N^+}\\
             (v_k,\varpi_k^{\mathsf C,m})_{k\in\mathcal I}}}
 &h\sum_{k\in\mathcal I}
       (J_k^{\mathsf C}+J_{{\rm release},k})+\Phi_{\mathsf C}^m,\\
 \text{subject to}\quad
 &y_0^{\mathsf C,m}=\hat y^{\mathsf C,m}
                         \quad\text{as substituted data},\\
 &y_{k+1}^{\mathsf C,m}
    =F_{\mathsf C,k}^m(y_k^{\mathsf C,m},v_k,\varpi_k^{\mathsf C,m}),
                                               &&k\in\mathcal I,\\
 &C_k^h=0,\quad C_{\mathsf C,k}^{\rm mode,m}=0,\quad
       b_{\mathsf C,k}^m(y_k^{\mathsf C,m},z_k)\geq0,
                                               &&k\in\mathcal N^+,\\
 &v_k\in\mathbb V_{\mathsf C}^m,\quad
       \text{wrench balance and force limits of layout }m,
                                               &&k\in\mathcal I,\\
 &y_H^{\mathsf C,m}\in\mathcal X_f^{\mathsf C,m}.
\end{aligned}
\end{equation}
In held placement, $F_{\mathsf C,k}^m$ is the loaded arm map
\eqref{eq:arm-discrete-dynamics}, $J_{{\rm release},k}=0$,
and $\Phi_{\mathsf C}^m=\Phi_{\mathsf C}^{\rm place}$.
After measured support acquisition, the same loaded map contains the
support torque and the force variables satisfy
\eqref{eq:arm-held-support-balance}; the independent support equations
are \eqref{eq:arm-release-support-mode}. In detached supported release
and retreat, $F_{\mathsf C,k}^m$ combines the unloaded arm map with
\eqref{eq:arm-release-object-dynamics}, the force variables are the
support normals $N_{i,k}^S$, and
$\Phi_{\mathsf C}^m$ is \eqref{eq:arm-terminal-cost} below.
The vector $b_{\mathsf C,k}^m$ contains the applicable future state
bounds, full-body collision and clearance constraints, the small-motion
support neighborhood, and the hold region
\eqref{eq:arm-release-before-retreat} until measured separation.
The control set $\mathbb V_{\mathsf C}^m$ contains the original torque
and aperture-rate limits, with $v^\gamma=0$ while the force servo is
active. Thus each symbol in \eqref{eq:arm-ocp-c} denotes one of the
explicit mechanical maps or constraint catalogues above; $m$ is
fixed supervisory data during an OCP solve.
The layout-dependent terminal set $\mathcal X_f^{\mathsf C,m}$ is
specified below; its retreat requirements become active after measured
separation.

\subsubsection{Assembled finite-horizon problems and terminal retreat}
\label{subsec:arm-mpc-implementation}
With $y_k^{\mathsf P}$ collecting the arm state and, in $\mathsf C$'s supported
release/retreat layout, the independent object state, all three
problems have the common form
\begin{equation}\label{eq:arm-nlp}
 \begin{aligned}
 \min\quad&h\sum_{k=0}^{H-1}\left(J_k^{\mathsf P}+
       \mathbf1_{\{\mathsf P=\mathsf C\}}J_{{\rm release},k}\right)+\Phi_{\mathsf P}\\
 \textnormal{s.t.}\quad&y_{k+1}^{\mathsf P}=F_{\mathsf P,k}(y_k^{\mathsf P},v_k,\varpi_k^{\mathsf P}),
                           &&k=0,\ldots,H-1,\\
 &z_k-h^z_{\mathsf P}(y_k^{\mathsf P})=0,
                           &&k=1,\ldots,H,\\
 &C_{\mathsf P,k}^{\rm mode}(y_k^{\mathsf P})=0,\quad
       b_{\mathsf P,k}(y_k^{\mathsf P},z_k)\geq0,
                           &&k=1,\ldots,H,\\
 &v_k\in\mathbb V_{\mathsf P},\quad
       \text{interval wrench constraints of the selected problem},
                           &&k=0,\ldots,H-1,\\
 &y_H^{\mathsf P}\in\mathcal X_f^{\mathsf P},\\
 &y_0^{\mathsf P}=\hat y^{\mathsf P}\quad\text{(substituted measured parameter)}.
 \end{aligned}
\end{equation}
The auxiliary primal variables $\varpi_k^{\mathsf P}$, $k\in\mathcal I$,
represent the displayed force allocations or support reactions. The
actuator controls are collected in the separate control block.
The force variables enter the costs stated so far only through the capacity
term of $\Phi_\mathsf A$, while the balance equations leave directions in
which they can change at no cost (Section~\ref{sec:arm-qualification-details}).
Along such directions the reduced Hessian is singular and no second-order
sufficient condition can hold. The objective of \eqref{eq:arm-nlp}
therefore also contains a small force-allocation regularization with a
weight $\rho_f>0$,
\begin{equation}\label{eq:arm-force-regularization}
 J^{\rm reg}=
 \begin{cases}
 \rho_f\displaystyle\sum_{\ell=1}^{N_w}\|f^\ell-f^{{\rm ref},\ell}\|^2,
      &\mathsf P=\mathsf A,\\[2mm]
 h\rho_f\displaystyle\sum_{k\in\mathcal I}\|f_k-f^{\rm ref}\|^2,
      &\text{held layouts of $\mathsf B$ and $\mathsf C$},\\[2mm]
 h\rho_f\displaystyle\sum_{k\in\mathcal I}\sum_i
      \bigl(N^S_{i,k}-\tfrac14m_bg\bigr)^2,
      &\text{detached layouts of $\mathsf C$}.
 \end{cases}
\end{equation}
The reference allocation $f^{\rm ref}$ is fixed data for the upright box:
a preload $N_{\rm pre}>0$ along each contact normal plus an equal share
$m_bg/N_c$ of the weight at each contact,
$f_j^{\rm ref}=N_{\rm pre}n_j+(m_bg/N_c)e_3$. In $f^{{\rm ref},\ell}$ the
weight $m_bg$ is replaced by the vertical force of scenario $w^\ell$. This
term makes the allocation unique.
The interval force block is empty in $\mathsf A$, whose certificate
forces $f^\ell$ are terminal variables in \eqref{eq:arm-acquisition-terminal-cost}.
The release term $J_{{\rm release},k}$ is zero in $\mathsf A$,
$\mathsf B$, and $\mathsf C$'s held placement layout. Supported release/retreat includes future
bounds on $\|\mathcal V_{b,k}\|$ and
$\|e_{b,k}\|$ defining the stated small-motion support neighborhood.
The maps and constraints are exactly \eqref{eq:arm-discrete-dynamics},
\eqref{eq:arm-explicit-z}, the collision catalogue, and
\eqref{eq:arm-transport-force-balance},
\eqref{eq:arm-held-support-balance}, or
\eqref{eq:arm-release-support-mode}--\eqref{eq:arm-release-object-dynamics}
for the selected problem. Problem $\mathsf A$ also has the terminal
variables and constraints in \eqref{eq:arm-acquisition-terminal-cost}.
The geometric rank representation
\eqref{eq:arm-rank-manifold}--\eqref{eq:arm-closure-section}
applies at acquisition and throughout the held layouts, with the
stated substitution of dependent factors/certificates.  The terminal
costs are $\Phi_\mathsf A$, $\Phi_\mathsf B$, and either
$\Phi_{\mathsf C}^{\rm place}$ or the supported-release/home cost below,
according to $\mathsf C$'s fixed measured internal layout.  Each such layout
is fixed during a solve; its change initiates a new local RTI episode
within the same third top-level OCP. Fixed future apertures and other structural identities are
substituted; only an independent equality set is assembled.

The regularity properties of this transcription, including mode-implied
constraint removal and the conditions on the complete shooting system,
are developed in Section~\ref{sec:arm-qualification-details}.

Let $q_{\rm home}$ be a reachable collision-free home configuration with
$h_{\rm home}=e_3^\top p_e(q_{\rm home})$. Once $b_{\rm sep}=1$, the
terminal cost is
\begin{equation}\label{eq:arm-terminal-cost}
 \begin{aligned}
 \Phi_\mathsf C={}&
 \|\operatorname{Log}_{q_{\rm home}}q_H\|_{Q_h}^2
 +\|x_H^v\|_{Q_{v,f}}^2
 +\rho_h(e_3^\top p_e(q_H)-h_{\rm home})^2\\
 &+\rho_{\gamma,f}|\gamma_H-\gamma_{\rm open}|^2
 +\|e_{b,H}\|_{Q_{b,f}}^2+\|\mathcal V_{b,H}\|_{Q_{b,vf}}^2.
 \end{aligned}
\end{equation}
In the supported-release layout before separation, replace the
home-position terms by their measured
release-entry references. The future retreat terminal set requires zero
joint rates, open aperture, full-body clearance, and admissible unloaded
holding torque. In particular, for the detached retreat layout,
\begin{equation}\label{eq:arm-retreat-terminal-set}
\mathcal X_f^{\mathsf C,\mathrm{ret}}=
\left\{y:\begin{array}{l}
 x^v=0,\quad\gamma=\gamma_{\rm open},\\
 \|\operatorname{Log}_{q_{\rm home}}q\|\leq\delta_{\rm retreat},\\
 b_{\mathsf C,H}^{\mathrm{ret}}(y,h^z_{\mathsf C}(y))\geq0,\\
 \underline v^q\leq g_0(q,\gamma_{\rm open})\leq\overline v^q
\end{array}\right\}.
\end{equation}
The margin vector includes full-body and hand--object clearance.
In \eqref{eq:arm-nlp}, $\mathcal X_f^{\mathsf C}$ denotes the set
for the currently fixed internal layout $m$.
For acquisition and the other $\mathsf C$ layouts, $\mathcal X_f$ is the
corresponding state domain with the terminal conditions already stated;
the optional hard pre-placement target in $\mathsf B$ is included in
$\mathcal X_f^{\mathsf B}$ when selected. A bound implied by a hard
terminal equality is represented once in the retained terminal catalogue,
as detailed in Section~\ref{sec:arm-qualification-details}.
A local home neighborhood may be chosen nonsingular by
certifying $\sigma_{\min}(\widetilde J_e(q_{\rm home}))\geq\sigma_0$,
where $\widetilde J_e=S_w^{-1}J_e^b$ is the body Jacobian of the gripper
frame with its angular rows multiplied by $\ell_b$, so that all six rows
have the units of a linear velocity,
and a Jacobian Lipschitz bound $L_J$, then restricting
$\|\operatorname{Log}_{q_{\rm home}}q_H\|\leq\delta_{\rm retreat}$ with
$\delta_{\rm retreat}\leq(\sigma_0-\varepsilon_J)/L_J$,
$\sigma_0>\varepsilon_J>0$; the constant-Jacobian case needs no division.
These are geometric task data to validate for the chosen robot and scene.

Each MPC uses smooth fixed data and a fixed contact catalogue during one
solve. Grasp acquisition, support acquisition, release and clearance
switches create new problems and require fresh state/dual initialization
and the local proximity/regularity checks of the pathfollowing analysis,
including Assumption~\ref{ass:offline-backup} for offline recovery.
The local theorem applies separately to each admitted regular segment.
The measured switches concatenate those segments through the specified
initialization and physical acceptance tests.
\section{Parametric Riemannian Nonlinear Program}
\label{sec:parametric-problem}

\paragraph{Indices.}
The integer $\hat t$ labels the current feedback measurement. During the solve
at $\hat t$, the measured state, environmental data and contact schedule are
frozen. The integer $j$ labels optional nonlinear corrector rounds at that
same measurement; $j=0$ is the single-round RTI choice. Within a frozen local
QP, $\tau\in[0,1]$ is an artificial homotopy parameter.

The entries of the trajectory $r$ indexed by prediction node describe
predicted plant motion. In the notation $f(r,t),C(r,t),D(r,t)$ below, $t$
is a shorthand for the current frozen problem data.

\paragraph{Original variable blocks and horizon.}
Write $H=H_{\hat t}$ for the number of control intervals and use relative state indices
$k=0,\ldots,H$. Let
\[
 \mathcal M_u=\prod_{m=1}^{M^u}\mathcal M_m^u,\qquad
 \mathcal M_q=\prod_{m=1}^{M^q}\mathcal M_m^q,\qquad
 n_u=\dim\mathcal M_u,\quad n_q=\dim\mathcal M_q.
\]
The primal trajectory is $r=(x,q,u,v,z)$, with
\begin{equation}\label{eq:part1-product-domain}
 \mathcal M_H=(\mathbb R^{d_x})^{H+1}\times\mathcal M_q^{H+1}
 \times\mathcal M_u^{H}\times(\mathbb R^{d_v})^{H}
 \times(\mathbb R^{d_z})^{H+1}.
\end{equation}
Its nominal intrinsic dimension is
$(H+1)(d_x+n_q+d_z)+H(n_u+d_v)$.
There is no terminal $u_{H}$ or $v_{H}$. The variables $z_0,\ldots,z_{H}$
remain explicit.

The node-zero quantities $(x_0,q_0)$ are the measured data and
$z_0=h^z_0(x_0,q_0)$ is evaluated from them. These three blocks are
substituted before the numerical tangent space is formed; their increments
are zero and no initial-state equality row is retained. The nominal product
above records the complete trajectory layout. Node-zero substitution reduces its nominal decision dimension by
$d_x+n_q+d_z$. The shooting and forward-kinematics equalities remain
constraints on the retained variables. Operational state bounds, collision margins, contact
equalities and stationary-phase requirements apply only at nodes
$k=1,\ldots,H$. The first derived state $z_0$ is unconstrained.

Coordinates prescribed as constants by the current OCP, including coordinates with coincident lower and upper bounds, are substituted before forming the tangent systems. Structural elimination is distinct from temporary fixing by a working set. A structurally eliminated coordinate is not a variable of the problem at all. A coordinate at a bound, in contrast, remains a free variable until that bound is selected in the operational working set, and it can be released again during the same QP; its inequality remains in the problem. Horizon-wide decision variables are represented once. For example, optimizing the remaining duration \(T_{\hat t}\) replaces \(r\) and \(\mathcal M_H\) by \((r,T_{\hat t})\) and \(\mathcal M_H\times\mathbb R\), respectively, increasing the nominal intrinsic dimension by one. The duration bounds are imposed once, while all stage and terminal dependence on \(T_{\hat t}\), including mixed derivatives, is assembled into its shared Jacobian column and Hessian row and column.

\paragraph{Manifold and Euclidean integrator equalities.}
The exact shooting and initial conditions retain their original meaning:
\begin{equation}\label{eq:intrinsic-shooting}
\begin{aligned}
 C^x_k(r,t)&:=x_{k+1}-f^x_k(x_k,q_k,u_k,v_k;t)=0,
                         &&0\leq k\leq H-1,\\
 q_{k+1}&=f^q_k(q_k,x_k,u_k,v_k;t)\quad\hbox{in }\mathcal M_q,
                         &&0\leq k\leq H-1,\\
 C^h_k(r,t)&:=z_k-h^z_k(x_k,q_k;t)=0,
                         &&1\leq k\leq H,\\
 x_0&=\hat x(\hat t),\qquad q_0=\hat q(\hat t).
\end{aligned}
\end{equation}
Geometric shooting is equivalently membership of
$(q_{k+1},f^q_k(r;t))$ in the diagonal submanifold
$\Delta_{\mathcal M_q}=\{(a,a):a\in\mathcal M_q\}$.
These equations compare points intrinsically on $\mathcal M_q$.

To write the geometric shooting condition as an equation $C=0$ with a
vector-valued residual, the mismatch between two points is measured by a
tangent vector. On a normal neighborhood, define the geometric defect
\[
 \operatorname{def}(a,b)=\operatorname{Log}_a(b)\in T_a\mathcal M_q.
\]
Its zero is exactly $a=b$. The defect lives in the tangent space at $a$,
which moves with the iterate. To compare and add residuals evaluated at
different trial points, they are transported to one fixed tangent space,
that of an \emph{anchor} $\bar q_k$: a fixed reference point near $q_k$,
for example its value when the local model is built. Fix residual anchors $\bar q_k$ during one local
model and each corresponding merit comparison, and set
\begin{equation}\label{eq:intrinsic-residual-representative}
\begin{aligned}
 C^q_k(r,t)&=\mathcal P_{q_{k+1}\to\bar q_{k+1}}
            \operatorname{Log}_{q_{k+1}}\!\bigl(f^q_k(r;t)\bigr),\\
 C^0(r,t)&=\operatorname{col}\left(x_0-\hat x(\hat t),\,
        \mathcal P_{q_0\to\bar q_0}
            \operatorname{Log}_{q_0}\hat q(\hat t)\right).
\end{aligned}
\end{equation}
Here $\mathcal P$ is parallel transport along the specified unique short
geodesic. A smooth equivalent Lie-group error or inverse-retraction defect may be
used, with its exact first and second differentials included in the SQP
model.
The displayed $C^0$ describes measurement discrepancy when comparing
prepared and measured problem data. It is omitted from the retained
constraint vector after the node-zero substitution.
An intrinsic logarithm here measures a mismatch between points. The variables
remain on $\mathcal M_H$, and numerical matrices represent linear operators
in tangent frames only \cite[Chapters~5,6,10]{boumal2023introduction}.

The residual representative in
\eqref{eq:intrinsic-residual-representative} is assumed to be
twice continuously differentiable on a neighborhood of the
reference trajectory. For the logarithmic construction, the
logarithm and transport are smooth and single valued on this
neighborhood. Derivatives are taken with respect to all geometric
arguments, including both endpoints of each shooting defect.

At a feasible geometric shooting pair $(a,a)$, the covariant first derivative
of the untransported defect is especially simple:
\[
 D\operatorname{def}(a,a)[\eta_a,\eta_b]=\eta_b-\eta_a\in T_a\mathcal M_q.
\]
Thus the familiar equation $D f^q_k[\eta]-\eta_{q_{k+1}}=0$ is legitimate
\emph{at a matching pair}, where the vectors share a tangent space; away from
matching, the full derivative of the defect is required. The full-horizon assembly retains the resulting two-node Hessian
couplings.

\paragraph{Primal problem and dual spaces.}
Let $C^{\rm event}$ collect the independent phase and terminal equalities. Set $C=\operatorname{col}(C^x,C^q,C^h,C^{\rm event})$, let $D$ contain
the operational inequalities, and let $b_{\rm smp}$ contain the nonstructural
simple margins, that is, the simple bounds written as $b_{\rm smp}(r,t)\geq0$. The current problem is
\begin{equation}\label{eq:priemnlp}
 \min_{r\in\mathcal M_H} f(r,t)
 \quad\text{subject to}\quad C(r,t)=0,\quad D(r,t)\geq0,
 \quad b_{\rm smp}(r,t)\geq0.
\end{equation}
This is \eqref{eq:riemnlp} in compact form: $f$ is the sum of the stage and
terminal costs, $C$ contains the shooting, feature and event equalities,
and the two inequality blocks contain the state constraints $c_k\geq0$
and the limits that define $\mathbb U$, $\mathbb V$ and $\mathbb Q$. The
parameter $t$ carries the measured initial state and all other data of the
current sample.

The geometric duals belong to the dual residual tangent spaces:
\begin{equation}\label{eq:dualspace}
\begin{aligned}
 \lambda^q_k&\in T^*_{\bar q_{k+1}}\mathcal M_q,\quad 0\leq k<H,\\
 \lambda^x&\in(\mathbb R^{Hd_x})^*,\qquad
 \lambda^h\in(\mathbb R^{Hd_z})^*,\\
 \mu&\leq0,\qquad \nu\leq0,
\end{aligned}
\end{equation}
and the multiplier $\lambda^{\rm event}$ of $C^{\rm event}$ is a covector on
the Euclidean event-residual space.

The geometric shooting condition compares two configurations on
$\mathcal M_q$. Its residual representative
$C_k^q(r,t)$ in \eqref{eq:intrinsic-residual-representative}, however,
is a tangent vector at the fixed residual anchor $\bar q_{k+1}$. Write
\[
    \mathcal V_k := T_{\bar q_{k+1}}\mathcal M_q,
    \qquad C_k^q(r,t)\in\mathcal V_k.
\]
The multiplier belongs to the dual vector space,
$\lambda_k^q\in\mathcal V_k^*$. This means that $\lambda_k^q$ is a
linear functional: it takes a residual vector as its argument and
returns a real number. The corresponding contribution to the
Lagrangian is therefore
\[
    \langle\lambda_k^q,C_k^q(r,t)\rangle
       := \lambda_k^q\!\left(C_k^q(r,t)\right)\in\mathbb R.
\]
Here the angle brackets denote evaluation of a covector on a vector.
Node zero is substituted, so the operational multiplier array has no
initial-condition multiplier block.

\paragraph{Numerical representation in a frame.}
A residual frame represents each covector by a column of coefficients in
the corresponding dual basis. Let $e_{k,1},\ldots,e_{k,n_q}$ be a basis
of $\mathcal V_k$. Its dual basis $e_k^1,\ldots,e_k^{n_q}$ consists
of the linear functionals satisfying
\[
    e_k^a(e_{k,b})=\delta^a_b,
\]
where $\delta^a_b$ is one when $a=b$ and zero otherwise. Write
\[
    C_k^q=\sum_{a=1}^{n_q}c_k^a e_{k,a},
    \qquad
    \lambda_k^q=\sum_{a=1}^{n_q}\lambda_{k,a}^q e_k^a.
\]
With the coefficient columns
$\boldsymbol c_k=(c_k^1,\ldots,c_k^{n_q})^\top$ and
$\boldsymbol\lambda_k^q=(\lambda_{k,1}^q,\ldots,
\lambda_{k,n_q}^q)^\top$, the dual pairing becomes
\[
    \langle\lambda_k^q,C_k^q\rangle
      =\sum_{a=1}^{n_q}\lambda_{k,a}^q c_k^a
      =(\boldsymbol\lambda_k^q)^\top\boldsymbol c_k.
\]
Thus the implementation uses the familiar multiplier--residual dot
product. The frame is a basis of the residual tangent space.

\paragraph{Changing the residual frame.}
The scalar pairing is invariant under a change of residual frame.
To specify the transformation convention unambiguously, suppose that
an invertible change of residual coordinates is written as
\[
    \widetilde{\boldsymbol c}_k=S_k\boldsymbol c_k.
\]
The multiplier coefficients in the new dual frame must then satisfy
\[
    \widetilde{\boldsymbol\lambda}_k^q
       =S_k^{-\top}\boldsymbol\lambda_k^q,
    \qquad S_k^{-\top}:=(S_k^{-1})^\top,
\]
because
\[
    (\widetilde{\boldsymbol\lambda}_k^q)^\top
       \widetilde{\boldsymbol c}_k
    =(S_k^{-\top}\boldsymbol\lambda_k^q)^\top
       S_k\boldsymbol c_k
    =(\boldsymbol\lambda_k^q)^\top\boldsymbol c_k.
\]
For example, doubling a residual coordinate requires halving its multiplier
coordinate, preserving both geometric objects and their pairing. This inverse-transpose rule is the
coordinate meaning of the inverse dual transformation. The local QP uses a fixed residual frame throughout its numerical
homotopy; the inverse-transpose identity specifies consistency between
representations.

\paragraph{The scalar Lagrangian.}
Collect all equality residuals, including the event residuals, in
$C(r,t)$, with corresponding multiplier $\lambda$. The notation
$\langle\lambda,C(r,t)\rangle$ denotes the sum of the dual pairings
for these blocks. With $D(r,t)$ the operational inequality margins and
$b_{\mathrm{smp}}(r,t)$ the simple-bound margins, define
\begin{equation}\label{eq:lags}
\begin{aligned}
 L(r,\lambda,\mu,\nu,t)
   ={}& f(r,t)+\langle\lambda,C(r,t)\rangle\\
     &+\mu^\top D(r,t)+\nu^\top b_{\mathrm{smp}}(r,t).
\end{aligned}
\end{equation}
In particular, the geometric shooting contribution is
\[
    \sum_{k=0}^{H-1}
       \langle\lambda_k^q,C_k^q(r,t)\rangle.
\]
Each multiplier acts on the tangent-valued shooting defect determined by
$q_{k+1}$ and $f_k^q(r;t)$. Thus configurations enter the scalar pairing
through their anchored residual.

The equality multipliers are unrestricted in sign. Since feasible
inequalities are written as $D\geq0$ and $b_{\mathrm{smp}}\geq0$,
the plus-sign convention in \eqref{eq:lags} gives
\[
    \mu\leq0,\qquad \nu\leq0,
    \qquad \mu_iD_i=0,\qquad \nu_j b_{\mathrm{smp},j}=0
\]
at a KKT point.

\paragraph{Differentiation for the local quadratic model.}
To form the model at a reference primal--dual point
$(\bar r,\bar\lambda,\bar\mu,\bar\nu)$, regard
\[
    F(r):=L(r,\bar\lambda,\bar\mu,\bar\nu,t)
\]
as a scalar function of $r$ alone. The current problem data $t$, the
residual anchors, and the multiplier coefficients in their fixed dual
frames are held constant in this partial differentiation. Because the
residuals already take values in fixed anchored vector spaces, no
additional state-dependent multiplier extension is needed in this
representation. Thus, for $\eta\in T_r\mathcal M_H$,
\[
\begin{aligned}
    D F(r)[\eta]={}&D_r f(r,t)[\eta]
       +\langle\bar\lambda,D_r C(r,t)[\eta]\rangle\\
       &+\bar\mu^\top D_r D(r,t)[\eta]
       +\bar\nu^\top D_r b_{\mathrm{smp}}(r,t)[\eta].
\end{aligned}
\]
Differentiation includes all residual dependence on $r$, including the
geometric endpoints and transport in the defect, with the multipliers
held fixed with respect to $r$.

The resulting operator is the Riemannian Hessian of the twice
continuously differentiable scalar function $F$, using the
Levi--Civita connection of the declared product metric. Consequently,
\[
    \langle\xi,\operatorname{Hess}F(\bar r)[\eta]\rangle_{\bar r}
    =\langle\eta,\operatorname{Hess}F(\bar r)[\xi]\rangle_{\bar r}
    \quad
    (\xi,\eta\in T_{\bar r}\mathcal M_H).
\]

In an orthonormal primal tangent frame its matrix $H_L$ is symmetric.

\paragraph{Regularity and Optimality Conditions.}
\begin{definition}\label{def:firstorderopt}
At a feasible point $r_*$, the first-order KKT conditions require duals
as in \eqref{eq:dualspace} satisfying
\[
 \operatorname{grad}_r L(r_*,\lambda_*,\mu_*,\nu_*,t)=0,\qquad
 \mu_*\odot D(r_*,t)=0,\qquad \nu_*\odot b_{\rm smp}(r_*,t)=0.
\]
All adjoints use the declared Riemannian metric.
\end{definition}
With $b=\operatorname{col}(D,b_{\rm smp})$ and
$\zeta_b=\operatorname{col}(\mu,\nu)$, let
$\mathcal A=\{i:b_i(r_*,t)=0\}$ and
$\mathcal A_- =\{i\in\mathcal A:(\zeta_b)_i<0\}$.
The critical cone and a larger supported subspace are
\begin{equation}\label{eq:part1-critical-support}
\begin{aligned}
 \mathcal C_*&=\{\eta\in T_{r_*}\mathcal M_H:DC_*[\eta]=0,\ Db_{i,*}[\eta]\geq0\ (i\in\mathcal A),
                         \ Df_*[\eta]=0\},\\
 \mathcal T_-&=\{\eta\in T_{r_*}\mathcal M_H:DC_*[\eta]=0,\ Db_{i,*}[\eta]=0\ (i\in\mathcal A_-)\}.
\end{aligned}
\end{equation}

The KKT conditions imply $\mathcal C_*\subseteq\mathcal T_-$.

Collect the dual variables as $\zeta=\operatorname{col}(\lambda,\zeta_b)$.
We distinguish the following second-order conditions.
\begin{definition}[Riemannian SOSC]\label{def:riem-sosc}
At a feasible KKT pair $(r_*,\zeta_*)$, the second-order sufficient
condition (SOSC) holds if there is $\alpha_{\zeta_*}>0$ such that
\[
 \langle\eta,\operatorname{Hess}_rL(r_*,\zeta_*,t)[\eta]\rangle_{r_*}
 \geq\alpha_{\zeta_*}\|\eta\|_{r_*}^2
 \quad\text{for every }\eta\in\mathcal C_* .
\]
The differential, adjoint, critical cone and Hessian all use the
declared product Riemannian metric.
\end{definition}
\begin{definition}[Riemannian SSOSC]\label{def:riem-ssosc}
The strong second-order sufficient condition (SSOSC) at a KKT pair
requires the same inequality for every $\eta\in\mathcal T_-(t,\zeta_*)$.
The generalized SSOSC requires this condition for every multiplier
$\zeta_*\in\Lambda(t)$. The uniform generalized condition on a compact
path uses one constant $\alpha>0$ independent of $t$ and of the
multiplier, as in \eqref{eq:uniform-path-ssosc}.
\end{definition}
SSOSC implies SOSC because $\mathcal C_*\subseteq\mathcal T_-$.
In words, SOSC requires positive curvature of the Lagrangian along the
critical directions. SSOSC requires it on the larger subspace obtained by
keeping only the inequalities with nonzero multiplier (the \emph{strongly
active} ones) as equalities and ignoring the \emph{weakly active} ones,
which are active with zero multiplier. This stronger condition is what
makes the KKT matrix of a working set nonsingular when the working set
contains the strongly active constraints but possibly none of the weakly
active ones.
Section~\ref{sec:pathfollowing-theory} uses uniform generalized SSOSC. These definitions extend the supported-subspace condition in
\cite{kungurtsev2017predictor} to the product manifold.

\begin{definition}[Riemannian MFCQ]\label{def:mfcq}
The Mangasarian--Fromovitz constraint qualification holds at $r_*$ if,
after structural constants and certified redundant equalities are removed,
$DC_*$ is surjective onto the residual space and some
$\eta\in T_{r_*}\mathcal M_H$ satisfies
$DC_*[\eta]=0$ and $Db_{i,*}[\eta]>0$ for all $i\in\mathcal A$.
\end{definition}

The operational working set $W\subseteq\mathcal A$ selects independent
rows and working-bound tangent restrictions for each linear solve. In the
terminology of active-set methods, the columns of its free-variable space
are then partitioned into basic variables, which are determined by the
working rows, and superbasic variables, which parametrize the null space of
those rows.

\begin{assumption}[Continuous solution path, MFCQ, and uniform SSOSC]
\label{as:sosc}
Fix a compact parameter interval $I=[t_0,t_f]$ on which the horizon $H$,
constraint catalogue, and structural substitutions are fixed. The
corresponding product manifold is denoted by $\mathcal M_H$. Write
\[
 b(r,t)=\operatorname{col}\bigl(D(r,t),b_{\rm smp}(r,t)\bigr),
 \qquad
 \zeta=(\lambda,\mu,\nu),
 \qquad \zeta_b=\operatorname{col}(\mu,\nu).
\]
The following conditions hold throughout $I$.

\begin{enumerate}
\item \textbf{Smooth data and a continuous primal solution path.}
There exists a continuous map
\[
 r_*:I\longrightarrow\mathcal M_H
\]
such that $r_*(t)$ is an isolated local minimizer of the current nonlinear
problem for every $t\in I$. In particular,
\[
 C(r_*(t),t)=0,\qquad b(r_*(t),t)\geq0.
\]
On one fixed neighborhood of the graph of this path, the objective and the
chosen residual representatives are twice continuously differentiable in
the primal variables and in the parameter. In common smooth tangent and
residual frames, their first and second derivatives are bounded and
Lipschitz with one common constant $L$ across the entire homotopy.
The same $L$ may be enlarged to bound the retraction, transport and frame
derivatives on that neighborhood. 

\item \textbf{MFCQ along the entire path.}
Let
\[
 \mathcal A(t)=\{i:b_i(r_*(t),t)=0\},
 \qquad A_C(t)=D_rC(r_*(t),t).
\]
After removal of structural constants and equality equations locally
implied by the retained equalities, $A_C(t)$ is surjective. More quantitatively, there exist
constants $\gamma_C,\gamma_b>0$, independent of $t$, such that
\[
 \|A_C(t)^*y\|_{r_*(t)}\geq\gamma_C\|y\|
 \quad\text{for every equality-residual covector }y,
\]
and for each $t\in I$ there exists $d_t\in T_{r_*(t)}\mathcal M_H$ with
\[
 \|d_t\|_{r_*(t)}\leq1,\qquad A_C(t)[d_t]=0,
 \qquad D_rb_i(r_*(t),t)[d_t]\geq\gamma_b
 \quad(i\in\mathcal A(t)).
\]
The corresponding condition is vacuous when its constraint block is
empty. In particular, MFCQ does not require the active inequality
gradients to be linearly independent.

\item \textbf{Multipliers and uniform generalized SSOSC.}
Define the full KKT multiplier set at the selected primal solution by
\[
\begin{split}
 \Lambda(t)=\bigl\{(\lambda,\mu,\nu):\;&
 \operatorname{grad}_rL(r_*(t),\lambda,\mu,\nu,t)=0,\\
 &\mu\leq0,\quad\nu\leq0,\\
 &\mu\odot D(r_*(t),t)=0,\quad
 \nu\odot b_{\rm smp}(r_*(t),t)=0\bigr\}.
\end{split}
\]
These sets are nonempty, and their elements are uniformly bounded:
\[
 \sup_{t\in I}\sup_{\zeta\in\Lambda(t)}\|\zeta\|\leq M_\Lambda<\infty.
\]
For $\zeta\in\Lambda(t)$ set
\[
 \mathcal A_-(t,\zeta)=\{i\in\mathcal A(t):(\zeta_b)_i<0\}
\]
and define the supported tangent subspace
\[
 \mathcal T_-(t,\zeta)=
 \left\{\eta\in T_{r_*(t)}\mathcal M_H:
 \begin{array}{l}
 A_C(t)[\eta]=0,\\
 D_rb_i(r_*(t),t)[\eta]=0\quad(i\in\mathcal A_-(t,\zeta))
 \end{array}\right\}.
\]
There exists $\alpha>0$, independent of $t$ and of the multiplier choice,
such that
\begin{equation}\label{eq:uniform-path-ssosc}
 \bigl\langle\eta,
   \operatorname{Hess}_rL(r_*(t),\zeta,t)[\eta]
 \bigr\rangle_{r_*(t)}
 \geq\alpha\|\eta\|_{r_*(t)}^2,
 \qquad
 \substack{t\in I,\ \zeta\in\Lambda(t),\\
           \eta\in\mathcal T_-(t,\zeta)}.
\end{equation}
Here $L(r,\zeta,t)$ abbreviates $L(r,\lambda,\mu,\nu,t)$.
This is the uniform, all-multiplier form of SSOSC used for a
possibly dual-degenerate solution path; compare
\cite{kungurtsev2017predictor}. A touching inequality with zero multiplier
contributes no equality to the supported tangent subspace.
The multiplier set may contain several elements, and strict
complementarity is optional.

\item \textbf{Retraction and local initialization.}
The product retraction is smooth and second order on a neighborhood of
the path, or its acceleration term is included in the pullback Hessian.
For a second-order retraction this means, for
$\gamma(s)=R_r(s\eta)$,
\[
 R_r(0)=r,\qquad DR_r(0)=\operatorname{id}_{T_r\mathcal M_H},
 \qquad
 \left.\nabla_{\dot\gamma}\dot\gamma\right|_{s=0}=0.
\]
The starting primal--dual estimate satisfies
\begin{equation}\label{eq:path-initial-distance}
 \delta(\bar r_0,\bar\zeta_0;t_0)
 :=\left(
 d_{\mathcal M_H}(\bar r_0,r_*(t_0))^2+
 \operatorname{dist}(\bar\zeta_0,\Lambda(t_0))^2
 \right)^{1/2}
 \leq\epsilon_0,
\end{equation}
where $\epsilon_0>0$ is sufficiently small. The dual distance is evaluated after expressing the
multipliers in common residual dual spaces. It is a distance to the
multiplier \emph{set}.
\end{enumerate}
\end{assumption}

Existence and continuity of $r_*(t)$ are additional path hypotheses.
The active set may change, and the multiplier set can have several
branches with differing regularity.
For a local theorem that uses ordinary Newton inverses uniformly, one
additionally needs a consistently selected collection of independent
working representations. In precise terms, at each $t$ such a
representation has a multiplier $\zeta^W\in\Lambda(t)$ with
\[
 \mathcal A_-(t,\zeta^W)\subseteq W\subseteq\mathcal A(t),\qquad
 A_W(t)=\begin{bmatrix}D_rC(r_*(t),t)\\D_rb_W(r_*(t),t)\end{bmatrix}
\]
where $A_W(t)$ has independent rows. The associated KKT operator
\[
 K_W(t)=
 \begin{bmatrix}
 \operatorname{Hess}_rL(r_*(t),\zeta^W,t)&A_W(t)^*\\
 A_W(t)&0
 \end{bmatrix}
\]
satisfies $\|K_W(t)^{-1}\|\leq\kappa$ on the selected pieces for a
uniform constant $\kappa$. The bound, when used, is an additional
working-representation hypothesis. Assumption~\ref{ass:regular-path} in
Section~\ref{sec:pathfollowing-theory} states it precisely; that
assumption strengthens Assumption~\ref{as:sosc} by the relaxed
constant-rank qualification, a finite partition of the parameter interval
into pieces with fixed active set, and these uniform bounds.

\paragraph{Real-time OCP update and prediction indices.}
The integer $\hat t$ counts feedback samples, whereas $k$ is a prediction
index. At sample $\hat t$, write the stored approximate trajectory as
\[
 \bar r_{\hat t}=
 \bigl((\bar x_{k|\hat t},\bar q_{k|\hat t},\bar z_{k|\hat t})_{k=0}^{H},
       (\bar u_{k|\hat t},\bar v_{k|\hat t})_{k=0}^{H-1}\bigr).
\]
There are $H+1$ state/feature nodes and $H$ controls. The first control is
applied to the plant. At the next sample, the measurement
$(\hat x(\hat t+1),\hat q(\hat t+1))$ defines the new initial condition.
The maps $f^x_{k|\hat t+1}$, $f^q_{k|\hat t+1}$, and
$h^z_{k|\hat t+1}$ below include the new stage data and physical mesh.
They describe the nominal predicted dynamics and features.

\paragraph{Fixed horizon: shift the tail and append a terminal rollout.}
Assume $H\geq1$ and that one applied control interval corresponds to one
prediction interval. First form the shifted tail
\begin{equation}\label{eq:fixed-horizon-tail-shift}
\begin{aligned}
 (\widetilde x_{k|\hat t+1},\widetilde q_{k|\hat t+1})
   &=(\bar x_{k+1|\hat t},\bar q_{k+1|\hat t}),
       &&k=0,\ldots,H-1,\\
 (\widetilde u_{k|\hat t+1},\widetilde v_{k|\hat t+1})
   &=(\bar u_{k+1|\hat t},\bar v_{k+1|\hat t}),
       &&k=0,\ldots,H-2.
\end{aligned}
\end{equation}
The second range is empty when $H=1$. The old terminal state now occupies
node $H-1$ of the new horizon. The complete warm start therefore includes one new control and one new state.

Choose a stated terminal extension rule $\kappa_{\rm ext}$ and set
\begin{equation}\label{eq:terminal-control-extension}
 (\widetilde u_{H-1|\hat t+1},\widetilde v_{H-1|\hat t+1})
 =\kappa_{\rm ext}
   (\widetilde x_{H-1|\hat t+1},\widetilde q_{H-1|\hat t+1};\hat t+1).
\end{equation}
For example, this rule can evaluate a terminal feedback policy or repeat
the previous last control, provided the latter remains admissible.
The new terminal-state guess is then the result of one nominal dynamics
step:
\begin{equation}\label{eq:terminal-state-extension}
\begin{aligned}
 \widetilde x_{H|\hat t+1}
 &=f^x_{H-1|\hat t+1}
   (\widetilde x_{H-1|\hat t+1},\widetilde q_{H-1|\hat t+1},
    \widetilde u_{H-1|\hat t+1},\widetilde v_{H-1|\hat t+1}),\\
 \widetilde q_{H|\hat t+1}
 &=f^q_{H-1|\hat t+1}
   (\widetilde q_{H-1|\hat t+1},\widetilde x_{H-1|\hat t+1},
    \widetilde u_{H-1|\hat t+1},\widetilde v_{H-1|\hat t+1}),\\
 \widetilde z_{H|\hat t+1}
 &=h^z_{H|\hat t+1}
   (\widetilde x_{H|\hat t+1},\widetilde q_{H|\hat t+1}).
\end{aligned}
\end{equation}
The manifold-valued prediction map generates the geometric terminal state.
At a suitable equilibrium this transition can reproduce the old terminal
state. In general the dynamics step determines the new state. The terminal
set is chosen to satisfy the required invariance, collision-feasibility and
stability conditions.

\paragraph{Incorporating the measurement.}
A shifted multiple-shooting initialization sets
\begin{equation}\label{eq:shift-measured-initial-state}
 \bar x^0_{0|\hat t+1}=\hat x(\hat t+1),\qquad
 \bar q^0_{0|\hat t+1}=\hat q(\hat t+1)
\end{equation}
and keeps the other shifted and extended state guesses. Initial feature
values may be recomputed at every node:
\[
 \bar z^0_{k|\hat t+1}
 =h^z_{k|\hat t+1}(\bar x^0_{k|\hat t+1},\bar q^0_{k|\hat t+1}),
 \qquad k=0,\ldots,H.
\]
Model--state mismatch generally produces shooting defects, even if the
old prediction was feasible. In a full-space SQP these are legitimate
initial residuals and are retained in the linearized equality equations;
the corrector acts on their measured values.

Alternatively, a measurement-consistent initialization rolls out the
shifted controls from \eqref{eq:shift-measured-initial-state}. For
$k=0,\ldots,H-2$, propagate $(\bar x^0_k,\bar q^0_k)$ with the copied
control $(\widetilde u_k,\widetilde v_k)$ and the new maps
$(f^x_{k|\hat t+1},f^q_{k|\hat t+1})$. Evaluate $\kappa_{\rm ext}$ at
the resulting node $H-1$ and propagate once more to obtain node $H$.
Set the feature guesses using $h^z$ as above. This gives an initial trajectory consistent with the nominal shooting
equations up to integration error. The path and terminal constraints are
retained for subsequent correction. Both initializers retain future states and $z$ variables as
independent multiple-shooting NLP variables during subsequent SQP iterations;
node zero is substituted from the measurement. Future path constraints begin
at node one in either initialization.

\paragraph{Shrinking horizon, shared duration, and warm-start metadata.}
For a shrinking horizon, $H_{\hat t+1}=H_{\hat t}-1$. Retain the old
state nodes $1,\ldots,H_{\hat t}$ and old controls
$1,\ldots,H_{\hat t}-1$, reindexing both from zero. The old terminal
state is already the new terminal state, so no additional terminal
control or state is appended. When the new horizon is zero, the problem
is terminal-only and contains no control. Measurement replacement or
nominal rollout is performed on the retained layout in the same manner.

A shared remaining-duration variable is also shifted consistently with
physical time. For a shrinking fixed-arrival plan, the unadjusted guess is
$T^0_{\hat t+1}=\bar T_{\hat t}-h^{\rm applied}_{\hat t}$, subject to
the new duration bounds. A fixed-node-count extension instead also adds
the duration of the appended interval. If the new mesh does not coincide
with the shifted old mesh, the trajectory must be resampled using a
manifold-consistent interpolation or rolled out again. Copied nodes retain their physical time labels.

The shift map $\mathscr S_{\hat t}$ includes these stated choices.
Multipliers and working-set entries are transferred only for retained
constraint identities, that is, for rows that exist in both problems as the
same constraint at the same physical node; a row is matched by a stable key
(its kind and its node), and a row without a match starts with a zero
multiplier. Geometric equality multipliers are
reexpressed in the new residual dual frames. Newly introduced terminal
and event rows enter the KKT-consistent auxiliary initialization of
Section~\ref{sec:pathfollowing}. Each independent working set then
determines its complete multiplier through its own KKT system. The
measured node remains substituted and has no operational multiplier.
Previously active bounds are reevaluated against the current problem.
Contact-mode, catalogue, or horizon changes require a new model version.
Assumption~\ref{as:sosc} applies on each fixed-layout parameter interval after this transfer, with its own decision space. To restart a local estimate
after a layout change, the transferred guess must again satisfy the
proximity requirement \eqref{eq:path-initial-distance} for the new family.

\paragraph{Real-time iteration (RTI).}
In regular RTI operation the solver makes one outer SQP correction per
measurement
\cite{diehl2005real}. Applied here, it builds the local models at the
shifted estimate, solves the admitted working-set QP homotopies, and
retracts and reevaluates their outputs on the current measured OCP.
With $j$ denoting the outer correction index, the geometric update is
\begin{equation}\label{eq:priemnlpcan}
 \bar r_{\hat t,j,w}\in\mathcal M_H,
 \qquad p\in T_{\bar r_{\hat t,j,w}}\mathcal M_H,
 \qquad r^{\rm trial}=R_{\bar r_{\hat t,j,w}}(p).
\end{equation}
A successful regular RTI round uses $j=0$ and applies only the first
control of the accepted trajectory. Multiple blocking events and
parallel working sets contribute inner QP work within each nonlinear SQP
round.

RTI computations may be divided into a preparation phase and a feedback
phase. Derivative evaluation, structural assembly, and admissible
factor preparation can be performed before the next measurement;
the feedback phase incorporates the measured-state dependence of the first shooting residual and
any required parameter-dependent updates, solves the prepared QP, and
returns the first control. Prepared matrices are reused only where their
dependence on the arriving data has been accounted for. In particular,
a prepared future-trajectory base may remain fixed while node zero is
substituted from the measurement and its influence on the first shooting
equation is updated. The corresponding explicit parameter derivatives are
included. If any free base block is reset as well, the affected model
quantities are reevaluated. 

\paragraph{Advanced-step NMPC.}
Advanced-step NMPC uses the interval before the next measurement to solve
a problem initialized at a \emph{predicted} next state
\cite{zavala2009advanced}. Given the measured current state and the control
actually applied, compute
\begin{equation}\label{eq:advanced-predicted-state}
\begin{aligned}
 \check x_{\hat t+1}
 &=f^x_{0|\hat t}
    (\hat x(\hat t),\hat q(\hat t),u^{\rm applied}_{\hat t},
       v^{\rm applied}_{\hat t}),\\
 \check q_{\hat t+1}
 &=f^q_{0|\hat t}
    (\hat q(\hat t),\hat x(\hat t),u^{\rm applied}_{\hat t},
       v^{\rm applied}_{\hat t}).
\end{aligned}
\end{equation}
The integration interval is the actual applied one. Using the shifted
and extended trajectory as a starting guess, solve the next OCP with
initial state $(\check x_{\hat t+1},\check q_{\hat t+1})$ and forecast
data, obtaining a prepared primal--dual estimate
$(r^{\rm prep},\zeta^{\rm prep})$. The classical advanced-step approach
seeks an accurately solved predicted OCP, whereas RTI deliberately
limits the number of outer corrections.

When the measurement arrives, form its local geometric discrepancy
\[
 \Delta\xi=
 \left(\hat x(\hat t+1)-\check x_{\hat t+1},\,
 R_{\check q_{\hat t+1}}^{-1}(\hat q(\hat t+1))\right).
\]
The inverse retraction is used only locally to express a perturbation
of the initial-condition data. The decision states remain manifold
points. Additional forecast-data discrepancies are included similarly.

Advanced-step RTI combines a prescribed
number of preparation iterations on the predicted OCP, then a short
measurement-dependent feedback correction \cite{frey2024advancedstep}.
Both approaches therefore admit preparation before measurement. They
differ in the accuracy of the predicted nonlinear solve and the
construction of the arriving-data correction. 
\section{Parallel Working-Set Parametric Pathfollowing}
\label{sec:pathfollowing}

We consider the parametric quadratic program defining one Riemannian SQP
step. The method follows the parametric active-set construction of
qpOASES~\cite{ferreau2014qpoases} and maintains a finite collection of
independent working sets in parallel. Each set determines a unique
primal--dual face solution, that is, the solution of the QP in which the
constraints of the working set are imposed as equalities and all other
inequalities are ignored. The resulting tasks follow their own paths. When
the active constraints are degenerate, different tasks carry different
multiplier vectors for the same primal point. The algorithm selects one
complete endpoint at synchronization, after all tasks have terminated. Following
\cite{gill2021numerical}, the $LBL^\top$ factorization supplies the
Newton and sensitivity solves, constraint additions and deletions, and
rank and inertia tests. We give these factorization operations as part of
the pathfollowing algorithm.

\subsection{Riemannian multiple shooting and the sparse tangent model}
\label{subsec:shooting-assembly}
The online QP is the tangent model of the full multiple-shooting OCP. At a feedback update the nodes
are $k=0,\ldots,H$, and controls are attached to the $H$ intervals
$k=0,\ldots,H-1$. Aggregate the dynamic states as $s_k=(x_k,q_k)\in\mathcal S
=\mathbb R^{d_x}\times\mathcal M_q$ and the actions as $a_k=(u_k,v_k)\in\mathcal A
=\mathcal M_u\times\mathbb R^{d_v}$. The initial state $s_0=\hat s$ is
measured data and is substituted before assembling the working KKT matrix.
Its associated fixed feature values are substituted as well. Thus there is
no node-zero state increment or initial-state multiplier in the operational
system. State, feature, collision and phase constraints are assembled at
future nodes $1,\ldots,H$; control constraints are assembled on intervals
$0,\ldots,H-1$. Horizon-wide parameters
such as a free remaining duration are collected in $\theta$ and stored once.

Choose orthonormal tangent frames $E^s_k$ and $E^a_k$ at the current state
and control nodes, and keep their versions fixed throughout this QP.
For each geometric shooting residual choose an anchor $\bar q_{k+1}$ and
an orthonormal residual frame $\bar E_{k+1}$ at that anchor. On the declared
normal neighborhood define the anchored, vector-valued map
\begin{equation}\label{eq:shooting-anchored-map}
 \Theta_{k+1}(a,b)=\bar E_{k+1}^{*}
       \mathcal P_{a\to\bar q_{k+1}}\operatorname{Log}_{a}(b),
 \qquad
 c_k=\begin{pmatrix}
 x_{k+1}-f^x_k(s_k,a_k;\theta)\\
  \Theta_{k+1}\bigl(q_{k+1},f^q_k(s_k,a_k;\theta)\bigr)
 \end{pmatrix}.
\end{equation}
The first row has the Euclidean sign used in
\eqref{eq:intrinsic-shooting}; the second has the geometric sign used in
\eqref{eq:intrinsic-residual-representative}. Both vanish at the same
shooting match. The anchor and its frame are held fixed when differentiating
and when comparing trial values in the current nonlinear model.

The exact geometric first derivative is
\begin{equation}\label{eq:shooting-exact-differential}
  Dc_k^q[\eta]=D_1\Theta_{k+1}[\eta_{q_{k+1}}]
       +D_2\Theta_{k+1}\!\left[D f^q_k[\eta_{s_k},\eta_{a_k},\eta_\theta]\right].
\end{equation}
Here $D_1\Theta$ includes the derivative of transport and of the base point
of the logarithm; $D_2\Theta$ differentiates its second point. At a matching
pair $(a,a)$ this reduces to
\begin{equation}\label{eq:shooting-match-differential}
 Dc_k^q[\eta]=\bar E_{k+1}^{*}\mathcal P_{a\to\bar q_{k+1}}
       \left(Df^q_k[\eta_{s_k},\eta_{a_k},\eta_\theta]
                       -\eta_{q_{k+1}}\right).
\end{equation}

The Hessian contribution is the covariant second derivative of the scalar
multiplier--residual pairing. With
$\Psi_k(r)=(q_{k+1},f^q_k(s_k,a_k;\theta))$ and
$\varphi_k=\langle\lambda_k^q,\Theta_{k+1}\rangle$, its exact chain rule is
\begin{equation}\label{eq:shooting-exact-second-differential}
 \nabla^2(\varphi_k\circ\Psi_k)[\eta,\xi]
 =\nabla^2\varphi_k[D\Psi_k\eta,D\Psi_k\xi]
       +D\varphi_k\!\left[(\nabla D\Psi_k)[\eta,\xi]\right].
\end{equation}
All connections in this formula are the product Levi--Civita connections
and the fixed vector-space connection on the anchored residual. The first
term includes mixed derivatives between $q_{k+1}$ and the predicted
configuration; the second includes the covariant second derivative of the
integrator. The endpoint and transport derivatives therefore enter the SQP
Hessian together with the dynamics derivatives. An equivalent smooth
Lie-group error may replace \eqref{eq:shooting-anchored-map} if both
\eqref{eq:shooting-exact-differential} and
\eqref{eq:shooting-exact-second-differential} are formed for that error.

Let $\delta s_k,\delta a_k$ be frame coordinates and set $\delta s_0=0$.
At the current model define
\begin{equation}\label{eq:shooting-stage-blocks}
 M_k=D_{s_{k+1}}c_k E^s_{k+1},\quad
 N_k=D_{s_k}c_k E^s_k,\quad
 U_k=D_{a_k}c_k E^a_k,\quad G_k=D_\theta c_k .
\end{equation}
The equality row for interval $k$ is exactly
\begin{equation}\label{eq:shooting-linearized-row}
 M_k\delta s_{k+1}+N_k\delta s_k
       +U_k\delta a_k+G_k\delta\theta=-c_k,
       \qquad k=0,\ldots,H-1.
\end{equation}
For $k=0$, the $N_k\delta s_k$ term is absent. A perturbation of the
measurement changes the target data in that first row and in the
stage-zero derivatives. Near a shooting match each $M_k$ is invertible, since its
Euclidean block is an identity and its geometric block is minus an
invertible transport/frame map. Invertibility is checked in the actual
off-match coordinates before using a condensed solve.

Feature rows $z_k-h_k^z(s_k;\theta)=0$ give
\begin{equation}\label{eq:shooting-feature-row}
 \delta z_k-D_{s_k}h_k^z E_k^s\delta s_k
                    -D_\theta h_k^z\delta\theta
       =-\bigl(z_k-h_k^z(s_k;\theta)\bigr),\qquad 1\leq k\leq H.
\end{equation}
Their $z_k$ block is an identity. A future-node smooth margin $b_{k,i}$
gives $Db_{k,i}E_rp\geq-b_{k,i}$; an interval control margin is treated
identically at its own interval. Constraints coupling several nodes retain
all their corresponding columns. The optional duration column is shared
by every dynamics, feature, cost and constraint block depending on duration.

In the ordering
$(\delta a_0,\delta s_1,\delta z_1,\delta a_1,\ldots,
\delta s_H,\delta z_H,\delta\theta)$,
\eqref{eq:shooting-linearized-row} has nonzeros only at adjacent states,
its interval control and the shared parameter. The feature row is local to
its node. The Lagrangian Hessian is assembled by adding, once per stable
row ID, the objective Hessian and the multiplier-weighted exact scalar
constraint Hessians, including \eqref{eq:shooting-exact-second-differential}.
For stage-local objectives and constraints this produces neighboring-node
blocks and a possible shared-parameter arrow; explicitly nonlocal event
constraints add their declared blocks. 
The uncondensed matrix retains the sparse multiple-shooting block structure
under the factor updates described below.

For completeness, dynamics condensing is an exact affine elimination of
\eqref{eq:shooting-linearized-row}: the state increments are expressed
through the control increments, which gives a smaller and denser QP. If every $M_k$ is certified invertible,
set
\begin{equation}\label{eq:shooting-condensing-recursion}
 \delta s_{k+1}=\mathsf A_k\delta s_k+
       \mathsf B_k\delta a_k+\mathsf G_k\delta\theta+d_k,
 \quad
 (\mathsf A_k,\mathsf B_k,\mathsf G_k,d_k)
       =-M_k^{-1}(N_k,U_k,G_k,c_k).
\end{equation}
Forward substitution from $\delta s_0=0$ defines $p=p_c+Vy$, where $y$
collects the retained controls, features and shared parameters. It yields
\begin{equation}\label{eq:shooting-condensed-data}
 H_c=V^\top HV,\quad g_c=V^\top(g+Hp_c),\quad
 A_c=AV,\quad b_c=b-Ap_c .
\end{equation}
The constant objective term is immaterial to the optimizer but retained if
objective values are compared. State bounds become general rows after
condensing. Selected
working bounds are interpreted in the coordinates actually
factored.
Partial condensing applies the same recursion within declared stage blocks.

The state increments are recovered by
\eqref{eq:shooting-condensing-recursion}. The dynamics multipliers are
recovered from the state part of full stationarity:
\begin{equation}\label{eq:shooting-dual-recovery}
 (A_{\mathrm{dyn},s})^\top\lambda_{\mathrm{dyn}}
       =-\bigl(Hp+g+A_{\mathrm{other}}^\top\zeta_{\mathrm{other}}\bigr)_s.
\end{equation}
The matrix on the left is block triangular with diagonal blocks $M_k^\top$;
backward substitution supplies its unique solution. The complete uncondensed
KKT residual is then checked. If $m_{\mathrm{dyn}}$ dynamics rows were
eliminated, exact congruence gives
\begin{equation}\label{eq:shooting-condensed-inertia}
 \operatorname{In}(K_{\mathrm{uncond}})
   =(m_{\mathrm{dyn}},m_{\mathrm{dyn}},0)
                 +\operatorname{In}(K_{\mathrm{cond}}).
\end{equation}
This identity requires that the retained selected rows remain independent
after dynamics elimination. The implementation may use uncondensed, partially
condensed or fully condensed matrices; its matrix key records that choice,
and its final certificate always expands back to the same OCP rows.

\subsection{Frozen tangent-space QP and data layout}
\label{subsec:target-qp}
Let $r\in\mathcal M_H$ be the current nonlinear trajectory and let
$E_r:\mathbb R^{n}\to T_r\mathcal M_H$ be the declared orthonormal product
tangent frame after structural substitutions. Use a fixed ordering of the
reduced tangent coordinates and a fixed ordering of every constraint ID. The
local QP is generically written as
\begin{equation}\label{eq:target-qp}
\begin{aligned}
  \min_{p\in\mathbb R^n}\quad &(g^1)^\top p+\tfrac12p^\top Hp,\\
 A_Ep&=b_E^1,\\
 A_Ip&\geq b_I^1,\\
 \ell^1&\leq S_b p\leq u^1.
\end{aligned}
\end{equation}
Here $g^1=E_r^*\operatorname{grad}f(r)$ and $H$ is the symmetric matrix of
the chosen Riemannian Lagrangian Hessian on the same frame. The superscript
$1$ marks target data, which the homotopy reaches at $\tau=1$. Let $b_{\rm gen}$ collect the inequality margins that are represented
as general rows: the operational inequalities $D$ and the geometric simple
margins. The assembled data are
\begin{equation}\label{eq:ocp-target-derivatives}
\begin{aligned}
 A_E&=DC(r,\hat t)E_r,& b_E^1&=-C(r,\hat t),\\
 A_I&=Db_{\rm gen}(r,\hat t)E_r,& b_I^1&=-b_{\rm gen}(r,\hat t),\\
 H&=E_r^*\operatorname{Hess}_rL(r,\zeta,\hat t)E_r .
\end{aligned}
\end{equation}
Component bounds are kept separately from $A_I$; geometric margins enter
$A_I$. For a Euclidean coordinate $y_i$, its target increment bounds are
$\ell_i^1=\underline y_i-y_{i,v}$ and $u_i^1=\overline y_i-y_{i,v}$.
The theorem uses the exact Hessian; a declared symmetric approximation may
be used with its measured error included in the derivative budget. Any
prescribed affine structural substitution is applied to the objective and
row right sides as in \eqref{eq:bound-affine-lift} before these data are
stored. Inequality rows use positive-feasible margins: each inequality is written so
that a nonnegative value means feasibility. The QP Lagrangian is
\begin{equation}\label{eq:qp-lagrangian}
\begin{aligned}
 L_Q(p,\lambda,\mu,\nu^\ell,\nu^u)
 ={}&(g^1)^\top p+\tfrac12p^\top Hp
 +\lambda^\top(A_Ep-b_E^1)+\mu^\top(A_Ip-b_I^1)\\
 &+(\nu^\ell)^\top(S_b p-\ell^1)+(\nu^u)^\top(u^1-S_b p),
\end{aligned}
\end{equation}
with $\mu,\nu^\ell,\nu^u\leq0$. Here $S_b$ selects the bounded
coordinates; the bounds are indexed by their original coordinate IDs.
Equality multipliers are free. The full dual arrays are indexed by the
original constraint IDs, and workers store selectors into those arrays.
These fixed IDs identify the multiplier components through an exchange or
a resumed task.

The bound catalogue is partitioned before any factorization. Coincident
bounds $\ell_i(\tau)=u_i(\tau)$ on the full path are structural coordinates
and are eliminated by substituting $p_i=c_i(\tau)$ into the objective and
all rows. On the remaining
(nonstructural) coordinate set $\mathcal I_{\rm ns}$, distinguish activity from working status:
\begin{equation}\label{eq:bound-catalogues}
\begin{aligned}
 \mathcal A_\ell&=\{i:p_i=\ell^1_i\},&
 \mathcal A_u&=\{i:p_i=u^1_i\},\\
 W_\ell&\subseteq\mathcal A_\ell,&
 W_u&\subseteq\mathcal A_u,\\
 N&=W_\ell\cup W_u,&
 F&=\mathcal I_{\rm ns}\setminus N .
\end{aligned}
\end{equation}
Only working bounds are fixed in the reduced KKT system. An active bound
omitted from $W_\ell$ or $W_u$ remains a free coordinate
and is monitored as an ordinary primal blocker. A nonstructural coordinate
is working at at most one end. Basic/superbasic column labels are assigned
only after this reduction: they are a pivot-column partition of the free
constraint Jacobian. Activity and working-set membership are recorded
separately.

\subsection{Task hypotheses and KKT-consistent auxiliary start}
\label{subsec:qp-start}
At feedback sample $\hat t$, let $w$ identify an initial active-set hypothesis
$\mathcal W_{\hat t,w}$. Denote the collection of parallel working sets by
$\mathbb W_{\hat t}=\{\mathcal W_{\hat t,w}:w\in\mathcal Q_{\hat t}\}$, where
$\mathcal Q_{\hat t}$ is the finite index set of the tasks at this sample.
Each hypothesis contains selected general inequality IDs and selected
geometric-bound IDs after the structural coordinates have been substituted.
Several hypotheses can describe the same primal point with
different valid multipliers. Such hypotheses are distinct QP tasks until
their complete checkpoint data are compared.

A QP warm start is accepted only if the supplied $(p^0,\zeta^0)$ satisfies
its stored start-QP KKT equations, primal feasibility, all dual signs,
complementarity, and the working-row rank test. Otherwise construct a start
QP that has the chosen checkpoint as an exact KKT point. The start QP
differs from the target QP only in its gradient and its right-hand sides;
the Hessian and all constraint normals are those of the target QP. Write the target
general rows as $A_Ep=b_E^1$ and $A_Ip\geq b_I^1$, and bounds as
$\ell^1\leq S_b p\leq u^1$. Select a primal vector $p^0$ and duals
$\lambda^0,\mu^0,\nu^{\ell,0},\nu^{u,0}$ with inequality and bound duals
nonpositive and supported on the proposed working set. Define
\begin{equation}\label{eq:start-qp-data}
\begin{aligned}
 b_E^0&=A_Ep^0,\\
 b_i^0&=\begin{cases}a_i^\top p^0,&i\in W_I^0,\\
        \min\{b_i^1,a_i^\top p^0\},&i\notin W_I^0,\end{cases}\\
 g^0&=-Hp^0-A_E^\top\lambda^0-A_I^\top\mu^0
       -S_b^\top(\nu^{\ell,0}-\nu^{u,0}),\\
 \ell_i^0&=\begin{cases}p_i^0,&i\in W_\ell^0,\\
                  \min\{\ell_i^1,p_i^0\},&i\notin W_\ell^0,\end{cases}\\
 u_i^0&=\begin{cases}p_i^0,&i\in W_u^0,\\
                  \max\{u_i^1,p_i^0\},&i\notin W_u^0.\end{cases}
\end{aligned}
\end{equation}
$S_b$ selects the bounded coordinates. For a lower side its normal is
$+e_i^\top$ and its multiplier is $\nu_i^\ell\leq0$; for an upper side the
oriented normal is $-e_i^\top$ and multiplier is $\nu_i^u\leq0$.  Direct substitution verifies zero start slack on
each working row/side, feasibility of every inactive row/side, dual signs,
and zero stationarity residual. 

The path data are affine in one artificial parameter $\tau\in[0,1]$:
\begin{equation}\label{eq:qp-path-data}
\begin{aligned}
  g(\tau)&=(1-\tau)g^0+\tau g^1,\\
 b_E(\tau)&=(1-\tau)b_E^0+\tau b_E^1,\\
 b_I(\tau)&=(1-\tau)b_I^0+\tau b_I^1,\\
 \ell(\tau)&=(1-\tau)\ell^0+\tau\ell^1,\\
 u(\tau)&=(1-\tau)u^0+\tau u^1.
\end{aligned}
\end{equation}
The Hessian $H$ and every row normal are frozen during this homotopy. Hence
target data at $\tau=1$ equal those of \eqref{eq:target-qp}. Because only
the linear term and the right-hand sides depend on $\tau$, and they do so
affinely, the solution on a fixed working set is an affine function of
$\tau$. This is what makes the path exactly computable. 

\subsection{Working selectors and reduced KKT equations}
\label{subsec:reduced-kkt}
For one operational task at homotopy state $(p,\zeta,\tau)$, let $W_E$
contain all retained equality rows and $W_I$ the independent selected
general inequality rows. Let $N$ contain working-bound coordinates and
$F$ the remaining free coordinates. With $p_N$ prescribed by its working
bounds, distinguish the full selected rows and their unreduced right side
from the free-coordinate rows:
\begin{equation}\label{eq:working-kkt}
\begin{aligned}
 B_W&=\begin{bmatrix}A_E\\A_{I,W_I}\end{bmatrix},&
 \beta_W(\tau)&=\begin{bmatrix}b_E(\tau)\\b_{I,W_I}(\tau)\end{bmatrix},\\
 A_W&=B_{W,F},&
 \widetilde\beta_W(\tau)&=\beta_W(\tau)-B_{W,N}p_N(\tau),\qquad
 g_F(\tau)=(g(\tau))_F .
\end{aligned}
\end{equation}
In particular, $g_F$ denotes only the free entries of the objective's linear
coefficient. The selected face equations
are
\begin{equation}\label{eq:face-kkt}
 F_W(p_F,\zeta_W,\tau)=
 \begin{pmatrix}
   H_{FF}p_F+H_{FN}p_N+g_F(\tau)+A_W^\top\zeta_W\\
   A_Wp_F-\widetilde\beta_W(\tau)
 \end{pmatrix}=0.
\end{equation}
Independence of the selected rows and positive definiteness of $H_{FF}$
on $\ker A_W$ make this bordered system nonsingular. Consequently the
working set, its bound values and the current path data determine one
face solution and one selected multiplier vector. In particular, if
$A_W^\top=Q_WR_W$ is the thin QR factorization, with $Q_W^\top Q_W=I_m$
and nonsingular triangular $R_W$, then
\begin{equation}\label{eq:working-dual-unique}
 R_W\zeta_W=-Q_W^\top(Hp+g(\tau))_F,\qquad
 \mu_i=0\ (i\notin W_I),\quad
 \nu_i^\ell=0\ (i\notin W_\ell),\quad
 \nu_i^u=0\ (i\notin W_u).
\end{equation}
Working-bound multipliers are recovered uniquely from
\eqref{eq:reduced-cost-map}; any condensed dynamics or geometric rows
are recovered by \eqref{eq:shooting-dual-recovery} or
\eqref{eq:geometric-reduced-cost}. Equation~\eqref{eq:working-dual-unique}
is a triangular solve and a zero extension by stable row ID. If $m=0$, the triangular solve is empty. With different independent working sets, the system may produce
different full multiplier arrays at the same primal point. Those alternatives
are the parallel tasks of this section.

The full residuals, before bound elimination, are
\begin{equation}\label{eq:full-face-residuals}
 r^d=Hp+g(\tau)+B_W^\top\zeta_W,\qquad
 r_W^p=B_Wp-\beta_W(\tau).
\end{equation}
The missing stationarity components at working coordinates are balanced by
their signed bound multipliers. With $d_N=p_N^{\rm prescribed}-p_N$, a
fixed-parameter correction of an inconsistent checkpoint solves
\begin{equation}\label{eq:qp-newton-correction}
 \begin{pmatrix}H_{FF}&A_W^\top\\A_W&0\end{pmatrix}
 \binom{d_F}{d_{\zeta_W}}
 =-\binom{r_F^d+H_{FN}d_N}{r_W^p+B_{W,N}d_N}.
\end{equation}
Here $r_W^p=A_Wp_F-\widetilde\beta_W$ whenever $\widetilde\beta_W$ is
formed using the current $p_N$. The correction then sets
$p^+=p+d$ and $\zeta_W^+=\zeta_W+d_{\zeta_W}$, recovers the bound
multipliers from \eqref{eq:reduced-cost-map}, and checks the full residual
\eqref{eq:full-face-residuals}. A sensitivity right side differentiates the
path data; a correction right side cancels the current residual.

When the free constraint matrix has $m$ independent rows, a column ordering
$P_F$ can expose an invertible basic block:
\begin{equation}\label{eq:basic-superbasic}
\begin{aligned}
 A_WP_F&=[A_B\ A_S],\qquad \det A_B\ne0,\\
 d_F&=P_F\binom{A_B^{-1}(r_A-A_Sd_S)}{d_S},\qquad
 Z_F=P_F\binom{-A_B^{-1}A_S}{I},\quad A_WZ_F=0.
\end{aligned}
\end{equation}
The vector $d_S$ parameterizes homogeneous null motion; $r_A$ supplies the
particular inhomogeneous correction. A deterministic rank-revealing pivot
rule selects the basic columns and records $P_F$. If no full-row-rank basic
block is found, the row-dependency handler in
Section~\ref{subsec:dependent-exchange} must resolve the representation
before an ordinary reduced solve. A dense KKT implementation can avoid
forming $Z_F$ explicitly and solve \eqref{eq:face-kkt} directly.

\subsection{Sensitivity, reduced costs, and exact event ratios}
\label{subsec:path-events}
On an interval with fixed $W_I,N,F$ and factorization coordinates, differentiate
\eqref{eq:face-kkt} with respect to $\tau$. A dot is derivative per unit
artificial $\tau$. Working bounds prescribe
$\dot p_N=\dot\ell_N$ or $\dot u_N$; structural coordinates use their
prescribed affine derivative. The predictor system is
\begin{equation}\label{eq:path-sensitivity}
 \underbrace{\begin{pmatrix}H_{FF}&A_W^\top\\
                              A_W&0\end{pmatrix}}_{K_W}
 \binom{\dot p_F}{\dot\zeta_W}
 =\binom{-\dot g_F-H_{FN}\dot p_N}
        {\dot\beta_W-B_{W,N}\dot p_N}.
\end{equation}
The lower right side is $\dot{\widetilde\beta}_W$; it contains the
working-bound velocity exactly once. The upper right side uses the explicit
data derivative $(\dot g)_F$, not the derivative of $r_F^d$.
Recover the complete primal derivative and then every reduced-gradient
derivative from
\begin{equation}\label{eq:reduced-gradient-sensitivity}
 \dot p=(\dot p_F,\dot p_N),\qquad
 \bar g=g(\tau)+Hp+B_W^\top\zeta_W,\qquad
 \dot{\bar g}=\dot g+H\dot p+B_W^\top\dot\zeta_W.
\end{equation}
With the plus-Lagrangian sign convention, a working general inequality
multiplier and bound multipliers satisfy $\mu_i,\nu_i^\ell,\nu_i^u\leq0$.
The coordinate reduced costs satisfy
\begin{equation}\label{eq:reduced-cost-signs}
 \bar g_i\geq0\quad(i\in W_\ell),\qquad
 \bar g_i\leq0\quad(i\in W_u),\qquad \bar g_F=0.
\end{equation}
More exactly, $\nu_i^\ell=-\bar g_i$ and $\nu_i^u=\bar g_i$ in the
full stationarity representation. General inequality rows use their signed
QP multipliers directly. Coupling between nonorthogonal geometric normals
is accounted for by the QR triangular solve in
Section~\ref{subsec:geometric-rows}.

For each inactive general row $a_i^\top p\geq b_i(\tau)$ define
\begin{equation}\label{eq:general-block-ratio}
 s_i=a_i^\top p-b_i(\tau),\qquad \dot s_i=a_i^\top\dot p-\dot b_i.
\end{equation}
If $\dot s_i<0$, the first primal blocking ratio is
$\alpha_i^p=-s_i/\dot s_i$. For an inactive lower and upper bound use,
respectively,
\begin{equation}\label{eq:bound-ratios}
 (s_i^\ell,\dot s_i^\ell)=(p_i-\ell_i,\dot p_i-\dot\ell_i),\qquad
 (s_i^u,\dot s_i^u)=(u_i-p_i,\dot u_i-\dot p_i),
\end{equation}
with a ratio $-s/\dot s$ only when $\dot s<0$. A working general inequality with $\dot\mu_i>0$
leaves when its multiplier reaches zero, at the step
$\alpha_i^d=-\mu_i/\dot\mu_i$. A working lower bound with
$\dot{\bar g}_i<0$, or a working upper bound with $\dot{\bar g}_i>0$,
leaves when its reduced cost $\bar g_i$ reaches zero, at the corresponding
zero-crossing ratio $-\bar g_i/\dot{\bar g}_i$. Equality multipliers have no sign event.

For every event candidate retain its type, stable constraint ID, ratio, and
one-sided derivative. A negative slack or a multiplier of the wrong sign at the current point
means that the checkpoint has failed. Zero ratios are
valid tied events. The step length is
\begin{equation}\label{eq:min-event-step}
 \Delta\tau=\min\left(1-\tau,
     \{\alpha_i^p\}_{i\in\mathcal I_p},
     \{\alpha_i^d\}_{i\in\mathcal I_d},
     \{\alpha_i^{\rm bound}\}_{i\in\mathcal I_b},
     \{\alpha_i^{\rm rc}\}_{i\in\mathcal I_{\rm rc}}\right),
\end{equation}
where empty candidate sets are omitted. Advance exactly once:
\begin{equation}\label{eq:advance-checkpoint}
 p^+=p+\Delta\tau\dot p,\qquad
 \zeta_W^+=\zeta_W+\Delta\tau\dot\zeta_W,\qquad
 \bar g^+=\bar g+\Delta\tau\dot{\bar g},\qquad
 \tau^+=\tau+\Delta\tau.
\end{equation}
At this checkpoint recompute all slacks, signs, residuals and
zero-distance ratios before the next pivot. These operations define the
event transition for the frozen affine-QP path.

\subsection{Relation to the qpOASES homotopy and path invariants}
\label{subsec:qpoases-homotopy}
The qpOASES homotopy moves the gradient and constraint offsets along an
affine segment and leaves the current face when a primal slack or dual
multiplier reaches zero~\cite{ferreau2014qpoases}. We apply this construction
to each parallel task in its fixed tangent frame. Write the fixed-Hessian
QP at parameter $\tau$ as
\begin{equation}\label{eq:parametric-qpoases-standard}
 \min_p\;\tfrac12p^\top Hp+g(\tau)^\top p,\qquad
 A_Ep=b_E(\tau),\quad A_Ip\geq b_I(\tau),\quad
 \ell(\tau)\leq S_b p\leq u(\tau).
\end{equation}
Let $\mathcal W$ contain all equality rows and the selected inequality and oriented bound rows, so that $A_\mathcal Wp=b_\mathcal W(\tau)$. The equality-constrained face problem is characterized by
\begin{equation}\label{eq:parametric-face-kkt-full}
 \mathcal F_\mathcal W(z,\tau)=
 \begin{pmatrix}Hp+g(\tau)+A_\mathcal W^\top\zeta_\mathcal W\\
                   A_\mathcal Wp-b_\mathcal W(\tau)\end{pmatrix}=0,
 \qquad z=(p,\zeta_\mathcal W).
\end{equation}
On an interval where $\mathcal W$ is fixed, differentiation gives the constant system
\begin{equation}\label{eq:parametric-face-slope}
 \begin{pmatrix}H&A_\mathcal W^\top\\A_\mathcal W&0\end{pmatrix}
 \begin{pmatrix}\dot p\\\dot\zeta_\mathcal W\end{pmatrix}
 =\begin{pmatrix}-\dot g\\\dot b_\mathcal W\end{pmatrix}.
\end{equation}
Consequently $z(\tau+\alpha)=z(\tau)+\alpha\dot z$ exactly until a face event; there is no ODE integrator and no step-size truncation error inside a frozen QP. Equations~\eqref{eq:path-sensitivity}--\eqref{eq:advance-checkpoint} are the same system after structural and working coordinates have been eliminated. Expanding a reduced solve must satisfy \eqref{eq:parametric-face-slope} to the declared residual tolerance; this is a direct implementation consistency check.

The face solution is a certified strict minimizer of its equality-constrained face problem when its projected Hessian is positive definite. Here $H$ and the face rows are constant during one QP homotopy, so $H_\mathcal W=Z_\mathcal W^\top HZ_\mathcal W$ is constant on the face. It is certified when the working set is factored; it cannot change continuously with $\tau$. A face failing the selected SSOSC/inertia filter is marked inadmissible and handled by another task or the backup.

Each affine segment of the parametric path ends at an event of one of three classes. A primal event in the form of a blocking constraint (previously strictly feasible inequality constraint becoming active) adds a row with zero multiplier. A dual event removes a selected inequality whose multiplier reaches zero. A representation event occurs when the incoming primal row is dependent on the current working rows and the multiplier can be redistributed; it triggers the exchange mechanism in Section~\ref{subsec:dependent-exchange}. For independent rows the first two event types are the usual qpOASES blocking and leaving constraints. For dependent rows, simply adding the incoming row would make $K_\mathcal W$ singular, so a dual-compatible exchange is required.

At every accepted event, the procedure checks that the previously selected face equations remain satisfied,
\begin{equation}\label{eq:event-primal-invariance}
 A_\mathcal Wp^+-b_\mathcal W(\tau^+)=0,
\end{equation}
and recomputes stationarity with the new representation,
\begin{equation}\label{eq:event-stationarity-invariance}
 r_d^+=Hp^++g(\tau^+)+A_{\mathcal W^+}^\top\zeta_{\mathcal W^+}^+,
 \qquad \|r_d^+\|_{\rm scaled}\leq\epsilon_d.
\end{equation}
Every inequality multiplier must retain its sign, every inactive margin must be nonnegative, and every omitted active row must be retested for a forward blocker. 

\paragraph{Simultaneous working-set changes.}
The next working-set change is determined by a ratio test.
Let $s_i\geq0$ denote an inactive constraint's slack and let
$\pi_j:=-\mu_j\geq0$ be the magnitude of a working-set inequality
multiplier, so that both quantities are nonnegative; working bounds use
their oriented multipliers in the same way. Dots denote derivatives
with respect to the homotopy parameter $\tau$.
On the current affine homotopy segment, the candidate step lengths are
\[
 \Delta\tau_i^{\mathrm p}
 =\frac{s_i}{-\dot s_i}
 \quad(\dot s_i<0),
 \qquad
 \Delta\tau_j^{\mathrm d}
 =\frac{\pi_j}{-\dot\pi_j}
 \quad(\dot\pi_j<0).
\]
These are the ratios $\alpha_i^p$ and $\alpha_j^d$ of
Section~\ref{subsec:path-events}. The first is the parameter increment at which an inactive
constraint becomes active; the second is the increment at which
a working-set multiplier reaches zero.

Advance by the smallest nonnegative candidate step, capped by
the remaining homotopy interval. Candidates within $\epsilon_\tau$
of the minimum are treated as possible simultaneous events;
recheck their event conditions at the new point.
At this fixed parameter, first process dual releases: remove an
inequality whose multiplier magnitude $\pi_j$ is zero and whose forward derivative
would make it negative, that is, would give the multiplier the wrong sign.
Retain a zero multiplier whose forward
derivative is nonnegative. Then process primal blockers through
the row additions of Section~\ref{subsec:add-working-row} or the exchanges of
Section~\ref{subsec:dependent-exchange}.
Within each event type, use increasing stable constraint IDs,
which remain unchanged under row permutations.

After every pivot or exchange, recompute the slacks, working-set
multipliers, their forward derivatives, and the candidate step
lengths before selecting the next event. Continue until no
zero-distance event remains, then resume forward motion.
A zero-distance event is processed at the current parameter
without introducing an artificial positive step.

If a checkpoint $(\mathcal W,\tau,p,\zeta)$ repeats, using the
checkpoint merge tolerance for comparison, return
\texttt{EVENT\_CYCLE}.
\paragraph{Start and endpoint residuals.}
At $\tau=0$, verify the primal--dual residual over all original rows and
bound sides. For fixed positive scale factors $s_p,s_d,s_c$ computed from the QP data, define
\begin{equation}\label{eq:qp-scaled-residual}
 \begin{aligned}
 r_p&=\max\{\|A_Ep-b_E\|_\infty,\|[b_I-A_Ip]_+\|_\infty,
                  \|[\ell-S_bp]_+\|_\infty,\|[S_bp-u]_+\|_\infty\},\\
 r_d&=\|Hp+g+A_E^\top\lambda+A_I^\top\mu+S_b^\top(\nu^\ell-\nu^u)\|_\infty,\\
 r_c&=\max\{\|\mu\odot(A_Ip-b_I)\|_\infty,
 \|\nu^\ell\odot(S_bp-\ell)\|_\infty,\|\nu^u\odot(u-S_bp)\|_\infty,\\
 &\qquad \|[\mu]_+\|_\infty,\|[\nu^\ell]_+\|_\infty,\|[\nu^u]_+\|_\infty\},\\
 r_{\rm KKT}&=\max\{r_p/s_p,r_d/s_d,r_c/s_c\}.
 \end{aligned}
\end{equation}
Separate lower and upper dual arrays are mapped to their oriented rows before evaluating stationarity. At $\tau=1$, return \texttt{TARGET\_QP} only when $r_{\rm KKT}\leq\epsilon_{\rm QP}$ and the final factor/inertia tests pass. Otherwise return the last checkpoint and the failed residual component.

\subsection{Factor data structure, pivot selection, and reproducible solves}
\label{subsec:factor-data-structure}
A stored factorization is valid only for the exact matrix from which it
was computed. To prevent a factor from being applied to a different matrix
of the same size, each factor is versioned by the exact key
\begin{equation}\label{eq:factor-key}
 \mathsf K=(\mathsf{id}_{\rm model},\mathsf{id}_{\rm frame},
             \mathsf{order}(F),\mathsf{condensing},\mathsf{IDs}(W_E,W_I,W_\ell,W_u),
             \mathsf{scaling},\mathsf{pivot\_policy}).
\end{equation}
The stored object contains $(\Pi,L,\{B_j\}_{j=1}^q,R,\mathsf{key})$, where $R=0$ for a complete factor and otherwise is the exact trailing Schur block together with its coupling to the factored variables (written $S_*$ and $L_{Rc}$ in Section~\ref{subsec:factor-completion}). It also records block sizes, signed pivot classifications, scaling vectors, reconstruction residual and factor version. A factor is applied only when every component of its key matches the current matrix. In particular, a frame change, bound substitution, or constraint-ID change requires a corresponding factor update even when the dimensions agree.

A symmetric elimination step on a remainder $S=\left(\begin{smallmatrix}S_{11}&S_{12}\\S_{12}^\top&S_{22}\end{smallmatrix}\right)$ with nonsingular pivot block $D=S_{11}$ stores
\begin{equation}\label{eq:lblt-schur-step}
 L_{21}=S_{12}^\top D^{-1},\qquad
 S_{\rm next}=S_{22}-L_{21}DL_{21}^\top
              =S_{22}-S_{12}^\top D^{-1}S_{12}.
\end{equation}
The step uses symmetric row-and-column permutations; a row permutation alone invalidates the factor. Pivot selection uses the declared thresholds $\delta_{\rm abs},\delta_{\rm rel}$ and a deterministic tie rule on original variable/row IDs. Accept a scalar pivot $d$ only when $|d|>\delta_{\rm abs}+\delta_{\rm rel}\|S\|_\infty$. Otherwise test the selected $2\times2$ block by its signed eigenvalues and scaled determinant. If neither is certified, stop and retain the unreduced Schur block as $R$; never divide by an uncertified pivot.

Diagonal equilibration may be used before factorization. For positive diagonal $D_s$, factor $D_sK_WD_s$ and transform the right side by $D_s$; recover $x=D_sy_s$ from
\begin{equation}\label{eq:scaled-kkt-solve}
 (D_sK_WD_s)y_s=D_sr,\qquad x=D_sy_s.
\end{equation}
The congruence preserves inertia. Its scaling and ordering are recorded in
\eqref{eq:factor-key}, and the matrix and right side are transformed
together. The final solve test uses the unscaled residual in
\eqref{eq:lblt-backward-error}.

For multiple right sides, reuse the triangular factors and solve each column independently. If $\mathrm{RHS}=[r^{(1)},\ldots,r^{(k)}]$, compute
\begin{equation}\label{eq:multi-rhs-solve}
 Y=L^{-1}\Pi^\top\,\mathrm{RHS},\quad Z=B^{-1}Y,\quad X=\Pi L^{-\top}Z,
\end{equation}
and evaluate each column's backward error. A sensitivity right-hand side, Newton correction and event-border right-hand side carry different semantic labels and residuals. Reusing one solution as another is valid only for an identical right-hand side or after recomputing its residual.

An update is attempted only when the old factor is complete, its key matches, the candidate matrix differs by exactly the specified border or coordinate insertion/deletion, and its old reconstruction residual is at most $\epsilon_{\rm fac}$. Accept the result only if the new reconstruction residual is at most $\epsilon_{\rm fac}$, the bordered solve residual is at most $\epsilon_{\rm solve}$, and the expected inertia matches. If any condition fails, perform one full symmetric refactorization of the new matrix and log \texttt{REFACTORIZED}. If refactorization fails, return a non-success status; no approximate factor continues to a later event.

\subsection{\texorpdfstring{$LBL^\top$}{LBLT} KKT factorization and solve operations}
\label{subsec:inertia-protocol}
\label{subsec:lblt}
For the current free variables and selected independent working rows, set
\begin{equation}\label{eq:working-kkt-matrix}
 K_W=\begin{pmatrix}H_{FF}&A_W^\top\\A_W&0\end{pmatrix},\qquad
 \Pi^\top K_W\Pi=LBL^\top,\qquad
 B=\operatorname{diag}(B_1,\ldots,B_q).
\end{equation}
$\Pi$ is a symmetric permutation, $L$ is unit lower triangular, and each
pivot $B_j$ is symmetric $1\times1$ or $2\times2$. A complete factorization
satisfies \eqref{eq:working-kkt-matrix}; if pivoting is stopped, store the
entire Schur remainder and its coupling, never replace it with zero. The
factorization is associated with a specific model ID, tangent frame version,
free-variable ordering, working-row IDs and bound selector.

To solve $K_Wx=r$, execute exactly these triangular/block operations:
\begin{equation}\label{eq:lblt-solve}
 r_p=\Pi^\top r,\quad Ly=r_p,\quad Bz=y,\quad L^\top w=z,\quad x=\Pi w.
\end{equation}
Each $2\times2$ diagonal block is solved as a $2\times2$ system. A scalar-pivot
representation requires the corresponding symmetric transformation. After
the solve compute the scaled backward error
\begin{equation}\label{eq:lblt-backward-error}
 \eta_{\rm bwd}=
 \frac{\|K_Wx-r\|_\infty}
 {\|K_W\|_\infty\|x\|_\infty+\|r\|_\infty}.
\end{equation}
For the zero right side and zero solution, define the backward error as
zero. An empty free KKT system has a complete empty factor, zero residuals
and inertia $(0,0,0)$; only prescribed bound motion and full reduced-cost
recovery remain. Retain the factor only if the reconstruction residual and
the solve residual meet the declared scaled tolerances. For a rank-revealing or partial factor,
check the residual against the complete stored representation including the
Schur remainder.

For an independent row matrix $A_W$ with $m$ rows and an $n_F$-dimensional
free primal space, positivity of $H_{FF}$ on $\ker A_W$ gives the expected
inertia
\begin{equation}\label{eq:expected-inertia}
 \operatorname{In}(K_W)=(n_F,m,0).
\end{equation}
This identity is the factorization check for the selected regular face. In an
orthonormal null basis $Z$ with $A_WZ=0$, the congruence relation is
\begin{equation}\label{eq:kkt-inertia-nullspace}
 \operatorname{In}(K_W)=(m,m,0)+\operatorname{In}(Z^\top H_{FF}Z).
\end{equation}
The constraint blocks thus contribute the expected saddle-point negative
pivots, and the inertia test is \eqref{eq:expected-inertia}. For a $1\times1$ pivot $[d]$, its inertia
is $(1,0,0)$ if $d>0$, $(0,1,0)$ if $d<0$. For a symmetric $2\times2$
pivot $B_j=\left(\begin{smallmatrix}a&b\\b&c\end{smallmatrix}\right)$, use
the signed eigenvalues
\begin{equation}\label{eq:two-by-two-inertia}
 \lambda_{\pm}=\frac{a+c\pm\sqrt{(a-c)^2+4b^2}}{2}
\end{equation}
or an equivalent stable signed-inertia test. Classify each signed eigenvalue against the scaled
threshold $\delta_{\rm abs}+\delta_{\rm rel}\|K_W\|$; a classification band
that straddles zero is recorded as \texttt{UNRESOLVED\_PIVOT} and requires
a resolved pivot before the face can be admitted.

A dense refactorization costs $O((n_F+m)^3)$. Once a validated complete
factorization is available, a single bordered row/coordinate update uses
triangular solves and a small pivot operation, with $O((n_F+m)^2)$ arithmetic
for dense factors. A deletion uses symmetric permutations and local
refactorization of the affected trailing block; its cost includes the explicit trailing Schur assembly and factorization
in \eqref{eq:delete-local-schur}. It is quadratic only when the selected
update implementation confines this work accordingly. If fill or pivot quality prevents this,
refactorize the same new matrix and record \texttt{REFACTORIZED}.
Sparse implementations replace these bounds by costs determined by
fill and elimination-tree changes, but must execute the same algebraic
identities and checks.

\subsection{Incomplete factors and admissible solve states}
\label{subsec:factor-completion}
A factorization stops early when no remaining pivot can be accepted
safely; what is stored in that case determines whether the inertia can
still be assessed. A numerical factor object has one of three states: \texttt{COMPLETE},
\texttt{STOPPED\_WITH\_REMAINDER}, or
\texttt{REFACTORIZED\_COMPLETE}. In the stopped state, with $k$ pivots
accepted, store the completed block diagonal factor $B_c$, leading unit
lower triangular factor $L_c$, coupling $L_{Rc}$, permutation $\Pi$, and
full trailing Schur block $S_*$:
\begin{equation}\label{eq:partial-lblt-certified}
 \Pi^\top K_W\Pi=
 \begin{pmatrix}L_c&0\\L_{Rc}&I\end{pmatrix}
 \begin{pmatrix}B_c&0\\0&S_*\end{pmatrix}
 \begin{pmatrix}L_c^\top&L_{Rc}^\top\\0&I\end{pmatrix}.
\end{equation}
The representation is complete only if its reconstructed residual meets
$\epsilon_{\rm fac}$ and $S_*$ is empty. A nonempty $S_*$ contributes its
own signed inertia whenever that inertia can be certified. The inertia of
$K_W$ includes both the completed pivots and the retained remainder.
The online path worker requires a nonsingular complete factor for its
sensitivity solve. It may request one full refactorization of $K_W$ with a
new deterministic pivot policy. If the remainder remains unresolved, the
worker returns \texttt{UNRESOLVED\_FACTOR} with the retained remainder.
A regular-face solve requires the completed factorization.

Completion uses the signed $1\times1$ and $2\times2$ pivots of
Section~\ref{subsec:lblt}. The remainder and the attempted pivot policy
are retained with the stopped checkpoint, so that a failed completion can
be reproduced. A refactorization of the full matrix is recorded separately
from an incremental update.

The reconstruction test is applied to the entire represented matrix. For a
complete factor it is \eqref{eq:factor-reconstruction}; for the stopped form
it is
\begin{equation}\label{eq:partial-reconstruction-error}
 \eta_{\rm partial}=
 \frac{\left\|\Pi^\top K_W\Pi-L_0\operatorname{diag}(B_c,S_*)L_0^\top\right\|_\infty}
      {1+\|K_W\|_\infty},\quad
 L_0=\begin{pmatrix}L_c&0\\L_{Rc}&I\end{pmatrix}.
\end{equation}
The solve residual is always computed in original coordinates. A solve
requiring a complete factor must pass both factor completion and solve
regularity, in addition to the bound on $\eta_{\rm partial}$.

\subsection{Bound and structural-coordinate event protocol}
\label{subsec:bound-coordinate-updates}
The bound catalogue is assembled from physical variable metadata before QP
construction. For every scalar coordinate $i$, the metadata contains a
stable coordinate ID, lower-side ID, upper-side ID, time functions
$\ell_i(\tau),u_i(\tau)$, and an enum in
$\{\texttt{STRUCTURAL},\texttt{LOWER},\texttt{UPPER},\texttt{FREE}\}$. If
$\ell_i(\tau)=u_i(\tau)$ for every $\tau\in[0,1]$, it is
\texttt{STRUCTURAL}; substitute $p_i=c_i(\tau)$ in all rows and remove the
coordinate before any rank or factor test. If endpoints coincide only at
isolated parameters, they are not structural: the event handler processes
both sides at that parameter and the next checkpoint must resolve whether a
nonzero-width interval opens. This prevents an isolated coincidence from
being mistaken for a permanently fixed degree of freedom.

A nonstructural scalar may have at most one working side. Store the oriented
rows and signed duals explicitly:
\begin{equation}\label{eq:oriented-bound-row}
 \begin{array}{c|c|c|c}
 \text{side}&a_i^\top&\text{face equation}&\text{dual sign}\\\hline
 \text{lower}&+e_i^\top&p_i=\ell_i(\tau)&\nu_i^\ell\leq0\\
 \text{upper}&-e_i^\top&-p_i=-u_i(\tau)&\nu_i^u\leq0
 \end{array}
\end{equation}
The signed multiplier for the upper row is the upper-side multiplier; the
full stationarity contribution is $e_i\nu_i^\ell-e_i\nu_i^u$. Therefore the
reduced-coordinate identity is $\bar g_i=-\nu_i^\ell$ at a lower face and
$\bar g_i=\nu_i^u$ at an upper face, as in
\eqref{eq:reduced-cost-map}. This sign table is authoritative when converting
between separate arrays and the combined general-row representation.

At a step start compute primal residuals independently for both sides:
\begin{equation}\label{eq:bound-side-residual}
 s_i^\ell=p_i-\ell_i,\qquad s_i^u=u_i-p_i,
 \qquad s_i^\ell\geq0,\quad s_i^u\geq0.
\end{equation}
For a nonworking side with negative slack beyond tolerance, reject the
checkpoint. For a zero or near-zero nonworking slack, calculate the forward
derivative from \eqref{eq:bound-ratios}: if it is negative, submit an
immediate zero-ratio primal event; if it is nonnegative, keep the coordinate
free and monitor it at the next event. A side can therefore be geometrically
active and still not be in the algebraic working set. The distinction is
important when the side is weakly active and its multiplier is zero.

When a lower-bound event is accepted, prescribe $p_i=\ell_i(\tau^+)$ and
remove coordinate $i$ from $F$. Its contribution to the reduced objective
and rows is included through the affine lift $p=E_Fx_F+p_B$. The new matrix
is the principal primal reduction plus the same working-row block; the
coordinate deletion operation is \eqref{eq:delete-local-schur} applied to
that primal index and must preserve symmetric ordering. When a
working lower bound releases, add coordinate $i$ back to $F$, prescribe no
coordinate value, and use its current bound value only as the event
checkpoint. The upper-side operations are identical with $u_i$ substituted
for $\ell_i$. In both cases, recalculate $r_F^d$ and $r_W^p$ from
\eqref{eq:bound-affine-rhs}; retaining old right sides would solve a different
face problem.

A moving working side has nonzero prescribed derivative. For a lower side,
$\dot p_i=\dot\ell_i$; for an upper side, $\dot p_i=\dot u_i$. In the
sensitivity solve, these values enter $H_{FN}\dot p_N$ and
$B_{W,N}\dot p_N$. In a Newton correction, use
$d_N=p_N^{\rm prescribed}-p_N$ instead. The two right-hand sides are
written separately in \eqref{eq:path-sensitivity} and
\eqref{eq:qp-newton-correction}; implementations must not substitute a bound
velocity for a correction displacement. These are the only moving-bound
terms in the frozen QP because $H$ and all row normals are constant in $\tau$.

\subsection{Candidate generation and working-set initialization}
\label{subsec:task-generation}
At each physical feedback update, build one primary hypothesis from the
previous accepted target working set and a finite configured set of alternate
hypotheses. Every hypothesis is a tuple of stable inequality and bound IDs;
equality IDs are always retained. The supported generation rules are:
\begin{enumerate}
 \item \texttt{WARM}: retain the previous accepted working IDs whose current
       primal margins are within a declared threshold
       $\epsilon_w\geq\epsilon_a$ and whose previous multipliers
       are below $-\epsilon_\mu$;
 \item \texttt{ACTIVE\_IN}: add to \texttt{WARM} one or more omitted rows
       whose margins are within $\epsilon_a$, in a declared deterministic
       order, only when an independent-row test succeeds;
 \item \texttt{ACTIVE\_OUT}: remove from \texttt{WARM} its weakly supported
       rows, those with $|\mu_i|\leq c_\mu\epsilon_\mu$ for a declared factor
       $c_\mu\geq1$; and
 \item \texttt{EXCHANGE}: use only the verified dual-compatible candidates
       from \eqref{eq:exchange-working-set}.
\end{enumerate}
The configured initial-hypothesis cap $N_{\max}\leq N_{\rm task}$ applies
after duplicate merging; $N_{\rm task}$ also counts exchange children.
A child is admitted after the following tests. Compute a rank-revealing
QR of $A_W^\top$ with fixed threshold, record its permutation and diagonal
rank profile, and require the selected rows to be independent at that
threshold. Construct and factor the selected KKT matrix, check its expected
inertia, and form the auxiliary data by \eqref{eq:start-qp-data} at the
proposed primal and sign-feasible dual. Verify its start-face equations,
using \eqref{eq:qp-newton-correction} only if a declared numerical correction
is needed. If a hypothesis cannot produce a verified start, return it with
a start failure status; changing its working rows creates a separate
hypothesis. A hypothesis whose proposed rows are linearly dependent may be
mapped, before these tests, to an independent sign-feasible representation:
its multiplier is redistributed onto an independent subset of the rows by
the family \eqref{eq:dependent-dual-family}. The redistribution leaves the
combination $\sum_i\zeta_ia_i$ unchanged, and therefore also the gradient
$g^0$ of the auxiliary start.

The warm start is formed from the preceding accepted manifold trajectory.
Apply the
horizon shift and terminal extension of Section~\ref{sec:parametric-problem}
to obtain $r^{\rm shift}$, substitute the current measured node-zero state,
and, on the declared inverse-retraction neighborhood, set
\begin{equation}\label{eq:warm-start-inverse-retraction}
 p^0=E_v^*R_{r_v}^{-1}(r^{\rm shift}).
\end{equation}
If $r_v=r^{\rm shift}$ this vector is zero. A transport of a stored tangent
vector is permitted only as a separately checked approximation to
\eqref{eq:warm-start-inverse-retraction}. The terminal extension rule $\kappa_{\rm ext}$ of
\eqref{eq:terminal-control-extension} supplies
the new terminal state and control interval before this inverse is taken.
Shift old dual arrays by retained row IDs and transform their residual dual
frames; new rows receive zero dual. Keep only supported nonpositive
inequality/bound multipliers for the proposed working selection. The pair
$(p^0,\zeta^0)$ is a candidate, and the exact auxiliary construction
\eqref{eq:start-qp-data} is always verified before its homotopy begins.
A shifted trajectory need not satisfy the new dynamics or active constraints;
those mismatches determine the new target row residuals.

For each generated task $w$, the registry key is
\begin{equation}\label{eq:task-id-key}
 \mathsf{key}(w)=\bigl(\mathsf{id}_{\rm model},\mathsf{id}_{\rm frame},
  \mathsf{id}_{\rm path},\mathsf{IDs}(W_w),p_w^0,\zeta_w^0,\mathsf{factor\_version}\bigr).
\end{equation}
Exact duplicates share one stable task ID. A tolerance-based merge is
permitted only for identical working-row/side selectors and immutable path
data, when the scaled norm of every primal and dual difference is below
its declared threshold. On a
merge, retain the lower stable ID and log the discarded ID, residuals, and
factor version. Factors are reusable across tasks only when their complete
matrix keys in \eqref{eq:factor-key} match exactly.

\subsection{Adding a working row: border, update, and inertia check}
\label{subsec:add-working-row}
Suppose $K_W$ has factorization \eqref{eq:working-kkt-matrix}. Adding one
independent general active row with free normal $a_F$ gives the bordered
matrix
\begin{equation}\label{eq:border-general-row}
 \widetilde K=\begin{pmatrix}K_W&c\\c^\top&0\end{pmatrix},
 \qquad c=\operatorname{col}(a_F,0_{m}).
\end{equation}
Set
\begin{equation}\label{eq:border-solves}
 \widehat c=\Pi^\top c,\qquad y=L^{-1}\widehat c,\qquad
 w=B^{-1}y,\qquad \theta=-y^\top B^{-1}y=-c^\top K_W^{-1}c.
\end{equation}
When the scalar pivot is acceptable, the exact update is
\begin{equation}\label{eq:addcon1}
 \begin{pmatrix}\Pi^\top&0\\0&1\end{pmatrix}\widetilde K
 \begin{pmatrix}\Pi&0\\0&1\end{pmatrix}
 =\begin{pmatrix}L&0\\w^\top&1\end{pmatrix}
   \begin{pmatrix}B&0\\0&\theta\end{pmatrix}
   \begin{pmatrix}L^\top&w\\0&1\end{pmatrix}.
\end{equation}
The resulting inertia update is
\begin{equation}\label{eq:addinertup1}
 \operatorname{In}(\widetilde K)=\operatorname{In}(K_W)+\operatorname{In}([\theta]).
\end{equation}
For a newly independent working inequality, the expected change is one
additional negative eigenvalue and no other change:
$(n_F,m,0)\mapsto(n_F,m+1,0)$. The event KKT residual is checked again with
new multiplier zero and the updated factors before the task resumes.

If a small scalar border is handled by combining the previous trailing
scalar pivot with the incoming coordinate, its local block is formed before
eliminating that old scalar. Write the old factors, after the declared
symmetric permutation, as
\begin{equation}\label{eq:border-trailing-partition}
 L=\begin{pmatrix}L_-&0\\l^\top&1\end{pmatrix},\quad
 B=\operatorname{diag}(B_-,[\eta]),\quad
 \Pi^\top c=\binom{c_-}{c_t}.
\end{equation}
For a border with diagonal entry $d$ ($d=0$ for a general row and $d=h_{ii}$
for an inserted primal coordinate), set
\begin{equation}\label{eq:addcon2}
 y_-=L_-^{-1}c_-,\quad w_-=B_-^{-1}y_-,\quad
 \rho=c_t-l^\top y_-,\quad \omega=d-y_-^\top B_-^{-1}y_-,\quad
 B_{\rm new}=\begin{pmatrix}\eta&\rho\\\rho&\omega\end{pmatrix}.
\end{equation}
The exact extended factor is
\begin{equation}\label{eq:border-promoted-factor}
 \begin{pmatrix}L_-&0&0\\l^\top&1&0\\w_-^\top&0&1\end{pmatrix}
 \operatorname{diag}(B_-,B_{\rm new})
 \begin{pmatrix}L_-^\top&l&w_-\\0&1&0\\0&0&1\end{pmatrix}.
\end{equation}
When $\eta\ne0$, the fully eliminated scalar Schur complement is
$\omega-\rho^2/\eta$; it is not the lower-right entry $\omega$ of the
promoted block. Thus the old scalar pivot is replaced once, giving
\begin{equation}\label{eq:addinertup2}
 \operatorname{In}(\widetilde K)=\operatorname{In}(K_W)
       -\operatorname{In}([\eta])+\operatorname{In}(B_{\rm new}).
\end{equation}
Compute its two signed eigenvalues by \eqref{eq:two-by-two-inertia}, or use
a stable determinant/trace classification with scaling. Never count the old
scalar pivot twice. The update is accepted only after verifying the full
reconstruction
\begin{equation}\label{eq:factor-reconstruction}
 \eta_{\rm fac}=
 \frac{\|\Pi^\top\widetilde K\Pi-\widetilde L\widetilde B\widetilde L^\top\|_\infty}
 {1+\|\widetilde K\|_\infty}
 \leq\epsilon_{\rm fac}
\end{equation}
and the expected inertia. In a dense implementation this residual need
not be recomputed by a cubic matrix product at every scalar border. Store
the old defect $\mathcal E_K=\Pi^\top K_W\Pi-LBL^\top$ with its exact matrix version. For
\eqref{eq:addcon1}, with the actual border diagonal $d$, calculate
\begin{equation}\label{eq:border-reconstruction-residual}
 e_c=\Pi^\top c-LBw,\qquad
 e_d=d-w^\top Bw-\theta,\qquad
 \mathcal E_{\widetilde K}=\begin{pmatrix}\mathcal E_K&e_c\\e_c^\top&e_d\end{pmatrix}.
\end{equation}
For the promoted factor \eqref{eq:border-promoted-factor}, replace $LBw$
in the column residual by
$\operatorname{col}(L_-B_-w_-,\,l^\top B_-w_-+\rho)$ and replace
$w^\top Bw+\theta$ in the diagonal residual by $w_-^\top B_-w_-+\omega$.
These identities evaluate every new residual entry using quadratic work,
retaining the old full defect and the declared roundoff bound. Every entry
is included in the residual test. A symmetric permutation also permutes
the cached defect. A factor change for which no validated incremental
residual identity is available requires a fresh reconstruction; its cubic
cost is entered separately in \eqref{eq:task-cost-ledger}. All matrices in
this check are the actual symmetrically permuted matrices. If the accepted pivot path cannot be
reconstructed, discard the copied factors and refactorize the same bordered
matrix. If refactorization cannot certify the selected independent face,
return that checkpoint with an explicit factorization status.

\subsection{Removing a working row and changing the factor order}
\label{subsec:remove-working-row}
Removing a working row deletes one row and one column of the KKT matrix.
The factors of the part that precedes the deleted coordinate in the pivot
order remain valid, and only the trailing part is recomputed.
To remove an inequality row, symmetrically move its dual coordinate last.
The right-hand side and original stable-ID permutation are moved by the same
permutation. A concrete local refactorization retains any unchanged leading
pivot blocks. If the first affected pivot begins after index $r$, partition
the symmetrically reordered matrix as
\begin{equation}\label{eq:delete-local-schur}
 K'=\begin{pmatrix}K_{11}&K_{12}'\\K_{21}'&K_{22}'\end{pmatrix},\quad
 K_{11}=L_{11}B_{11}L_{11}^\top,\quad
 L_{21}'=K_{21}'L_{11}^{-\top}B_{11}^{-1},\quad
 S'=K_{22}'-L_{21}'B_{11}L_{21}'^\top.
\end{equation}
Delete the last coordinate from the trailing rows of $L_{21}'$ and the
corresponding row and column of $S'$, factor the retained trailing Schur
matrix with \eqref{eq:lblt-schur-step}, and concatenate it with the unchanged
prefix. Equivalently, update the affected trailing factors by symmetric
pivot interchanges. Both operations reconstruct the same principal
submatrix; the implementation records the number of refactored trailing
coordinates. For $r=0$, this operation is a full refactorization and is
counted as such. No row-only permutation or unrecorded pivot deletion is
permitted.

If the resulting
factorization has the exact trailing scalar form
\begin{equation}\label{eq:subcon1}
 \Pi^\top K_W\Pi=
 \begin{pmatrix}L_-&0\\w^\top&1\end{pmatrix}
 \operatorname{diag}(B_-,[\theta])
 \begin{pmatrix}L_-^\top&w\\0&1\end{pmatrix},
\end{equation}
then deletion returns $K_-=L_-B_-L_-^\top$ and
\begin{equation}\label{eq:subinertup1}
 \operatorname{In}(K_-)=\operatorname{In}(K_W)-\operatorname{In}([\theta]).
\end{equation}
For a trailing $2\times2$ block
$E=\left(\begin{smallmatrix}\eta&\rho\\\rho&\theta\end{smallmatrix}\right)$,
if the second coordinate is the one deleted, retain the first coordinate
block $[\eta]$ and update inertia by
\begin{equation}\label{eq:subcon2}
 \operatorname{In}(K_-)=\operatorname{In}(K_W)-\operatorname{In}(E)
                                     +\operatorname{In}([\eta]).
\end{equation}
The formula is a block-count identity for every $\eta$ and requires no
division by $\eta$. Use the scalar Schur expression
$\theta-\rho^2/\eta$ only when $\eta$ has a certified nonzero bound and the
corresponding symmetric block elimination was actually performed.

After deletion, recompute the full expected inertia for the retained matrix;
for a regular face with $m-1$ independent working rows it is
$(n_F,m-1,0)$. Reconstruct the retained factors and compute \eqref{eq:lblt-backward-error}
before resuming the worker. A zero pivot or incorrect inertia returns a
non-success face status for the same unchanged matrix. Equality rows remain
selected throughout this operation. When the
deleted object is a working bound, apply the corresponding primal-coordinate
operation in
Section~\ref{subsec:bound-coordinate-updates}.

\subsection{Bound-coordinate updates, basic variables, and reduced costs}
\label{subsec:euclidean-bounds}
The free KKT matrix is built using only $H_{FF}$ and $A_W=B_{W,F}$.
Structural coordinates have already been substituted in the stored QP.
Let $E_F$ inject the free coordinates and let $p_B$ be the full-length vector
whose entries on $N$ equal the prescribed working-bound values and whose
entries on $F$ are zero. Thus
\begin{equation}\label{eq:bound-affine-lift}
 p=E_Fx_F+p_B,\qquad
 H_F=E_F^\top HE_F,\quad g_F^{\rm red}=E_F^\top(g+Hp_B),\quad
 A_W=B_WE_F,\quad \widetilde\beta_W=\beta_W-B_Wp_B .
\end{equation}
Here $H_F$ is the matrix written $H_{FF}$ elsewhere. Substitution gives the complete reduced linear system
\begin{equation}\label{eq:bound-affine-rhs}
 \begin{pmatrix}H_{FF}&A_W^\top\\A_W&0\end{pmatrix}
 \binom{x_F}{\zeta_W}
   =\binom{-g_F-H_{FN}p_N}{\beta_W-B_{W,N}p_N}.
\end{equation}
The multiplier term $A_W^\top\zeta_W$ remains on the left. The symbol $\beta_W$ here is always the unreduced selected
right side from \eqref{eq:working-kkt}; $\widetilde\beta_W$ denotes its
reduced value. This convention also applies after any bound insertion,
release or task exchange.

For a moving working bound, use its prescribed
$\dot p_N=\dot\ell_N$ or $\dot u_N$. Differentiating the affine lift gives the sensitivity system \eqref{eq:path-sensitivity} again,
\begin{equation}\label{eq:free-bound-kkt}
 \begin{pmatrix}H_{FF}&A_W^\top\\A_W&0\end{pmatrix}
 \binom{\dot p_F}{\dot\zeta_W}
 =\binom{-(\dot g)_F-H_{FN}\dot p_N}
        {\dot\beta_W-B_{W,N}\dot p_N}.
\end{equation}
For a fixed-parameter correction, the right side is instead
\eqref{eq:qp-newton-correction}, using $d_N=p_N^{\rm prescribed}-p_N$.
If a pre-substitution model contains structural coordinates $C$ with
prescribed affine values $p_C(\tau)$, first replace its objective coefficient
by $g_{\mathcal I_{\rm ns}}+H_{\mathcal I_{\rm ns},C}p_C(\tau)$ and its row right sides by
$b-A_Cp_C(\tau)$, and then apply \eqref{eq:bound-affine-lift} on
$\mathcal I_{\rm ns}$. Their derivatives enter this first substitution exactly once.
They are not reintroduced as working bounds or as extra terms in
\eqref{eq:free-bound-kkt}.

The full signed reduced gradient is
\begin{equation}\label{eq:reduced-cost}
 \bar g=g+Hp+B_W^\top\zeta_W.
\end{equation}
For the Lagrangian with lower normals $+e_i$ and upper normals $-e_i$,
\begin{equation}\label{eq:reduced-cost-map}
 \bar g_F=0,\quad \bar g_i=-\nu_i^\ell\geq0\ (i\in W_\ell),\quad
 \bar g_i=\nu_i^u\leq0\ (i\in W_u).
\end{equation}
These are signed coordinate reduced costs. A bound release test is evaluated
from both $\bar g_i$ and $\dot{\bar g}_i$ in
Section~\ref{subsec:path-events}; no unsigned complementarity residual can
replace that sign test.

When one free coordinate $i$ is fixed at a new working bound, delete its
primal row and column from the KKT matrix. When a working bound releases,
insert that primal variable. In the current factor ordering the new coordinate
border is
\begin{equation}\label{eq:coordinate-border}
 \widetilde K=\begin{pmatrix}K_W&c_i\\c_i^\top&h_{ii}\end{pmatrix},\qquad
 c_i=\operatorname{col}(H_{F,i},A_{W,i}),\qquad
 \theta_i=h_{ii}-y_i^\top B^{-1}y_i,
 \quad y_i=L^{-1}\Pi^\top c_i.
\end{equation}
Use \eqref{eq:addcon1} with the lower-right scalar $0$ replaced by $h_{ii}$
and the computed $\theta_i$. The reverse operation permutes the primal
coordinate last and removes its actual trailing pivot block with the same
block rules as \eqref{eq:subcon1}--\eqref{eq:subcon2}. Expected inertia is
$(n_F-1,m,0)$ after fixing one independent coordinate and
$(n_F+1,m,0)$ after releasing it, provided the working row rank remains $m$.
If it does not, return to the row-dependency handler before solving.

The full lower and upper slacks and derivatives are those in
\eqref{eq:bound-ratios}. Coincident bounds were eliminated in
\eqref{eq:bound-catalogues}; their prescribed motion affects the affine
right sides, but they never enter either a primal release event or a reduced-
cost release event. A coordinate at a bound that is not selected as working
is still free, and an outward derivative creates a zero-distance event. This
is required to distinguish an active-but-weak bound from a fixed coordinate.

\subsection{Riemannian simple-bound rows and tangent-frame updates}
\label{subsec:riemannian-bounds}
\label{subsec:geometric-rows}
A smooth geometric margin $b_s(q,t)\geq0$ contributes, in an orthonormal
frame $E_q:\mathbb R^{n_q}\to T_q\mathcal M$, the row
\begin{equation}\label{eq:geometric-bound-row}
 a_s(q,t)=D b_s(q,t)E_q,\qquad
 b_s(q,t)+a_s(q,t)\xi\geq0.
\end{equation}
Its stable ID is kept whether it is selected in the working rows or only
monitored for blocking. For the independent working geometric margins
$W_b$, let $N_q$ stack their rows. The local affine tangent restriction is
\begin{equation}\label{eq:geometric-affine-tangent}
 N_q\xi=r_b(\tau),\qquad r_b(\tau)=(1-\tau)\beta_{W_b}^0+\tau\beta_{W_b}^1,\quad
 \beta_{W_b}^1=-b_{W_b}(q,t).
\end{equation}
The start offsets $\beta_{W_b}^0$ are those of
\eqref{eq:start-qp-data}. At a feasible target boundary the target value
is $r_b(1)=0$. Solve for a particular coordinate $\xi_0$ and choose
orthonormal $U_F$ spanning
$\ker N_q$:
\begin{equation}\label{eq:geometric-free-basis}
 \xi=\xi_0+U_F\xi_F,\quad N_q\xi_0=r_b,\quad N_qU_F=0,\quad
 U_F^\top U_F=I.
\end{equation}
This parameterization supplies the face nullspace and, if used for a
condensed face solve, reduces the objective and remaining rows by the exact
affine substitution in \eqref{eq:shooting-condensed-data}. The affine term
$E_q\xi_0$ contributes to both right sides as in
\eqref{eq:bound-affine-rhs}. The canonical event worker retains geometric
rows in $B_W$ and updates them as general rows, retaining their normal
coordinates in the KKT matrix. A backend using the
face elimination expands its primal and recovers the eliminated geometric
duals before applying the full-row event and residual tests. A QR factorization with column pivoting supplies the rank-revealing
orthonormal basis; the task records the rank, permutation, verified row-rank
bound from its triangular factor, and tangent-frame version.

For independent working normals $N_q$, write the thin QR factorization
$N_q^\top=Q_NR_N$. The eliminated geometric multipliers are recovered by
the triangular system
\begin{equation}\label{eq:geometric-reduced-cost}
 R_N\nu=-Q_N^\top g_0,\qquad
 (I-Q_NQ_N^\top)g_0=0,\qquad g_0+N_q^\top\nu=0,
\end{equation}
where $g_0$ is the $q$-block of $g+Hp$ plus all other selected
constraint-gradient terms, excluding precisely the geometric rows in $N_q$. The computed components must
satisfy $\nu\leq0$. The final two identities are checked in the original
tangent coordinates. The triangular factor includes the coupling between
nonorthogonal normals and determines their multipliers uniquely.
If normals are dependent, construct an independent working representation
and retain the other admissible representations as tasks using
Section~\ref{subsec:dependent-exchange}.

A manifold retraction preserves membership in $\mathcal M$. Exact membership
in an active inequality boundary requires an additional construction: the
retraction of a QP solution that lies on the linearized boundary of an
active constraint need not lie on its nonlinear boundary. The canonical
worker does not correct this defect. It accepts the retracted point only if
the original nonlinear return tolerances hold there. An implementation additionally requiring exact geometric
face feasibility uses a known face retraction or the joint local correction
\begin{equation}\label{eq:boundary-retraction-equations}
 b_{W_b}\bigl(R_q(\eta+E_qN_q^\top a),t\bigr)=0
\end{equation}
for normal coefficients $a$, and verifies the joint residual and local
solution. The extra displacement and the difference from the canonical
product retraction are recorded in the retraction-error budget used in
Section~\ref{sec:pathfollowing-theory}. A
failed face correction rejects that nonlinear candidate at synchronization,
while the frozen QP predictor remains fixed. For $SO(3)$
log-coordinate bounds, the row includes the differential of the log map,
for example $e_i^\top J_r(\phi)^{-1}$ under the right-update convention
$Y^+=Y\exp(\widehat\eta)$. Thus a bound on $\phi_i$ is not generally a
coordinate bound on $\eta_i$. At a tangent-frame change $E_F^+=E_FQ$, rotate
all reduced data by
\begin{equation}\label{eq:tangent-congruence}
 H_F^+=Q^\top H_FQ,\qquad A_F^+=A_FQ,\qquad
 K_F^+=\operatorname{diag}(Q,I)^\top K_F\operatorname{diag}(Q,I).
\end{equation}
If dimensions change, perform this congruence on the common coordinates
first, then add/delete the row or variable and refactor the affected blocks.
A factor from a previous frame cannot be applied to new coordinates without
\eqref{eq:tangent-congruence}.

\subsection{Dependent blocking rows and exact dual-compatible exchanges}
\label{subsec:dependent-exchange}
Both general rows and component bounds use the same full-row dependence
identity. A working bound on coordinate $i\in N$ has orientation
$\sigma_i=+1$ at a lower bound and $\sigma_i=-1$ at an upper bound. Define
\begin{equation}\label{eq:full-working-catalogue}
\begin{aligned}
 D_N&=\operatorname{diag}(\sigma_i)_{i\in N},&
 T_N&=D_NE_N^\top,\\
 \mathscr B_W&=\begin{bmatrix}B_W\\T_N\end{bmatrix},&
 \mathscr\beta_W&=\binom{\beta_W}{D_Np_N},\qquad
 \mathscr\zeta_W=\binom{\zeta_W}{\nu_N}.
\end{aligned}
\end{equation}
Here $E_N$ is the coordinate injection, $\nu_N$ stores the signed lower or
upper multiplier selected at each coordinate, and every omitted multiplier
is zero. The full stationarity equation is
$Hp+g+\mathscr B_W^\top\mathscr\zeta_W=0$. Independence of $A_W=B_{W,F}$
is equivalent to independence of the combined general/equality/bound rows
$\mathscr B_W$.

An incoming general row uses its full normal $a_l$ and right side
$\beta_l(\tau)=b_l(\tau)$. An incoming lower or upper component bound uses
$(a_l,\beta_l)=(e_i,\ell_i)$ or $(-e_i,-u_i)$, respectively. For an
independent general row, use the border
\eqref{eq:border-general-row}. For an independent incoming bound, use the
coordinate deletion in \eqref{eq:coordinate-border} in reverse, set its
prescribed value, and recover its multiplier from
\eqref{eq:reduced-cost-map}. These two updates represent the same full-row
addition before and after bound elimination.

Suppose now that the normal of a blocking row is a linear combination of
the working rows. Adding the row would make the KKT matrix singular. The
same linear relation, however, shows how multiplier weight can be moved from
the working rows to the incoming row without changing stationarity, and the
first old inequality multiplier that reaches zero identifies a row that can
leave in exchange for the incoming one.
For such a dependent blocker first solve the free-space relation and then lift
it to the working-bound coordinates:
\begin{equation}\label{eq:dependent-row-identities}
\begin{aligned}
 a_{l,F}&=A_W^\top\gamma,&
 \gamma_N&=D_N\bigl(a_{l,N}-B_{W,N}^\top\gamma\bigr),\\
 \Gamma&=\binom{\gamma}{\gamma_N},&
 a_l&=\mathscr B_W^\top\Gamma,\qquad
 \beta_l(\tau)=\Gamma^\top\mathscr\beta_W(\tau).
\end{aligned}
\end{equation}
The second identity on the last line is verified at the actual event
checkpoint. Equivalently, its reduced form is
$\beta_l-a_{l,N}^\top p_N=\gamma^\top\widetilde\beta_W$.
The full normal residual is scaled by
$\|a_l\|+\|\mathscr B_W\|\|\Gamma\|$, and the right-side residual by
the corresponding right-side scale. The coefficients in $\gamma_N$ are
essential to preserve the bound multipliers and stationarity at fixed
coordinates when lifting the free-column dependence. The right-side slope
need not obey the same dependence away from the event.

With $\delta\leq0$ the incoming multiplier, every full multiplier vector
preserving stationarity along this one-dimensional redistribution is
\begin{equation}\label{eq:dependent-dual-family}
 \mathscr\zeta_W^+(\delta)=\mathscr\zeta_W-\delta\Gamma,
 \qquad \mu_l^+=\delta.
\end{equation}
Equality components are unrestricted. For each removable general inequality
or oriented bound row $j$ with $\Gamma_j\ne0$, calculate
\begin{equation}\label{eq:exchange-ratio}
 \delta_j=\frac{(\mathscr\zeta_W)_j}{\Gamma_j}.
\end{equation}
Retain the candidate precisely when $\delta_j\leq0$ and
\begin{equation}\label{eq:exchange-sign-test}
 (\mathscr\zeta_W)_i-\delta_j\Gamma_i\leq0
       \quad\hbox{for every old inequality or bound row }i.
\end{equation}
The equalities are never candidates for removal. Zero ratios are retained;
for $\Gamma_j=0$ there is no exchange ratio. At each retained value the old
multiplier of row $j$ is exactly zero, all other old inequality multipliers
are sign-feasible, and the incoming multiplier is sign-feasible.

Let $\mathcal W$ denote the full selected ID catalogue, including oriented
working bounds. The exact exchange is
\begin{equation}\label{eq:exchange-working-set}
 \mathcal W_j^+=(\mathcal W\cup\{l\})\setminus\{j\},\qquad
 \mathscr\zeta_W^+=\mathscr\zeta_W-\delta_j\Gamma,\quad
 \mu_l^+=\delta_j,
\end{equation}
with the zero component $j$ removed from the selected storage and kept zero
in the full stable-ID array. The primal $p$ and parameter $\tau$ are
unchanged. If $j$ is a general inequality, remove its row from $W_I$; if it
is a bound, release its coordinate from $N$ to $F$. If $l$ is a general
inequality, insert its row in $W_I$; if it is a bound, select its side and
move its coordinate from $F$ to $N$. Rebuild $B_W,\beta_W,A_W$ and
$\widetilde\beta_W$ by \eqref{eq:working-kkt}, and recompute the reduced
costs and their oriented bound multipliers. The computed reduced costs must
agree with the lifted multiplier vector in
\eqref{eq:exchange-working-set}; otherwise the child is rejected.
The family in \eqref{eq:dependent-dual-family} is used only to identify
these working-set exchanges. After row $j$ is removed, independence of
$\mathcal W_j^+$ makes its displayed multiplier the unique multiplier
for that new face. Each admissible exchange therefore defines a distinct
primal--dual task with its own independent working set.

For each candidate, factor the final matrix $K_{\mathcal W_j^+}$ directly
or perform the corresponding row and coordinate updates on copied factors.
Intermediate factor matrices must reconstruct their declared matrix, but no
intermediate row selection is published as a checkpoint and no parameter
step is taken until the final exchange is complete. If an intermediate pivot
is unresolved, refactor the final matrix once. Check full row independence,
expected inertia, reconstruction and backward error, then solve the new
sensitivity system. A single exchange can expose another tied zero-distance event. The
private closure below processes those events before any child is
published or a positive parameter step is taken. Its final outgoing
working sets satisfy all forward tests as well as the full primal,
stationarity, sign and support checks. If no outgoing set survives,
retain the failed branch statuses with the parent checkpoint. This
local event closure explores only the exchange graph reached at the
current event; exhaustive enumeration of unrelated independent subsets
is a separate, explicitly budgeted mode.

As a small bound-lifting example, consider
$\min \tfrac12\|p\|^2+0.1p_1+p_2$ with $p_1\geq0$, $p_2\geq0$ and
$p_1+p_2\geq\tau$. At $\tau=0$ choose the first margin as a working
bound and the second as a working general row. Their multipliers are
$-0.1$ and $-1$. The incoming sum row has $\gamma=1$ in the sole free
column, but also $\gamma_N=1$ at the bound. Equations
\eqref{eq:exchange-ratio}--\eqref{eq:exchange-sign-test} select the bound
release at $\delta=-0.1$, producing $p=(\tau,0)$ for small positive
$\tau$. Omitting $\gamma_N$ would miss this valid continuation.

\paragraph{Simultaneous events and zero-distance closure.}
At a fixed event point $(p,\tau)$, several constraints can be blocking.
A valid outgoing working set can require more than one zero-distance
pivot, that is, more than one change of the working set made at the current
point without any step in $\tau$. The forward conditions are imposed on the completed event
transition, after these intermediate exchanges have been performed.
For an independent, sign-feasible candidate checkpoint $X$, let
$\mathcal W_{\rm ineq}(X)=\mathsf{IDs}(W_I,W_\ell,W_u)$ be the
selected general-inequality and oriented bound IDs. Write $\zeta_j$ for
the corresponding coefficient in its complete stable-ID multiplier
array. Recompute \eqref{eq:path-sensitivity} and define
\begin{equation}\label{eq:zero-distance-sets}
\begin{aligned}
 \mathcal D_0(X)&=\{j\in\mathcal W_{\rm ineq}(X):
               \zeta_j=0,\ \dot\zeta_j>0\},\\
 \mathcal P_0(X)&=\{i\in\mathcal A^0\setminus\mathcal W_{\rm ineq}(X):
               \dot s_i<0\}.
\end{aligned}
\end{equation}
Bounds use their oriented rows and multipliers; the zero classifications
and numerical tolerances are those of the ratio scan. Equality rows
never belong to $\mathcal D_0$. When $\mathcal D_0(X)\ne\varnothing$, let $j=\min\mathcal D_0(X)$.
Otherwise, if $\mathcal P_0(X)\ne\varnothing$, let
$l=\min\mathcal P_0(X)$. For a dependent $l$, let $\mathcal J_l(X)$
contain precisely the removable old inequality and bound IDs with
$\Gamma_j\ne0$, $\delta_j\leq0$ and all signs in
\eqref{eq:exchange-ratio}--\eqref{eq:exchange-sign-test} satisfied.
The set-valued one-pivot map is
\begin{equation}\label{eq:zero-distance-successors}
\mathsf{Next}_0(X)=
\begin{cases}
 \{\mathsf{Release}_j(X)\},
    &\mathcal D_0(X)\ne\varnothing,\\
 \{\mathsf{Add}_l(X)\},
    &\mathcal D_0(X)=\varnothing,\quad a_l\text{ independent},\\
 \{\mathsf{Exchange}_{lj}(X):j\in\mathcal J_l(X)\},
    &\mathcal D_0(X)=\varnothing,\quad a_l\text{ dependent}.
\end{cases}
\end{equation}
The map is evaluated only when $\mathcal D_0(X)\cup\mathcal P_0(X)$
is nonempty. Every minimum uses stable row IDs, and dual releases
precede primal blockers. The release operation removes the zero-multiplier general row
by \eqref{eq:subcon1}--\eqref{eq:subcon2}, or releases its bound coordinate
by \eqref{eq:coordinate-border}. The addition operation uses
\eqref{eq:border-general-row}--\eqref{eq:addcon2}, or the reverse
coordinate operation for a bound, with zero incoming multiplier.
The exchange is exactly \eqref{eq:exchange-working-set}. Each operation
rebuilds its selectors and affine lift, expands its unique multiplier,
and checks the final KKT factor and residuals. Thus every retained
intermediate selector is linearly independent, even though its new
sensitivity can identify a further zero-distance blocker.

Algorithm~\ref{alg:zero-distance-closure} applies these maps on private
copies at the same $(p,\tau)$. Only its completed outgoing states enter
the parallel family. Distinct intermediate states are keyed by their
full selector, primal and unique dual array at the fixed immutable
path. A repeated state is merged with its retained copy and recorded.
For queue ordering, append $(0,j,-1)$ to the lineage after a release,
$(1,l,-1)$ after an independent addition, and $(2,l,j)$ after an exchange,
and use lexicographic order. Only a state whose factor, KKT equations
and sensitivity solve have passed their tests enters the visited set.
There are finitely many independent selectors at fixed data. All private
successors charge one shared remaining operation counter, which is
preserved when a checkpoint is copied. Budget exhaustion retains pending
states and has its own stopping status.
No artificial positive step is introduced to leave a tie.

\begin{algorithm}[H]
\caption{Private zero-distance closure of an event transition}
\label{alg:zero-distance-closure}
\begin{algorithmic}[1]
\Require A finite seed collection $\mathcal S$ at the same event $(p,\tau)$ and immutable path, with the caller's remaining operation budget.
\State Set the pending queue $\mathcal P\gets\mathcal S$, visited keys $\mathcal V\gets\varnothing$, completed outputs $\mathcal C\gets\varnothing$ and failed records $\mathcal F\gets\varnothing$.
\While{$\mathcal P\ne\varnothing$ and the original budget remains}
 \State Remove its lowest-lineage checkpoint $X$ from $\mathcal P$.
 \If{the full key of $X$ is not in $\mathcal V$}
  \State Verify \eqref{eq:working-set-relevance}.
  \If{the relevance test passes}
   \State Solve \eqref{eq:path-sensitivity} with the full solve-residual and backward-error tests. A failed solve appends its specific stopped record to $\mathcal F$ and ends this branch. Otherwise insert the key in $\mathcal V$, recompute all zero-distance ratios, and form \eqref{eq:zero-distance-sets}.
   \If{$\mathcal D_0(X)=\mathcal P_0(X)=\varnothing$}
    \State Append the complete forward-relevant $X$ to $\mathcal C$.
   \Else
    \State Form every successor in \eqref{eq:zero-distance-successors} on a private copy; append each independently verified successor to $\mathcal P$ in stable lineage order, and every failed factor/rank/sign/residual record to $\mathcal F$.
   \EndIf
  \Else
   \State Append the specific stopped relevance-test record to $\mathcal F$.
  \EndIf
 \Else
  \State Merge the duplicate with its retained exact key and record the repeated zero-distance state.
 \EndIf
\EndWhile
\State Return $(\mathcal C,\mathcal P,\mathcal F)$; a nonempty $\mathcal P$ records \texttt{BUDGET\_EXHAUSTED} or deadline cancellation. Only $\mathcal C$ contains eligible outgoing states.
\end{algorithmic}
\end{algorithm}

\subsection{Parallel task scheduler and event-processing algorithms}
\label{subsec:parallel-tasks}
A task checkpoint is the tuple
\begin{equation}\label{eq:task-checkpoint}
\begin{aligned}
 \mathsf X_w=(&\mathsf{id}_{\rm model},\mathsf{id}_{\rm frame},
 \mathsf{id}_{\rm path},\tau_w,p_w,
 \lambda_w,\mu_w,\nu_w^\ell,\nu_w^u,\bar g_w,\\
 &W_{E,w},W_{I,w},W_{\ell,w},W_{u,w},F_w,N_w,
 \Pi_w,L_w,B_w,\eta_{{\rm fac},w},\eta_{{\rm bwd},w},\\
 &\dot p_w,\dot\zeta_w,\dot{\bar g}_w,
 \mathcal R_w,B_{{\rm op},w}^{\rm remaining},\mathsf{status}_w).
\end{aligned}
\end{equation}
Here $\mathcal R_w$ is the task's ordered event list. Its entries are
$(\mathsf{type},\mathsf{rowID},\alpha,\dot s)$ for primal events or
$(\mathsf{type},\mathsf{rowID},\alpha,\dot\mu)$ for dual events, with
the bound orientations of \eqref{eq:bound-ratios} and
\eqref{eq:reduced-cost-signs}. A new sensitivity solve replaces this list;
it is never shared across differing working sets. The checkpoint
stores the complete arrays and their selectors. The task
registry key uses model/frame/path IDs, parameter, primal vector, all dual values,
selectors and factorization version. Exact identity reuses a task ID. A
numerical merge is allowed only for the same immutable path data and under
documented tolerances for every continuous entry; it records which checkpoint
is retained. Duplicate detection compares the complete immutable path and
checkpoint data, including the artificial path identity. Before an operation
consumes its budget, retain the last verified checkpoint as immutable data;
perform solves, factor changes and full-array updates on a private copy.
Publish the new checkpoint only after all post-operation tests pass. A
budget or deadline return therefore contains no half-completed exchange.

A working set remains in the running family only while it represents the
current QP at its own path parameter: its face solution must be feasible
for every constraint of that QP, its multipliers must have the correct
signs, and its rows must be independent. Let $\mathcal A_w^0$ be the catalogue
of general rows and oriented bound sides with zero slack at that
checkpoint. The algebraic relevance predicate is
\begin{equation}\label{eq:working-set-relevance}
\begin{gathered}
 \mathsf{Relevant}_v(w)=1\quad\Longleftrightarrow\\
 (\mathsf{id}_{\rm model},\mathsf{id}_{\rm frame})_w=(v,\mathsf{id}_{E_v}),\\
 \text{the immutable path data and factor key match},\\
 \|F_{W_w}\|_{\rm scaled}\leq\epsilon_{\rm face},\quad
   r_{{\rm KKT},w}\leq\epsilon_{\rm QP},\\
 \mathsf{IDs}(W_{I,w},W_{\ell,w},W_{u,w})\subseteq\mathcal A_w^0,\\
 \operatorname{rank}A_{W_w}=m_w,\quad
   \operatorname{In}(K_{W_w})=(n_{F,w},m_w,0).
\end{gathered}
\end{equation}
Activity, rank and residual identities in this display use their declared
verified tolerances in floating point; $\epsilon_{\rm face}$ is the
scaled tolerance for \eqref{eq:face-kkt}. The selected multiplier is always
the unique array from \eqref{eq:working-dual-unique}. Before a positive
homotopy step, its one-sided forward relevance additionally requires
\begin{equation}\label{eq:working-set-forward-relevance}
 \dot s_{i,w}\geq0\quad
       (i\in\mathcal A_w^0\setminus\mathcal W_{\rm ineq}(X_w)),\qquad
 \dot\zeta_{j,w}\leq0\quad
       (j\in\mathcal W_{\rm ineq}(X_w),\ \zeta_{j,w}=0).
\end{equation}
Here bound sides use their oriented multipliers; equivalently their tests
are the reduced-cost tests in \eqref{eq:reduced-cost-signs}. The
zero-slack/zero-multiplier classifications are exactly those used in the
zero-distance event scan. At $\tau_w=1$ the forward test is unnecessary.

An active-set change can invalidate
\eqref{eq:working-set-forward-relevance} even though the old face equations
still hold at the event. The event handler then retires that incoming
selector from the running family and commits only its verified outgoing
selector or exchange children. The saved parent is retained as history and
as an offline seed. If no outgoing
representation passes the relevance tests, the task is returned as
\texttt{RETIRED\_IRRELEVANT} with the failed condition. An earlier rank,
factor or residual failure retains its more specific terminal status.
Retirement is evaluated at the task's own immutable path and parameter;
other relevant auxiliary paths continue independently.

The parallel family at logical generation $g$ is
\begin{equation}\label{eq:parallel-working-family}
  \mathcal T_v^{(g)}=\{\mathsf X_w:w\in\mathcal Q_v^{(g)}\},\qquad
 \mathbb W_v^{(g)}=\{(w,W_w,\tau_w):w\in\mathcal Q_v^{(g)}\}.
\end{equation}
Here $\mathcal Q_v^{(g)}$ is the index set of the tasks of generation $g$;
the set $\mathcal Q_v$ of all admitted tasks, used below, is their union over
the generations. Each running worker owns its selector, dual array, parameter $\tau_w$,
factor updates and event list. Workers share only the immutable nonlinear
model and constraint catalogue; immutable matching factors may be shared
until the first update makes a private copy. No worker's pivot changes
another worker's selection or dual. In particular, two tasks with the same
primal point and different relevant working sets are retained in
\eqref{eq:parallel-working-family} and continue independently. The selected
nonlinear endpoint is chosen only after their barrier tests.

The scheduler uses logical generations for reproducible admission. Initial
hypotheses form generation zero. Workers in one generation run concurrently
until they return an endpoint or stopped record and may publish certified
child proposals immediately. At the generation barrier, merge proposals
only under the path/checkpoint rule above, sort them by
$(\mathsf{parentID},\mathsf{eventID},\mathsf{removedID},\mathsf{incomingID})$,
and admit the first proposals fitting the remaining global task cap. Their
ordered lineage defines their stable IDs. Run the admitted next generation
in parallel. Repeat until the queue is empty or the deadline closes the
barrier. A rejected proposal is logged as \texttt{CAP\_REACHED}. This policy
makes the cap and pivot ordering independent of thread completion order;
wall-clock deadline cancellation is recorded separately.
Equivalently, with $\mathsf{Child}(w)$ the certified proposals produced by
Algorithm~\ref{alg:event-handler}, the next index set is the stable-ID
ordering of
\begin{equation}\label{eq:parallel-generation-update}
 \mathsf{Admit}_{N_{\rm task}-|\mathcal Q_v|}
 \left(\mathsf{Dedup}\!\left(
              \bigcup_{w\in\mathcal Q_v^{(g)}}\mathsf{Child}(w)\right)\right).
\end{equation}
Here $\mathcal Q_v$ contains every task already admitted for model $v$;
deduplication uses the full identity rule above, and $\mathsf{Admit}_k$
takes the first $\max\{0,k\}$ proposals in the declared lineage ordering.
Terminal endpoints remain in $\mathcal Q_v$ for nonlinear comparison.

\begin{algorithm}[H]
\caption{One parallel working-set QP homotopy worker}
\label{alg:one-path-worker}
\begin{algorithmic}[1]
\Require A verified start checkpoint from \eqref{eq:start-qp-data}, target $\tau_f=1$,
         factorization, stable row IDs, registry and per-task operation budget.
\State Verify \eqref{eq:face-kkt}, signs, feasibility, reconstruction, backward error and expected inertia.
\If{a check fails} \State Return the full checkpoint with its failed-test status. \EndIf
\While{$\tau<\tau_f$ and the worker budget remains}
 \State Solve the sensitivity system \eqref{eq:path-sensitivity} using \eqref{eq:lblt-solve}.
 \State Compute all general-row, bound, multiplier and reduced-cost ratios in \eqref{eq:general-block-ratio}--\eqref{eq:bound-ratios}.
 \State Test \eqref{eq:working-set-relevance}; return the saved checkpoint with its specific failed-test status if the predicate fails.
 \If{\eqref{eq:working-set-forward-relevance} fails}
  \State Call Algorithm~\ref{alg:zero-distance-closure} with the single seed $\{\mathsf X_w\}$ and the same remaining budget.
  \State If there is one completed output, retire the incoming selector, commit that output and restart this while-iteration at unchanged $\tau$. If there are several completed outputs, publish them as distinct child proposals, retire the parent and return \texttt{EXCHANGE\_HANDOFF}. If there is no completed output, return the saved parent with its budget/deadline or specific failed-test status. Retain every stopped/pending branch record.
 \EndIf
 \State Compute $\Delta\tau$ by \eqref{eq:min-event-step}; advance once by \eqref{eq:advance-checkpoint}.
 \If{$\tau=\tau_f$}
  \State Verify target QP feasibility, stationarity, dual signs, complementarity, reconstruction, backward error and expected inertia.
  \State Return \texttt{TARGET\_QP} only if every test passes and $r_{\rm KKT}\leq\epsilon_{\rm QP}$; otherwise return the last verified checkpoint with \texttt{STOPPED\_RESIDUAL} or \texttt{STOPPED\_FACTOR} and the failed component.
 \Else
  \State Process tied zero-ratio dual releases from Section~\ref{subsec:path-events} before primal blockers, in stable ID order; after each pivot recompute the sensitivity and every zero-distance event.
  \State For a blocker, call Algorithm~\ref{alg:event-handler}. If it returns \texttt{EXCHANGE\_HANDOFF} or another terminal status, return the parent checkpoint and stop this worker; otherwise retire the old selector, commit the outgoing selector and recompute every zero-distance event.
 \EndIf
\EndWhile
\State Return the last full checkpoint as \texttt{BUDGET\_EXHAUSTED} if the operation budget expires.
\end{algorithmic}
\end{algorithm}

\begin{algorithm}[H]
\caption{Blocking-row event handler with factor update and parallel exchange}
\label{alg:event-handler}
\begin{algorithmic}[1]
\setlength{\itemsep}{2pt}
\Require Event checkpoint $\mathsf X$, incoming blocker $l$, shared task registry and the remaining operation budget.
\State Save the parent; form \eqref{eq:full-working-catalogue}; set the private seed collection $\mathcal S\gets\varnothing$.
\State Test independence of $a_{l,F}$ against $A_W$ using the declared rank threshold.
\If{the incoming row is certified independent}
 \If{$l$ is a general or geometric row}
  \State On a private copy apply \eqref{eq:border-general-row}--\eqref{eq:addcon2}; insert its ID and zero multiplier.
 \Else
  \State On a private copy select the incoming bound side; remove its coordinate from $F$ by \eqref{eq:coordinate-border} in reverse and recover its zero event multiplier by \eqref{eq:reduced-cost-map}.
 \EndIf
 \State Rebuild \eqref{eq:working-kkt}; verify rank, full KKT residual, reconstruction and expected inertia, with one full refactorization if required. Put a verified copy in $\mathcal S$; otherwise retain its specific failed-test record.
\Else
 \State Solve for $\gamma$ and lift to $\Gamma$ by \eqref{eq:dependent-row-identities}; verify both full-normal and event-right-side identities. If dependence is unresolved, return the parent with \texttt{STOPPED\_RANK}.
 \For{each removable inequality or bound ID $j$ with $\Gamma_j\ne0$, in stable ID order}
  \State Calculate \eqref{eq:exchange-ratio}; retain $j$ only when $\delta_j\leq0$ and every test \eqref{eq:exchange-sign-test} passes.
  \State For a retained $j$, copy the parent, form \eqref{eq:exchange-working-set}, rebuild the affine lift and verify the final factor, independence, reduced-cost recovery and full KKT residual. Put a verified copy in $\mathcal S$; otherwise retain its specific failed-test record.
 \EndFor
\EndIf
\State Compute $(\mathcal C,\mathcal P,\mathcal F)$ from Algorithm~\ref{alg:zero-distance-closure} with seeds $\mathcal S$ and the same remaining budget. No intermediate selector is published.
\If{the initial addition was independent and $\mathcal C$ contains exactly one state}
 \State Retire the old selector, commit that completed outgoing state and return \texttt{CONTINUE}, retaining any stopped sibling records.
\Else
 \State Publish each state in $\mathcal C$ as a child with the parent's path ID, or attach its registered exact equivalent, in stable lineage order; retain $\mathcal P$ and $\mathcal F$ as stopped event history.
 \State Retire the parent selector. If $\mathcal C\ne\varnothing$, return \texttt{EXCHANGE\_HANDOFF}; if $\mathcal P\ne\varnothing$, return its budget/deadline status; otherwise return \texttt{RETIRED\_IRRELEVANT} with reason \texttt{NO\_CERTIFIED\_EXCHANGE} or the specific failed certificate.
\EndIf
\end{algorithmic}
\end{algorithm}

A worker may publish a valid child proposal before its siblings finish;
the scheduler admits the next generation and coordinates nonlinear selection. For a simultaneous event, the
worker first advances once to the common ratio; it then applies the tied zero-ratio dual releases in stable ID order,
followed by incoming primal blockers in stable ID order. A zero multiplier
whose forward derivative preserves its sign is retained. After each change it recomputes the factorization, all KKT residuals,
and the complete zero-ratio list. This deterministic rule defines one
reproducible traversal. If another pivot ordering is desired, each admissible
ordering is a separate registered child. A branch without a forward-feasible
pivot returns an explicit stopped checkpoint at $\Delta\tau=0$.

\subsection{Synchronized nonlinear reevaluation and output selection}
\label{subsec:nonlinear-recheck}
A synchronization barrier is the point at which all tasks of the current
model have terminated. At a synchronization barrier only the target-QP endpoints in
$\mathcal E_v$ of \eqref{eq:eligible-qps} are eligible for the online
nonlinear acceptance test. A stopped intermediate checkpoint may supply an offline seed. Every endpoint stores the
same physical model base $r_v$ and frame $E_v$, and forms
\begin{equation}\label{eq:checkpoint-retraction}
 r_e=R_{r_v}(E_vp_e).
\end{equation}
The endpoint is admitted when the retraction is in its declared domain and
$\|E_vp_e\|\leq\Delta_{\rm loc}$. The worker reevaluates the original objective,
constraints $C,D,b_{\rm smp}$ and their differentials at $(r_e,\hat t)$.
The original physical data remain fixed. The residual anchors of model $v$
also remain fixed for this comparison; any change of residual basis uses
the inverse-transpose multiplier transformation in Section~\ref{sec:parametric-problem}.
A new residual anchor is introduced only with a new nonlinear model and
fresh derivatives.

Write $b=\operatorname{col}(D,b_{\rm smp})$. For endpoint task $e$, first
recover every eliminated selected multiplier and expand the complete array
by stable row ID. Denote this operation by $\mathsf{Expand}_{W_e}$; it uses
\eqref{eq:working-dual-unique}, \eqref{eq:reduced-cost-map}, and the exact
condensing/normal recovery equations when those eliminations were used.
Let $\mathsf T_e$ be the declared change from the model's residual
coordinates to the coordinates used for reporting at $r_e$. The endpoint
multiplier is
\begin{equation}\label{eq:endpoint-working-dual}
 \zeta_e=\zeta_{Q,e}
 :=\mathsf T_e^{-\top}\mathsf{Expand}_{W_e}
                  (\zeta_{W_e},\nu_{N_e})
 =\operatorname{col}(\lambda_e,\zeta_{b,e}).
\end{equation}
For the fixed residual anchors and frames of the current model,
$\mathsf T_e=I$. A different constant residual basis uses precisely the
inverse-transpose transformation in
Section~\ref{sec:parametric-problem}. Changing a tangent frame acts on the
primal tangent vector and derivative matrices; a scalar constraint's
multiplier is unchanged. All row scalings and bound orientations are undone
before assembling the original signed array. Nonworking inequalities have
zero multiplier, and an exchange supplies the unique array of its new
working set. Thus \eqref{eq:endpoint-working-dual} carries the task's dual
through retraction without selecting another multiplier at the nonlinear
point.

Define the near-active catalogue
$\mathcal A_e=\{i:b_i(r_e,\hat t)\leq\epsilon_a\}$ with the same fixed
threshold for every candidate. Its dual admission test is
\begin{equation}\label{eq:endpoint-dual-admission}
 \bigl\|[\zeta_{b,e}]_+\bigr\|_\infty\leq\epsilon_{\rm dual},\qquad
 \bigl\|(\zeta_{b,e})_{\mathcal A_e^c}\bigr\|_\infty
                                      \leq\epsilon_{\rm dual}.
\end{equation}
The exact-arithmetic convention is $\epsilon_{\rm dual}=0$; a declared
floating-point tolerance is part of the common numerical budget. No entry
is clipped or reassigned to satisfy this test. Failure rejects this
endpoint, while other working-set tasks remain available. The stationarity
residual is evaluated directly with this same array,
\begin{equation}\label{eq:endpoint-working-stationarity}
 r_{s,e}=\operatorname{grad}f(r_e,\hat t)
       +DC(r_e,\hat t)^*\lambda_e
       +Db(r_e,\hat t)^*\zeta_{b,e}.
\end{equation}
Its value enters the original nonlinear return test and the local
inexact-SQP estimate in Section~\ref{sec:pathfollowing-theory}. Numerical
error in expanding or changing the residual frame is charged to the
corresponding transport/assembly budget. There is no additional multiplier
optimization at this barrier.

For the retained task dual, let
$S_e=\{i\in\mathcal A_e:(\zeta_{b,e})_i< -\epsilon_\mu\}$, using a
fixed sign threshold $\epsilon_\mu\geq\epsilon_{\rm dual}$. Define the
fresh supported tangent
space
\begin{equation}\label{eq:fresh-curvature-filter}
 \mathcal T_e=\{\eta\in T_{r_e}\mathcal M_H:
 DC_e[\eta]=0,\ Db_{i,e}[\eta]=0\ (i\in S_e)\}.
\end{equation}
A rank-revealing QR supplies an orthonormal basis $Z_e$ for
$\mathcal T_e$ and its verified row-rank certificate. The curvature test is
\begin{equation}\label{eq:fresh-curvature-test}
 \lambda_{\min}\bigl(Z_e^\top\operatorname{Hess}L(r_e,\zeta_e,\hat t)Z_e\bigr)
 \geq\sigma_{\min}>0,
\end{equation}
with a verified eigenvalue lower bound. When the rows defining
$\mathcal T_e$ are linearly independent, the strict form of
\eqref{eq:fresh-curvature-test} is equivalent to an inertia condition: with
$A_{S_e}$ the matrix of these rows in the tangent frame and $m_{S_e}$ their
number,
\begin{equation}\label{eq:fresh-curvature-inertia}
 \operatorname{In}
 \begin{pmatrix}H_e-\sigma_{\min}I&A_{S_e}^\top\\A_{S_e}&0\end{pmatrix}
 =(n,m_{S_e},0),
\end{equation}
where $H_e$ is the matrix of $\operatorname{Hess}L(r_e,\zeta_e,\hat t)$ in
the frame. This follows from \eqref{eq:kkt-inertia-nullspace} applied to
$H_e-\sigma_{\min}I$. The test can therefore be certified by one $LBL^\top$
factorization, without computing $Z_e$ or an eigenvalue.
Ambiguous multipliers are treated as
weak and omitted from $S_e$, making the tested subspace larger. Weak active inequalities remain in the primal feasibility tests; their
normals are omitted from the equality rows defining this curvature space. When
$\mathcal T_e=\{0\}$, \eqref{eq:fresh-curvature-test} is vacuous, but all
first-order, primal, dual and physical admissibility tests still apply.
The successful-tracking assumptions in Section~\ref{sec:pathfollowing-theory}
include uniform positivity on this same supported space, with a margin for
the declared thresholds and numerical errors.

A checkpoint is eligible only if the scaled residual vector
\begin{equation}\label{eq:original-kkt-residuals}
\begin{aligned}
 \chi_f&=\max\{\|C_e\|_\infty,\|[-b_e]_+\|_\infty\},\\
 \chi_s&=\|\operatorname{grad}f_e+DC_e^*\lambda_e+Db_e^*\zeta_{b,e}\|_\infty,\\
 \chi_c&=\max\{\|b_e\odot\zeta_{b,e}\|_\infty,\|[\zeta_{b,e}]_+\|_\infty\}
\end{aligned}
\end{equation}
meets its declared tolerances, the fresh test \eqref{eq:fresh-curvature-test}
passes, and the first control satisfies the independent physical admissibility
checks. The scaled infinity norms use declared orthonormal frames and units.
Eligible points are compared by the scaled score
\eqref{eq:nonlinear-candidate-score} after all individual tests pass.

A new nonlinear correction round at the same measurement, if configured,
uses the complete selected primal--dual point as its new base, reevaluates
derivatives, and restarts each artificial QP at $\tau=0$. It defines a new
model version for the same physical data. A further round may also be
started when no endpoint was accepted: the eligible endpoint with the
smallest first component of \eqref{eq:nonlinear-candidate-score} then
becomes the new base, provided that this component has decreased by a
declared factor relative to the previous base, and the offline procedure is
invoked after a declared maximal number of rounds. Single-round RTI uses only $j=0$. The inexact SQP
recursion in Section~\ref{sec:pathfollowing-theory} states the local
conditions under which that one round meets the return tolerances.

\subsection{Synchronization barrier and nonlinear model refresh}
\label{subsec:synchronization-barrier}
The task queue operates on one immutable nonlinear model snapshot. Let its
version be $v$ and the associated manifold base be $r_v$. Every generated
QP stores $(v,r_v,E_v)$, where $E_v$ is the frame and tangent-coordinate
ordering. Workers may run concurrently only while these IDs remain fixed.
A measurement arriving during a homotopy is queued as the next physical
update; it cannot change the target model or any in-flight path data. Each
auxiliary path has its own scalar $\tau$ and immutable start differences;
its exchange children share that path. The task family shares the physical
target at $\tau=1$; intermediate QPs are specific to each auxiliary path.
In other words, different hypotheses start from different start QPs and
pass through different intermediate QPs, but all of them end at the same
target QP.

A task has one of the statuses
\texttt{QUEUED}, \texttt{RUNNING}, \texttt{TARGET\_QP},
\texttt{STOPPED\_START}, \texttt{STOPPED\_FACTOR}, \texttt{STOPPED\_RANK},
\texttt{STOPPED\_RESIDUAL},
\texttt{STOPPED\_SIGN}, \texttt{EXCHANGE\_HANDOFF},
\texttt{RETIRED\_IRRELEVANT},
\texttt{NO\_CERTIFIED\_EXCHANGE}, \texttt{EVENT\_CYCLE},
\texttt{BUDGET\_EXHAUSTED}, or \texttt{CANCELLED}. Only
\texttt{TARGET\_QP} supplies a QP endpoint. Every stopped status retains its
last primal-dual checkpoint; a worker never returns only a status string.
The scheduler enforces the operation budget as an integer count of KKT solves
and factor changes per task, plus a separate global cap $N_{\rm task}$ on
admitted task IDs. Budget exhaustion is deterministic for a fixed event
sequence. It returns the last verified checkpoint without converting an
incomplete path into a target solution.

At a barrier, let $\mathcal Q_v$ contain every admitted task with model
version $v$. The barrier closes only when every task is terminal or when the
predeclared global deadline expires. The barrier includes all logical
generations admitted under the preceding rule. The eligible QP set is
\begin{equation}\label{eq:eligible-qps}
\begin{aligned}
 \mathcal E_v=\{e\in\mathcal Q_v:\;&
   \mathsf{status}(e)=\texttt{TARGET\_QP},\quad
   \mathsf{Relevant}_v(e)=1,\\
 &\eta_{\rm fac}(e)\leq\epsilon_{\rm fac},\quad
   \operatorname{In}(K_e)=(n_{F,e},m_e,0)\}.
\end{aligned}
\end{equation}
Only these endpoints are passed to the online nonlinear acceptance tests.
Retired working sets have no endpoint in $\mathcal E_v$; their immutable
records remain available as history and offline seeds. A deadline cancels
uncompleted tasks after their last immutable verified checkpoints have been
copied. Only completed target endpoints remain eligible for selection.

For every endpoint $p_e$ form the retracted candidate
$r_e=R_{r_v}(E_vp_e)$. Before nonlinear reevaluation, verify that the
retraction returned a point in the declared manifold and that the tangent
step obeys $\|E_vp_e\|\leq\Delta_{\rm loc}$. Then evaluate all original
nonlinear constraints and the objective at $r_e$. The accepted feedback
candidate is the eligible point minimizing
\begin{equation}\label{eq:nonlinear-candidate-score}
 J(e)=\left(\max\left\{
 \frac{\chi_f(e)}{\epsilon_f},
 \frac{\chi_s(e)}{\epsilon_s},
 \frac{\chi_c(e)}{\epsilon_c},
 \frac{\max(0,\sigma_{\min}-\lambda_{\min}(Z_e^\top\operatorname{Hess}L_eZ_e))}{\epsilon_H}
 \right\},\ f(r_e),\ \mathsf{taskID}(e)\right)
\end{equation}
lexicographically, provided the first component is at most one,
\eqref{eq:fresh-curvature-test} passes, and the physical control
admissibility test passes. The first control block is extracted only after
this selection. Its applied value is
\begin{equation}\label{eq:verified-control-output}
 u^{\rm applied}_{\hat t}=\mathcal A_{\rm phys}
                    (u_{0,e},v_{0,e}),
\end{equation}
where the physical actuator map is the same map used by the OCP's first
interval dynamics. Check the applied command, including its physical limits,
before returning it. If an additional saturation or actuator transformation
would change the modeled command beyond its declared tolerance, reject the
candidate and rebuild the corresponding model. Log both the optimizer
coordinates and the applied command.

If no endpoint passes, set \texttt{NO\_ACCEPTED\_ENDPOINT} and choose an
offline seed from complete verified checkpoint records. For each record
whose retraction is defined in the declared recovery domain, evaluate its
original nonlinear residuals with its complete transported multiplier array
and assign
\begin{equation}\label{eq:backup-seed-score}
 J_{\rm seed}(e)=\left(
 \max\{\chi_f(e)/\epsilon_f,\chi_s(e)/\epsilon_s,\chi_c(e)/\epsilon_c\},
 \operatorname{dist}(r_e,r_v),\mathsf{taskID}(e)\right).
\end{equation}
Choose its lexicographic minimum among finite scores, using the original
OCP residuals to rank recovery seeds. If no such checkpoint exists, use the
shifted incoming manifold point when it belongs to the admitted recovery
seed set; otherwise use the designated admitted initialization seed of
Assumption~\ref{ass:offline-backup}. The whole primal--dual seed is passed
together.

Invoke the offline procedure of Appendix~\ref{app:pdsqp} under
Assumption~\ref{ass:offline-backup}. At initialization this computation
finishes before the first control handoff. During recovery the retained
certified plan, that is, the last accepted trajectory, whose remaining
controls may still be applied for a limited time, supplies admissible
controls over its declared validity interval, and the backup finishes
within that interval. Each pending command
is checked against the retained plan's state/prediction uncertainty
certificate and expiry. Its successful return is verified
against the original residual, supported-curvature and physical-output
tests above before it supplies a command. A reported failure or an invalid
return produces \texttt{OFFLINE\_RECOVERY\_FAILURE} and terminates this
feedback invocation before the new command-output and plant-advance
operations; its failure record retains the existing certified plan and its
expiry. Pending recovery commands continue only within that certified
interval.
Under the reliable-backup assumption a recovery call returns verified
success. Thus every successful feedback invocation has a single explicit
applied command, obtained either from an accepted online endpoint or from a
verified offline return.

After that command is applied and the physical state advances, increment
the model version $v$, shift the retained complete solution, substitute the
new measured node-zero state, build the new tangent frame, and reevaluate
derivatives. Retire every old-version working set from the online family.
Transfer its stable row IDs only as a proposed new hypothesis, and admit
it again only after \eqref{eq:start-qp-data} and
\eqref{eq:working-set-relevance} pass for the new model. Start new QP
homotopies at $\tau=0$. No factor, sensitivity,
event ratio or working selector is reused across model versions without
validating its exact key and coordinate transformation.

\begin{algorithm}[H]
\caption{Barrier selection and next-model transition}
\label{alg:barrier-selection}
\begin{algorithmic}[1]
\Require Immutable model version $v$, terminal task records $\mathcal Q_v$,
         declared nonlinear tolerances, and the reliable offline interface of Assumption~\ref{ass:offline-backup}.
\State Form $\mathcal E_v$ using every condition in \eqref{eq:eligible-qps}.
\For{each $e\in\mathcal E_v$}
 \State Retract by \eqref{eq:checkpoint-retraction}; reject if the domain or local-radius test fails.
 \State Recover and expand the unique working-set dual by \eqref{eq:endpoint-working-dual}; test \eqref{eq:endpoint-dual-admission} without changing it.
 \State Evaluate \eqref{eq:original-kkt-residuals}, \eqref{eq:fresh-curvature-test}, and the physically applied command \eqref{eq:verified-control-output}.
\EndFor
\If{at least one endpoint passes every test}
 \State Retain the entire endpoint minimizing \eqref{eq:nonlinear-candidate-score}.
\Else
 \State Select the complete seed minimizing \eqref{eq:backup-seed-score}, or the declared shifted incoming seed; check its recovery-domain membership.
 \State During recovery use the retained certified plan within its validity interval; at initialization finish before the first handoff. Invoke Appendix~\ref{app:pdsqp} under Assumption~\ref{ass:offline-backup} and verify its complete return by the same original tests.
 \If{recovery fails or a return test fails}
  \State Return \texttt{OFFLINE\_RECOVERY\_FAILURE}, the retained certified plan and its expiry; issue pending controls only while that plan's validity certificate holds.
 \EndIf
 \State Retain the complete verified offline solution.
\EndIf
\State Extract and check \eqref{eq:verified-control-output}; return/log the command and apply it once.
\State After the physical advance, increment $v$, shift the complete solution, substitute the new measurement, and rebuild $r_v,E_v,g,H,A,b$.
\State Retire all old-version selectors; invalidate their event ratios and nonmatching factors. Admit transferred working sets only after the new start and relevance tests; initialize their paths at $\tau=0$.
\end{algorithmic}
\end{algorithm}

\subsection{Per-event cost, operation budget, and termination statuses}
\label{subsec:operation-costs}
Let $n_F$ be the number of free tangent coordinates, $m_W$ the number of
selected equality/general rows, $N_K=n_F+m_W$, $m_I$ the number of monitored
general inequalities, and $n_b$ the number of bounded nonstructural
coordinates. There are at most $2n_b$ component-bound sides. With dense
storage, a complete symmetric factorization costs $O(N_K^3)$ arithmetic
and $O(N_K^2)$ memory; one solve with the factors costs $O(N_K^2)$, and
$k$ right sides cost $O(kN_K^2)$. The event scan evaluates $m_I+2n_b$
primal slopes and at most $m_W+|N|$ inequality/reduced-cost slopes.
For dense rows it costs $O((m_I+n_b)N_K)$ after the sensitivity solve;
for sparse rows it costs $O(\operatorname{nnz}(A_I)+n_b)$ plus that solve.

A scalar border and the cached defect certificate
\eqref{eq:border-reconstruction-residual} cost $O(N_K^2)$. A fresh product
$LBL^\top$ for a full reconstruction costs $O(N_K^3)$ and is counted as such.
For the concrete deletion in \eqref{eq:delete-local-schur}, let $r$ be the
retained prefix size and $q$ the affected trailing size. Triangular solves,
Schur assembly, trailing factorization and reconstruction have the bound
\begin{equation}\label{eq:trailing-refactor-cost}
 C_{\rm tail}(r,q)\leq
       c_1r^2q+c_2rq^2+c_3q^3+c_4N_K^2.
\end{equation}
A deletion is recorded as a quadratic update only when the actual update
scheme and affected sizes justify that bound. A refactorization of the
whole retained matrix is logged as a full refactorization. Sparse costs
are determined by the actual fill and elimination tree; the same full
residual and inertia conditions apply.

For a worker with $N_e$ accepted events, $N_u$ scalar updates, $N_f$ full
factorizations, and $N_R$ standalone full reconstruction products, a dense
linear-algebra ledger is
\begin{equation}\label{eq:task-cost-ledger}
\begin{aligned}
 C_{\rm lin,w}\leq{}&N_fc_fN_K^3+N_Rc_RN_K^3
       +(N_e+N_u)c_uN_K^2\\
       &+\sum_{d\in\mathcal D_w}C_{\rm tail}(r_d,q_d)
          +N_ec_s(m_I+n_b)N_K .
\end{aligned}
\end{equation}
Here $\mathcal D_w$ contains proper trailing refactorizations; full ones
are already counted in $N_f$. A scalar-update certificate using
\eqref{eq:border-reconstruction-residual} is included in $c_u$ and is not
also counted in $N_R$. The constants are backend operation bounds or
measured accounting constants. Each count is obtained directly from the
executed operations.

Barrier work is additional to \eqref{eq:task-cost-ledger}. In particular,
log its retractions, original residual/derivative/Hessian evaluations,
working-dual expansion/frame transformations, nullspace constructions and
curvature solves:
\begin{equation}\label{eq:barrier-cost-ledger}
 C_{\rm barrier}=\sum_{e\in\mathcal E_v}
  \bigl(C_{R,e}+C_{\rm assembly,e}+C_{\rm dual\,map,e}
                     +C_{Z,e}+C_{H,e}\bigr).
\end{equation}
A dense eigensolve on the supported space has cubic cost in its dimension;
its work is not represented as a triangular KKT solve. Offline recovery
has its own ledger. With at most $N_{\rm task}$ admitted workers and a
per-worker integer budget $B_{\rm op}$ counting solves, factor changes
and standalone certificates, the dense factor/certificate work has the
loose cubic cap $N_{\rm task}B_{\rm op}c_{\max}N_{K,\max}^3$.
Barrier assembly and recovery costs are accounted for separately.
Parallel execution distributes the same event equations and arithmetic
among the available workers.

A worker terminates with a certified target endpoint, a certified exchange
handoff, retirement of an irrelevant working set, a failed
start/rank/factor/residual/sign test, no certified exchange,
a detected zero-distance cycle, budget exhaustion, or deadline cancellation.
Each terminal status retains its latest immutable verified checkpoint and
the failed component or outgoing child references. Target success requires
$\tau=1$ and every target certificate. An exchange handoff records the
parent transition and its children. The scheduler stops admitting children when
$|\mathcal Q_v|=N_{\rm task}$; later candidate exchanges are recorded as
\texttt{CAP\_REACHED} with their parent ID and are not silently dropped.
The finite cap means the queue may omit some admissible multiplier
representations. Accordingly, Corollary~\ref{cor:parallel-theory} applies to
accepted tasks and not to completeness of the configured search.

\subsection{Operation ledger and implementation obligations}
\label{subsec:pathfollowing-ledger}
Every task records the following per event: model, frame and path IDs; pre/post
$\tau$; event ID and event type; primal and all dual arrays; primal slack and
reduced-cost arrays; working selectors; permutation and $LBL^\top$ block
structure; cached or full factor reconstruction and solve backward errors;
signed inertia; row-dependence residuals and all $\gamma,\gamma_N,\Gamma,\delta_j$ exchange data; child IDs;
and final status. The target-QP record also includes the direct KKT residual
and the count of factor updates versus full refactorizations. These records
make an update reproducible and distinguish an actual factor update from a
factor rebuilt from the current matrix.

The concrete procedure is therefore: (i) identify the fixed tangent frame and
structural substitutions; (ii) construct and verify a start QP; (iii) factor
the selected KKT system; (iv) solve \eqref{eq:path-sensitivity}; (v) compute
every primal and dual event ratio; (vi) advance once; (vii) update or downdate
$LBL^\top$ and verify it; (viii) when a row is dependent, form each signed
exchange candidate from \eqref{eq:dependent-dual-family}; (ix) schedule all
valid distinct children; and (x) synchronize, retract and reevaluate before
using a candidate at the nonlinear level. Each operation uses the displayed formulas, its immutable model/path data,
and its declared algebraic and nonlinear return certificates.
\section{Riemannian Pathfollowing Theory Under Degeneracy}
\label{sec:pathfollowing-theory}

We analyze the shift, frozen-QP solve, retraction, and working-set
multiplier transport of Section~\ref{sec:pathfollowing}. The local analysis combines the
strict-working-set perspective of Wright~\cite{wright2002modifying}, the
dual-degenerate parametric analysis of Kungurtsev and
J\"aschke~\cite{kungurtsev2017predictor}, and the parametric SQP framework of
Kungurtsev and Diehl~\cite{kungurtsev2014sequential}. The assumptions address
the regularity of the nonlinear solution path, solvability of the auxiliary
QP homotopies, and availability of a verified recovery solution.
All derivatives use the product Levi--Civita connection. Tangent vectors
and residual covectors are compared in fixed smooth local identifications.

\subsection{Problem parameter, qualifications, and second-order condition}
Consider the fixed-layout problem
\begin{equation}\label{eq:theory-problem}
 \min_{r\in\mathcal M}f(r,t)
 \quad\text{subject to}\quad C(r,t)=0,\qquad b(r,t)\geq0,
 \qquad t\in I=[t_0,t_f].
\end{equation}
Here \(t\) indexes physical problem data, including the measurement and
reference data on the interval under consideration. The artificial
parameter \(\tau\in[0,1]\) of \eqref{eq:qp-path-data} parametrizes a frozen
QP and is a different variable. Scalar notation simplifies the sensitivity
formulas. For a measurement parameter \(\vartheta\in\mathbb R^q\), replace
\(\partial_t\mathscr F\,h\) by
\(D_\vartheta\mathscr F[\Delta\vartheta]\) and \(|h|\) by
\(\|\Delta\vartheta\|\) in the estimates below.

Use the plus Lagrangian \(L=f+\lambda^\top C+\mu^\top b\), with \(\mu\leq0\).
In this section \(b\) denotes all inequality margins and \(\mu\) all their
multipliers; they are \(\operatorname{col}(D,b_{\rm smp})\) and
\(\zeta_b=\operatorname{col}(\mu,\nu)\) in the notation of
Section~\ref{sec:parametric-problem}.
For a selected continuous primal solution \(r_*(t)\), define
\[
 \mathcal A(t)=\{i:b_i(r_*(t),t)=0\},\qquad
 \mathcal A_-(t,\zeta)=\{i\in\mathcal A(t):\mu_i<0\},
 \qquad \zeta=(\lambda,\mu)\in\Lambda(t).
\]
The full multiplier set at that primal solution is
\begin{equation}\label{eq:multiplier-set}
 \Lambda(t)=\left\{(\lambda,\mu):
 \begin{array}{l}
 \operatorname{grad}f+DC^*\lambda+Db^*\mu=0,\\
 \mu\leq0,\quad \mu_i=0\quad(i\notin\mathcal A(t))
 \end{array}\right\}.
\end{equation}
For a fixed independent working representation, the adjoint of its
selected row operator is injective. At a selected KKT root, its
stationarity equation therefore determines a unique multiplier on those
rows; extending by zero on omitted inequalities gives the task's complete
multiplier. The set
\(\Lambda(t)\) may nevertheless contain several such working-set
representations. Parallel tasks retain their own selected multipliers,
as in \eqref{eq:working-dual-unique}, including during nonlinear reevaluation.

The critical cone and multiplier-dependent supported tangent space are
\begin{align}
 \mathcal C(t)&=\{d:DC[d]=0,\ Db_i[d]\geq0\ (i\in\mathcal A(t)),
                         \ \langle\operatorname{grad}f,d\rangle=0\},
 \label{eq:theory-critical-cone}\\
 \mathcal T_-(t,\zeta)&=\{d:DC[d]=0,\
                  Db_i[d]=0\ (i\in\mathcal A_-(t,\zeta))\}.
 \label{eq:theory-supported-space}
\end{align}
The KKT equation gives
\(\mathcal C(t)\subseteq\mathcal T_-(t,\zeta)\) for every
\(\zeta\in\Lambda(t)\).

\begin{definition}[Riemannian MFCQ]
\label{def:theory-mfcq}
At \((r_*,t)\), \(DC(r_*,t)\) is onto its retained equality residual
space and there exists \(d\in T_{r_*}\mathcal M\) such that
\(DC[d]=0\) and \(Db_i[d]>0\) for every active inequality. Uniform MFCQ
means that the equality differential has a uniformly bounded right inverse
and such a direction can be chosen with \(\|d\|\leq1\) and
\(Db_i[d]\geq\gamma_b>0\).
\end{definition}

\begin{definition}[Riemannian relaxed constant-rank qualification]
\label{def:rcrq}
RCRQ, also called RCRCQ, holds at \((r_*,t)\) if, for every
\(J\subseteq\mathcal A(t)\), the rank of
\([DC(r,t);(Db_i(r,t))_{i\in J}]\) is constant as \(r\) varies in a
neighborhood of \(r_*\), with \(t\) fixed. Constraint residual
representatives are fixed on that neighborhood~\cite{minchenko2011relaxed}.
The uniform path assumption below additionally specifies common parameter
neighborhoods and singular-value bounds for the selected independent rows.
\end{definition}
In words, RCRQ requires that any subset of the active inequality gradients,
taken together with the equality gradients, keeps the same rank at all
nearby points: gradients that are linearly dependent at the solution remain
dependent in a neighborhood. MFCQ permits such dependence, and RCRQ
controls how it varies.

The uniform SSOSC used here is precisely the supported-space condition in
Assumption~\ref{as:sosc}: for some \(\alpha>0\),
\begin{equation}\label{eq:uniform-theory-ssosc}
 \langle d,\operatorname{Hess}_rL(r_*(t),\zeta,t)[d]\rangle
 \geq\alpha\|d\|^2
 \quad\text{for all }t\in I,\ \zeta\in\Lambda(t),\
 d\in\mathcal T_-(t,\zeta).
\end{equation}
It implies the critical-cone SOSC
\begin{equation}\label{eq:uniform-critical-ssosc}
 \langle d,\operatorname{Hess}_rL(r_*(t),\zeta,t)[d]\rangle
 \geq\alpha\|d\|^2,\qquad d\in\mathcal C(t).
\end{equation}
The supported-space condition is the one used by
\eqref{eq:fresh-curvature-test}. Weakly active inequalities retain their one-sided meaning in
\eqref{eq:theory-critical-cone}.

\subsection{Uniform neighborhoods and finite regular continuation}
\begin{assumption}[Uniform regular path and finite event partition]
\label{ass:regular-path}
The following conditions hold on the fixed-layout interval \(I\).
\begin{enumerate}
\item There is a continuous path of isolated local minimizers \(r_*(t)\).
Uniform MFCQ and RCRQ hold along this path, and its nonempty multiplier
sets are uniformly bounded. The all-multiplier SSOSC
\eqref{eq:uniform-theory-ssosc} holds.

\item A compact tube \(\mathcal U\) contains the graph of the path,
the admitted local warm starts, and every retraction and frozen-QP path
used in the local analysis. A finite atlas, smooth residual
identifications, and second-order retractions are fixed on this tube.
For \(a=f,C,b\), the joint covariant derivatives in \((r,t)\) satisfy
\begin{equation}\label{eq:uniform-data-lipschitz}
 \|\nabla^j a(r,t)\|\leq M,\qquad
 \|\nabla^j a(r,t)-\mathcal P_{r'\to r}\nabla^j a(r',s)\|
 \leq L\bigl(d_{\mathcal M}(r,r')+|t-s|\bigr),
 \quad j=1,2.
\end{equation}
These comparisons use points in a common chart. The same constants
\(M,L\) apply throughout each entire homotopy
\(\tau\in[0,1]\), and uniformly across the finitely many pieces below.
The transport of each tensor slot is included in \(\mathcal P\).
Retraction and transport derivatives are uniformly bounded. Increasing
\(L\) if necessary gives a common Lipschitz bound for the selected KKT
derivatives on the tube with bounded multipliers.

\item There is an explicitly finite partition
\[
 t_0=s_0<s_1<\cdots<s_J=t_f
\]
such that each open piece has fixed active catalogue and admits a finite
collection of relevant strict working representations. Relevance means
that the selected KKT section satisfies the full constraint and sign
conditions on that piece. A representation \(W\subseteq\mathcal A(t)\)
contains the negative support of its selected multiplier
\(\zeta_*^W(t)\), and its equality and selected inequality rows are
independent. The selected multiplier is continuous on the piece and
has finite one-sided limits that remain valid KKT multipliers in that
representation. Its primal component is \(r_*(t)\), with zero
multipliers on omitted rows. At every endpoint at least one specified
outgoing representation is valid. After an active-set change the
retained family is updated: a representation that is no longer relevant
is removed from the parallel flow, while valid outgoing representations
provide the local continuation cover.

\item In the common local frames, the selected row operators have
nonzero singular values bounded below by \(\gamma_A>0\), and the
bordered operators
\begin{equation}\label{eq:theory-working-operator}
 \mathscr K_W(t)=
 \begin{pmatrix}H_W(t)&A_W(t)^*\\A_W(t)&0\end{pmatrix},
 \qquad H_W(t)=\operatorname{Hess}_rL(r_*(t),\zeta_*^W(t),t)
\end{equation}
satisfy \(\|\mathscr K_W(t)^{-1}\|\leq\kappa\), including one-sided
limits. Structural coordinates are eliminated before these tests.
The same bounds, with fixed enlarged constants, hold for the admitted
working systems along each auxiliary QP homotopy.
\end{enumerate}
\end{assumption}

In words, the third item requires that the parameter interval splits into
finitely many pieces, on each of which the active set is constant and the
solution is described by finitely many independent working sets, each with
its own continuous multiplier.
The finite partition is an additional path hypothesis. A finite row
catalogue and a uniform Lipschitz constant alone allow event accumulation.
For example, \(\min_{x\geq0}(x-a(t))^2/2\), with
\(a(t)=e^{-1/t^2}\sin(1/t)\) for \(t>0\) and \(a(0)=0\), has smooth
data, Hessian one, and the two regular working representations
\(\varnothing\) and \(\{x=0\}\), but infinitely many active-set changes
accumulate at zero. Assumption~\ref{ass:regular-path} imposes
nonaccumulation explicitly. Horizon, contact-catalogue, and structural
changes use the transfer and restart rules of
Section~\ref{sec:parametric-problem}, on separate fixed-layout intervals.

For an independent representation with
\(\mathcal A_-(t,\zeta_*^W)\subseteq W\), one has
\(\ker A_W\subseteq\mathcal T_-(t,\zeta_*^W)\).
Thus \eqref{eq:uniform-theory-ssosc} supplies positive reduced curvature.
The explicit row and inverse bounds also control the constants in the
numerical event and factorization tests. MFCQ permits dependent active
inequalities; independence concerns the selected rows.

\subsection{Sensitivity of a selected intrinsic KKT section}
\label{subsec:riemannian-sensitivity}
Fix one piece and one representation \(W\). Let
\(c_W=\operatorname{col}(C,b_W)\), let \(\zeta_W\) contain only the
corresponding multipliers, and let \(A_W=D_rc_W\).
The selected intrinsic KKT section is
\begin{equation}\label{eq:intrinsic-face-kkt}
 \mathscr F_W(r,\zeta_W,t)=
 \begin{pmatrix}\operatorname{grad}_r
       [f(r,t)+\zeta_W^\top c_W(r,t)]\\c_W(r,t)\end{pmatrix}.
\end{equation}
At a root, its covariant derivative is
\begin{equation}\label{eq:riem-bordered-derivative}
 D_{(r,\zeta_W)}\mathscr F_W=
 \begin{pmatrix}\operatorname{Hess}_rL&A_W^*\\A_W&0\end{pmatrix}
 =\mathscr K_W .
\end{equation}
The residual-coordinate and tangent-frame transformations act on both
blocks. In an orthonormal frame, this is the symmetric KKT matrix of
Section~\ref{sec:pathfollowing}.

Writing \(z_*^W=(r_*,\zeta_*^W)\), differentiation gives
\begin{equation}\label{eq:intrinsic-path-derivative}
 \mathscr K_W(t)\dot z_*^W(t)
 =-\partial_t\mathscr F_W(z_*^W(t),t).
\end{equation}
The lower right-hand side contains the explicit derivative of the
selected constraints, and the upper one contains parameter derivatives
of their differentials contracted with the multipliers. Let \(B_t\)
bound \(\|\partial_t\mathscr F_W\|\) on the common tube. Then
\begin{equation}\label{eq:path-derivative-bound}
 \|\dot r_*(t)\|+\|\dot\zeta_*^W(t)\|
 \leq B_*:=c_0\kappa B_t ,
\end{equation}
where \(c_0\) converts the product norm to the displayed sum.
This parameter-derivative bound enters the predictor constants.

\begin{lemma}[Uniform local face sensitivity]
\label{lem:face-sensitivity}
Under Assumption~\ref{ass:regular-path}, each selected KKT root has a
unique local \(C^1\) continuation in its declared representation. Its
derivative is locally Lipschitz, uniformly on the finite collection of
pieces. There is \(C_{\rm pred}<\infty\), depending on
\(M,L,\kappa,B_t\) and the chart, retraction, and transport bounds,
such that, whenever \(t,t+h\) lie in one piece,
\begin{equation}\label{eq:predictor-remainder-bound}
 d_{\mathcal M}\bigl(r_*(t+h),R_{r_*(t)}(h\dot r_*(t))\bigr)
 +\|\zeta_*^W(t+h)-\zeta_*^W(t)-h\dot\zeta_*^W(t)\|
 \leq C_{\rm pred}h^2 .
\end{equation}
The estimate holds one-sided at an endpoint of that piece.
\end{lemma}
\begin{proof}
Apply the implicit-function theorem to
\eqref{eq:intrinsic-face-kkt} in the fixed local identifications.
The inverse bound and \eqref{eq:uniform-data-lipschitz} give
\(\dot z_*^W=-\mathscr K_W^{-1}\partial_t\mathscr F_W\).
Both factors are Lipschitz along the root, with constants controlled
by \(L,\kappa\), and \(B_*\). Taylor expansion in these identifications
has a uniform quadratic remainder. The retraction and transport bounds
give \eqref{eq:predictor-remainder-bound}. Finitely many pieces and
charts permit a common constant.
\end{proof}

\subsection{One frozen QP is an inexact Riemannian Newton step}
\label{subsec:actual-sqp-map}
Let \(\bar z=(\bar r,\bar\zeta_W)\) be a shifted estimate at the new
physical parameter \(t^+\). At this fixed parameter write
\[
 F_0=\mathscr F_W(\bar z,t^+),\qquad
 K_0=
 \begin{pmatrix}
  \operatorname{Hess}_rL(\bar r,\bar\zeta_W,t^+)&A_W(\bar r,t^+)^*\\
  A_W(\bar r,t^+)&0
 \end{pmatrix}.
\]
The tangent QP with the exact Lagrangian Hessian, linearized constraints,
and final multiplier \(\zeta_Q\) has selected equations
\begin{equation}\label{eq:sqp-newton-equivalence}
 K_0\binom{p}{\zeta_Q-\bar\zeta_W}=-F_0.
\end{equation}
Its upper equation is
\(Hp+\operatorname{grad}f+A_W^*\zeta_Q=0\), and its lower equation is
\(A_Wp+c_W=0\). Structural and working-bound elimination gives the
reduced systems in Section~\ref{subsec:reduced-kkt}. A complete
frozen-QP homotopy computes this endpoint while maintaining all QP
inequalities and multiplier signs. Independently initialized tasks
may have different auxiliary paths to this same physical target QP.

For numerical assembly, let
\(\widetilde F_0=F_0+\Delta F\) and
\(\widetilde K_0=K_0+\Delta K\). If the full step
\(u=(p,\zeta_Q-\bar\zeta_W)\) obeys
\(\|\widetilde K_0u+\widetilde F_0\|\leq\epsilon_{\rm lin}\), then
\begin{equation}\label{eq:theory-residual-budget}
 \|K_0u+F_0\|\leq
 \eta_{\rm N}:=\epsilon_{\rm lin}
       +\|\Delta F\|+\|\Delta K\|\,\|u\|.
\end{equation}
Here \(\epsilon_{\rm lin}\) is evaluated against the assembled full
system and includes factor-solve backward error. The assembly terms
include derivative, prepared-model, frame-transport, and bound-substitution
defects in these same coordinates. A Hessian approximation contributes
\[
 \bigl\|[H_{\rm used}
   -\operatorname{Hess}_rL(\bar r,\bar\zeta_W,t^+)]p\bigr\|
\]
to \(\|\Delta K\,u\|\). This term is also retained if a dual exchange
changes the selected multiplier while the QP Hessian remains frozen.
Direct residuals are measured quantities; an assembly-error bound
requires the declared derivative or model accuracy certificate.

Use
\[
 d_Z((r,\zeta),(r',\zeta'))=
 d_{\mathcal M}(r,r')+\|\zeta-\zeta'\|
\]
in common residual coordinates, and set
\(\delta=d_Z(\bar z,z_*^W(t^+))\).

\begin{lemma}[Local error of the implemented QP--retraction map]
\label{lem:one-sqp-map}
There are uniform constants \(\rho,C_{\rm N}>0\) such that, if
\(\delta\leq\rho\), \(\eta_{\rm N}\leq\rho\), the returned endpoint uses
the declared representation, and its step remains in the common tube,
then \(z_Q^+=(R_{\bar r}(p),\zeta_Q)\) satisfies
\begin{equation}\label{eq:one-sqp-error}
 d_Z(z_Q^+,z_*^W(t^+))
 \leq C_{\rm N}(\delta^2+\eta_{\rm N})
       +\epsilon_{\rm rt}.
\end{equation}
Here \(\zeta_Q\) is expressed in the common residual dual coordinates
by \eqref{eq:endpoint-working-dual}. The term
\(\epsilon_{\rm rt}\) bounds numerical retraction and residual-frame
transport error, including the corresponding error in the represented
multiplier. The original nonlinear primal, stationarity, and
complementarity residuals satisfy
\begin{equation}\label{eq:one-sqp-nonlinear-residual}
 \max\{\chi_f,\chi_s,\chi_c\}_{Q}
 \leq C_{\rm res}(\delta^2+\eta_{\rm N}+\epsilon_{\rm rt}),
\end{equation}
provided the endpoint satisfies the full QP inequalities, signs, and
complementarity to their declared residual tolerance, included in
\(\epsilon_{\rm lin}\).
\end{lemma}
\begin{proof}
For this Taylor estimate, transport tangent gradients parallel along
the radial retraction curves using the Levi--Civita connection.
The derivative at zero is then the covariant derivative \(K_0\),
including away from a root. An implemented transport with a different
first-order action contributes its defect to \(\Delta F,\Delta K\).
In retraction coordinates at \(\bar r\), let \(u_*\) be the displacement
from \(\bar z\) to \(z_*^W(t^+)\). Uniform local equivalence with \(d_Z\)
gives \(\|u_*\|\leq c_1\delta\).
Taylor expansion of the transported KKT section gives
\[
 0=F_0+K_0u_*+q,\qquad
 \|q\|\leq c_2 L\|u_*\|^2 .
\]
For sufficiently small \(\rho\), perturbation from the root gives
\(\|K_0^{-1}\|\leq2\kappa\). Subtraction and
\eqref{eq:theory-residual-budget} yield
\(\|u-u_*\|\leq2\kappa(\eta_{\rm N}+c_2Lc_1^2\delta^2)\).
Retraction and transport bounds prove \eqref{eq:one-sqp-error}.

For the residual estimate, each original constraint differs from its QP
linearization by at most \(c_3\|p\|^2\). Thus nonlinear negative slack
has that order plus the QP feasibility residual. QP complementarity and
bounded multipliers give the same bound for the nonlinear
complementarity product. After transport to the base,
\[
 \operatorname{grad}L(R_{\bar r}(p),\zeta_Q,t^+)
 =\operatorname{grad}f(\bar r,t^+)+Hp+A_W^*\zeta_Q
   +O(\|p\|^2+\|p\|\,\|\zeta_Q-\bar\zeta_W\|).
\]
QP stationarity and the bound on \(u-u_*\) prove
\eqref{eq:one-sqp-nonlinear-residual}, enlarging constants and shrinking
\(\rho\) to absorb \(\eta_{\rm N}^2\).
\end{proof}

The endpoint multiplier remains the one determined by its final working
set. In the notation of Section~\ref{sec:pathfollowing},
\(\zeta_e=\zeta_{Q,e}\) after the prescribed residual-frame change.
Zero extension, inverse-transpose coordinate changes and stable-ID
transfers preserve this multiplier in the new representation. Combining these operations with
Lemma~\ref{lem:one-sqp-map} gives
\begin{equation}\label{eq:inexact-face-bound}
 d_Z(z_Q^+,z_*^W(t^+))
 \leq C_{\rm SQP}
   \bigl(\delta^2+\eta_{\rm N}+\epsilon_{\rm rt}\bigr).
\end{equation}
The original nonlinear tests are evaluated with this same complete
primal--dual endpoint. A candidate failing the sign, support or residual
conditions of \eqref{eq:endpoint-dual-admission} is ineligible; another
parallel working set or offline recovery supplies the next candidate.

\subsection{Task admission, certification margins, and event resets}
The next assumption states what the finite task configuration must provide.
In short, at least one admitted task reaches the target QP through regular
working sets (coverage), the endpoints that are compared lie close to the
solution, the numerical certificates hold with positive margins, and a
change of working set at an event resets the multiplier without moving the
primal point.
\begin{assumption}[Local task coverage and robust certification]
\label{ass:task-coverage}
The task configuration and tolerances have the following properties on
the tube of Assumption~\ref{ass:regular-path}.
\begin{enumerate}
\item For each admitted local estimate there is a strict-representation
route from a verified auxiliary start to a target-QP endpoint. Each
such route has finitely many events without accumulation in its own
artificial parameter \(\tau\), uses the independent working systems
of Assumption~\ref{ass:regular-path}, and satisfies complete QP slack
and multiplier sign conditions. The configured generator admits at
least one such route, including its required tied-event alternatives
and dual-compatible exchanges. Operation and task budgets permit
that route to finish. The dispatcher allocates the required work
fairly, and the global wall-clock deadline allows at least one
covered route to return a verified target endpoint before the barrier
closes. The required time allowance is included in the scheduling hypothesis.
The frozen-QP route conditions can be checked from its affine
segments and event ratios. Within each auxiliary path, an active-set
event applies \eqref{eq:working-set-relevance} and, for \(\tau<1\),
\eqref{eq:working-set-forward-relevance} at that path's own parameter
value. An irrelevant working set is retired from that path's further
scheduling and endpoint selection. The private zero-distance closure of
Algorithm~\ref{alg:zero-distance-closure} processes all remaining tied
blockers and releases at the same event point. Its completed outgoing
sets include at least one forward-relevant covered continuation, so
retirement preserves coverage. Other live auxiliary paths update
their working sets at their own events.

\item Every eligible endpoint considered by the local selection theorem
lies in the common neighborhood of a covered representation.
For every such representation \(W\), its selected initialization
\(\bar z_k^W\) is obtained by its declared shift and root-preserving
reset, and obeys the quantitative bound
\begin{equation}\label{eq:all-task-prestep-bound}
 d_Z(\bar z_k^W,z_{*,k+1}^W)
 \leq B_Se_k+B_*h_k+\sigma_k+\omega_k=D_k
\end{equation}
with the uniform constants and defects defined in
\eqref{eq:theory-shift-defect}--\eqref{eq:theory-prestep-error}.
The reset precedes the selected Newton comparison. All tasks retain
the same frozen physical target model; any difference between its
Hessian and the Hessian at the reset multiplier is included in
\eqref{eq:theory-residual-budget}.
The barrier retains complete primal--dual endpoints only from the current
target model and relevant working sets, including each task/path identity,
working-set identity and uniquely determined multiplier. An endpoint
from a retired task state is ineligible after its relevance update.
At a physical model refresh all old-version states are retired, and
new-version starts are admitted through the declared transfer and
revalidation procedure.
For each covered endpoint, its transported QP dual passes
\eqref{eq:endpoint-dual-admission}: omitted multipliers remain zero,
retained multipliers have the prescribed signs, and their nonzero
components belong to the fresh near-active catalogue. These are
admission hypotheses for the unchanged task dual.

\item Rank, factor, inertia, and curvature certificates have positive
numerical margins. In each local neighborhood there is a certificate
subset \(S^\circ\) of the exact negative support with
\(\mu_{*,i}\leq-4\epsilon_\mu\) for \(i\in S^\circ\).
The relevant row singular values are bounded below, and the
Lagrangian Hessian has lower bound \(2\sigma_{\min}\) on the nullspace
of \(DC\) and \(Db_{S^\circ}\). The active-row tolerance includes
these rows after retraction. Derivative, multiplier-transport, and
eigenvalue perturbations are small enough that the task dual's
classified support contains
\(S^\circ\) and \eqref{eq:fresh-curvature-test} retains lower bound
\(\sigma_{\min}\). At a zero-multiplier event a vanishing row can
be omitted from \(S^\circ\); the larger tangent space must have this
curvature margin.

\item The independent physical admissibility tests have a local
certificate margin or a certified feasibility-preserving control map.
Residual tolerances \(\epsilon_f,\epsilon_s,\epsilon_c\) dominate
the bound from \eqref{eq:one-sqp-nonlinear-residual}. Retraction,
residual-frame transport, assembly and factor-solve budgets also lie
below the rank, sign-classification and curvature certificate margins.

\item A specified working-set exchange and coordinate transfer preserves
the event primal state and is locally Lipschitz in the complete stable-ID
variables. The outgoing multiplier is the unique multiplier of its
selected independent working representation, obtained by the declared
KKT solve or exact dual-compatible exchange. At the exact event its reset
satisfies
\begin{equation}\label{eq:theory-event-reset}
 \mathcal R_e(r_*(s_j),\zeta_*^-(s_j))
       =(r_*(s_j),\zeta_*^+(s_j)),\qquad
 \operatorname{Lip}(\mathcal R_e)\leq B_e.
\end{equation}
Numerical reset error is bounded by \(\epsilon_e\). A new physical
model is assembled when needed to respect
\eqref{eq:theory-residual-budget}; retaining a frozen Hessian retains
its explicit Hessian defect in that budget. The active-set transition
updates the relevant working-set family before the next continuation
step or barrier selection.
\end{enumerate}
\end{assumption}

These conditions allow weak activity and a nonunique full multiplier
set, while each retained working set supplies a unique multiplier. They
specify the local coverage and numerical margins required of the finite
task configuration. At a physical event the outgoing section is
identified in the same primal neighborhood as the incoming section.
Equation \eqref{eq:theory-event-reset} supplies their common primal
value, and uniqueness of each selected KKT section supplies the
one-sided identification. A jump between valid multiplier
representations is removed by the reset before measuring the outgoing
dual error.

\subsection{Tracking recurrence for the actual RTI map}
Let \(z_k\) be the accepted iterate and let \(z_{*,k}^{W_k}\) be its
selected reference root. Its error is
\(e_k=d_Z(z_k,z_{*,k}^{W_k})\), expanding omitted multipliers by stable
ID. For each eligible outgoing representation \(W\), let
\(\mathscr S_k^W\) be its complete horizon shift, measurement insertion,
and representation reset. The maps have the same primal base and a
uniform local Lipschitz bound \(B_S\), and their exact-reference
consistency defects satisfy
\begin{equation}\label{eq:theory-shift-defect}
 d_Z(\mathscr S_k^W(z_{*,k}^{W_k}),z_{*,k+1}^{W})
 \leq B_*h_k+\sigma_k+\omega_k,\qquad
 h_k=|t_{k+1}-t_k|.
\end{equation}
Here \(\sigma_k\) bounds shift, resampling, and numerical reset error;
\(\omega_k\) bounds prediction error used to build the warm start.
For vector measurement data, \(h_k\) is their normed increment.
A stale or approximate model used in the new QP instead contributes
to \(\Delta F,\Delta K\) in \eqref{eq:theory-residual-budget}.
Thus every eligible task's shifted estimate obeys
\begin{equation}\label{eq:theory-prestep-error}
 \delta_k^W:=d_Z(\mathscr S_k^W(z_k),z_{*,k+1}^W)
 \leq D_k:=
       B_Se_k+B_*h_k+\sigma_k+\omega_k .
\end{equation}
For a layout change, \(\mathscr S_k^W\) explicitly maps between the two
specified problem spaces, and \eqref{eq:theory-shift-defect} is a
restart condition between their regular neighborhoods.
For a barrier with several eligible local endpoints, use the largest
of their certified \(\eta_{{\rm N},k}\) and
\(\epsilon_{{\rm rt},k}\) as the respective budgets in the following
theorem.

\begin{theorem}[Local Riemannian strict-working-set RTI tracking]
\label{thm:riemannian-pathfollowing}
Under Assumptions~\ref{ass:regular-path} and
\ref{ass:task-coverage}, there are uniform constants
\(\rho,C_{\rm RTI},C_{\rm KKT}>0\) with the following property.
If \(D_k\leq\rho\) and the declared numerical budgets are sufficiently
small, the covered route of Section~\ref{sec:pathfollowing} returns a
QP endpoint whose retraction and transported working-set multiplier pass
the original nonlinear residual and curvature tests. Each eligible local
endpoint
satisfies
\begin{equation}\label{eq:actual-rti-recurrence}
 e_{k+1}\leq C_{\rm RTI}
 \left[
  (B_Se_k+B_*h_k+\sigma_k+\omega_k)^2
  +\eta_{{\rm N},k}+\epsilon_{{\rm rt},k}
 \right].
\end{equation}
Its original nonlinear KKT residual is bounded by the same expression
with constant \(C_{\rm KKT}\). Selection by
\eqref{eq:nonlinear-candidate-score} preserves this bound in the
selected representation because its constants are uniform across the
eligible local tasks.

If
\[
 a_k:=B_*h_k+\sigma_k+\omega_k\leq\bar a,\qquad
 \eta_{{\rm N},k}+\epsilon_{{\rm rt},k}\leq\bar\eta,
\]
then a sufficiently small \(E>0\) satisfying
\begin{equation}\label{eq:rti-invariant-neighborhood}
 B_SE+\bar a\leq\rho,\qquad
 C_{\rm RTI}\bigl[(B_SE+\bar a)^2+\bar\eta\bigr]\leq E
\end{equation}
defines an invariant local error neighborhood for \(e_k\leq E\).
At zero parameter and
shift defect, the exact single-round map is locally quadratically
convergent. For small bounded increments it tracks in the neighborhood
specified by \eqref{eq:rti-invariant-neighborhood}.
\end{theorem}
\begin{proof}
The shift estimate follows by its Lipschitz bound, the triangle
inequality, and \eqref{eq:theory-shift-defect}. On a covered route
the worker returns a complete QP endpoint. Its selected equations are
\eqref{eq:sqp-newton-equivalence}, with perturbation bounded by
\eqref{eq:theory-residual-budget}. Apply
Lemma~\ref{lem:one-sqp-map} with \(\delta=\delta_k^W\), in the
common residual dual coordinates, to obtain
\eqref{eq:actual-rti-recurrence}. The residual estimate follows from
the lemma and bounded derivatives. Classification and curvature
margins preserve the fresh certificates; the residual and physical
admissibility margins imply acceptance. All eligible endpoints remain
in the declared local cover, so the estimate applies to the barrier
selection.

At a representation event, \(\mathcal R_e\) preserves the exact primal
root, has bounded local amplification, and contributes numerical
defect to \(\sigma_k\). Restart with the outgoing multiplier and its
own Hessian assembly budget. The finite event partition permits
uniform enlargement of \(B_S\) and the constants across these resets.
Induction proves the invariant-neighborhood claim.
\end{proof}

In words, the error after one update is bounded by the square of the error
of the shifted estimate, which is composed of the previous error, the
change of the problem data over one sample and the shift defects, plus the
numerical budgets of the linear solve and of the retraction.
The quadratic term in \eqref{eq:actual-rti-recurrence} concerns one
SQP/RTI update of the original nonlinear problem. Each frozen QP path
is piecewise affine and is followed exactly between events, so there
is no \(\tau\)-integration truncation term. An approximate Hessian
or an unrefreshed multiplier reset can make \(\eta_{{\rm N},k}\)
first order; the displayed budget records that accuracy directly.

\begin{corollary}[Parallel task interpretation]
\label{cor:parallel-theory}
If the generator admits the covered routes in
Assumption~\ref{ass:task-coverage}, retires each path's working sets when
they lose relevance at that path's active-set changes, and the barrier retains complete
endpoints only from the relevant sets and the same immutable physical
target model, parallel
execution and the specified selection preserve
\eqref{eq:actual-rti-recurrence}. Each task follows its own certified
auxiliary path. Different dual representations are compared with
their selected roots and reset maps.
\end{corollary}
\begin{proof}
Lemma~\ref{lem:one-sqp-map} and the theorem constants are uniform over
the finite local family. Coverage survives each relevance update by
assumption. The barrier selects one complete endpoint from the retained
family.
\end{proof}

\subsection{Recovery and availability of a verified control}
\begin{assumption}[Reliable offline backup and control handoff]
\label{ass:offline-backup}
Fix a declared compact feasible operating domain of problem data and
an admitted set of initialization or recovery seeds. The seed chosen
by the algorithm belongs to this set; otherwise a designated
initialization seed in the set is supplied. Each invocation of
Appendix~\ref{app:pdsqp} in that domain satisfies this contract.
\begin{enumerate}
\item It returns a point of the original manifold problem with verified
primal, signed-dual, stationarity, and complementarity residuals within
the declared backup tolerances. Its supported-space curvature and
independent physical control-admissibility certificates also pass.
The solution lies in a declared regular neighborhood for the next
local segment, with selected-representation initial error at most
\(\epsilon_0\).

\item The solve and certification complete within an available
preparation interval \(T_{\rm avail}\). At initialization this interval
precedes the first control handoff. During recovery, a retained
certified plan supplies admissible controls over an explicit validity
interval, with the required state/prediction uncertainty certificate;
the backup completes before this interval expires.

\item The handoff uses either an accepted RTI control or the verified
backup control. The plant is advanced only with an available control
covered by one of these certificates. Pending computation uses the
retained plan within its certified validity interval.
\end{enumerate}
\end{assumption}

Assumption~\ref{ass:offline-backup} specifies the feasible domain,
initialization, numerical accuracy, and preparation time required for
recovery. Appendix~\ref{app:pdsqp} gives the Riemannian pdSQP construction
and its convergence hypotheses. The implementation supplies the stated
availability interval.

\begin{corollary}[Verified continuation with recovery]
\label{cor:backup-handoff}
Under Assumptions~\ref{ass:regular-path},
\ref{ass:task-coverage}, and \ref{ass:offline-backup}, each control
handoff uses a certified RTI or backup endpoint. A failed local
certificate triggers recovery, whose verified output initializes the
next local estimate with error at most \(\epsilon_0\).
Whenever \eqref{eq:rti-invariant-neighborhood} holds with
\(\epsilon_0\leq E\), subsequent RTI
iterates remain in that neighborhood until the next declared layout
or operating-domain transition.
\end{corollary}
\begin{proof}
Theorem~\ref{thm:riemannian-pathfollowing} supplies an accepted endpoint
on the covered local route. Failure outside that guarantee uses the
handoff and preparation contract of
Assumption~\ref{ass:offline-backup}. Its certified output supplies
the initial error for the next invariant-neighborhood argument.
\end{proof}
\section{Numerical examples}
\label{sec:numerical-results}
\label{sec:results1019B}

We assess numerical tracking of optimality for the drone and arm
OCPs of Section~\ref{sec:geometric-case-studies} in closed loop. The
controller is the real-time iteration of Section~\ref{sec:pathfollowing};
the stabilized primal--dual SQP (pdSQP) method of Appendix~\ref{app:pdsqp}
computes every cold initialization and serves as the backup. The plant is
simulated in MuJoCo~\cite{todorov2012mujoco} (version 3.10.0). Four runs
are reported: the drone on its nominal plant and on a perturbed plant, and
the arm pick-and-place task at two discretizations. For each run we report
the KKT and feasibility residuals of the original nonlinear problem at
every feedback sample and along the QP homotopies, the outcome and wall
time of every feedback step, and MuJoCo frames of the motion.
Appendix~\ref{app:supplementary-figures} contains the remaining figures.

\subsection{Reported quantities and acceptance test}
\label{subsec:results-measures}

All nonlinear residuals refer to the \emph{original} problem
\eqref{eq:priemnlp} at a primal--dual point on the manifold. They are not
residuals of a linearized model.
\begin{itemize}
 \item $r_K$ is the projected KKT residual
       \eqref{eq:pdsqp-projected-residual} of the equality--slack form
       \eqref{eq:pdsqp-reformulation}. It is a two-norm over feasibility,
       stationarity and complementarity, and it vanishes exactly at a KKT
       point.
 \item $\chi_f$, $\chi_s$ and $\chi_c$ are the feasibility, stationarity
       and complementarity components in the infinity norm. At a
       retracted homotopy endpoint they are given by
       \eqref{eq:original-kkt-residuals}. At a pdSQP point, with
       independent slacks $s\geq0$ and signed inequality multipliers
       $\zeta_b\leq0$, they are
       \begin{equation}\label{eq:results-residuals1028}
       \begin{aligned}
        \chi_f&=\max\{\|C\|_\infty,\|b-s\|_\infty,
                            \|[-b]_+\|_\infty,\|[-s]_+\|_\infty\},\\
        \chi_s&=\|\operatorname{grad}f+DC^*\lambda+Db^*\zeta_b\|_\infty,\\
        \chi_c&=\max\{\|s\odot\zeta_b\|_\infty,
                                             \|[\zeta_b]_+\|_\infty\},
       \end{aligned}
       \end{equation}
       which is \eqref{eq:original-kkt-return} with slack nonnegativity
       included in $\chi_f$.
 \item The \emph{parametric-QP residual} is the scaled KKT residual
       $r_{\rm KKT}$ of \eqref{eq:qp-scaled-residual} for the frozen
       tangent QP at homotopy parameter $\tau$. It measures how accurately
       the piecewise-affine QP path is followed. Its target tolerance is
       $\epsilon_{\rm QP}=10^{-7}$.
\end{itemize}
Section~\ref{sec:pathfollowing-theory} bounds the distance $e_k$ of the
iterate to the solution path, which cannot be measured without the
solution. Theorem~\ref{thm:riemannian-pathfollowing} bounds the original
KKT residual by the same recurrence, so $r_K$ and $\chi_f$ are the
computable tracking errors used here.

Each feedback sample has one of three outcomes. \emph{Online}: a retracted
homotopy endpoint passes the tests of
Section~\ref{subsec:synchronization-barrier} and its first control is
applied. \emph{Backup}: no endpoint passes and the pdSQP method supplies a
verified point. \emph{Retained}: neither produces a verified point; the
shifted seed is returned unchanged and the supervisor decides how to
continue. In the drone runs the supervisor replans at once, from the same
measured state, with a fresh horizon of six intervals. In the arm runs the
controller would continue from the seed and the run would be counted as
failed; no arm sample is retained. A point is
accepted when each of $\chi_f$, $\chi_s$
and $\chi_c$ is at most the tolerance $\epsilon$, its multipliers pass the
admission test \eqref{eq:endpoint-dual-admission}, and the curvature test
\eqref{eq:fresh-curvature-test} holds. The tolerance is
$\epsilon=10^{-6}$ for the nominal drone and $\epsilon=10^{-4}$ for the
perturbed drone and the arm. The residual test is per component. Since $r_K$ is a
two-norm over all components, it can exceed $\epsilon$ slightly at an
accepted point; the largest such value in the four runs is
$1.7\times10^{-4}$. At most three correction rounds $j$
(Section~\ref{subsec:nonlinear-recheck}) are run per sample. In these runs
no round contained more than two working-set tasks.

For the homotopy figures, a separate run with identical outcomes and
residuals records the nonlinear residuals along each path: the path point
$p(\tau)$ is retracted and evaluated with the path multipliers
$\zeta(\tau)$ at $21$ (drone) or $11$ (arm) equally spaced values of
$\tau$ and at every active-set event.

\subsection{Implementation}
\label{subsec:results-implementation}
The method is implemented in Rust, and the runs reported here were made on
an eight-core x86-64 laptop with 32~GB of memory. This subsection records the choices
made where Sections~\ref{sec:geometric-case-studies}--\ref{sec:pathfollowing}
and Appendix~\ref{app:pdsqp} leave one open, and the places where the
implementation covers only part of the stated formulation.

\paragraph{Derivatives and linear algebra.}
A model lists its cost, equality and inequality terms once, each as a
function of a few tangent coordinates. Gradients, Jacobian rows and the
exact Lagrangian Hessian are obtained by forward-mode automatic
differentiation with second-order jets in the coordinates of the
retraction; the terms are evaluated in parallel. The KKT matrices are not
condensed. The KKT matrix of the first hypothesis of a round is factored
once by a block-tridiagonal elimination over the stages, from the terminal
stage backward, with the dense Bunch--Kaufman $LBL^\top$ factorization of
Section~\ref{sec:background} in each block; the free duration of the drone
is a border variable. This base factor is shared by all tasks. A task
represents its own row additions \eqref{eq:border-general-row} and
releases by a bordered Schur complement on the shared factor, tracks the
inertia through the border, and refactors after $40$ border updates. Simple
bounds are carried as the oriented unit rows \eqref{eq:oriented-bound-row}
in the same matrix, so that $N=\varnothing$ in
Section~\ref{subsec:dependent-exchange} and bounds and general rows are
treated uniformly. Pivots smaller than $10^{-11}$ times the largest matrix
entry are rejected ($\delta_{\rm abs}=0$), no equilibration is used, and
row independence is tested by a bordered solve with relative threshold
$10^{-9}$ in place of a rank-revealing QR factorization.

\paragraph{Feedback step.}
The auxiliary start \eqref{eq:start-qp-data} is always constructed, with
$p^0=0$ at the shifted trajectory. The hypotheses of
Section~\ref{subsec:task-generation} are \texttt{WARM} with
$\epsilon_w=\max\{\epsilon_a,10^{-3}\}$, \texttt{ACTIVE\_OUT} with
$c_\mu=10^3$, and \texttt{ACTIVE\_IN} with at most six added rows, taken in
order of increasing margin. Every hypothesis is first mapped to an
independent representation by \eqref{eq:dependent-dual-family}. This is
needed for the arm, where the four coplanar points of a pad and the four
bottom vertices of the box give linearly dependent active normals. The
tolerances are $\epsilon_f=\epsilon_s=\epsilon_c=\epsilon_a=\epsilon$ with
the value of $\epsilon$ given above, $\epsilon_\mu=10^{-7}$,
$\epsilon_{\rm dual}=10^{-9}$, $\sigma_{\min}=10^{-8}$,
$\Delta_{\rm loc}=10^{3}$ and $\epsilon_{\rm QP}=10^{-7}$; zero slacks and
multipliers are classified with relative tolerance $10^{-9}$. The scale
factors of \eqref{eq:qp-scaled-residual} are
$s_p=\max\{1,\|b(\tau)\|_\infty\}$,
$s_d=\max\{1,\|g^1\|_\infty,\|\zeta\|_\infty\}$ and
$s_c=s_p\max\{1,\|\zeta\|_\infty\}$. A round admits at most
$N_{\rm task}=16$ tasks, each with a budget of $2000$ KKT solves and factor
changes. Algorithms~\ref{alg:one-path-worker} and~\ref{alg:event-handler}
are realized as one loop that runs
Algorithm~\ref{alg:zero-distance-closure} at every checkpoint, including
the start. Of the statuses of Section~\ref{subsec:synchronization-barrier}
the implementation reports \texttt{TARGET\_QP}, \texttt{STOPPED\_START},
\texttt{STOPPED\_FACTOR}, \texttt{STOPPED\_RANK},
\texttt{STOPPED\_RESIDUAL}, \texttt{EXCHANGE\_HANDOFF},
\texttt{RETIRED\_IRRELEVANT}, \texttt{NO\_CERTIFIED\_EXCHANGE} and
\texttt{BUDGET\_EXHAUSTED}. A working set that repeats within a closure is
discarded, exchange children beyond $N_{\rm task}$ are not admitted, and a
sign failure appears as a failed residual test at $\tau=1$. At the barrier
the curvature test is the inertia condition
\eqref{eq:fresh-curvature-inertia}, with $S_e$ restricted to the working
rows of the task. The second test of \eqref{eq:endpoint-dual-admission} is
applied in the product form $(\zeta_{b,e})_ib_i\geq-\epsilon_c$ to the
rows with $b_i>\max\{\epsilon_a,\epsilon_f\}$. Ties in
\eqref{eq:nonlinear-candidate-score} are broken by the objective value and
the task ID. No separate physical-output test is applied, since the control
bounds are constraints of the problems. At most three correction rounds are
run per sample: when no endpoint is accepted, the endpoint with the
smallest score becomes the new base if its score is less than half that of
the previous base (in the first round, if it is finite). The backup
starts from the best such endpoint, or from the shifted seed, and is
verified with the same tolerances. There is no wall-clock deadline and no
separate preparation phase: all derivatives are evaluated after the
measurement, and the simulation waits for the solver.

\paragraph{Initialization and backup.}
The pdSQP method uses $\aleph_0=10^{-4}$, $\gamma=1/2$, $\beta=10^{-5}$,
$\eta_S=10^{-3}$, $\bar\lambda=0.2$, $\epsilon_a=0.1$, $y_{\max}=10^6$ and
$\tau_0=10^{-2}$, with at most $200$ iterations ($150$ for the arm
initializations). The slack variables are
eliminated analytically from its linear systems. The convexification
\eqref{eq:pdsqp-convexification} is $\widehat H_y=H_y+\delta I$ with
$\delta$ increased by factors of ten until the working matrix
\eqref{eq:pdsqp-stabilized-kkt} has the correct inertia, and the direction
of negative curvature is computed by shifted inverse iteration. The
return tolerance on $\chi_f$, $\chi_s$ and $\chi_c$ is $10^{-9}$ for the
drone initializations and $10^{-8}$ for the arm initializations. No
restoration phase is implemented: a failed backup is reported and gives the
outcome ``retained'' defined above.

\paragraph{Problem data.}
The drone problem \eqref{eq:ex1nlp} uses $m=1$~kg, $\omega=0.4$~rad/s,
$v_{\max}=20$~N, $\rho_T=1$, $\rho_v=0.01$, $h_{\min}=0.05$~s and
$h_{\max}=0.5$~s. The drone starts at rest at $(0,0,0.6)$~m, the goal is
$(1,0.4,1)$~m, and the obstacle is a sphere of radius $0.15$~m centered at
$(0.5,0.05,0.6)$~m; $a_d=0.2$~m and $d_{\rm safe}=d_{\rm ground}=0.05$~m.
Its attitude defect is the Lie-group error given by the logarithm of
$(q^R_{k+1})^\top q^R_ke^{h\Omega_k}$, which
Section~\ref{subsec:shooting-assembly} allows in place of the anchored map
\eqref{eq:shooting-anchored-map}. The commanded body rates are imposed
exactly on the simulated drone, which is the ideal inner loop assumed in
Section~\ref{sec:drone-example}.

The arm has the UR5 geometry of Section~\ref{sec:arm-example} with
torque-controlled joints and a parallel-jaw hand limited to $6$~N per jaw.
The box has mass $0.12$~kg and half extents $(0.035,0.05,0.035)$~m and
rests on four posts; the goal posts are $0.27$~m away and rotated by
$0.35$~rad. The problems use the inverse-dynamics form of the rate
equation, $36$ explicit features per node (tool position, tool attitude
and eight points), $\mu_{\rm safe}=0.5$, normal forces between $0.2$ and
$3$~N per contact, $N_w=2$ load scenarios and $\ell_b=0.05$~m. Two terms
differ from a literal reading of
\eqref{eq:arm-ocp-a}--\eqref{eq:arm-ocp-c}. The force-allocation
regularization \eqref{eq:arm-force-regularization} is used with
$\rho_f=10^{-2}$ and $N_{\rm pre}=1.2$~N, and the support-friction law
\eqref{eq:arm-release-support-friction} is used with $\mu_S=0.6$ and
$v_S=37.5\,h$, as explained there. The retreat neighborhood of
\eqref{eq:arm-retreat-terminal-set} is imposed with
$\delta_{\rm retreat}=0.2$~rad as task data; its Jacobian certificate is
not evaluated.

The implemented collision catalogue is a subset of the one in
Section~\ref{sec:arm-example}. It contains pad--box nonpenetration and a
hand-height bound in $\mathsf A$, the payload clearance
\eqref{eq:arm-payload-clearance} in $\mathsf B$ and $\mathsf C$, and
hand--object separation, the hold region and the small-motion neighborhood
in the detached layouts, together with all state and actuator bounds. It
omits the link sphere covers, the self-collision pairs and the
continuous-time clearance certificate
\eqref{eq:continuous-clearance-certificate}; constraints are therefore
imposed at the nodes only. In the detached layouts the residual pad wrench
is taken as zero after the measured unloading, the terminal holding-torque
bound is not imposed, and the support equalities
\eqref{eq:arm-release-support-mode} are taken relative to the measured
supported pose, which lies about $1$~mm above the nominal plane. The
transitions between problems are made by joint PD servos with gravity
compensation and by threshold tests on the measured contact forces and
poses; force closure is not recertified numerically at a transition. The
jaw position servo has gain $4000$~N/m, and the scene uses elliptic
friction cones and MuJoCo's no-slip pass. The arm experiments use the
nominal MuJoCo plant only; a perturbed plant is tested for the drone.

\subsection{Drone navigation}
\label{sec:drone-numerics}

\begin{figure}[tbp]
 \centering
 \begin{minipage}{.49\linewidth}\centering
  \includegraphics[width=\linewidth]{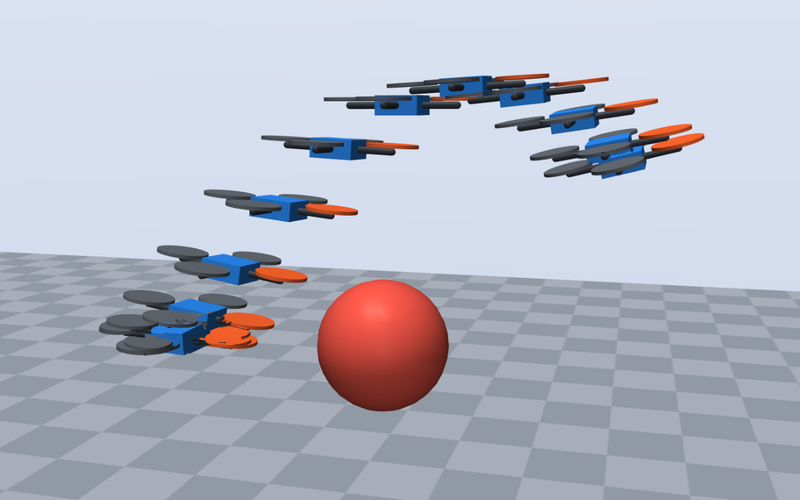}\\[2pt]
  {\small (a) nominal plant, all 13 sample frames}
 \end{minipage}\hfill
 \begin{minipage}{.49\linewidth}\centering
  \includegraphics[width=\linewidth]{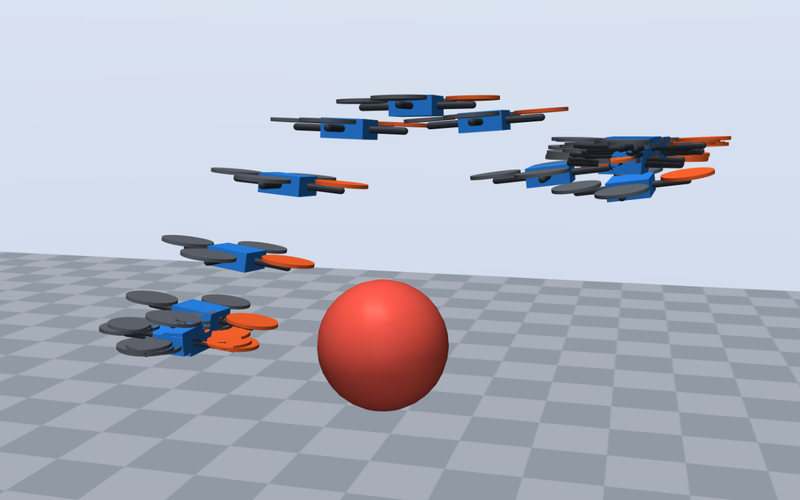}\\[2pt]
  {\small (b) perturbed plant, all 20 sample frames}
 \end{minipage}\\[8pt]
 \begin{minipage}[t]{.24\linewidth}\centering
  \includegraphics[width=\linewidth]{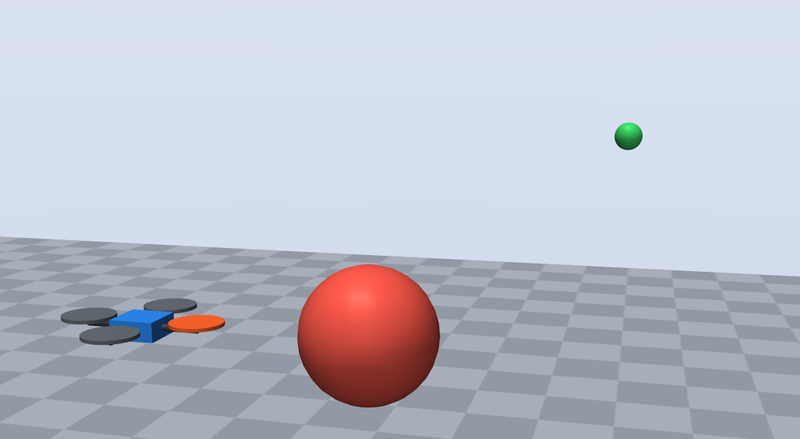}\\[-1pt]
  {\scriptsize $t=0.00$ s\par}%
 \end{minipage}\hfill
 \begin{minipage}[t]{.24\linewidth}\centering
  \includegraphics[width=\linewidth]{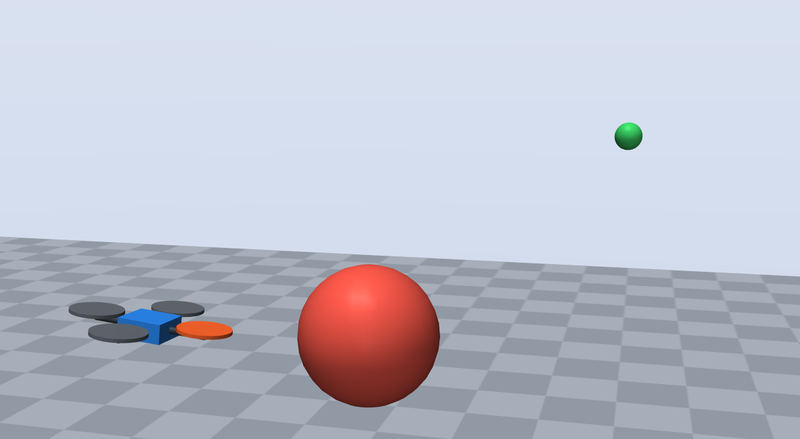}\\[-1pt]
  {\scriptsize $t=0.28$ s\par}%
 \end{minipage}\hfill
 \begin{minipage}[t]{.24\linewidth}\centering
  \includegraphics[width=\linewidth]{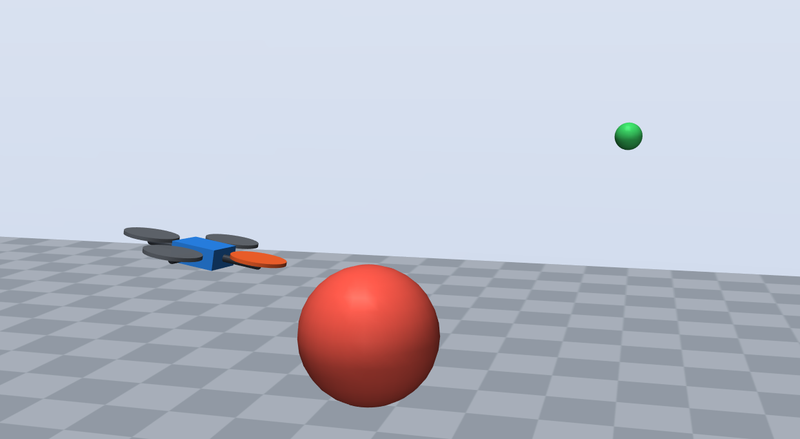}\\[-1pt]
  {\scriptsize $t=0.56$ s\par}%
 \end{minipage}\hfill
 \begin{minipage}[t]{.24\linewidth}\centering
  \includegraphics[width=\linewidth]{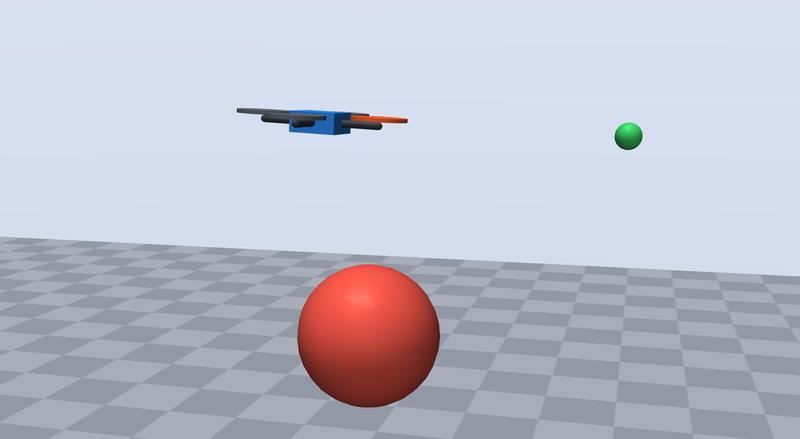}\\[-1pt]
  {\scriptsize $t=0.85$ s\par}%
 \end{minipage}\\[3pt]
 \begin{minipage}[t]{.24\linewidth}\centering
  \includegraphics[width=\linewidth]{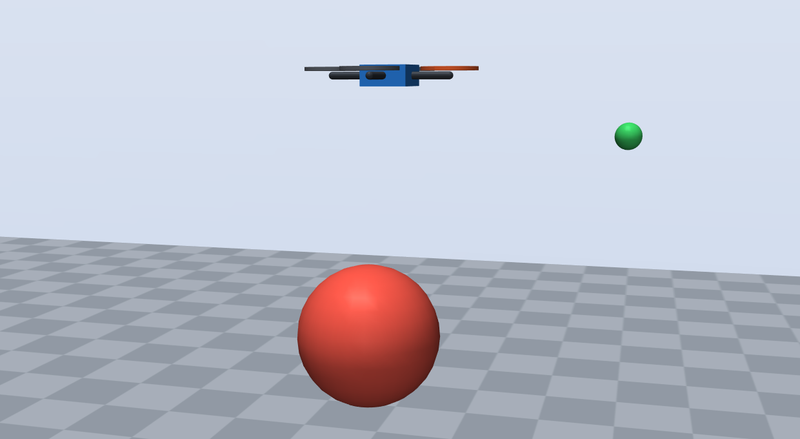}\\[-1pt]
  {\scriptsize $t=0.99$ s\par}%
 \end{minipage}\hfill
 \begin{minipage}[t]{.24\linewidth}\centering
  \includegraphics[width=\linewidth]{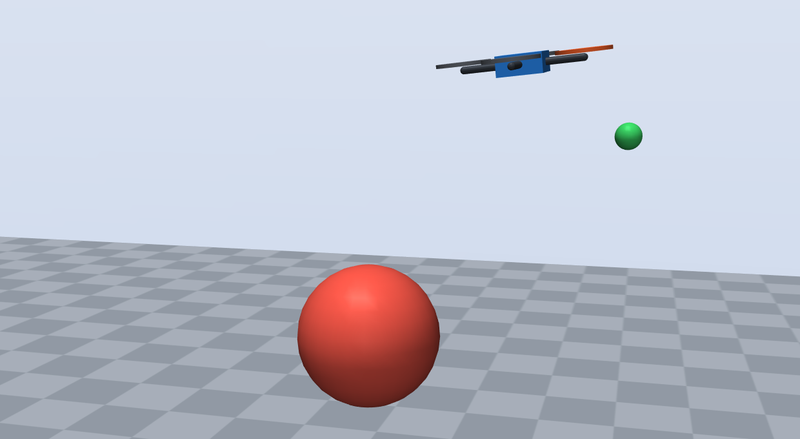}\\[-1pt]
  {\scriptsize $t=1.27$ s\par}%
 \end{minipage}\hfill
 \begin{minipage}[t]{.24\linewidth}\centering
  \includegraphics[width=\linewidth]{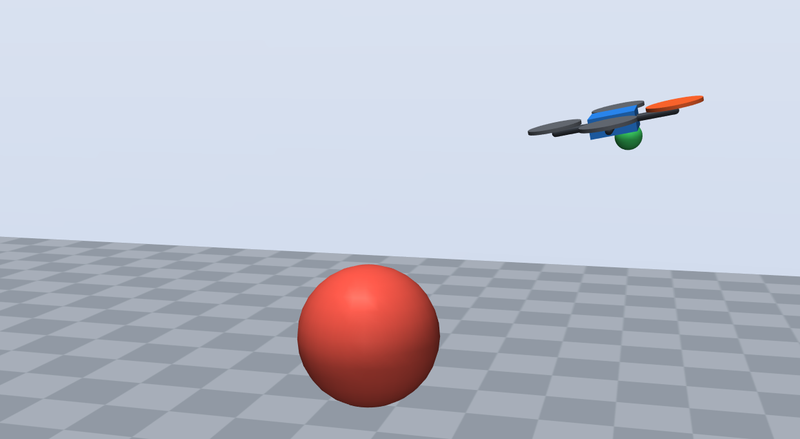}\\[-1pt]
  {\scriptsize $t=1.55$ s\par}%
 \end{minipage}\hfill
 \begin{minipage}[t]{.24\linewidth}\centering
  \includegraphics[width=\linewidth]{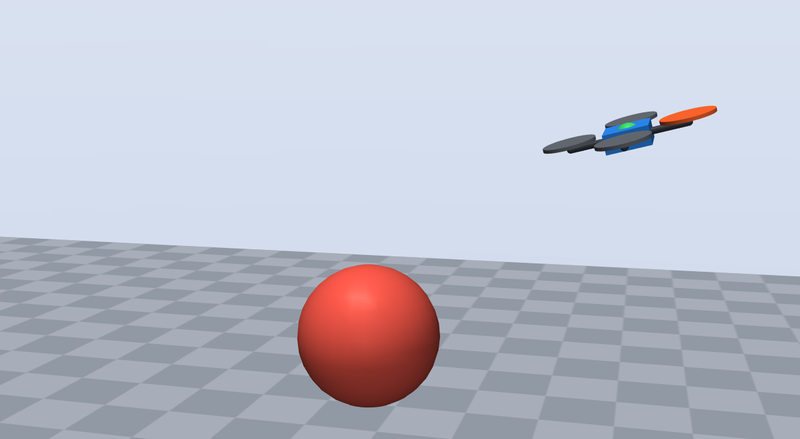}\\[-1pt]
  {\scriptsize $t=1.69$ s\par}%
 \end{minipage}\\[5pt]
 {\small (c) nominal plant}\\[9pt]
 \begin{minipage}[t]{.24\linewidth}\centering
  \includegraphics[width=\linewidth]{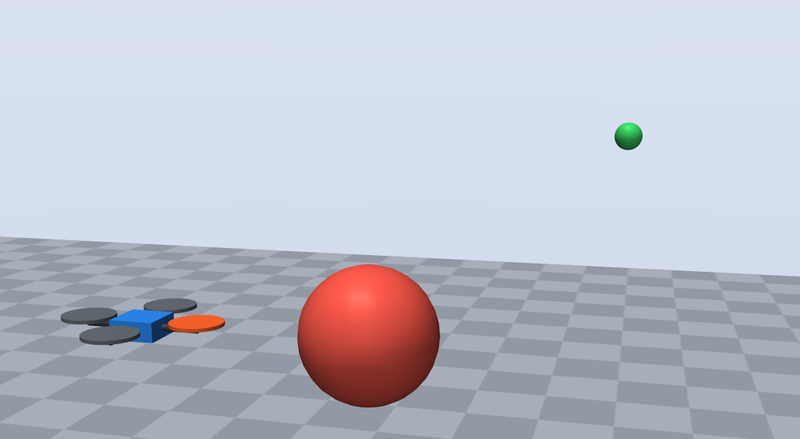}\\[-1pt]
  {\scriptsize $t=0.00$ s\par}%
 \end{minipage}\hfill
 \begin{minipage}[t]{.24\linewidth}\centering
  \includegraphics[width=\linewidth]{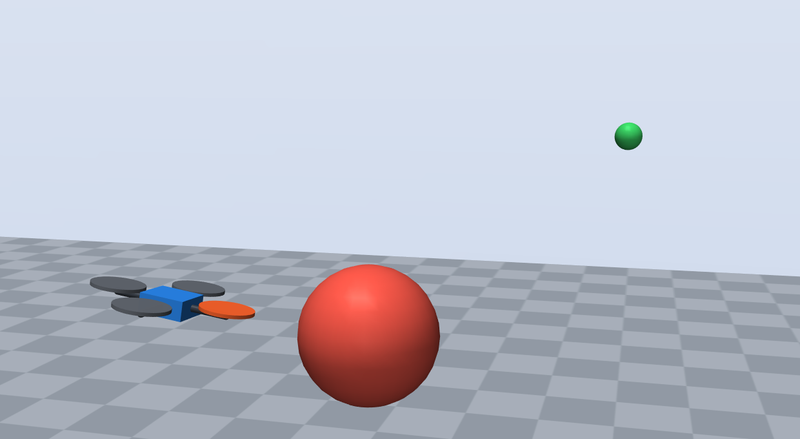}\\[-1pt]
  {\scriptsize $t=0.42$ s\par}%
 \end{minipage}\hfill
 \begin{minipage}[t]{.24\linewidth}\centering
  \includegraphics[width=\linewidth]{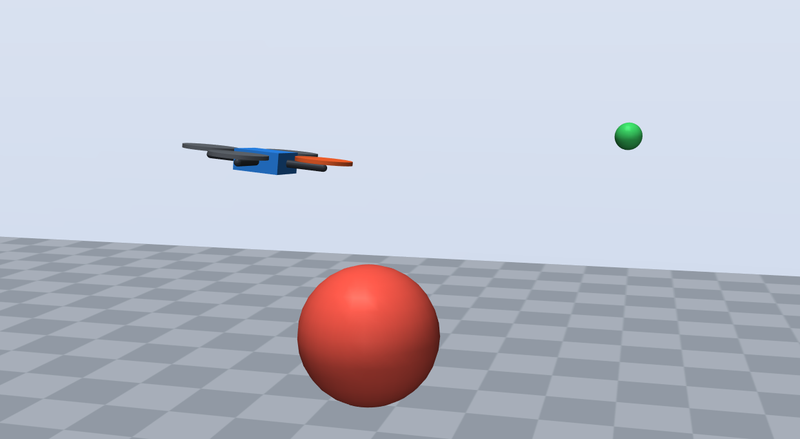}\\[-1pt]
  {\scriptsize $t=0.69$ s\par}%
 \end{minipage}\hfill
 \begin{minipage}[t]{.24\linewidth}\centering
  \includegraphics[width=\linewidth]{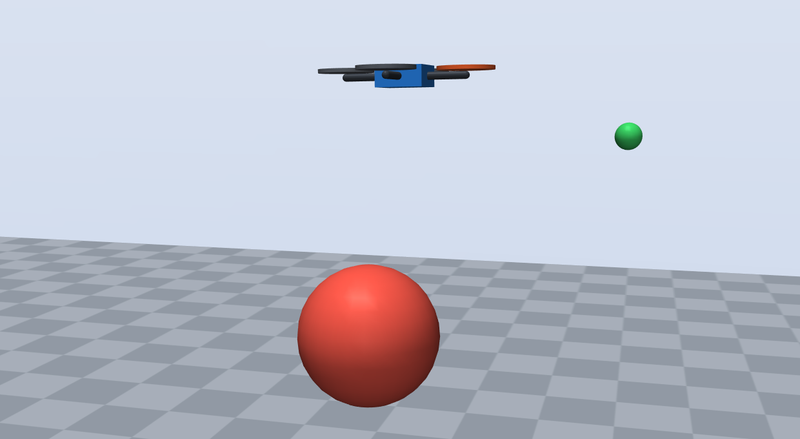}\\[-1pt]
  {\scriptsize $t=0.97$ s\par}%
 \end{minipage}\\[3pt]
 \begin{minipage}[t]{.24\linewidth}\centering
  \includegraphics[width=\linewidth]{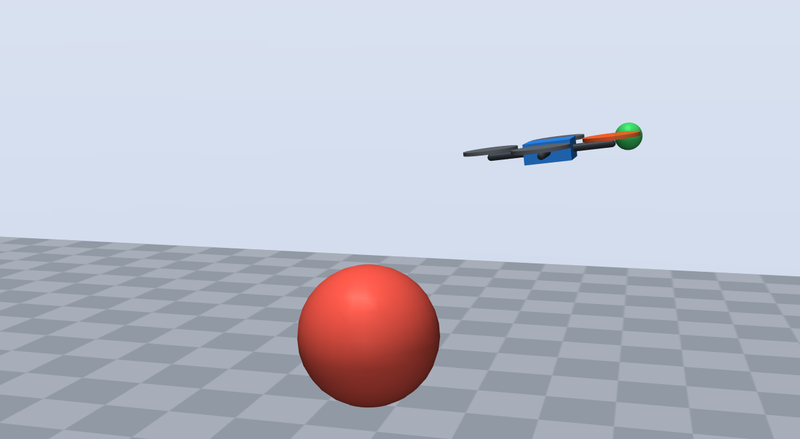}\\[-1pt]
  {\scriptsize $t=1.27$ s\par}%
 \end{minipage}\hfill
 \begin{minipage}[t]{.24\linewidth}\centering
  \includegraphics[width=\linewidth]{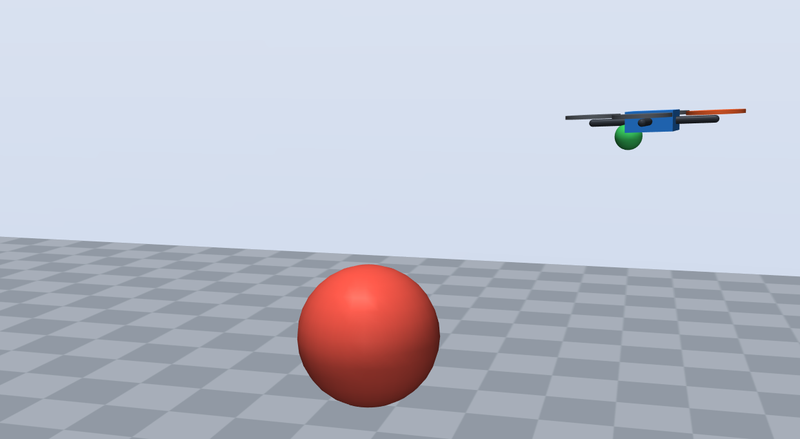}\\[-1pt]
  {\scriptsize $t=1.98$ s\par}%
 \end{minipage}\hfill
 \begin{minipage}[t]{.24\linewidth}\centering
  \includegraphics[width=\linewidth]{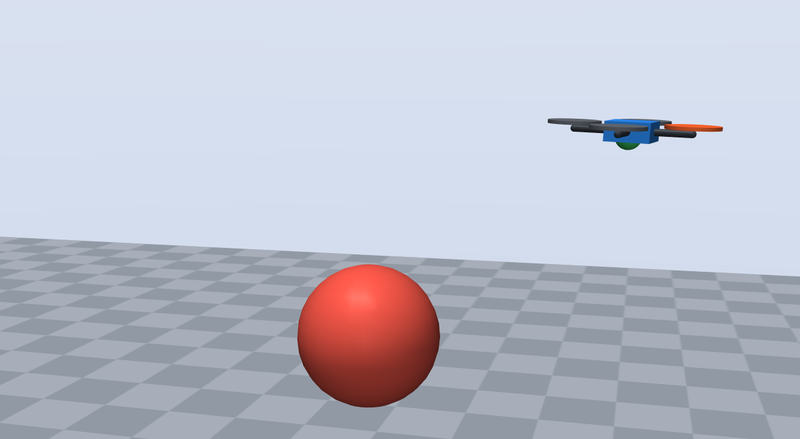}\\[-1pt]
  {\scriptsize $t=2.49$ s\par}%
 \end{minipage}\hfill
 \begin{minipage}[t]{.24\linewidth}\centering
  \includegraphics[width=\linewidth]{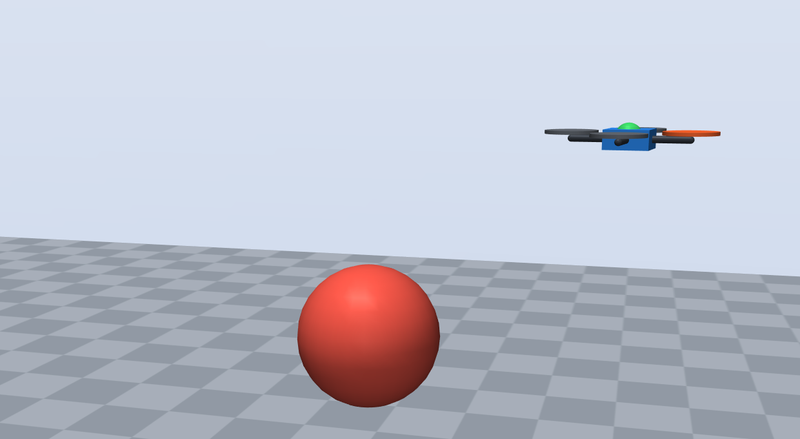}\\[-1pt]
  {\scriptsize $t=2.64$ s\par}%
 \end{minipage}\\[5pt]
 {\small (d) perturbed plant; the plans are recomputed at $t=1.98$ s and $t=2.49$ s}
 \caption{Drone navigation in MuJoCo. The red sphere is the obstacle and
 the green sphere marks the goal position. (a,\,b) All frames recorded at
 the feedback samples of one flight, composited into one image; the flight
 starts at the lower left. (c,\,d) Individual frames at the indicated
 simulation times.}
 \label{fig:drone-mujoco}
\end{figure}

The drone problem is the shrinking-horizon, free-final-time problem
\eqref{eq:ex1nlp} with $H=12$ control intervals. The drone starts at rest
near the ground, passes a spherical obstacle, and must arrive at rest at
the goal. Figure~\ref{fig:drone-mujoco} shows the two flights. The
perturbed plant has $3\%$ more mass than the prediction model and is subject to a constant wind force of $(0.15,-0.1,0)$~N, of magnitude $0.18$~N; the prediction model is unchanged.

\begin{figure}[tbp]
 \centering
 \includegraphics[width=\linewidth]{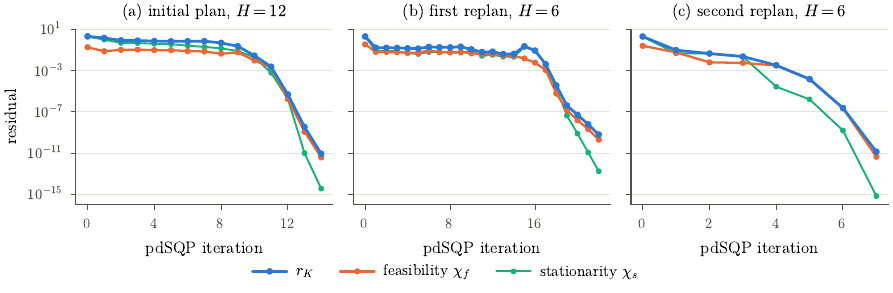}
 \caption{Convergence of the pdSQP initializations of the drone problem.
 Iteration $0$ is the seed. (a) The initial plan, identical for the two
 runs. (b,\,c) The two plans recomputed during the perturbed run.}
 \label{fig:drone-pdsqp}
\end{figure}

\paragraph{Initialization.}
The cold start is computed by pdSQP (Figure~\ref{fig:drone-pdsqp}a). It
converges in $14$ iterations and $0.09$~s to $r_K=8.7\times10^{-12}$, with
$\chi_f=3.9\times10^{-12}$, $\chi_s=3.6\times10^{-15}$ and
$\chi_c=6.0\times10^{-17}$. The optimized arrival time is
$T=1.6914$~s, the supported-space curvature test
\eqref{eq:fresh-curvature-test} passes, and one inequality row is strongly
active (at node $k=4$). All fourteen iterations are V-iterates in the
classification of Appendix~\ref{app:pdsqp}; the first eight use step
lengths below one and the last six take full steps.

\begin{figure}[tbp]
 \centering
 \includegraphics[width=\linewidth]{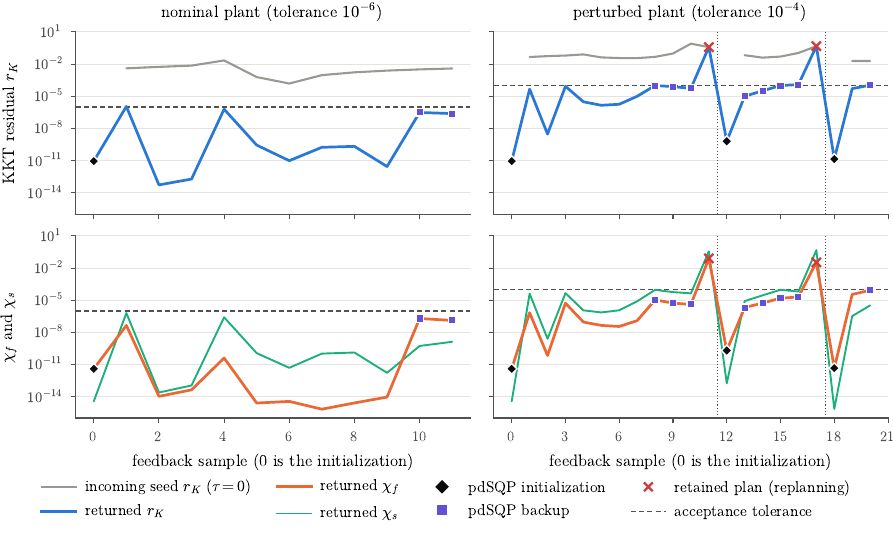}
 \caption{Drone: residuals of the original nonlinear problem at every
 feedback sample. Top: KKT residual $r_K$ of the incoming shifted seed
 ($\tau=0$, gray) and of the returned point (blue). Bottom: feasibility
 $\chi_f$ and stationarity $\chi_s$ of the returned point. Unmarked points
 are accepted online. In the perturbed run the vertical dotted lines
 separate the three plans.}
 \label{fig:drone-closed-loop}
\end{figure}

\paragraph{Nominal plant.}
Table~\ref{tab:drone-runs} and Figure~\ref{fig:drone-closed-loop} (left)
summarize the closed loop. The first nine of the eleven feedback samples
are accepted online, seven after one correction round and two after two.
The shifted seed has $r_K$ between $1.5\times10^{-4}$ and
$2.2\times10^{-2}$; the accepted point has $r_K$ between
$5.3\times10^{-14}$ and $1.1\times10^{-6}$. The last two samples, with
$H=2$ and $H=1$ remaining intervals, are resolved by the backup, which
returns $r_K=3.0\times10^{-7}$ and $2.4\times10^{-7}$. The end of the
shrinking horizon is the regime in which the regularity assumptions are
lost. For $H=1$ the retained equality system of \eqref{eq:ex1nlp} has
fifteen rows (nine shooting rows and six terminal rows) in the thirteen
unknowns $(T,u_0,v_0,x_1,q^R_1)$, so its Jacobian cannot be surjective and
MFCQ (Definition~\ref{def:mfcq}) fails. For $H=2$ it has $24$ rows in $25$
unknowns: MFCQ is not excluded by this count, but a single degree of
freedom remains for the free arrival time.
The drone arrives at $t=1.6910$~s with position error
$1.8\times10^{-7}$~m and velocity error $2.6\times10^{-7}$~m/s. A feedback
step takes between $2.0$ and $10.8$~ms, against a sampling interval of
$141$~ms.

The collision constraint \eqref{eq:ex1ineqcon} is imposed at the nodes.
Between nodes its margin reaches $-7.4\times10^{-3}$~m, that is, the
trajectory enters the safety margin by $7.4$~mm, while the physical
clearance to the obstacle remains $4.3\times10^{-2}$~m.

\begin{table}[tbp]
\centering\small
\setlength{\tabcolsep}{5pt}
\begin{tabular}{lcc}
\toprule
 & Nominal plant & Perturbed plant\\
\midrule
Plans (pdSQP initializations) & $1$ ($H=12$) & $3$ ($H=12,6,6$)\\
pdSQP iterations per plan & $14$ & $14$, $22$, $7$\\
Feedback samples & $11$ & $18$\\
Online / backup / retained & $9$ / $2$ / $0$ & $8$ / $8$ / $2$\\
Mean correction rounds per sample & $1.55$ & $2.83$\\
Acceptance tolerance $\epsilon$ & $10^{-6}$ & $10^{-4}$\\
Largest accepted component of $(\chi_f,\chi_s,\chi_c)$ & $6.3\times10^{-7}$ & $9.3\times10^{-5}$\\
Largest $r_K$ at an accepted point & $1.1\times10^{-6}$ & $1.2\times10^{-4}$\\
Arrival time & $1.691$ s & $2.635$ s\\
Arrival position error & $1.8\times10^{-7}$ m & $5.7\times10^{-3}$ m\\
Arrival velocity error & $2.6\times10^{-7}$ m/s & $1.2\times10^{-2}$ m/s\\
Minimum margin between nodes & $-7.4\times10^{-3}$ m & $-3.0\times10^{-3}$ m\\
Minimum physical clearance & $4.3\times10^{-2}$ m & $4.7\times10^{-2}$ m\\
Largest step time: online / backup / retained & $10.8$ / $3.4$ / -- ms & $13.7$ / $88.0$ / $659$ ms\\
\bottomrule
\end{tabular}
\caption{Closed-loop drone runs. Accepted points are those returned by an
online or a backup sample.}
\label{tab:drone-runs}
\end{table}

\paragraph{Perturbed plant.}
The model mismatch raises the residual of the incoming seed to between
$2\times10^{-2}$ and $0.8$, roughly one to two orders of magnitude above the
nominal run (Figure~\ref{fig:drone-closed-loop}, right). The first seven
samples are still accepted online, now after two or three rounds, and the
next three by the backup. At the last sample of the plan ($H=1$) neither
produces a verified point: with the disturbance, the hard terminal
equalities cannot be met from the measured state within the one remaining
interval. The supervisor then computes a new plan with $H=6$ from the same
measured state by pdSQP, in $22$ iterations
(Figure~\ref{fig:drone-pdsqp}b). The second plan is tracked by the backup
alone and again ends without a verified point at $H=1$; a third plan
($H=6$, seven pdSQP iterations) is followed for two samples until arrival.
The drone arrives at $t=2.635$~s with position error
$5.7\times10^{-3}$~m and velocity error $1.2\times10^{-2}$~m/s, having
kept a physical clearance of at least $4.7\times10^{-2}$~m. Online steps
take $5.8$--$13.7$~ms and backup steps $11$--$88$~ms; the two samples that
end in replanning take $127$~ms and $659$~ms. The sampling intervals of
this run range from $50$~ms to $387$~ms, so the replanning samples exceed
their interval.

\begin{figure}[tbp]
 \centering
 \includegraphics[width=\linewidth]{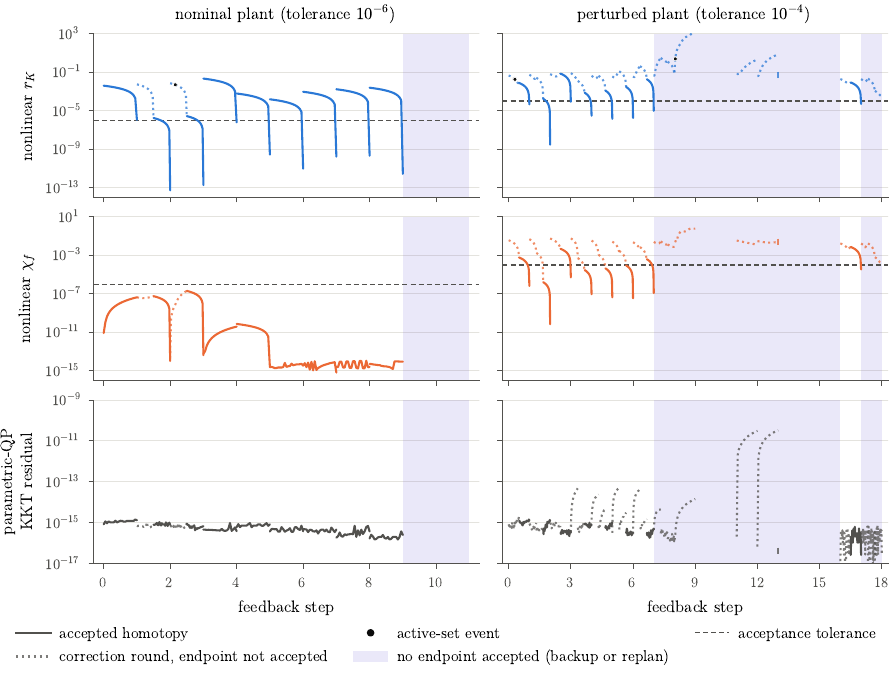}
 \caption{Drone: residuals along every traced homotopy. Feedback step $s$
 occupies the interval $[s-1,s]$; within it each correction round runs
 from $\tau=0$ to $\tau=1$, and successive rounds share the interval
 equally. Top and middle: nonlinear KKT residual $r_K$ and feasibility
 $\chi_f$ of the retracted path point. Bottom: KKT residual of the
 parametric QP. Shaded samples without a curve had no task that advanced
 beyond its start checkpoint.}
 \label{fig:drone-unrolled}
\end{figure}

\paragraph{Residuals along the homotopies.}
Figure~\ref{fig:drone-unrolled} lays every traced homotopy of the two runs
end to end. Three observations hold for all accepted paths.

First, the parametric-QP residual stays at rounding level: at most
$1.4\times10^{-15}$ on accepted paths and $3.3\times10^{-11}$ on rounds
whose endpoint was not accepted, against the tolerance $10^{-7}$. The QP
path is followed essentially exactly, as expected for a piecewise-affine
path advanced by exact event ratios (Section~\ref{subsec:qpoases-homotopy}).
The nonlinear residual left at $\tau=1$ is therefore the linearization
error of the frozen model, not a pathfollowing error.

Second, the nonlinear KKT residual decreases over an accepted path by
$3.6$ to $9.0$ orders of magnitude in the nominal run (median $6.9$) and
by $2.2$ to $4.8$ orders in the perturbed run (median $2.9$). Almost all
of the decrease occurs at the end of the path: for $\tau<1$ the residual
follows $r_K(\tau)\approx(1-\tau)\,r_K(0)$, because the interpolated QP
data carry the share $1-\tau$ of the start data until that share falls
below the linearization error
(Figure~\ref{fig:drone-tau} in Appendix~\ref{app:supplementary-figures}).

Third, in the nominal run the seed is feasible to rounding, so along four
of the nine accepted paths $\chi_f$ does not decrease: it starts at
rounding level, at most $8.5\times10^{-12}$, and ends at the size of the
linearization error, at most $4.4\times10^{-8}$, while stationarity
decreases. In the perturbed run the measured state makes the
seed infeasible, and feasibility decreases along every accepted path by
$2.3$ to $4.4$ orders of magnitude.

No accepted drone path contains an active-set event. The four recorded
events lie on correction rounds whose endpoint was not accepted.

\FloatBarrier
\subsection{Arm pick-and-place}
\label{sec:arm-numerics}

\begin{figure}[tbp]
 \centering
 \begin{minipage}[t]{.325\linewidth}\centering
  \includegraphics[width=\linewidth]{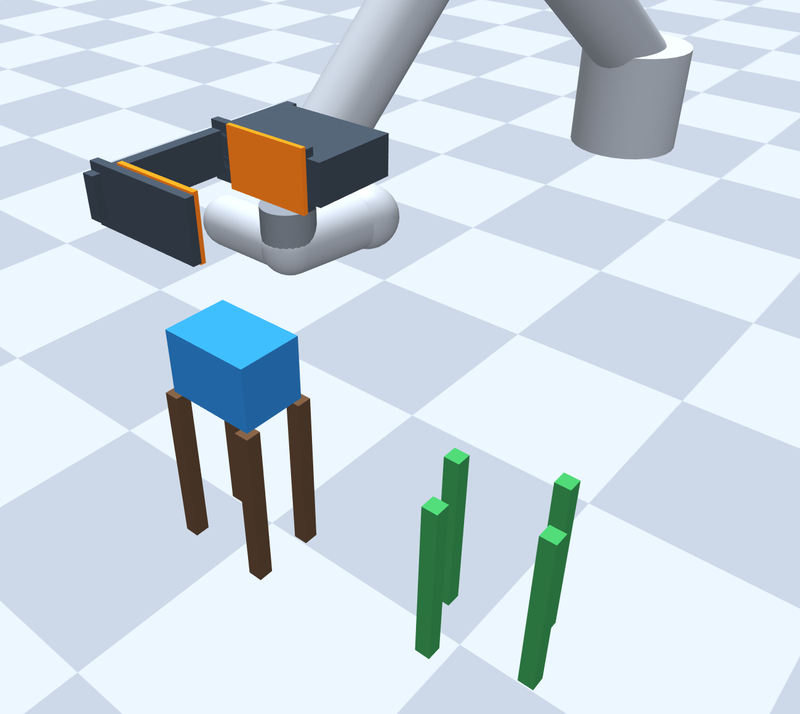}\\[-1pt]
  {\scriptsize 1. start, $t=0.30$ s: jaws open above the box\par}%
 \end{minipage}\hfill
 \begin{minipage}[t]{.325\linewidth}\centering
  \includegraphics[width=\linewidth]{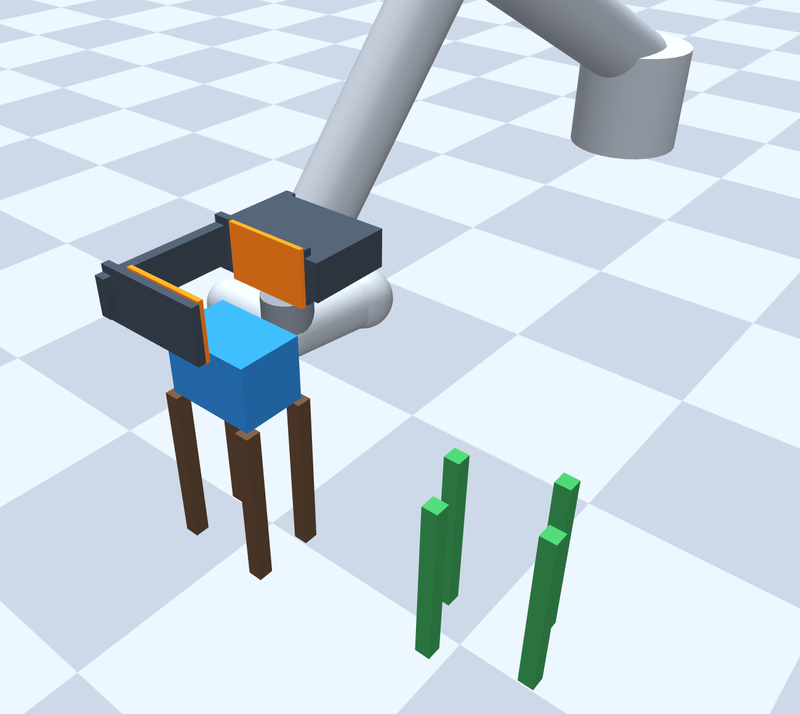}\\[-1pt]
  {\scriptsize 2. $\mathsf A$, $t=3.11$ s: descending over the box\par}%
 \end{minipage}\hfill
 \begin{minipage}[t]{.325\linewidth}\centering
  \includegraphics[width=\linewidth]{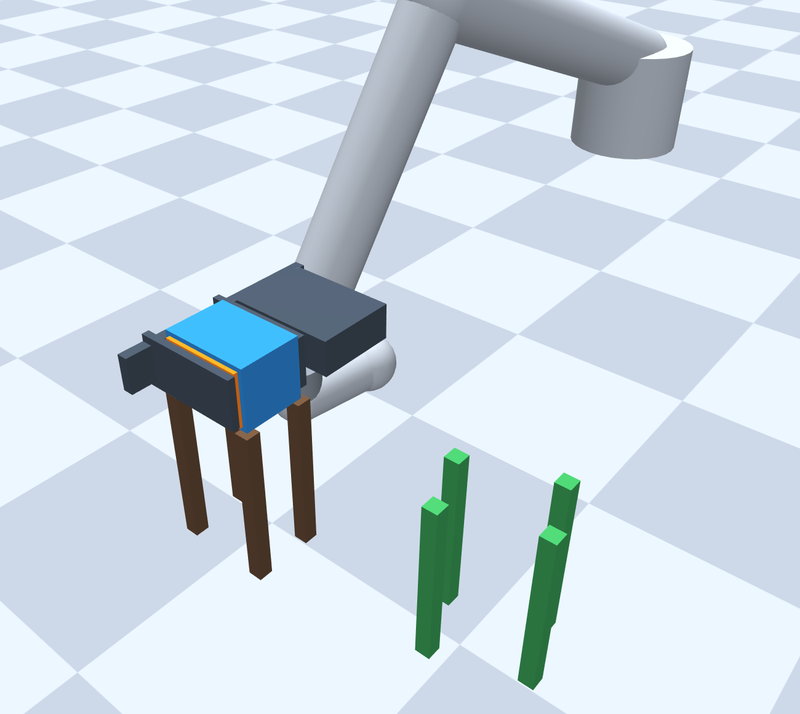}\\[-1pt]
  {\scriptsize 3. end of $\mathsf A$, $t=6.20$ s: grasp loaded\par}%
 \end{minipage}\\[4pt]
 \begin{minipage}[t]{.325\linewidth}\centering
  \includegraphics[width=\linewidth]{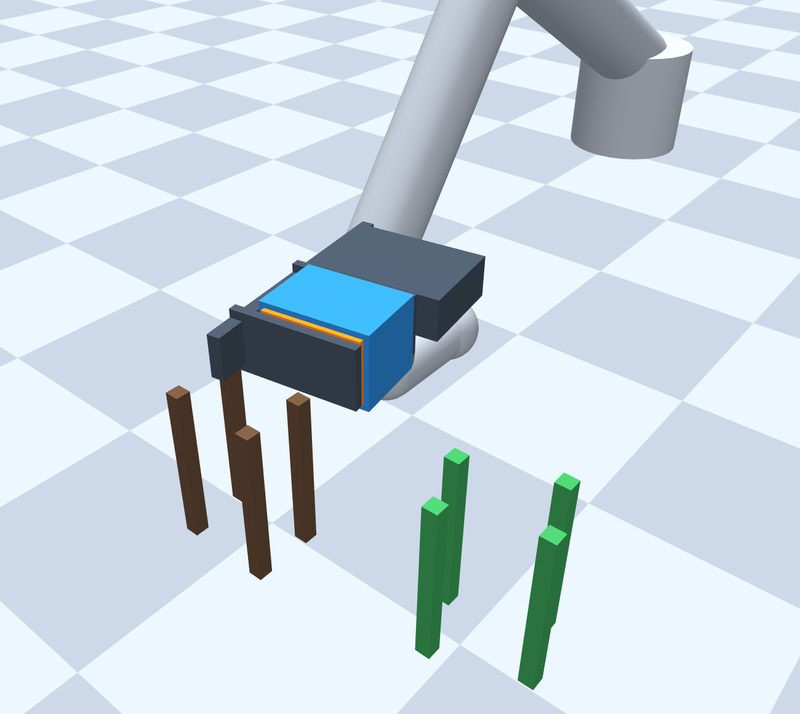}\\[-1pt]
  {\scriptsize 4. $\mathsf B$, $t=11.51$ s: lifted off the supports\par}%
 \end{minipage}\hfill
 \begin{minipage}[t]{.325\linewidth}\centering
  \includegraphics[width=\linewidth]{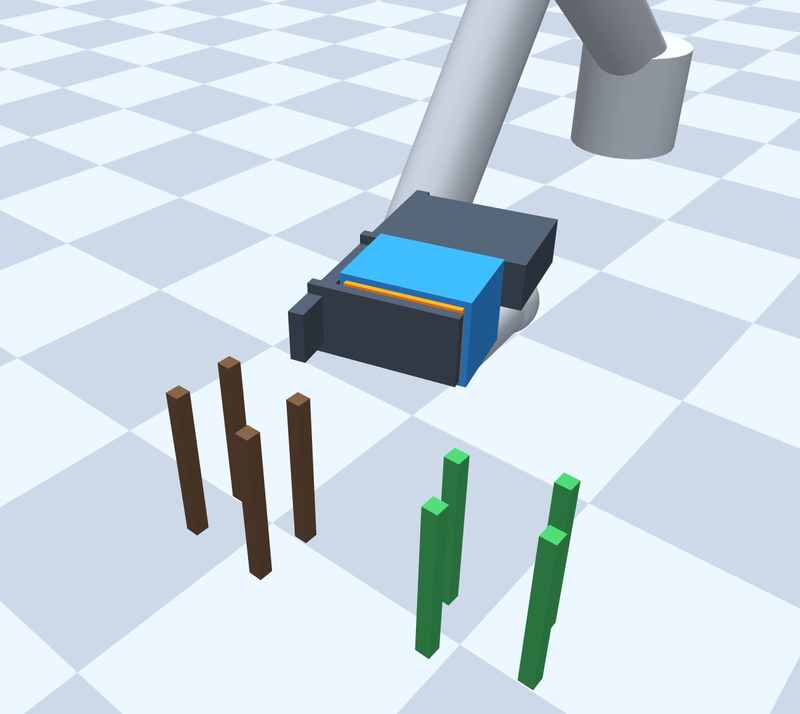}\\[-1pt]
  {\scriptsize 5. $\mathsf B$, $t=13.11$ s: carried toward the goal\par}%
 \end{minipage}\hfill
 \begin{minipage}[t]{.325\linewidth}\centering
  \includegraphics[width=\linewidth]{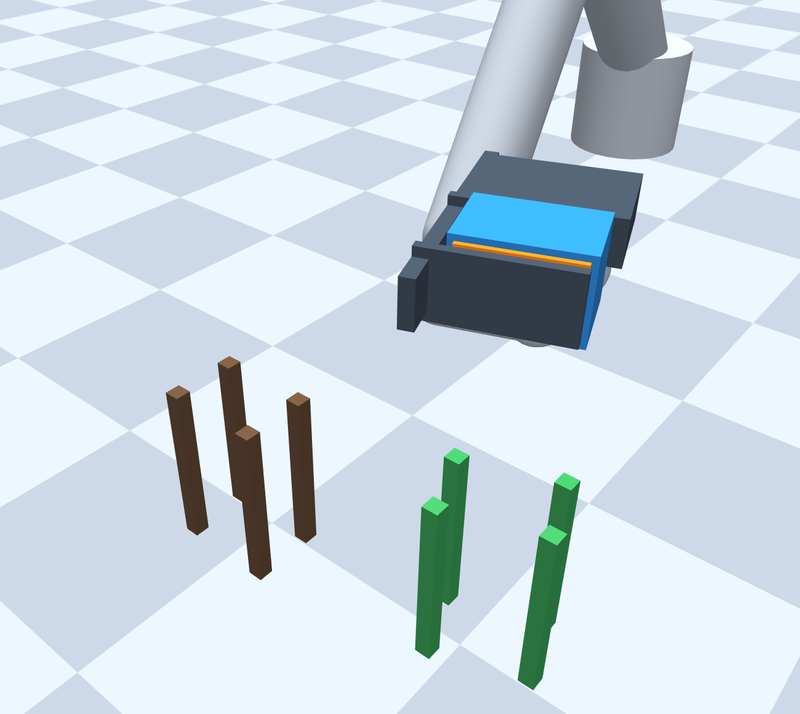}\\[-1pt]
  {\scriptsize 6. end of $\mathsf B$, $t=17.20$ s: above the goal supports\par}%
 \end{minipage}\\[4pt]
 \begin{minipage}[t]{.325\linewidth}\centering
  \includegraphics[width=\linewidth]{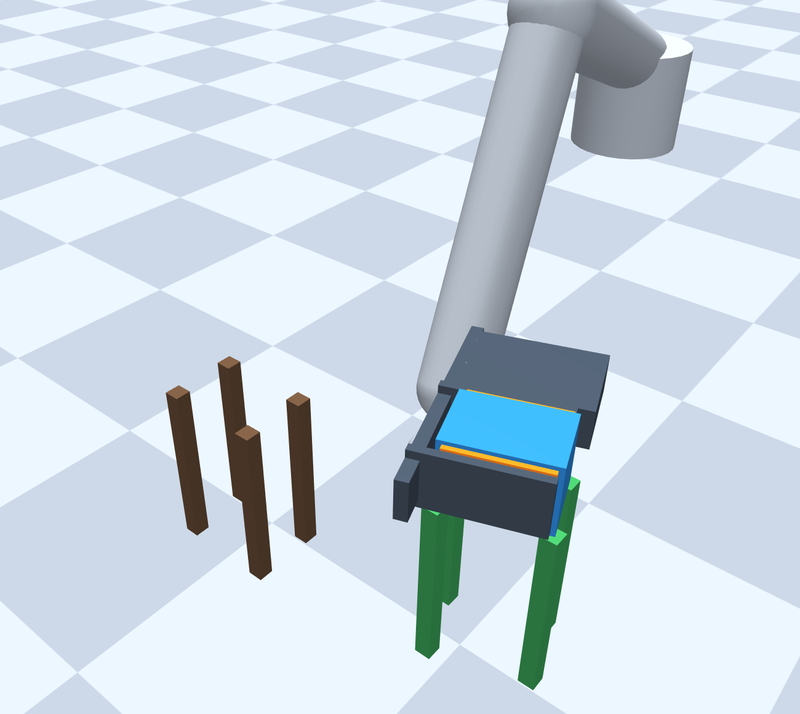}\\[-1pt]
  {\scriptsize 7. $\mathsf C$-place, $t=20.20$ s: set down, pads unloaded\par}%
 \end{minipage}\hfill
 \begin{minipage}[t]{.325\linewidth}\centering
  \includegraphics[width=\linewidth]{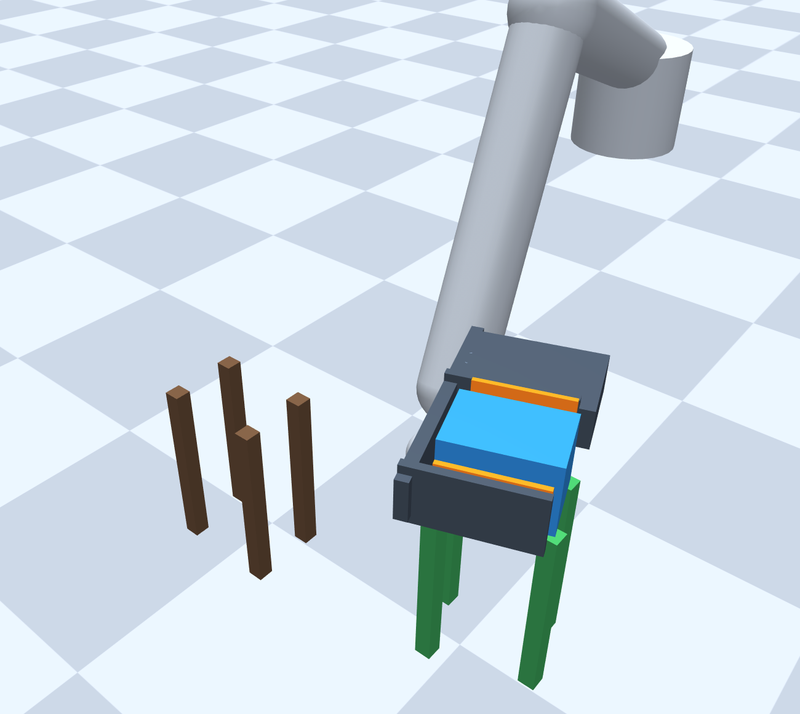}\\[-1pt]
  {\scriptsize 8. $\mathsf C$-open, $t=20.70$ s: jaws open\par}%
 \end{minipage}\hfill
 \begin{minipage}[t]{.325\linewidth}\centering
  \includegraphics[width=\linewidth]{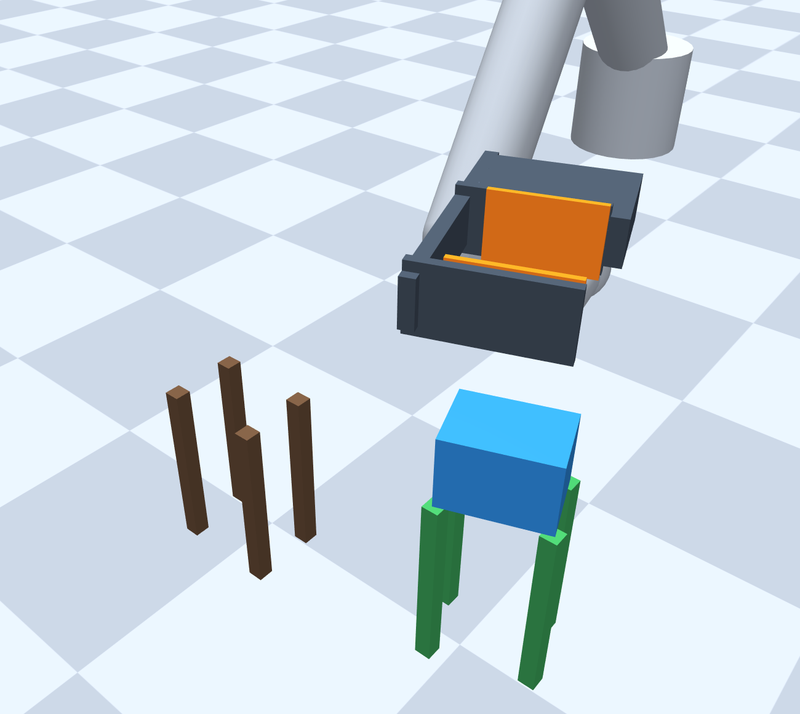}\\[-1pt]
  {\scriptsize 9. end, $t=23.00$ s: hand raised, box at the goal\par}%
 \end{minipage}
 \caption{Arm pick-and-place in MuJoCo, full-horizon run, close camera.
 Frames 1, 3, 6, 7, 8 and 9 are the plant states at the measured
 transitions between the problems. The blue box is the payload and the
 orange plates are the two gripper pads. The brown posts are the initial
 supports and the green posts the goal supports. The jaw rail, finger
 bodies and wrist flange are drawn for visibility only; they have no mass
 and take part in no contact, so that the contact model consists of the
 palm and the two pads.}
 \label{fig:arm-mujoco-key}
\end{figure}

\begin{figure}[tbp]
 \centering
 \begin{minipage}[t]{.24\linewidth}\centering
  \includegraphics[width=\linewidth]{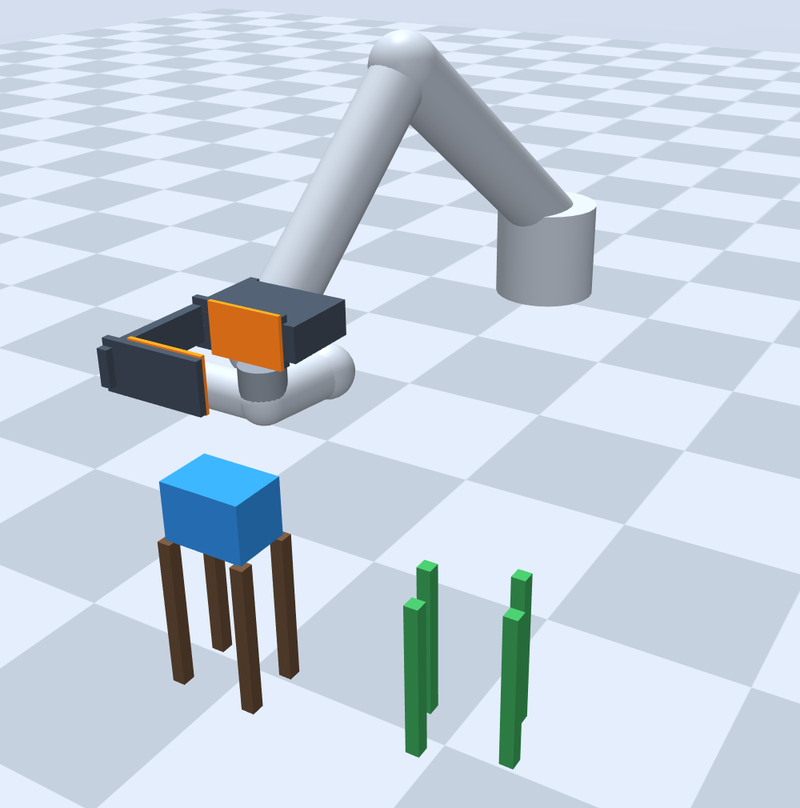}\\[-1pt]
  {\scriptsize $\mathsf A$, $t=0.31$ s\par}%
 \end{minipage}\hfill
 \begin{minipage}[t]{.24\linewidth}\centering
  \includegraphics[width=\linewidth]{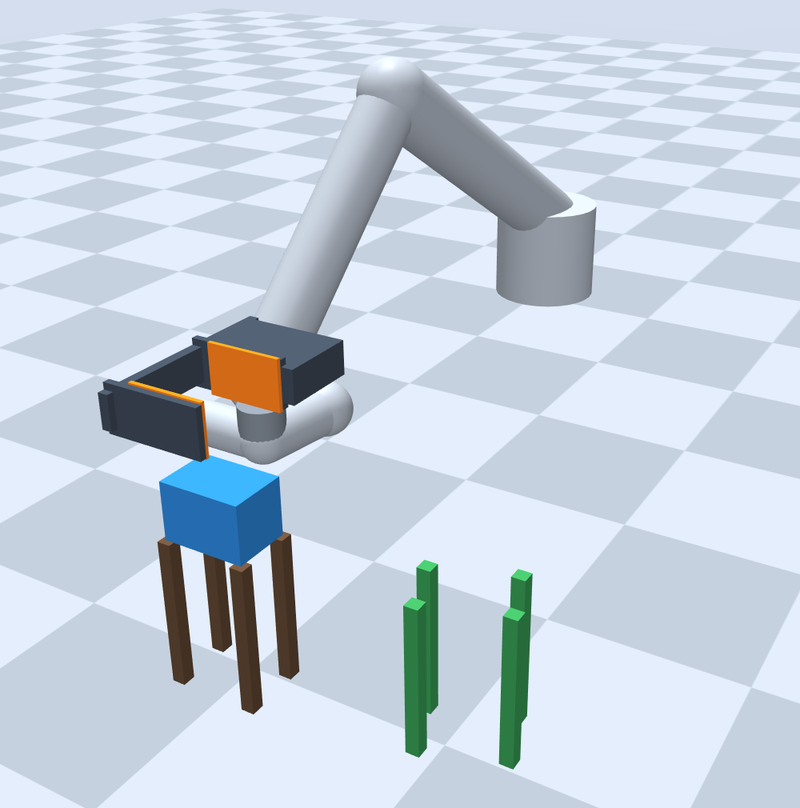}\\[-1pt]
  {\scriptsize $\mathsf A$, $t=2.31$ s\par}%
 \end{minipage}\hfill
 \begin{minipage}[t]{.24\linewidth}\centering
  \includegraphics[width=\linewidth]{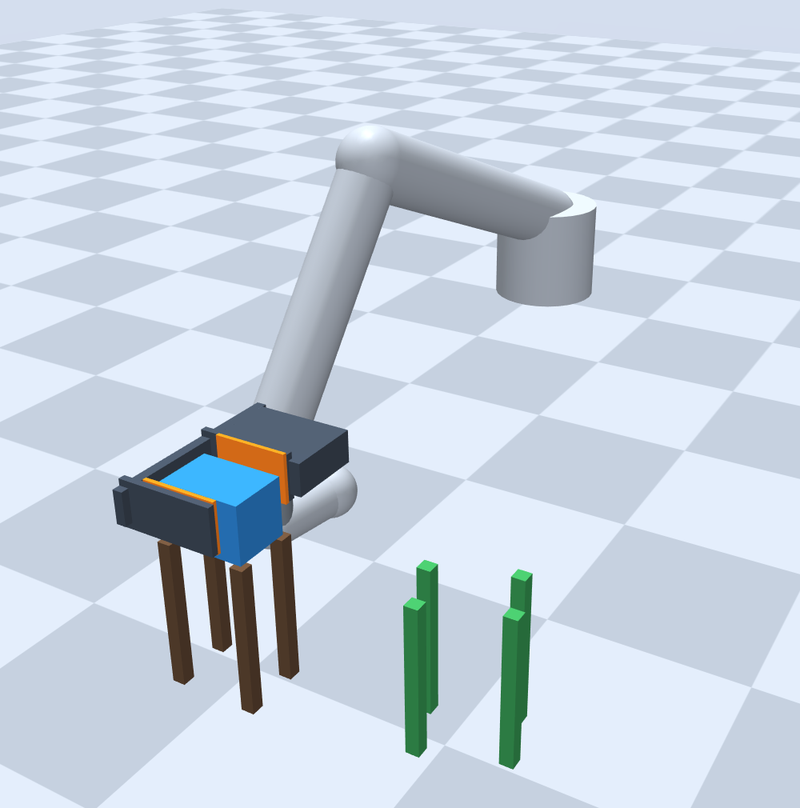}\\[-1pt]
  {\scriptsize $\mathsf A$, $t=4.31$ s\par}%
 \end{minipage}\hfill
 \begin{minipage}[t]{.24\linewidth}\centering
  \includegraphics[width=\linewidth]{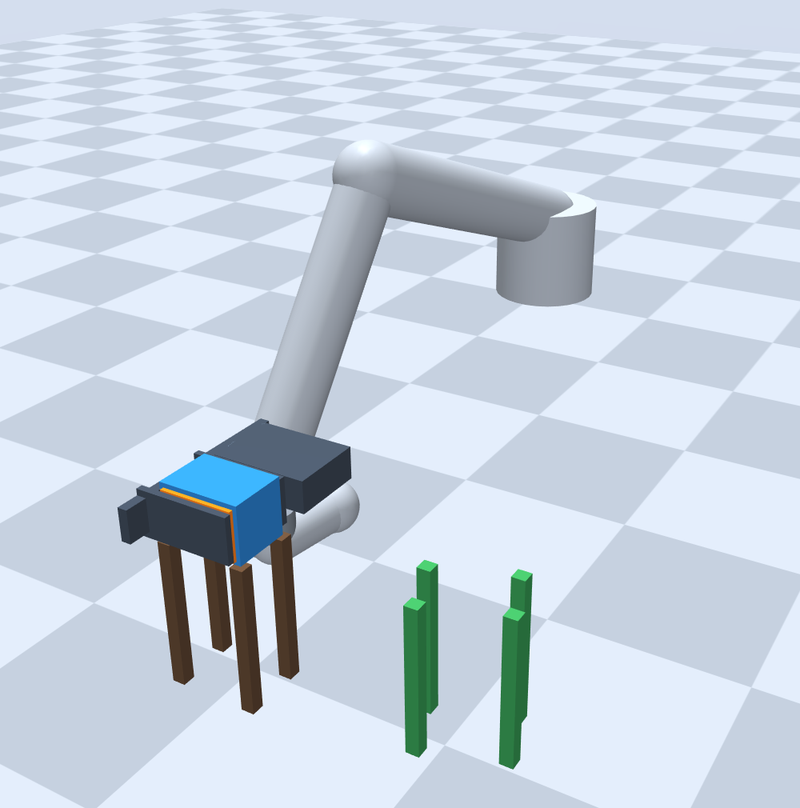}\\[-1pt]
  {\scriptsize $\mathsf A$, $t=5.71$ s\par}%
 \end{minipage}\\[3pt]
 \begin{minipage}[t]{.24\linewidth}\centering
  \includegraphics[width=\linewidth]{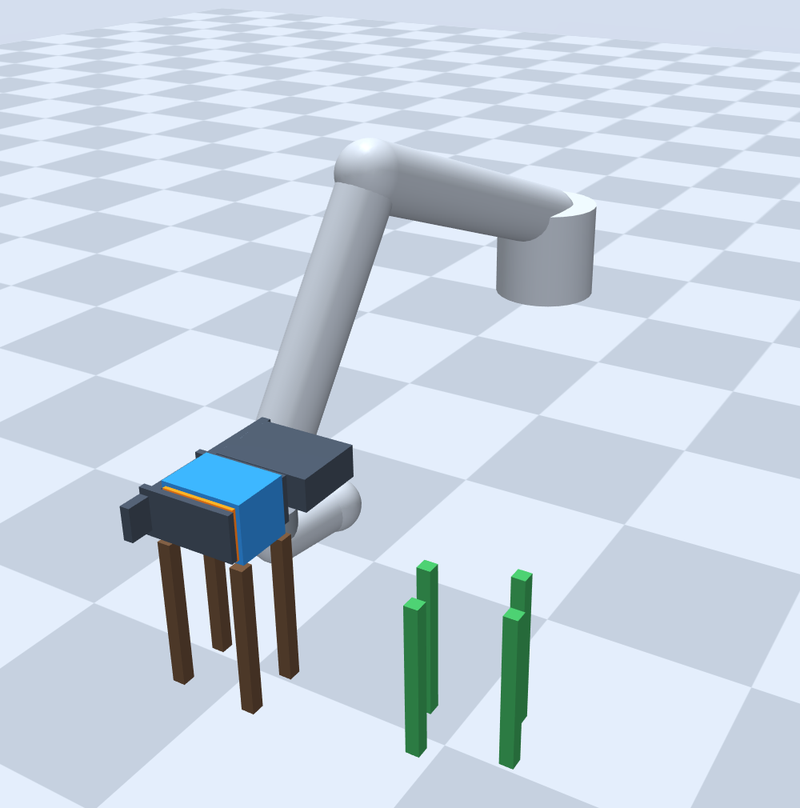}\\[-1pt]
  {\scriptsize $\mathsf B$, $t=7.51$ s\par}%
 \end{minipage}\hfill
 \begin{minipage}[t]{.24\linewidth}\centering
  \includegraphics[width=\linewidth]{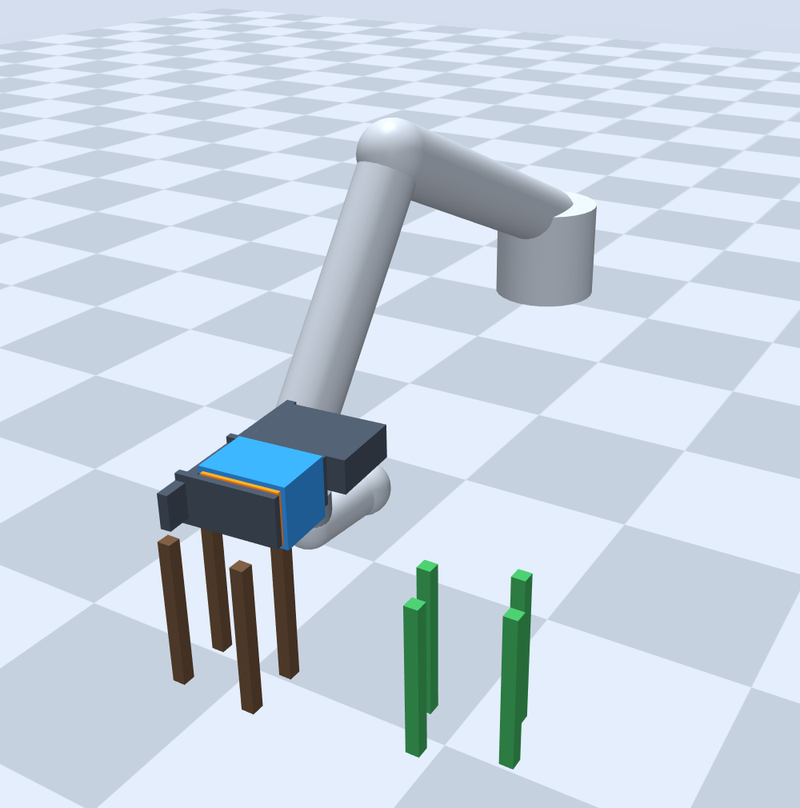}\\[-1pt]
  {\scriptsize $\mathsf B$, $t=9.91$ s\par}%
 \end{minipage}\hfill
 \begin{minipage}[t]{.24\linewidth}\centering
  \includegraphics[width=\linewidth]{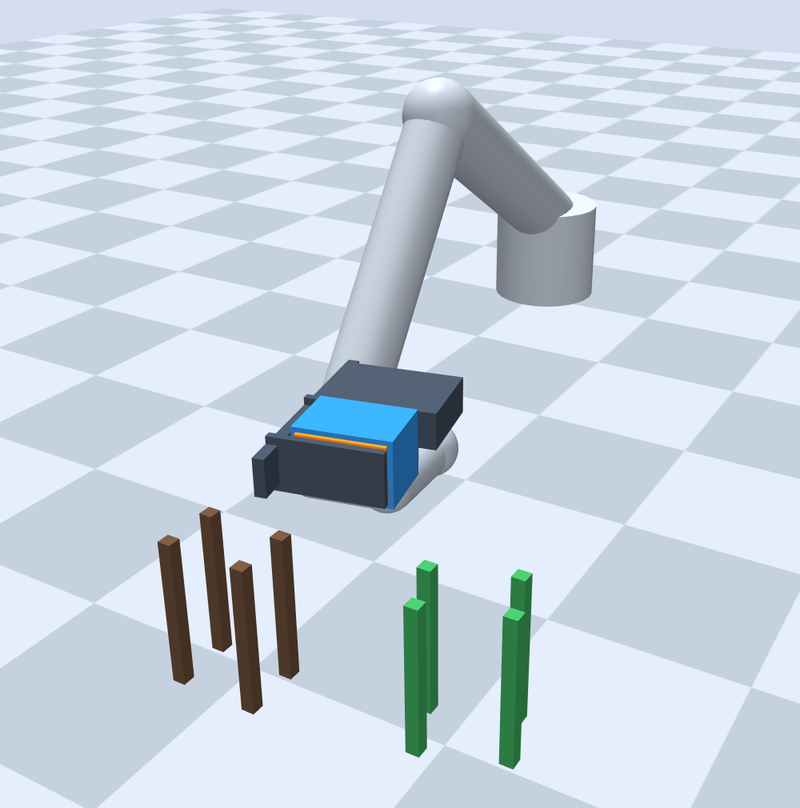}\\[-1pt]
  {\scriptsize $\mathsf B$, $t=12.31$ s\par}%
 \end{minipage}\hfill
 \begin{minipage}[t]{.24\linewidth}\centering
  \includegraphics[width=\linewidth]{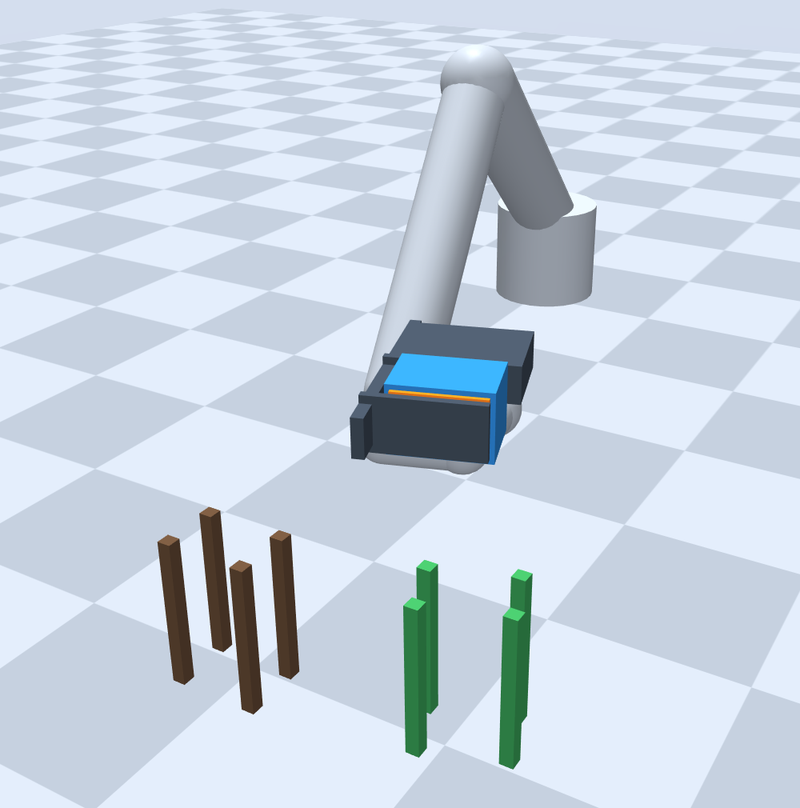}\\[-1pt]
  {\scriptsize $\mathsf B$, $t=14.71$ s\par}%
 \end{minipage}\\[3pt]
 \begin{minipage}[t]{.24\linewidth}\centering
  \includegraphics[width=\linewidth]{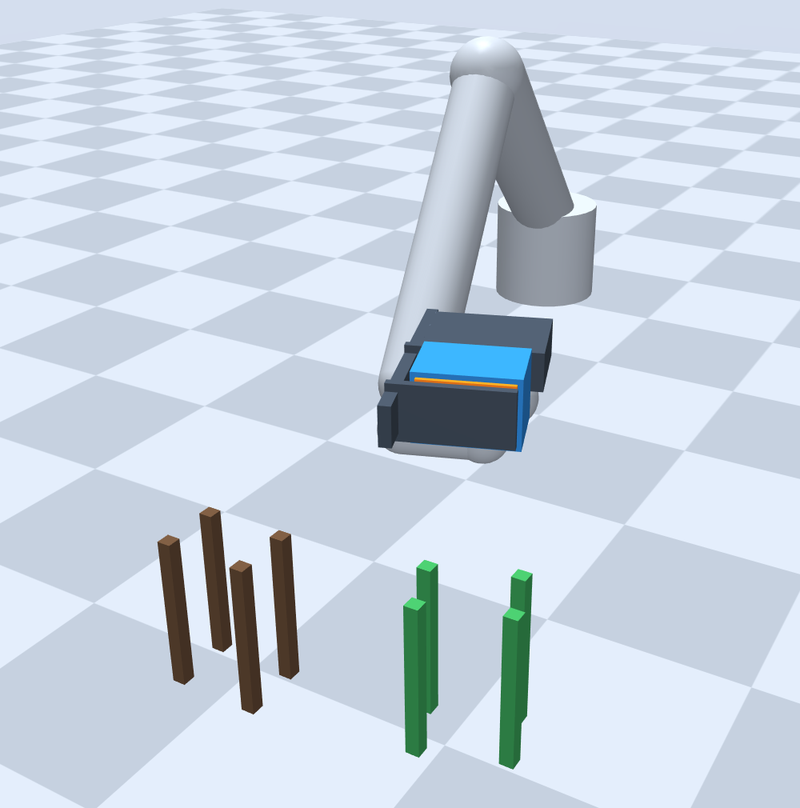}\\[-1pt]
  {\scriptsize $\mathsf B$, $t=17.11$ s\par}%
 \end{minipage}\hfill
 \begin{minipage}[t]{.24\linewidth}\centering
  \includegraphics[width=\linewidth]{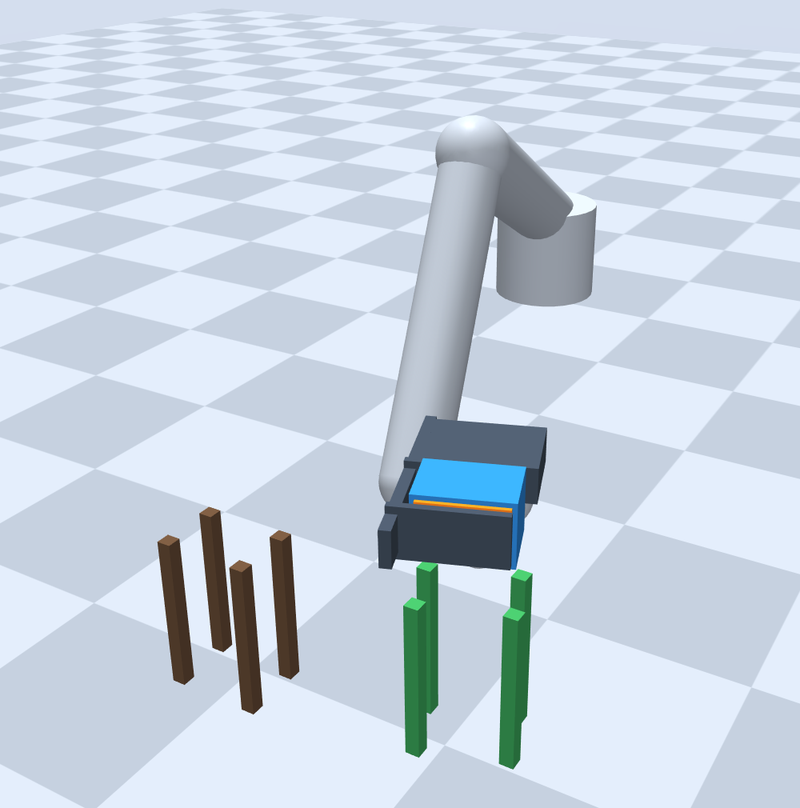}\\[-1pt]
  {\scriptsize $\mathsf C$-place, $t=18.71$ s\par}%
 \end{minipage}\hfill
 \begin{minipage}[t]{.24\linewidth}\centering
  \includegraphics[width=\linewidth]{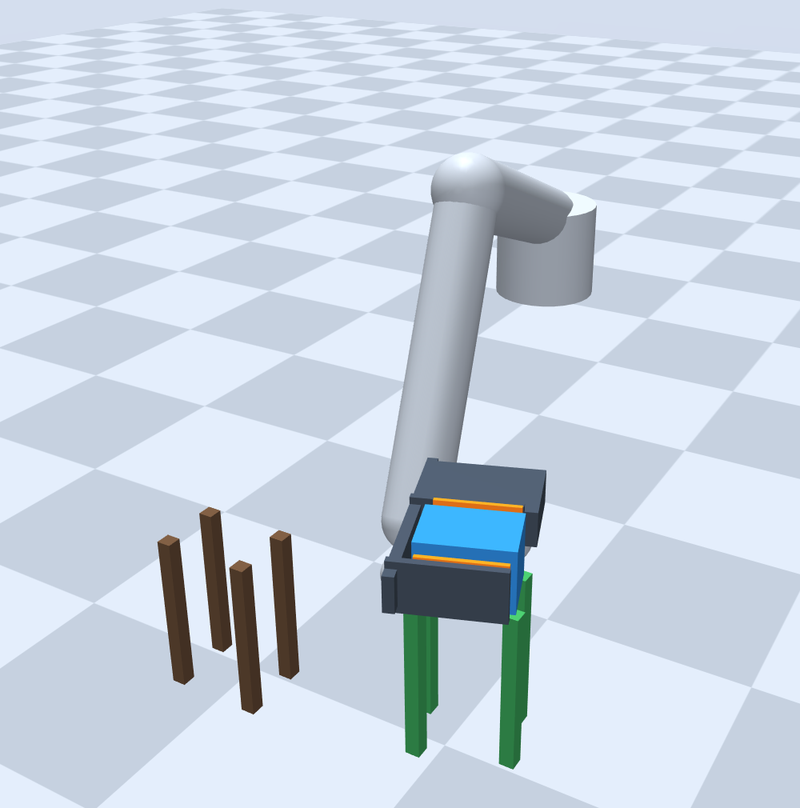}\\[-1pt]
  {\scriptsize $\mathsf C$-open, $t=20.51$ s\par}%
 \end{minipage}\hfill
 \begin{minipage}[t]{.24\linewidth}\centering
  \includegraphics[width=\linewidth]{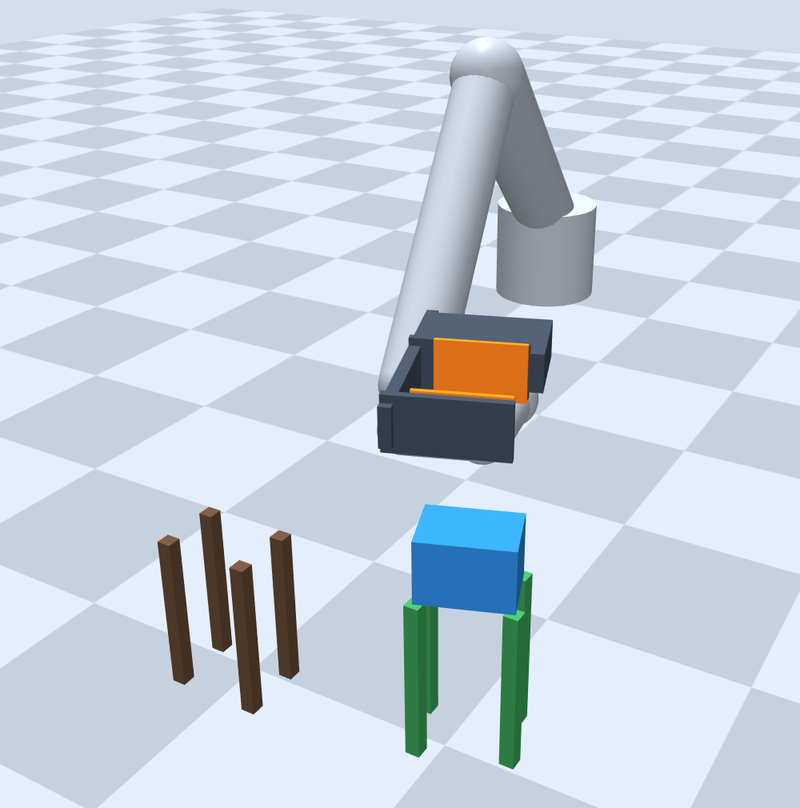}\\[-1pt]
  {\scriptsize $\mathsf C$-retreat, $t=22.51$ s\par}%
 \end{minipage}
 \caption{Arm pick-and-place in MuJoCo, full-horizon run, wide camera:
 twelve frames spread over the task. Problem $\mathsf A$ runs from $t=0.30$~s to
 $5.79$~s, $\mathsf B$ from $6.20$~s to $17.19$~s, $\mathsf C$-place from
 $17.20$~s to $19.69$~s, $\mathsf C$-open from $20.20$~s to $20.69$~s and
 $\mathsf C$-retreat from $20.70$~s to $22.69$~s.}
 \label{fig:arm-mujoco-film}
\end{figure}

The arm task consists of the three problems of
Section~\ref{subsec:arm-three-mpc}, solved in sequence with measured
transitions. The implementation treats the three layouts of Problem
$\mathsf C$ as separate phases: \emph{$\mathsf C$-place} is the held
placement, \emph{$\mathsf C$-open} is the supported release during which
the hand holds its position and the jaws open ($b_{\rm sep}=0$), and
\emph{$\mathsf C$-retreat} is the retreat after measured separation
($b_{\rm sep}=1$). Each phase is a shrinking-horizon problem with its own
pdSQP initialization. Two runs are reported:
\begin{itemize}
 \item the \emph{full-horizon run}, with sampling interval $h=10$~ms and
       horizons $H=550$, $1100$ and $250+50+200$ for $\mathsf A$,
       $\mathsf B$ and the three layouts of $\mathsf C$, traced at every
       25th sample;
 \item the \emph{reduced run}, with $h=20$~ms and horizons $H=100$, $150$
       and $100+30+100$, traced at every sample.
\end{itemize}
Figures~\ref{fig:arm-mujoco-key} and~\ref{fig:arm-mujoco-film} show the
full-horizon run in MuJoCo, from a camera close to the hand and from a
wider camera. The frames are rendered after the run from the recorded
plant states. The box is a free body that initially rests on
four posts. The arm approaches it, closes the two gripper pads, lifts and
carries it to a pose above a second set of posts, sets it down, opens the
jaws and returns to its home configuration.

\paragraph{Problem sizes and initialization.}
Table~\ref{tab:arm-init} lists, for every phase, the number $n$ of tangent
coordinates and the numbers $m_{\rm eq}$ and $m_{\rm in}$ of equality and
inequality rows of the full-space problem, together with the pdSQP
initialization. The largest problem, transport at $H=1100$, has
$85\,800$ variables, $59\,400$ equality rows and $77\,000$ inequality
rows. Every initialization converges to $r_K\leq2.3\times10^{-8}$ and
passes the supported-space curvature test
\eqref{eq:fresh-curvature-test}. The initializations of $\mathsf B$ and
$\mathsf C$-retreat need two iterations and that of $\mathsf C$-open
three. The initialization of $\mathsf C$-place is the slowest: after a
fast initial decrease the residual falls slowly, at a nearly constant
rate, over $68$ iterations at $h=10$~ms and $44$ at $h=20$~ms (Figure~\ref{fig:arm-pdsqp} in
Appendix~\ref{app:supplementary-figures}); $33$ and $9$ of these are
O-iterates, and all other initialization iterates of the arm are
V-iterates taken with unit step length.

\begin{table}[tbp]
\centering\small
\setlength{\tabcolsep}{4.5pt}
\begin{tabular}{llrrrrrrcr}
\toprule
Run & Phase & $H$ & $n$ & $m_{\rm eq}$ & $m_{\rm in}$ & iter. & time (s) & $r_K$ & act.\\
\midrule
$h=10$ ms & $\mathsf A$ & 550 & 30\,860 & 26\,962 & 27\,002 & 19 & 10.50 & $2.2\times10^{-8}$ & 15\\
 & $\mathsf B$ & 1100 & 85\,800 & 59\,400 & 77\,000 & 2 & 3.77 & $5.7\times10^{-11}$ & 0\\
 & $\mathsf C$-place & 250 & 19\,500 & 13\,500 & 17\,500 & 68 & 27.01 & $2.2\times10^{-8}$ & 6\\
 & $\mathsf C$-open & 50 & 3\,600 & 3\,200 & 3\,000 & 3 & 0.74 & $6.1\times10^{-9}$ & 2\\
 & $\mathsf C$-retreat & 200 & 14\,400 & 12\,807 & 11\,601 & 2 & 2.28 & $1.3\times10^{-9}$ & 0\\
\midrule
$h=20$ ms & $\mathsf A$ & 100 & 5\,660 & 4\,912 & 4\,952 & 19 & 2.57 & $2.3\times10^{-8}$ & 4\\
 & $\mathsf B$ & 150 & 11\,700 & 8\,100 & 10\,500 & 2 & 0.75 & $1.3\times10^{-10}$ & 0\\
 & $\mathsf C$-place & 100 & 7\,800 & 5\,400 & 7\,000 & 44 & 9.30 & $1.8\times10^{-8}$ & 5\\
 & $\mathsf C$-open & 30 & 2\,160 & 1\,920 & 1\,800 & 3 & 0.44 & $2.1\times10^{-9}$ & 1\\
 & $\mathsf C$-retreat & 100 & 7\,200 & 6\,407 & 5\,801 & 2 & 1.86 & $5.6\times10^{-9}$ & 0\\
\bottomrule
\end{tabular}
\caption{Arm: dimensions of the full-space problems and their pdSQP
initializations. ``iter.'' is the number of pdSQP iterations, ``time'' the
wall time of the initialization, $r_K$ the residual of the returned
point, and ``act.'' the number of strongly active inequality rows at that
point.}
\label{tab:arm-init}
\end{table}

\begin{figure}[tbp]
 \centering
 \includegraphics[width=\linewidth]{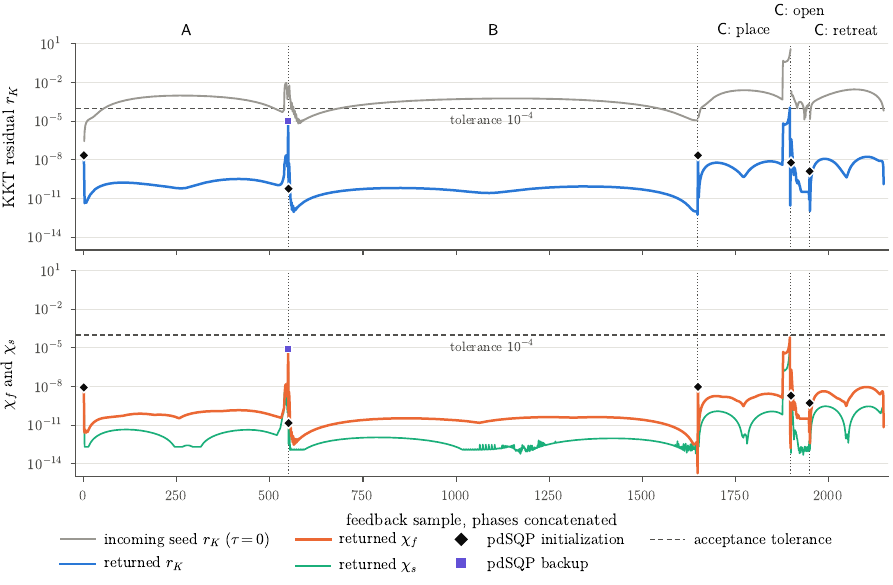}
 \caption{Arm, full-horizon run ($h=10$~ms): residuals of the original
 nonlinear problem at every feedback sample, phases concatenated. Top: KKT
 residual $r_K$ of the incoming shifted seed ($\tau=0$, gray) and of the
 returned point (blue). Bottom: feasibility $\chi_f$ and stationarity
 $\chi_s$ of the returned point. Diamonds mark the pdSQP initialization of
 each phase.}
 \label{fig:arm-full-closed-loop}
\end{figure}

\paragraph{Tracking over the closed loop.}
Table~\ref{tab:arm-closed-loop} and
Figure~\ref{fig:arm-full-closed-loop} give the residuals at every feedback
sample of the full-horizon run; Figure~\ref{fig:arm-reduced-closed-loop}
in Appendix~\ref{app:supplementary-figures} is the corresponding figure
for the reduced run.

In the full-horizon run, $2144$ of the $2145$ feedback samples are
accepted online, all but two of them after a single correction round. The
shifted seed has a median KKT residual of $4.7\times10^{-4}$, and it
exceeds the tolerance $10^{-4}$ at $84\%$ of the online samples. One
correction reduces it to a median of $7.7\times10^{-11}$, with median
feasibility $\chi_f=4.1\times10^{-11}$ and stationarity
$\chi_s=9.4\times10^{-13}$; $95\%$ of the accepted points have
$r_K<10^{-8}$. The reduction factor per sample has a median of $6.5$
orders of magnitude and is never smaller than $3.6$ orders. The single
backup occurs at the last sample of $\mathsf A$ ($H=1$) and returns
$r_K=1.1\times10^{-5}$.

The largest residuals occur at the end of the held placement. In the last
$22$ samples of $\mathsf C$-place the residual of the seed rises above
$0.05$ and reaches $4.0$. The corrections still return points that pass
the acceptance test, with $r_K$ at most $1.1\times10^{-4}$ (median
$6.3\times10^{-6}$); the last two samples need a second round. The
initialization of the next phase restores $r_K=6.1\times10^{-9}$.

The reduced run behaves in the same way at a coarser level: $473$ of $475$
samples are accepted online, the median seed residual is
$3.0\times10^{-3}$ and the median accepted residual
$7.1\times10^{-8}$. Its two backups are at $H=3$ and $H=1$ of
$\mathsf A$; the task log attributes them to stopped rank or factor tests
of the working-set tasks. Relative to the reduced run, the full-horizon run
halves the sampling interval and also lengthens the approach and transport
phases, so the problem data change less between consecutive samples. The seed residual is then
smaller by a factor of about six and the accepted residual by about three
orders of magnitude.

\begin{table}[tbp]
\centering\footnotesize
\setlength{\tabcolsep}{3pt}
\begin{tabular}{@{}llrcccccrr@{}}
\toprule
 & & & & \multicolumn{2}{c}{largest accepted} & \multicolumn{2}{c}{median} & \multicolumn{2}{c}{step time (ms)}\\
\cmidrule(lr){5-6}\cmidrule(lr){7-8}\cmidrule(lr){9-10}
Run & Phase & samples & on./back. & component & $r_K$ & seed $r_K$ & $r_K$ & mean & max\\
\midrule
$h=10$ ms & $\mathsf A$ & 549 & 548 / 1 & $8.6\times10^{-6}$ & $1.1\times10^{-5}$ & $5.9\times10^{-4}$ & $1.5\times10^{-10}$ & 734 & 4457\\
 & $\mathsf B$ & 1099 & 1099 / 0 & $4.2\times10^{-9}$ & $5.9\times10^{-9}$ & $3.5\times10^{-4}$ & $5.0\times10^{-11}$ & 1167 & 2319\\
 & $\mathsf C$-place & 249 & 249 / 0 & $6.6\times10^{-5}$ & $1.1\times10^{-4}$ & $1.8\times10^{-3}$ & $4.6\times10^{-9}$ & 515 & 2222\\
 & $\mathsf C$-open & 49 & 49 / 0 & $2.0\times10^{-7}$ & $4.0\times10^{-7}$ & $2.2\times10^{-4}$ & $5.7\times10^{-11}$ & 111 & 224\\
 & $\mathsf C$-retreat & 199 & 199 / 0 & $9.3\times10^{-9}$ & $1.7\times10^{-8}$ & $1.6\times10^{-3}$ & $8.3\times10^{-9}$ & 321 & 787\\
\midrule
$h=20$ ms & $\mathsf A$ & 99 & 97 / 2 & $7.1\times10^{-5}$ & $1.7\times10^{-4}$ & $3.6\times10^{-3}$ & $9.8\times10^{-8}$ & 162 & 320\\
 & $\mathsf B$ & 149 & 149 / 0 & $4.3\times10^{-8}$ & $1.0\times10^{-7}$ & $2.7\times10^{-3}$ & $5.5\times10^{-8}$ & 243 & 546\\
 & $\mathsf C$-place & 99 & 99 / 0 & $2.5\times10^{-5}$ & $3.8\times10^{-5}$ & $4.4\times10^{-3}$ & $1.1\times10^{-7}$ & 224 & 469\\
 & $\mathsf C$-open & 29 & 29 / 0 & $8.6\times10^{-7}$ & $1.6\times10^{-6}$ & $8.8\times10^{-4}$ & $1.7\times10^{-9}$ & 85 & 154\\
 & $\mathsf C$-retreat & 99 & 99 / 0 & $1.4\times10^{-7}$ & $1.8\times10^{-7}$ & $3.2\times10^{-3}$ & $9.4\times10^{-8}$ & 260 & 497\\
\bottomrule
\end{tabular}
\caption{Arm: closed-loop tracking by phase. ``on./back.'' counts the
samples accepted online and those resolved by the pdSQP backup; no sample
was retained. ``Component'' is the largest of $\chi_f,\chi_s,\chi_c$ over
the accepted points of the phase, which the acceptance test bounds by
$\epsilon=10^{-4}$. The medians are over all samples of the phase.}
\label{tab:arm-closed-loop}
\end{table}

\begin{figure}[tbp]
 \centering
 \includegraphics[width=\linewidth]{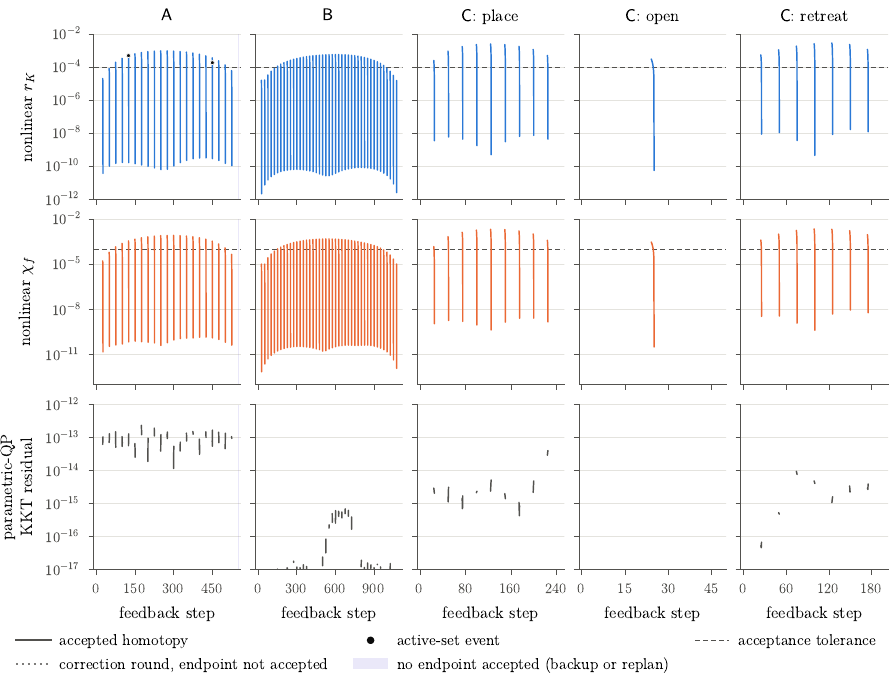}
 \caption{Arm, full-horizon run: residuals along the homotopies traced at
 every 25th feedback step. Each curve occupies the unit interval of its
 step and runs from $\tau=0$ to $\tau=1$, so it appears as a nearly
 vertical segment from the residual of the seed down to that of the
 endpoint. Top and middle: nonlinear KKT residual $r_K$ and feasibility
 $\chi_f$. Bottom: KKT residual of the parametric QP; values below
 $10^{-17}$ are drawn at $10^{-17}$.}
 \label{fig:arm-full-unrolled}
\end{figure}

\paragraph{Residuals along the homotopies.}
Figure~\ref{fig:arm-full-unrolled} shows the traced homotopies of the
full-horizon run, and
Figures~\ref{fig:arm-reduced-unrolled}--\ref{fig:arm-reduced-tau} in
Appendix~\ref{app:supplementary-figures} show every homotopy of the reduced
run. The three observations made for the drone hold again.

The parametric-QP residual is at most $2.4\times10^{-13}$ on the $81$
traced accepted paths of the full-horizon run. In the reduced run it is at
most $2.7\times10^{-14}$ on $472$ of the $473$ accepted paths and
$6.4\times10^{-10}$ on the remaining one, at the end of $\mathsf A$.

Over an accepted path the nonlinear KKT residual decreases by $4.8$ to
$7.3$ orders of magnitude in the full-horizon run (median $6.6$) and by
$3.1$ to $8.9$ orders in the reduced run (median $4.5$). The feasibility
residual decreases along every accepted path, by $5.1$ to $7.5$ orders
(median $6.8$) and by $3.3$ to $9.4$ orders (median $4.8$) respectively.
At each of the nine interior checkpoints $\tau=0.1,\ldots,0.9$ of all
$554$ traced accepted paths the residual equals $(1-\tau)\,r_K(0)$ to
within a relative deviation of $2.1\times10^{-5}$, so the whole decrease
below that line occurs between $\tau=0.9$ and $\tau=1$.

Active-set events are rare. In the reduced run $26$ of the $473$ accepted
paths contain at least one event ($19$ in $\mathsf A$, $3$ in
$\mathsf C$-place and $4$ in $\mathsf C$-open), with at most four events
on one path; all other accepted paths go from $\tau=0$ to $\tau=1$ on a
single working set. In $\mathsf A$ and $\mathsf C$-place a second
task with an added active row (\texttt{ACTIVE\_IN}) runs beside the
warm-started task, and at $22$ of the $29$ samples of $\mathsf C$-open a
second task with a weakly active row removed (\texttt{ACTIVE\_OUT}).
Whenever both tasks reach the target QP they return the same endpoint
($218$ of the $221$ two-task rounds; in the other three a task stopped
early), and the warm-started task is selected in $472$ of the $473$
accepted samples.

\paragraph{Measured transitions and task outcome.}
The supervisor advances between problems only on measured tests
(Section~\ref{subsec:arm-three-mpc}). Table~\ref{tab:arm-gates} lists the
recorded values; every test passed in both runs. The tests are the
following. After $\mathsf A$, each pad carries more than $1$~N and the
largest component of the box twist is below $0.05$. After $\mathsf B$, the
box is within $10$~mm and $0.05$~rad of the pre-placement pose, the grasp
has slipped by less than $5$~mm, and both pads are still loaded. After
$\mathsf C$-place, at least three posts carry the box, its twist is below
$0.02$ and the pads carry less than $0.2$~N. After $\mathsf C$-open, both
pad forces are zero and the aperture exceeds the grasp aperture by more
than $10$~mm. At the end, the box is within $10$~mm and $0.05$~rad of the
goal pose and at rest, and the arm is within $0.05$~rad of its home
configuration. The final placement error of
the box is $1.2$~mm and $1.4$~mm, the grasp slips by at most
$3.7\times10^{-5}$~m during transport, and no contact outside the
designated pad and support contacts is recorded.

\begin{table}[tbp]
\centering\footnotesize
\setlength{\tabcolsep}{4pt}
\begin{tabular}{@{}llcc@{}}
\toprule
Transition & Measured quantity & $h=10$ ms & $h=20$ ms\\
\midrule
End of $\mathsf A$ & pad normal forces, left / right (N) & $5.40$ / $5.41$ & $5.40$ / $5.42$\\
 & jaw aperture (m) & $0.0687$ & $0.0687$\\
End of $\mathsf B$ & grasp slip (m) & $1.02\times10^{-5}$ & $3.71\times10^{-5}$\\
  & box position error (m) & $1.39\times10^{-5}$ & $5.24\times10^{-5}$\\
  & box attitude error (rad) & $8.72\times10^{-6}$ & $5.54\times10^{-5}$\\
 & pad normal forces, left / right (N) & $5.41$ / $5.41$ & $5.41$ / $5.41$\\
End of $\mathsf C$-place & support force on the four posts (N) & $1.18$ & $1.18$\\
 & residual pad load (N) & $0.00$ & $0.00$\\
  & placement error (m, rad) & $1.25\times10^{-3}$, $3.82\times10^{-4}$ & $1.47\times10^{-3}$, $1.09\times10^{-3}$\\
End of $\mathsf C$-open & jaw aperture (m) & $0.1099$ & $0.1117$\\
 & pad normal forces (N) & $0.00$ / $0.00$ & $0.00$ / $0.00$\\
End of task & box placement error (m, rad) & $1.22\times10^{-3}$, $1.10\times10^{-5}$ & $1.36\times10^{-3}$, $4.78\times10^{-4}$\\
 & arm home error (rad) & $2.98\times10^{-6}$ & $6.26\times10^{-6}$\\
 & undesignated contacts & $0$ & $0$\\
 & final simulation time (s) & $23.00$ & $11.10$\\
\bottomrule
\end{tabular}
\caption{Arm: quantities measured in MuJoCo at the transitions between
problems and at the end of the task.}
\label{tab:arm-gates}
\end{table}

\begin{figure}[tbp]
 \centering
 \includegraphics[width=\linewidth]{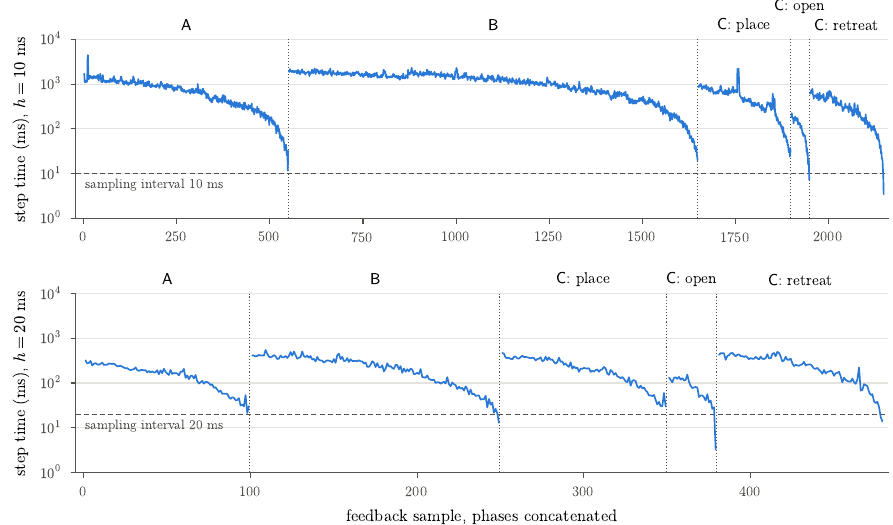}
 \caption{Arm: wall time of every feedback step for the full-horizon run
 (top) and the reduced run (bottom), with the sampling interval for
 comparison.}
 \label{fig:arm-step-time}
\end{figure}

\paragraph{Computation time.}
The arm controller does not run in real time in this implementation
(Figure~\ref{fig:arm-step-time}). In the full-horizon run a feedback step
takes $0.88$~s on average, against the sampling interval of $10$~ms, and
the simulation waits for the controller: the $22.4$~s of simulated motion
take $1883$~s of computation. The cost of a step is dominated by the exact
second-order differentiation of the arm model, not by the linear algebra. The step time decreases with the remaining
horizon within each phase, from about $1$--$2$~s at the start of
$\mathsf A$ and $\mathsf B$ to less than $40$~ms at the end, and it is
below the sampling interval at three samples only. The reduced run takes
$0.22$~s per step on average, against $20$~ms. Its timings were recorded
while unrelated jobs occupied three of the eight cores; an earlier run on
the idle machine gave mean step times between $67$ and $179$~ms per phase
with identical residuals. The drone controller, in contrast, runs well
inside its sampling interval except at the two replanning samples of the
perturbed run.

\subsection{Summary of the experiments}
\label{subsec:results-summary}
On the nominal plants one correction per sample suffices almost
everywhere: $2626$ of the $2631$ feedback samples of the three nominal
runs are accepted online, and the five exceptions lie within the last
three samples of a shrinking horizon. One correction reduces the KKT
residual of the original nonlinear problem by a median of $4.5$ to $6.9$
orders of magnitude, and the feasibility residual of the arm by a median
of $4.8$ to $6.8$ orders. The QP homotopies are followed to rounding
accuracy, so the residual that remains is the linearization error of one
SQP step, which is the quantity bounded in
Theorem~\ref{thm:riemannian-pathfollowing}. Under a plant perturbation the
pre-correction error grows by one to two orders of magnitude; the method
then needs more correction rounds, resolves eight of the eighteen samples
by the pdSQP backup and replans twice, and still completes the flight with
a positive clearance. The two limitations observed are the computation time
of the arm controller and the reliance on the backup near the end of a
shrinking horizon, where the regularity assumptions of
Section~\ref{sec:pathfollowing-theory} are not expected to hold.
\FloatBarrier
\subsection{Geometric properties of the arm formulation}
\label{sec:arm-qualification-details}

The finite-pad model of Section~\ref{sec:arm-example} provides explicit
conditions for force closure and for the constraint qualifications
used in Section~\ref{sec:pathfollowing-theory}.
The following propositions describe these properties and the physical
force-allocation directions.
\paragraph{Analytic force closure of the finite-pad grasp.}
The finite-pad geometry and wrench scaling are those in
\eqref{eq:arm-contact-geometry}--\eqref{eq:arm-rank-matching}.
The following calculation supplies the rank and positivity margins used
for the nominal grasp and its nearby measured geometries.

\begin{proposition}[Verified force closure of the finite-pad model]
\label{prop:arm-analytic-force-closure}
For \eqref{eq:arm-contact-geometry}, $\mu_f>0$, and $L\geq3$, the
matrix $G_{\rm ray}$ has rank six, and
$\alpha_*=(8L)^{-1}\mathbf1$ is a strictly positive normalized null
certificate.  In particular $\varepsilon_*=(8L)^{-1}$ in
\eqref{eq:arm-force-closure}.  With $c_f=\mu_f^2/2$,
\begin{equation}\label{eq:arm-ray-gram}
 \begin{aligned}
 G_{\rm ray}G_{\rm ray}^\top
 &=8L\operatorname{diag}(d_1,d_2',d_3',1,c_f,c_f),\\
 d_1&=c_f(d_2^2+d_3^2),\quad
 d_2'=c_fa_1^2+d_3^2,\quad d_3'=c_fa_1^2+d_2^2.
 \end{aligned}
\end{equation}
Consequently the scaled map satisfies $WW^\top\succeq\beta_g^2I_6$, where
\begin{equation}\label{eq:arm-rank-margin}
 \beta_g^2=8L\min\{d_1/\ell_b^2,d_2'/\ell_b^2,
                       d_3'/\ell_b^2,1,c_f\}>0.
\end{equation}
\end{proposition}
\begin{proof}
Equally spaced angles have zero first trigonometric sums and second
moments $L/2$ in each tangent direction.  Hence
$\sum_\ell d_{j\ell}=Ln_j$ and
$\sum_\ell d_{j\ell}d_{j\ell}^\top
 =L\operatorname{diag}(1,c_f,c_f)$, since every face normal is
$\pm e_1$.  Summing over the independent signs in
\eqref{eq:arm-contact-geometry} cancels all first and mixed spatial
moments.  It gives both $G_{\rm ray}\mathbf1=0$ and the diagonal
Gram matrix \eqref{eq:arm-ray-gram}.  Every displayed diagonal entry
is positive.  Normalization implies
$1=\sum_i\alpha_i\geq8L\varepsilon$, and the uniform certificate
attains that upper bound.  Scaling the moment rows gives
\eqref{eq:arm-rank-margin}.
\end{proof}
The calculation verifies force closure analytically for the specified finite
pads. A two-point pinch or a small acquisition cost alone provides insufficient
information for this verification. For measured perturbations, the same regular
stratum and the positive-certificate section in \eqref{eq:arm-closure-section} apply.  A bound
$\|W-W_0\|_2<\beta_g$ preserves full rank of the calibrated nominal
map $W_0$; positivity must also remain inside its declared section.

\paragraph{Qualifications of the geometric closure block.}
The normalized closure equality is \eqref{eq:arm-geometric-closure}.
Its positive-certificate domain is open. The factor and certificate
columns therefore permit the following rank and lifting argument.

\begin{proposition}[Regular rank and force-closure manifold]
\label{prop:arm-closure-cq}
On \eqref{eq:arm-rank-manifold}--\eqref{eq:arm-geometric-closure},
the matching and normalized-certificate equalities have independent
differentials, of joint rank $6m_r+7$.  They therefore satisfy
Riemannian MFCQ and RCRQ.  If a base OCP system $C_0=0$, $b_0\geq0$
is independent of the auxiliary factors and certificate, adjoining
these equalities preserves its MFCQ and RCRQ.  More precisely, for
each subset $J$ of its active inequalities,
\begin{equation}\label{eq:arm-closure-rank-inheritance}
 \operatorname{rank}D(C_0,b_{0,J},K_g,H_g)
 =\operatorname{rank}D(C_0,b_{0,J})+6m_r+7.
\end{equation}
\end{proposition}
\begin{proof}
For tangent increments the factor differential is
$\dot U_gT_g+U_g\dot T_g$, where
$U_g^\top\dot U_g$ is skew symmetric and $\dot T_g$ is upper
triangular.  If this differential is zero, then
$U_g^\top\dot U_g=-\dot T_gT_g^{-1}$ is also upper triangular,
and hence is zero.  It follows that both increments vanish.
Equal dimensions show that the differential is an isomorphism.
Thus $K_g$ is a regular graph constraint.  Next
$D_\alpha H_g=[U_g^\top;\mathbf1^\top]$ has rank seven: a row
relation $u^\top U_g^\top+\lambda\mathbf1^\top=0$, multiplied by
$\alpha$, gives $\lambda=0$, and then $u=0$.
The separate auxiliary columns prove
\eqref{eq:arm-closure-rank-inheritance} on a neighborhood.
Every base MFCQ direction lifts by first solving the factor
matching differential and then solving
$D_\alpha H_g[\dot\alpha]=-D_{U_g}H_g[\dot U_g]$.
These lifts do not change any base inequality derivative.  There
are no active positivity rows on this domain.  Equality
surjectivity proves MFCQ for the block itself; its locally constant
rank proves RCRQ, and \eqref{eq:arm-closure-rank-inheritance}
proves inheritance for every base active subset as required by
Definition~\ref{def:rcrq}.
\end{proof}

\paragraph{Qualifications of the complete transcription.}
The three-problem construction removes a qualification obstruction
arising in a monolithic contact transcription. If an active gap inequality
$b_i\geq0$ is already an identity of the mode equalities, then
$DC\,d=0$ implies $Db_i\,d=0$; hence no MFCQ direction can strictly
increase that gap. Retaining such an inequality invalidates MFCQ even
when redundant rows satisfy RCRQ. The transcription omits mode-implied pad/support gaps and bounds on substituted
constants before testing the remaining qualifications.
The measured geometry is still checked physically, and a mismatch rejects
the corresponding mode.  Propositions~\ref{prop:arm-analytic-force-closure}
and~\ref{prop:arm-closure-cq} verify the force-closure contribution to
MFCQ and RCRQ throughout every admitted held segment.  The following
criterion makes the remaining complete-OCP qualification check explicit.

\begin{proposition}[Complete-OCP regularity on an admitted tube]
\label{prop:arm-complete-cq}
Use a fixed local layout of $\mathsf A$, $\mathsf B$, or $\mathsf C$,
with measured node zero substituted and mode-implied rows removed.
Let $C_0$ denote its retained equalities before the regular geometric
lift.  Suppose a tangent-coordinate QR factorization at a reference
point gives
\[
 D_rC_0(\bar r,\bar t)^\top=Q_0R_0,\qquad
 \|R_0^{-1}\|_2\leq\kappa_0,
 \qquad \|D_rC_0(r,t)-D_rC_0(\bar r,\bar t)\|_2\leq L_0\rho
\]
in common tangent/residual coordinates throughout a joint
parameter/chart tube of radius $\rho$, with
$L_0\rho<\kappa_0^{-1}$.  Suppose every remaining operational
inequality has strictly positive margin on that tube.  Then the
complete OCP has MFCQ and RCRQ there, uniformly with equality
right-inverse bound at most
$(\kappa_0^{-1}-L_0\rho)^{-1}$ before the geometric lift.
The lift preserves these qualifications; its uniformly bounded
local inverse gives a uniform bound for the lifted system as well.
\end{proposition}
\begin{proof}
The QR bound and perturbation inequality preserve full row rank of
$DC_0$ and give the displayed right-inverse bound.  With no active
operational inequalities, MFCQ is equality surjectivity and RCRQ
reduces to this neighborhood-constant equality rank.  The graph
and certificate derivatives in Proposition~\ref{prop:arm-closure-cq}
extend the conclusion, with positivity and inverse bounds on a
compact subdomain of their charts.
\end{proof}
This operating regime uses the original physical limits and the
continuous-clearance tests to certify strict margins.  Shooting and
feature equalities have their recursive identity blocks.  Hard terminal
pose/rate conditions and support equations also enter the joint rank
bound.

With active operational inequalities, the general conditions of
Definitions~\ref{def:theory-mfcq} and~\ref{def:rcrq} apply: a common
strict inward direction and constant rank of $[DC_0;Db_{0,J}]$ for
every active subset $J$, on the common parameter neighborhoods of
Assumption~\ref{ass:regular-path}.  Equation~\eqref{eq:arm-closure-rank-inheritance}
then establishes preservation by the geometric force-closure block.
The argument uses neighborhood rank identities and bounds. Sampled ranks
provide numerical diagnostics of these conditions.

SSOSC is checked separately on the physical reduced objective.
Unobserved directions in terminal force allocations or statically
indeterminate support reactions require the corresponding curvature
check, an optimization-equivalent force-image reduction, or a
regularization of the allocation such as
\eqref{eq:arm-force-regularization}.  The
analytic force-closure calculation concerns the stated grasp geometry;
physical second-order conditions are assessed for the chosen force-allocation
formulation.

\paragraph{Physical force allocation and strict curvature.}
The contact representation separates geometric force closure from uniqueness
of a physical load allocation. Force closure guarantees that a load can be
balanced; it does not make the balancing contact forces unique, and a
non-unique allocation is a direction of zero curvature for the NLP.
As an example, a two-interval acquisition transcription with two load
scenarios contains
$40$ geometry/control coordinates and $48$ physical force coordinates
(two scenarios, eight contacts and three force components per contact). With the geometry fixed, each scenario enters
the dynamics through a six-dimensional wrench. If all $52$ force
guards are strictly inactive and the objective is independent of
allocation along a fixed wrench fiber, at a KKT point there are at
least $48-12=36$ zero-curvature allocation directions. This conclusion
holds for any smooth positive Riemannian metric: the feasible fiber
has constant objective, and at a stationary point its second
variation equals the Lagrangian Hessian on its tangent. Thus the force-closure certificate and the SSOSC of
Section~\ref{sec:pathfollowing-theory} impose separate conditions on this
enlarged physical formulation.

The supported-release problem has an explicit instance of the same
phenomenon. Denote its four first-interval support normals by
$N_1,\ldots,N_4$ in the corner order
\[
 (-a_1,-a_2),\quad(-a_1,+a_2),\quad
 (+a_1,-a_2),\quad(+a_1,+a_2).
\]

\begin{proposition}[First-interval support-allocation fiber]
\label{prop:release-normal-fiber}
Suppose the first-interval box pose is fixed measured data, its measured
twist is zero, and all four normals are strictly between their individual
lower and upper bounds. Suppose the objective and nonlinear equations
depend on these normals through the support wrench in
\eqref{eq:arm-release-support-friction}, with only affine individual
normal bounds in addition. Then the direction
$d_N=(1,-1,-1,1)$ leaves the objective and every nonlinear equation
constant locally. Its Lagrangian second variation is zero for every
finite multiplier vector in the flat force coordinates. At a KKT point,
strict positive curvature on the full supported tangent is impossible
while this fiber is retained.
\end{proposition}
\begin{proof}
Zero twist makes all tangential contact velocities zero. With
$\nu=R_b^\top e_3$, the wrench is
$\sum_i\operatorname{col}(\rho_i^S\times\nu,\nu)N_i$.
The signed sums of the four coefficients and of their corner positions
under $d_N$ both vanish. Hence the wrench, the nonlinear equations and
the objective are unchanged along $N+t d_N$. Strict interior bounds
give a two-sided feasible interval for $t$. The remaining normal-bound
functions are affine and contribute no Hessian in this direction.
\end{proof}

The argument applies to the fixed measured first interval. Future
intervals can have nonzero predicted twists and different contact
Jacobian ranks. A reduced load representation must preserve the finite
feasible image and reconstruct every physical normal and bound.
For example, the linear coordinates
\[
 S=N_1+N_2+N_3+N_4,\quad X=-N_1-N_2+N_3+N_4,
\]
\[
 Y=-N_1+N_2-N_3+N_4,\quad T=N_1-N_2-N_3+N_4
\]
separate wrench coordinates $(S,X,Y)$ from the allocation coordinate
$T$. The original normal bounds give an interval $L(S,X,Y)\leq
T\leq U(S,X,Y)$; its projection is the set $L\leq U$.
All original force limits, multiplier stationarity and the induced
metric remain part of an equivalent reduction. The OCPs of
Section~\ref{sec:arm-example} retain the physical allocation variables and
remove these zero-curvature directions by the regularization
\eqref{eq:arm-force-regularization}. With $\rho_f>0$ the objective is a
strictly convex quadratic function of the forces along every fiber of
constant wrench, so the allocation is unique and the supported-space
curvature test \eqref{eq:fresh-curvature-test} can hold. In the experiments
of Section~\ref{sec:results1019B} it passed at every initialization.

\FloatBarrier
\FloatBarrier
\section{Conclusion}
This paper develops a Riemannian parametric-QP pathfollowing method that
coordinates parallel strict-working-set tasks and includes explicit
factorization, event, synchronization, and status rules. Its local analysis
uses RCRQ, MFCQ, uniform generalized Riemannian SSOSC, and a regular
strict-working-set cover on finitely many path segments. The analysis
relates the frozen-QP correction, retraction and numerical defects to the
tracking error. The arm formulation uses three sequential MPC problems
for grasp acquisition, transport, and placement followed by supported release and retreat. Riemannian
stabilized primal--dual SQP is developed in the appendix for reliable
offline initialization and recovery under the stated operating-domain
assumption. The numerical examples illustrate drone navigation and arm
pick-and-place in closed loop. On the nominal plants, one correction per
measurement reduces the KKT residual of the original nonlinear problem by a
median of four to seven orders of magnitude, the QP homotopies are followed
to rounding accuracy, and the backup is needed only within the last three
samples of a shrinking horizon; the MuJoCo sequences show the associated
physical motion. Three limitations remain: the arm controller does not yet
run in real time, the perturbed drone run relies on the backup or on
replanning at ten of its eighteen samples, and the implemented arm
collision catalogue is a subset of the stated one. The geometric and
force-allocation analyses explain
the role of the representation in the regularity and curvature
assumptions of the method.

\section*{Acknowledgments}
The author would like to thank Johannes J\"aschke for
helpful discussions during the preparation of this manuscript.
The author used OpenAI's ChatGPT and Codex to assist with revision of the
exposition, preparation of figures from the recorded data, and consistency
checks. The author also used Anthropic's Claude (Claude Fable~5.2) to
integrate the closed-loop numerical results, to prepare the corresponding
figures and tables from the recorded data, to revise the exposition, and
to check the manuscript against the implementation and for internal
consistency. The author is responsible for the mathematical content and
interpretation of the results.
\appendix
\section{Riemannian PDSQP2 Initialization and Backup}
\label{app:pdsqp}
This appendix develops Riemannian stabilized primal--dual SQP for offline
initialization and recovery when a pathfollowing task requires a new regular
warm start. It specifies the tangent subproblems, merit function, line search,
parameter updates, and geometric hypotheses supporting the reliable-backup
contract in Assumption~\ref{ass:offline-backup}.
The stabilized subproblem, its primal--dual KKT equations, inner bound
events, merit function and outer parameter updates are given explicitly.
\subsection{Riemannian stabilized primal--dual SQP}
\label{subsec:global-handoff}
Cold initialization and recovery after a synchronized no-survivor exit,
that is, when no task has returned an acceptable endpoint at the
synchronization barrier, use stabilized primal--dual SQP (pdSQP). The method
is called stabilized because its QP subproblem contains a regularization
term in the multipliers. The subproblem therefore remains well posed when
the constraint gradients are linearly dependent, which is the situation in
which the online method is most likely to have stopped. A stabilized SQP
method for degenerate problems on Riemannian manifolds is also developed
in~\cite{yamakawaoka2026stabilized}. The construction here follows the primal--dual
augmented Lagrangian of Gill and Robinson
\cite{gillrobinson2012pdal}, their globally convergent stabilized method
\cite{gillrobinson2013stabilized}, and the global and superlinear analyses of
Gill, Kungurtsev and Robinson
\cite{gillkungurtsevrobinson2017global,gillkungurtsevrobinson2017superlinear}.
We give the geometric problem, stabilized subproblem, globalization and
recovery conditions. The current measurement, horizon and contact mode remain
fixed throughout one backup call. All parallel pathfollowing tasks are joined
or cancelled before this call starts.

\paragraph{Equality--slack problem and signs.}
Collect the nonlinear and simple margins into
$b(r)=\operatorname{col}(D(r),b_{\rm smp}(r))$. Eliminate structural constants,
including the measured initial state, before constructing derivatives. With
$m_C=\dim C$, $m_b=\dim b$, and independent slacks
$s\in\mathbb R^{m_b}$, write
\begin{equation}\label{eq:pdsqp-reformulation}
 \min_{a=(r,s)\in\mathcal N=\mathcal M_H\times\mathbb R^{m_b}}f(r),
 \qquad c(a)=\begin{pmatrix}C(r)\\b(r)-s\end{pmatrix}=0,
 \qquad s\geq0.
\end{equation}
Write $m=m_C+m_b$ for the equality-residual dimension. Unless a return
test explicitly specifies another norm, vector norms are the product
two-norm and operator norms are induced by the stated metric.
The source Lagrangian is $L^-(a,y)=f(r)-y^\top c(a)$, with
\begin{equation}\label{eq:pdsqp-signs}
 y_C=-\lambda,\qquad y_b=-\operatorname{col}(\mu,\nu).
\end{equation}
Thus $y$ is unrestricted during stabilized iterations, while an original KKT
point satisfies $y_b\geq0$ and $s\odot y_b=0$. Initialize
$s_1=[b(r_1)]_+$; subsequent slack steps are independent primal variables.
The product metric is the declared Riemannian metric on $\mathcal M_H$ and
the Euclidean metric on slacks. Let $g=\operatorname{grad}_a f$,
$J=Dc(a)$, and $H_y=\operatorname{Hess}_a L^-(a,y)$.
Adjoints below are metric adjoints; their matrices are transposes in
orthonormal tangent frames.
More explicitly, with $J_C=DC(r)$ and $J_b=Db(r)$,
\begin{equation}\label{eq:pdsqp-slack-derivatives}
 g=\binom{\operatorname{grad}f(r)}{0},\qquad
 J=\begin{pmatrix}J_C&0\\J_b&-I\end{pmatrix},\qquad
 H_y=\begin{pmatrix}
 \operatorname{Hess}_r(f-y_C^\top C-y_b^\top b)&0\\0&0
 \end{pmatrix}.
\end{equation}
Thus the slack coordinates have their own independent steps and bounds;
the geometric Hessian differentiates the complete signed scalar
Lagrangian, including its connection terms.

For a tangent vector $w=(w_r,w_s)$ define the bound KKT residual
\begin{equation}\label{eq:pdsqp-projected-residual}
 \mathcal P_s(w)=\operatorname{col}(w_r,\min(s,w_s)),\qquad
 r_K(a,y)=\|\operatorname{col}(c(a),\mathcal P_s(g-J^*y))\|.
\end{equation}
The minimum is componentwise. In particular, its slack block is
$\min(s,y_b)$, so $r_K=0$ is equivalent to feasibility, stationarity,
nonnegative slack multipliers and complementarity in
\eqref{eq:pdsqp-reformulation}.

\paragraph{Primal--dual augmented Lagrangian.}
Keep the stabilization $\aleph>0$ separate from the line-search penalty
$\rho>0$: the parameter $\aleph$ regularizes the QP subproblem, whereas
$\rho$ weights the constraint terms of the merit function used in the line
search. For a fixed multiplier estimate $y^E$, use
\begin{equation}\label{eq:pdsqp-merit}
 M_\rho(a,y;y^E)=f(r)-(y^E)^\top c(a)
       +\frac\rho2\left(\|c(a)\|^2+
             \|c(a)+\rho^{-1}(y-y^E)\|^2\right).
\end{equation}
The positive primal--dual weighting in the reference merit is set to one.
The reference merit has $\rho=1/\aleph$. With
$\pi=y^E-c(a)/\aleph$,
\begin{equation}\label{eq:pdsqp-gradient}
 G_\aleph=\operatorname{grad}_{(a,y)}M_{1/\aleph}
 =\begin{pmatrix}g-J^*(2\pi-y)\\c(a)+\aleph(y-y^E)\end{pmatrix}.
\end{equation}
All parameters and $y^E$ stay fixed within a line search.

\paragraph{The stabilized SQP subproblem.}
At a current primal--dual point $(a,y)$, let $p_a\in T_a\mathcal N$
be the primal step and let $\widehat y\in\mathbb R^{m_C+m_b}$ be the
new equality-multiplier trial. The canonical stabilized SQP subproblem is
\begin{equation}\label{eq:pdsqp-canonical-stabilized}
 \begin{aligned}
 \min_{p_a,\widehat y}\quad&
  \langle g,p_a\rangle+\tfrac12\langle p_a,H_yp_a\rangle
                        +\tfrac{\aleph}{2}\|\widehat y\|^2,\\
 \text{subject to}\quad&
  c(a)+Jp_a+\aleph(\widehat y-y^E)=0,\\
 &s+(p_a)_s\geq0.
 \end{aligned}
\end{equation}
For a V or O iteration, $y^E=y$, so its affine equation contains the
stabilized multiplier increment $\aleph(\widehat y-y)$. At
$\aleph=0$, its primal equations reduce formally to the ordinary SQP
linearization. For $\aleph>0$, every primal step satisfying the slack
bounds has a unique affine-feasible multiplier trial. In particular,
\begin{equation}\label{eq:pdsqp-known-feasible-point}
 p_a=0,\qquad \widehat y=y^E-\aleph^{-1}c(a)
\end{equation}
is a known feasible point, even when the current nonlinear equalities
are infeasible or $J$ has dependent rows.

The unmodified Hessian in
\eqref{eq:pdsqp-canonical-stabilized} can yield an indefinite reduced
quadratic away from a solution. The global direction uses the following
convexified subproblem; the relaxed local direction below retains $H_y$
when its stated tests hold. Throughout this offline method the equality
multiplier is an iterated primal--dual variable: its trial value is
determined together with the primal step by the stabilized equations.

\paragraph{Convexification and the global subproblem.}
Choose a symmetric modification $\widehat H_y$ satisfying
\begin{equation}\label{eq:pdsqp-convexification}
 \|\widehat H_y\|\leq M_H,\qquad
 \widehat H_y+\aleph^{-1}J^*J\succeq\epsilon_Q I,
 \qquad M_H<\infty,\quad\epsilon_Q>0.
\end{equation}
A uniformly bounded positive-definite choice is
$\widehat H_y=\epsilon_Q I$. When $H_y$ is uniformly bounded and a
uniformly bounded certified number $\ell_y$ satisfies
$\ell_y\leq\lambda_{\min}(H_y)$, the closer modification
\begin{equation}\label{eq:pdsqp-explicit-convexification}
 \widehat H_y=H_y+\delta_y I,\qquad
 \delta_y=\max\{0,\epsilon_Q-\ell_y\}
\end{equation}
also satisfies \eqref{eq:pdsqp-convexification}. Structured modifications satisfying
\eqref{eq:pdsqp-convexification} can reduce the change on the relevant
free subspace. These modifications affect the QP model; the objective,
constraint values, signed derivatives and return tests remain those of
the original OCP.
The global descent step is obtained from the strictly convex reduced problem
associated with
\begin{equation}\label{eq:pdsqp-bound-qp}
 \begin{aligned}
 \min_{p_a,\widehat y}\quad&
   \langle g,p_a\rangle+\tfrac12\langle p_a,\widehat H_yp_a\rangle
                   +\tfrac{\aleph}{2}\|\widehat y\|^2,\\
 \text{subject to}\quad&c(a)+Jp_a+\aleph(\widehat y-y^E)=0,
                  \qquad (p_a)_s\geq-s.
 \end{aligned}
\end{equation}
Eliminating its affine equation gives
\begin{equation}\label{eq:pdsqp-reduced-qp}
 \begin{aligned}
 \min_{p_a:\,s+(p_a)_s\geq0}\quad&
 q_{\rm red}(p_a)
 =\langle g,p_a\rangle+\tfrac12\langle p_a,\widehat H_yp_a\rangle
   +\tfrac{\aleph}{2}
        \|y^E-\aleph^{-1}(c+Jp_a)\|^2,\\
 \widehat y&=y^E-\aleph^{-1}(c+Jp_a),
 \end{aligned}
\end{equation}
and therefore
\begin{equation}\label{eq:pdsqp-reduced-gradient}
 \nabla q_{\rm red}(p_a)
 =g-J^*y^E+\aleph^{-1}J^*c+
       (\widehat H_y+\aleph^{-1}J^*J)p_a.
\end{equation}
Thus the reduced QP is always feasible, and
\eqref{eq:pdsqp-convexification} gives a unique solution. Its complete step
is $d=(p_a,\widehat y-y)$. To state its KKT equations, let
$E_s z=\operatorname{col}(0,z)$ inject slack-bound multipliers into
$T_a\mathcal N$. The subproblem Lagrangian uses minus signs for both
its affine equality multiplier $w$ and its nonnegative-bound multiplier
$z$. Differentiation in $\widehat y$ gives $w=\widehat y$, and the
complete KKT conditions are
\begin{equation}\label{eq:pdsqp-full-subproblem-kkt}
 \begin{aligned}
 g+\widehat H_yp_a-J^*\widehat y-E_s z&=0,\\
 c+Jp_a+\aleph(\widehat y-y^E)&=0,\\
 s+(p_a)_s\geq0,\qquad z\geq0,\qquad
 z\odot(s+(p_a)_s)&=0.
 \end{aligned}
\end{equation}
The stabilized QP returns the equality multiplier appearing in these KKT
conditions.

\paragraph{Regularized KKT solves and slack-bound events.}
Let $A$ be the working slack-bound indices, fixed at
$(p_a)_{s,A}=-s_A$, and let $F$ contain every geometric primal
coordinate and all remaining slack coordinates. At an affine-feasible
QP iterate, an inner correction obeys
\begin{equation}\label{eq:pdsqp-stabilized-kkt}
 \begin{pmatrix}(\widehat H_y)_{FF}&J_F^*\\J_F&-\aleph I\end{pmatrix}
 \begin{pmatrix}\delta p_F\\-\delta\widehat y\end{pmatrix}
 =\begin{pmatrix}-\bar g_F\\0\end{pmatrix},\qquad
 \bar g=g+\widehat H_yp_a-J^*\widehat y.
\end{equation}
Eliminated slack-bound multipliers are retained as reduced costs. Inner
termination checks the affine equation, slack bounds and projected reduced
stationarity. The regularized dual block is $-\aleph I$. Its Schur
complement in \eqref{eq:pdsqp-stabilized-kkt} is
$(\widehat H_y)_{FF}+\aleph^{-1}J_F^*J_F\succ0$, so the symmetric
matrix has inertia $(|F|,m_C+m_b,0)$ regardless of the row rank of
$J_F$. It can be solved with the signed $LBL^\top$ operations of
Section~\ref{sec:pathfollowing}.

For completeness, a primal bound-active-set solution of
\eqref{eq:pdsqp-reduced-qp} starts at
\eqref{eq:pdsqp-known-feasible-point}, or at a feasible local-step
candidate. Solve \eqref{eq:pdsqp-stabilized-kkt} with zero correction
on $A$. If the resulting primal correction is nonzero, take
\begin{equation}\label{eq:pdsqp-inner-bound-ratio}
 \alpha_Q=\min\left\{1,
   \min_{\substack{j\notin A\\(\delta p_s)_j<0}}
       \frac{s_j+(p_a)_{s,j}}{-(\delta p_s)_j}\right\},\qquad
 (p_a,\widehat y)\gets
       (p_a+\alpha_Q\delta p,\widehat y+\alpha_Q\delta\widehat y).
\end{equation}
An empty inner minimum is $+\infty$. Add the blocking slack bound
when $\alpha_Q<1$; an exact tied bound at $\alpha_Q=1$ is also
eligible to enter the working set. If the primal correction is zero,
the free reduced costs vanish; the working costs are
$z_A=(\bar g_s)_A$. Nonnegative working costs complete the KKT test.
Otherwise release a negative-cost bound and solve again. A consistent
lexicographic anti-cycling rule on the fixed slack IDs resolves
degenerate ties. The inner problem is a bound-constrained QP, and the
outer nonlinear iterates are updated by the manifold retraction.

The numerical stopping quantities can be written as
\begin{equation}\label{eq:pdsqp-inner-residuals}
 \begin{aligned}
 r_{Q,\mathrm{aff}}&=\|c+Jp_a+\aleph(\widehat y-y^E)\|_\infty,\\
 r_{Q,\mathrm{bnd}}&=\|[-s-(p_a)_s]_+\|_\infty,\\
 r_{Q,\mathrm{stat}}&=
   \max\{\|\bar g_r\|_\infty,
     \|\min(s+(p_a)_s,\bar g_s)\|_\infty\}.
 \end{aligned}
\end{equation}
These are measured against the complete assembled subproblem, with
eliminated costs expanded before testing.

\paragraph{Equivalent primal--dual quadratic model and descent.}
An alternative realization minimizes a quadratic in the complete step
$d=(p_a,p_y)$:
\begin{equation}\label{eq:pdsqp-complete-step-qp}
 \min_{d:\,s+(p_a)_s\geq0}
   \langle G_\aleph,d\rangle
       +\tfrac12\langle d,\widehat B_\aleph d\rangle,
 \qquad p_y=\widehat y-y,
\end{equation}
where
\begin{equation}\label{eq:pdsqp-model-hessian}
 \widehat B_\aleph=
 \begin{pmatrix}
  \widehat H_y+2\aleph^{-1}J^*J&J^*\\J&\aleph I
 \end{pmatrix}.
\end{equation}
Its $p_y$ equation is
$c+\aleph(y-y^E)+Jp_a+\aleph p_y=0$. Substituting it into the
primal equation recovers
\eqref{eq:pdsqp-full-subproblem-kkt}, proving equivalence to the
global stabilized SQP subproblem.
Its Schur complement is precisely
$\widehat H_y+\aleph^{-1}J^*J$. Hence
$\widehat B_\aleph\succ0$, and the exact QP minimizer satisfies
$\langle G_\aleph,d\rangle\leq
 -\langle d,\widehat B_\aleph d\rangle<0$ unless $d=0$.

\paragraph{Active-set estimate and local descent step.}
The superlinear analysis uses the following local step in addition to the
modified global QP. Let
$0<\gamma<1$, $0<\epsilon_a<1$, and
\begin{equation}\label{eq:pdsqp-epsilon-active}
 \begin{aligned}
 \epsilon(a,y;\delta)&=\min\{\epsilon_a,
           \max(\delta,r_K(a,y))^\gamma\},\\
 A_\epsilon(a,y;\delta)&=\{j:s_j\leq\epsilon(a,y;\delta)\}.
 \end{aligned}
\end{equation}
For the local step set $A_\epsilon=A_\epsilon(a_k,y_k;\aleph_k)$,
after the current outer parameters have been updated.
All geometric primal coordinates and all dual coordinates are free; only
slack coordinates in $A_\epsilon$ are fixed. Denote the unmodified matrix
\eqref{eq:pdsqp-model-hessian}, with $\widehat H_y$ replaced by $H_y$, by
$B_\aleph$. On a V or O iteration defined below, if its free principal
submatrix is positive definite, form $d^0$ by
$(d^0)_{s,A_\epsilon}=-(s_k)_{A_\epsilon}$ and zero elsewhere, and solve
\begin{equation}\label{eq:pdsqp-local-step}
 (\delta d)_{A_\epsilon}=0,\qquad
 (B_\aleph)_{FF}(\delta d)_F
      =-[G_\aleph+B_\aleph d^0]_F,\qquad d^{\rm loc}=d^0+\delta d.
\end{equation}
With an exponent $0<\bar\lambda<\min\{\gamma,1-\gamma\}$, accept this candidate when
\begin{equation}\label{eq:pdsqp-local-step-tests}
 s_k+(d^{\rm loc})_s\geq0,\qquad
 [G_\aleph+B_\aleph d^{\rm loc}]_{s,A_\epsilon}
       \geq-r_K(a_k,y_k)^{\bar\lambda}\mathbf1,\qquad
 \langle G_\aleph,d^{\rm loc}\rangle<0.
\end{equation}
Otherwise use \eqref{eq:pdsqp-bound-qp}. These are the relaxed active-bound
optimality tests in the descent-direction construction of
\cite[Section~2.3]{gillkungurtsevrobinson2017global} and
\cite[Section~2.2]{gillkungurtsevrobinson2017superlinear}.

\paragraph{Curvature safeguard and merit pullbacks.}
For the second-order variant, compute the least-curvature estimate before
the outer update, using the previous regularization
$\aleph=\aleph_{k-1}$. Let $P_\epsilon$ inject the primal coordinates
outside $A_\epsilon(a_k,y_k;\aleph_{k-1})$ into $T_{a_k}\mathcal N$;
all geometric coordinates are included. In orthonormal frames form
\begin{equation}\label{eq:pdsqp-curvature-matrix}
 S_{k,-}=P_\epsilon^*(H_y+\aleph^{-1}J^*J)P_\epsilon.
\end{equation}
If $S_{k,-}\succeq0$, set $\xi_k=0$ and $q_k=0$. Otherwise choose a
nonzero $v$ with
\[
 \frac{v^\top S_{k,-}v}{\|v\|^2}
       \leq\theta\lambda_{\min}(S_{k,-})<0,
 \qquad 0<\theta\leq1,
\]
where an exact least eigenvector is an admissible choice with $\theta=1$.
Set $u=P_\epsilon v$ and
\begin{equation}\label{eq:pdsqp-curvature-pair}
 q^-=(u,-Ju/\aleph),\qquad
 \xi_k=-\frac{\langle u,(H_y+\aleph^{-1}J^*J)u\rangle}{\|u\|^2}>0.
\end{equation}
Then $\langle q^-,B_\aleph q^-\rangle=-\xi_k\|u\|^2$ and the
$\epsilon$-active slack components of $u$ vanish. These formulas give the
least-curvature construction of
\cite[Section~2.1, Algorithm~1]{gillkungurtsevrobinson2017global}.

After updating the parameters and computing $d_k$, if
$\langle G_{\aleph_k},q^-\rangle>0$, replace
$(u,q^-)\gets(-u,-q^-)$. Thus the curvature pair has nonpositive
reference-merit slope. For $\xi_k>0$, define
\begin{equation}\label{eq:pdsqp-curvature-scale}
 \begin{aligned}
 \sigma_k=\min\left\{
   \frac{\max\{\xi_k,\|(d_k)_a\|\}}{\|u\|},\quad
   \min_{j:(u_s)_j<0}
      \frac{(s_k)_j+((d_k)_s)_j}{-(u_s)_j}\right\},\\
 q_k=\sigma_k q^- .
 \end{aligned}
\end{equation}
The inner minimum is $+\infty$ when no slack component decreases.
The primal descent step is slack-feasible, so $\sigma_k\geq0$.
Equivalently, $\sigma_k$ is the largest nonnegative scale for which
$s_k+(d_k)_s+\sigma_k u_s\geq0$ and
$\|\sigma_k u\|\leq\max\{\xi_k,\|(d_k)_a\|\}$.
Set $\Delta_k=d_k+q_k$. This retains the feasible curvature scaling of
\cite[Section~2.3, Algorithm~3]{gillkungurtsevrobinson2017global}.

Use the product retraction
\begin{equation}\label{eq:part1-riemannian-backup-retraction}
 \mathscr R_{(r,s,y)}(\alpha\Delta)=
 \bigl(R_r(\alpha\Delta_r),s+\alpha\Delta_s,
                         y+\alpha\Delta_y\bigr).
\end{equation}
For every smooth scalar merit $M$,
\begin{equation}\label{eq:merit-retraction-hessian}
 D^2(M\circ\mathscr R)(0)[\Delta,\Delta]
 =\langle\Delta,\operatorname{Hess}M[\Delta]\rangle
  +\langle\operatorname{grad}_rM,
                   \mathcal A_r(\Delta_r,\Delta_r)\rangle.
\end{equation}
Here $\mathcal A_r(\eta,\eta)=
\left.\nabla_{\dot\gamma}\dot\gamma\right|_{t=0}$ for
$\gamma(t)=R_r(t\eta)$ is the retraction acceleration.
A second-order retraction has $\mathcal A_r=0$. In particular, true reference
merit curvature uses $H_{2\pi-y}$ in its upper-left block, whereas the SQP
model uses $H_y$. The two are evaluated separately in the exceptional
zero-slope curvature test.

\paragraph{Flexible line search and its independent penalty.}
Write
$\Psi_k(\alpha;\rho)=
M_\rho(\mathscr R_{v_k}(\alpha\Delta_k);y_k^E)$,
where $v_k=(a_k,y_k)$. Let $\ell_k=1$ if $d_k\ne0$, $q_k=0$ and the
iteration is V or O, and let $\ell_k=2$ otherwise. With
$B_{k,-}=B_{\aleph_{k-1}}(a_k,y_k)$ define the reference predicted decrease
\begin{equation}\label{eq:pdsqp-line-model}
 \mathcal D_k(\alpha)=
 -\alpha\langle G_{\aleph_k},\Delta_k\rangle
 -\tfrac12(\ell_k-1)\alpha^2
             \min\{0,\langle\Delta_k,B_{k,-}\Delta_k\rangle\}.
\end{equation}
First set $\rho_k\gets\min\{\rho_k,1/\aleph_k\}$. For
$0<\eta_S<1/2$, test $\alpha=1,1/2,1/4,\ldots$ and accept the first valid
retraction-domain trial satisfying either inequality
\begin{equation}\label{eq:pdsqp-flexible-search}
 \Psi_k(0;\rho)-\Psi_k(\alpha;\rho)
       \geq\eta_S\mathcal D_k(\alpha),\qquad
 \rho\in\{\rho_k,1/\aleph_k\}.
\end{equation}
Both tests compare actual decrease with the same reference predicted decrease.
The accepted primal slack segment is feasible by its construction.

The exceptional cases follow the flexible-line-search procedure of
\cite[Section~2.4, Algorithm~4]{gillkungurtsevrobinson2017global}.
If $d_k=q_k=0$, keep the point. If $d_k=0$, $q_k\ne0$, its reference slope
is zero and $\aleph_k\ne\aleph_{k-1}$, first require
\[
 -\frac{D^2(M_{1/\aleph_k}\circ\mathscr R)(0)[q_k,q_k]}
              {\|(q_k)_a\|^2}>\eta_S\xi_k.
\]
If this fails, take $\alpha_k=0$ and continue the outer parameter updates
at the same point. Otherwise perform \eqref{eq:pdsqp-flexible-search}.
After acceptance, the penalty update is
\begin{equation}\label{eq:pdsqp-penalty-update}
 \rho_{k+1}=\begin{cases}
 \rho_k,&\Psi_k(0;\rho_k)-\Psi_k(\alpha_k;\rho_k)
                       \geq\eta_S\mathcal D_k(\alpha_k),
                       \text{ or }\alpha_k=0,\text{ or }\Delta_k=0,\\
 \min\{2\rho_k,1/\aleph_k\},&\text{otherwise}.
 \end{cases}
\end{equation}
This is the reciprocal-penalty form of the source update.

\begin{proposition}[Finite local backtracking]\label{prop:pdsqp-backtracking}
At fixed outer parameters, assume the chosen retraction segment is defined
locally. If $\langle G_{\aleph_k},\Delta_k\rangle<0$, sufficiently small
trials satisfy the reference test in \eqref{eq:pdsqp-flexible-search}.
The same conclusion holds for a pure zero-slope curvature step passing the
preceding true-curvature test.
\end{proposition}
\begin{proof}
The first derivative of the reference pullback at zero is
$\langle G_{\aleph_k},\Delta_k\rangle$, because $DR_r(0)$ is the identity.
For negative slope its linear term dominates the quadratic term in
\eqref{eq:pdsqp-line-model}. For the pure curvature case, the tested strict
inequality compares its actual negative quadratic Taylor coefficient with
$\eta_S$ times the model coefficient. Taylor expansion therefore proves the
reference decrease for all sufficiently small positive steps. Slack feasibility
holds along the entire segment.
\end{proof}

\paragraph{Complete V/O/M/F update policy.}
At iteration $k$, first compute $\mathfrak r_k=r_K(a_k,y_k)$ and the curvature estimate
$\xi_k$ using $\aleph_{k-1}$. For a fixed $0<\beta<1$, let
\begin{equation}\label{eq:pdsqp-vo-measures}
 \begin{aligned}
 \eta_k&=\|c(a_k)\|,\qquad
 \omega_k=\max\{\|\mathcal P_{s_k}(g_k-J_k^*y_k)\|,\xi_k\},\\
 \phi_{V,k}&=\eta_k+\beta\omega_k,\qquad
 \phi_{O,k}=\beta\eta_k+\omega_k,\\
 e_k&=\max\{\mathfrak r_k,\xi_k\}.
 \end{aligned}
\end{equation}
Start at $k=1$ with the supplied seed $r_1$,
$a_1=(r_1,[b(r_1)]_+)$ and the converted multiplier seed $y_1$.
Choose positive bounds $\Phi_{V,0},\Phi_{O,0}$, tolerance $\tau_0$,
$\aleph_0>0$, and a finite clipping bound $y_{\max}$. Set
$y_0^E=\operatorname{clip}(y_1,-y_{\max},y_{\max})$ and choose
$0<\rho_1\leq1/\aleph_0$. Thus the first curvature and M tests use the
initialized previous outer state $(\aleph_0,y_0^E,\tau_0)$, whereas the
first line search uses the subsequently updated state at $k=1$. Apply the
following ordered rules, retaining quantities not explicitly changed:
\begin{enumerate}
\item If $\phi_{V,k}\leq\Phi_{V,k-1}/2$, declare V,
set $\Phi_{V,k}=\Phi_{V,k-1}/2$ and $y_k^E=y_k$.
\item Otherwise, if $\phi_{O,k}\leq\Phi_{O,k-1}/2$, declare O,
set $\Phi_{O,k}=\Phi_{O,k-1}/2$ and $y_k^E=y_k$.
For either V or O use
\begin{equation}\label{eq:pdsqp-vo-update}
 \aleph_k=\begin{cases}
       \min\{\aleph_0,e_k^\gamma\},&e_k>0,\\
       \aleph_{k-1}/2,&e_k=0.
       \end{cases}
\end{equation}
The V/O tolerance convention is $\tau_k=\tau_{k-1}/2$, as in the parameter
definition in \cite[Section~2.1]{gillkungurtsevrobinson2017superlinear}.
\item Otherwise evaluate $G_{\aleph_{k-1}}$ with $y^E=y^E_{k-1}$.
Declare M if
\begin{equation}\label{eq:pdsqp-m-test}
 \|\mathcal P_{s_k}((G_{\aleph_{k-1}})_a)\|\leq\tau_{k-1},\quad
 \|(G_{\aleph_{k-1}})_y\|\leq\tau_{k-1}\aleph_{k-1},\quad
 \xi_k\leq\tau_{k-1}.
\end{equation}
Then set
\begin{equation}\label{eq:pdsqp-m-update}
 \begin{aligned}
 y_k^E&=\operatorname{clip}(y_k,-y_{\max},y_{\max}),\qquad
 \tau_k=\tau_{k-1}/2,\\
 \aleph_k&=\begin{cases}
       \min\{\aleph_{k-1}/2,e_k^\gamma\},&e_k>0,\\
       \aleph_{k-1}/2,&e_k=0.
       \end{cases}
 \end{aligned}
\end{equation}
\item If none of these tests holds, declare F and retain
$y_k^E=y_{k-1}^E$, $\aleph_k=\aleph_{k-1}$ and $\tau_k=\tau_{k-1}$.
\end{enumerate}
These rules retain the source distinction between improvement in the original
problem (V and O iterations, which reduce a measure weighted toward
constraint violation and toward optimality, respectively), approximate
minimization of a fixed augmented Lagrangian (M iterations), and
continued minimization with fixed outer parameters (F iterations). In particular, the V/O
regularization need not decrease monotonically. The clipping operation
changes the auxiliary estimate $y^E$ on initialization or an M iteration;
the iterated multiplier $y$ is updated by the primal--dual step.
Numerical stabilization floors define the attainable accuracy of these updates.

\paragraph{Return tests and recovery interface.}
Using the original signed multipliers
$\zeta_b=\operatorname{col}(\mu,\nu)$ of Section~\ref{sec:parametric-problem}, reevaluate
\begin{equation}\label{eq:original-kkt-return}
 \begin{aligned}
 \chi_f&=\max\{\|C\|_\infty,\|b-s\|_\infty,
                                      \|[-b]_+\|_\infty\},\\
 \chi_s&=\|\operatorname{grad}f+DC^*\lambda+Db^*\zeta_b\|_\infty,\\
 \chi_c&=\max\{\|s\odot\zeta_b\|_\infty,\|[\zeta_b]_+\|_\infty\}.
 \end{aligned}
\end{equation}
Norms refer to declared scaled orthonormal frames. Return a recovered OCP
point only when these residuals meet the requested tolerances and the fresh
SSOSC/working-system tests required for pathfollowing reentry pass. Slack
nonnegativity is checked separately in floating-point arithmetic. To
relate the equality--slack complementarity test directly to the original
margin, use
\begin{equation}\label{eq:pdsqp-original-complementarity}
 \|b\odot\zeta_b\|_\infty
 \leq\|s\odot\zeta_b\|_\infty
      +\|\zeta_b\|_\infty\|b-s\|_\infty.
\end{equation}
The original-margin product is reevaluated, or this upper bound is
required to meet the original complementarity tolerance.

The reference second-order termination test is
\begin{equation}\label{eq:pdsqp-source-stopping}
 \mathfrak r_k\leq\tau_{\rm stop},\qquad
 \xi_k\leq\tau_{\rm stop},\qquad
 \aleph_{k-1}\leq\tau_{\rm stop}.
\end{equation}
It records approximate second-order stationarity of the equality--slack
problem. Online reentry additionally requires the supported-space SSOSC
and regular working representation specified in
Section~\ref{sec:pathfollowing-theory}.
A candidate
stationary point of constraint violation is detected by
\[
 \min\{\|c(a_k)\|,\varepsilon_{\rm ISP}\}>\aleph_k,\qquad
 \|\mathcal P_{s_k}(J_k^*c(a_k))\|\leq\varepsilon_{\rm ISP},
 \qquad k\text{ is an M iteration},
\]
as in \cite[Section~2.5]{gillkungurtsevrobinson2017global}.
This test triggers restoration or a new initialization seed. Reseeding
also applies when the original residual tests pass, the SSOSC reentry gate
fails, and $\alpha_k\Delta_k=0$. If a numerical line search cannot
produce a valid descent trial at its declared arithmetic accuracy, control
returns to the same restoration/reseeding mechanism. The backup
continues until its recovery contract in Assumption~\ref{ass:offline-backup}
is satisfied. A completed backup return therefore includes the current
primal state, independent slacks, signed multipliers, residual frames,
working representation and reentry certificates. An infeasibility-stationary
or arithmetic-limited attempt is retained with its stopping reason and
passed to restoration. Assumption~\ref{ass:offline-backup} requires
successful completion of this restoration process on the stated operating
domain, together with the numerical accuracy and preparation interval
needed for the controller handoff.

\begin{algorithm}[htbp]
\caption{Riemannian pdSQP initialization and offline recovery}\label{alg:pdsqp2}
\begin{algorithmic}[1]
\Require Fixed current OCP, manifold seed, signed multiplier seed, product metric and retraction, requested residual and reentry tolerances, and parameters satisfying the bounds above.
\State Form \eqref{eq:pdsqp-reformulation}; initialize $k=1$, $a_1,y_1,y_0^E,\aleph_0,\rho_1,\tau_0,\Phi_{V,0},\Phi_{O,0}$ as specified above.
\While{the recovery tests \eqref{eq:original-kkt-return}--\eqref{eq:pdsqp-original-complementarity}, slack nonnegativity and the SSOSC reentry gate have not passed}
 \State Evaluate $g,J,H_y$ by \eqref{eq:pdsqp-slack-derivatives}, $\mathfrak r_k$ by \eqref{eq:pdsqp-projected-residual}, and the curvature pair by \eqref{eq:pdsqp-curvature-matrix}--\eqref{eq:pdsqp-curvature-pair}.
 \State Apply \eqref{eq:pdsqp-vo-measures}--\eqref{eq:pdsqp-m-update}; set $\rho_k\gets\min\{\rho_k,1/\aleph_k\}$ and freeze $y_k^E,\aleph_k,\rho_k$.
 \State Assemble $G_{\aleph_k}$; on a V/O iteration test the unmodified local step \eqref{eq:pdsqp-local-step}--\eqref{eq:pdsqp-local-step-tests} using $A_\epsilon(a_k,y_k;\aleph_k)$, and set $d_k=d^{\rm loc}$ when its tests pass.
 \State If that step is unavailable or fails a test, solve the modified stabilized SQP subproblem \eqref{eq:pdsqp-bound-qp} with \eqref{eq:pdsqp-convexification}, using \eqref{eq:pdsqp-stabilized-kkt}--\eqref{eq:pdsqp-inner-residuals}; set $d_k=(p_a,\widehat y-y_k)$.
 \State Orient and scale the curvature pair by \eqref{eq:pdsqp-curvature-scale}; set $q_k=0$ if $\xi_k=0$, and form $\Delta_k=d_k+q_k$.
 \State Compute $\alpha_k$ by \eqref{eq:pdsqp-line-model}--\eqref{eq:pdsqp-flexible-search}, including the exceptional true-curvature test; unresolved numerical search invokes restoration/reseeding under Assumption~\ref{ass:offline-backup}.
 \State Set $v_{k+1}=\mathscr R_{v_k}(\alpha_k\Delta_k)$ and update $\rho_{k+1}$ by \eqref{eq:pdsqp-penalty-update}.
 \State If infeasibility stationarity persists, or $\alpha_k\Delta_k=0$ at a residual-qualified point failing SSOSC, restore or reseed $v_{k+1}$ under Assumption~\ref{ass:offline-backup}.
 \State $k\gets k+1$.
\EndWhile
\State Return the reevaluated primal--dual solution, original-margin complementarity check \eqref{eq:pdsqp-original-complementarity}, and residual and curvature certificates.
\end{algorithmic}
\end{algorithm}

\paragraph{Geometric assumptions and transfer of the reference analysis.}
For the intrinsic backup, assume the following uniformly on all current OCPs
in the declared recovery domain: primal iterates lie in a compact subset of
$\mathcal N$; the retraction has a common domain radius and bounded first three
derivatives; the metric and its inverse have uniform chart bounds; and the
objective, constraints and their first two covariant derivatives have the
common smoothness bounds specified in Section~\ref{sec:pathfollowing-theory}. Use a
second-order retraction and smooth isometric transport when tangent quantities
from different iterates are compared. Require
\eqref{eq:pdsqp-convexification}, the stated source parameter and search rules,
and accurately solved inner systems. In the exact analysis no positive
regularization floor or work truncation is imposed. The convergence transfer below concerns exact derivatives and the specified
subproblem solutions. An inexact local implementation may use a computed step
$\widetilde d_k$ provided its error after the linear solve obeys
$\|\widetilde d_k-d_k\|=o(e_k^{1+\gamma})$ and the reference-merit decrease
is verified with freshly evaluated function values. The original residuals
in \eqref{eq:original-kkt-return} determine the delivered accuracy.

The geometric replacement is explicit. For $h=f$ or a component of $c$, a
second-order retraction gives
\[
 D(h\circ R_r)(0)[p]=Dh(r)[p],\qquad
 D^2(h\circ R_r)(0)[p,p]
       =\langle p,\operatorname{Hess}h(r)[p]\rangle.
\]
Consequently the source quadratic models, bound KKT equations and merit
slopes are exactly the tangent pullback models above. Uniform Lipschitz
second derivatives give the same Taylor remainder estimates, with constants
independent of the base point and current homotopy parameter. In orthonormal
frames the QP, inertia and multiplier-update arguments are the Euclidean
arguments; the compactness and line-search arguments use these uniform
pullback estimates. These observations transfer the reference global
convergence alternatives to the intrinsic formulation under the stated
hypotheses. Feasibility of the limit and a constraint qualification are the
additional conditions selecting the KKT alternative from an
infeasibility-stationary alternative.

For the local superlinear conclusion, also require a V/O subsequence
approaching a KKT pair, the compact multiplier-neighborhood SOSC assumption
of \cite[Section~3]{gillkungurtsevrobinson2017superlinear}, and its relaxed
local QP step. The uniform SSOSC imposed in Section~\ref{sec:pathfollowing-theory} supplies
the required stronger curvature on the selected solution branch. In a smooth
local chart, the additional terms from the second-order retraction are of
order $O(\|p\|^3)$ in function values and $O(\|p\|^2)$ in first-order
residuals. For $0<\gamma<1$, these terms are higher order than the stabilized
local error $O(e_k^{1+\gamma})$ when $\|p_k\|=O(e_k)$; the merit Taylor
remainder is $O((1+\aleph_k^{-1})\|p_k\|^3)$, as in the reference analysis.
Thus the source unit-step and superlinear argument applies with the same
local assumptions and controlled inexactness. Assumption~\ref{ass:offline-backup}
specifies the reliable recovery domain, restoration mechanism and returned
accuracy needed by the online controller.

\FloatBarrier
\section{Supplementary Numerical Figures}
\label{app:supplementary-figures}

This appendix collects the figures of the numerical study that are not
shown in Section~\ref{sec:numerical-results}. All of them are drawn from
the same recorded runs, and they use the conventions of
Section~\ref{subsec:results-measures}.

\begin{figure}[H]
 \centering
 \includegraphics[width=\linewidth]{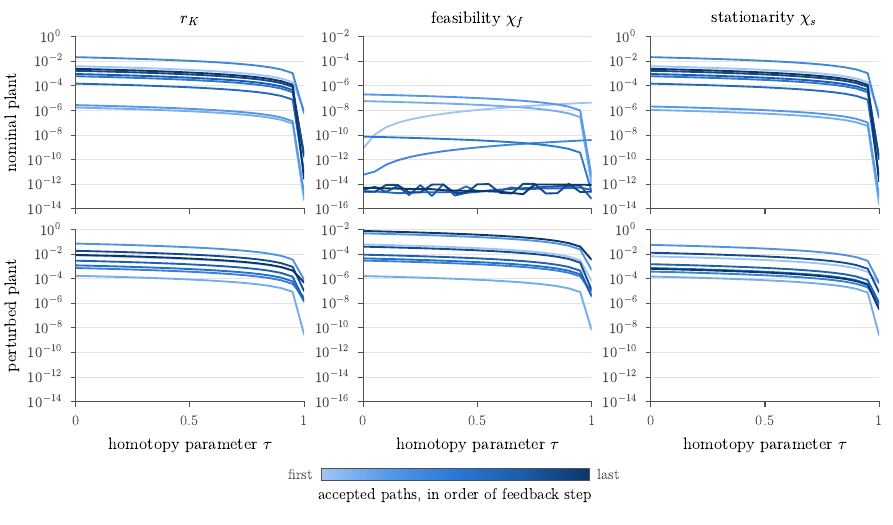}
 \caption{Drone: residuals of the original nonlinear problem against the
 homotopy parameter for every accepted path. On a logarithmic scale the
 linear decrease $r_K(\tau)\approx(1-\tau)\,r_K(0)$ appears as a curve
 that bends down towards $\tau=1$; the final drop is the passage from the
 interpolated data to the target QP. In the nominal run the feasibility
 residual of several paths starts at rounding level and rises to the
 linearization error.}
 \label{fig:drone-tau}
\end{figure}

\begin{figure}[H]
 \centering
 \includegraphics[width=\linewidth]{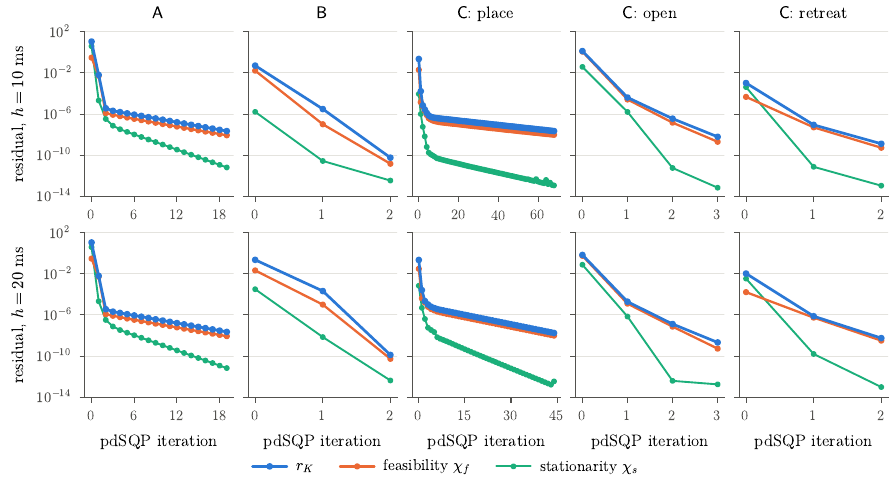}
 \caption{Arm: convergence of the pdSQP initialization of each phase for
 the full-horizon run (top) and the reduced run (bottom). Iteration $0$ is
 the seed.}
 \label{fig:arm-pdsqp}
\end{figure}

\begin{figure}[H]
 \centering
 \includegraphics[width=\linewidth]{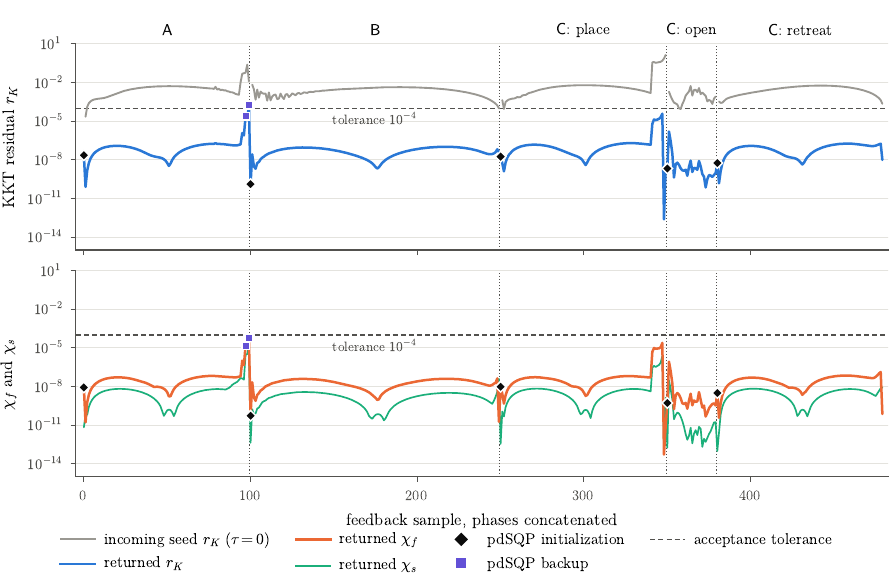}
 \caption{Arm, reduced run ($h=20$~ms, $H=100$, $150$, $100+30+100$):
 residuals of the original nonlinear problem at every feedback sample,
 phases concatenated. The layout is that of
 Figure~\ref{fig:arm-full-closed-loop}. The two backup samples are at the
 end of $\mathsf A$.}
 \label{fig:arm-reduced-closed-loop}
\end{figure}

\begin{figure}[H]
 \centering
 \includegraphics[width=\linewidth]{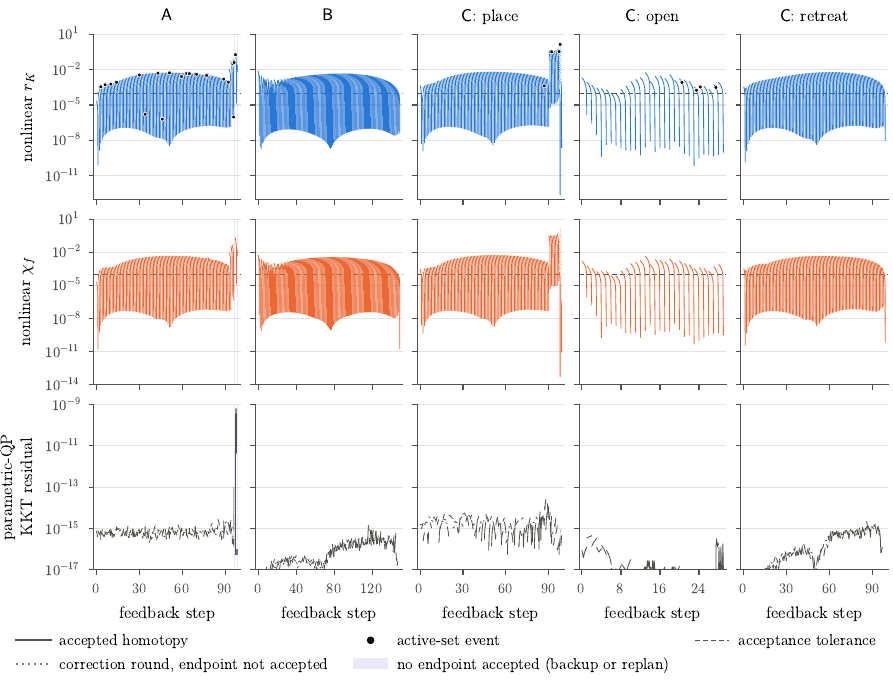}
 \caption{Arm, reduced run: residuals along the homotopies of all $475$
 feedback steps. Each curve occupies the unit interval of its step. The
 layout is that of Figure~\ref{fig:arm-full-unrolled}; black dots are
 active-set events.}
 \label{fig:arm-reduced-unrolled}
\end{figure}

\begin{figure}[H]
 \centering
 \includegraphics[width=\linewidth]{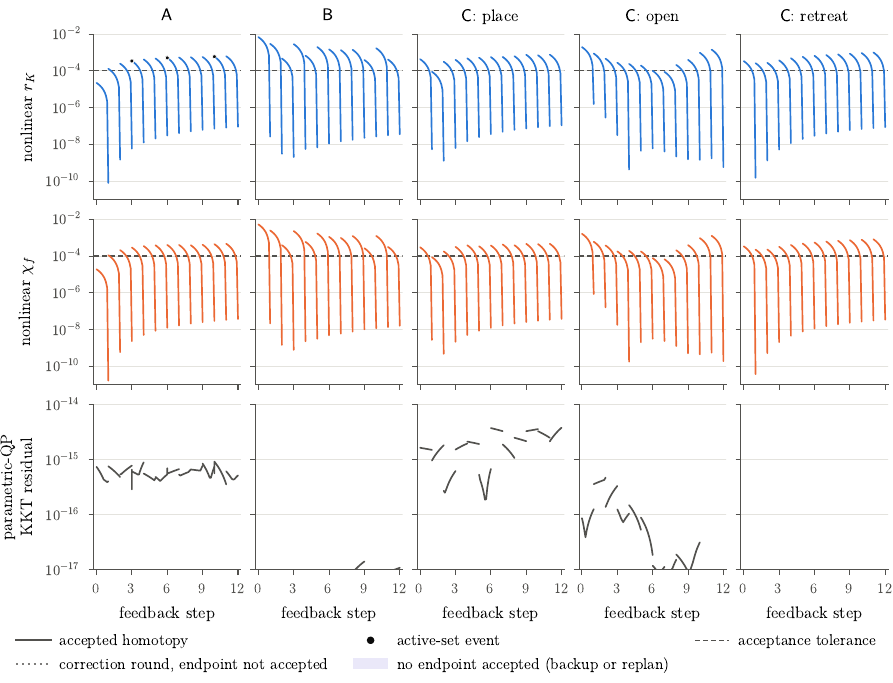}
 \caption{Arm, reduced run: the first twelve feedback steps of each phase,
 where the individual paths are visible. Within each step the nonlinear
 residuals decrease slowly until the last checkpoint before $\tau=1$ and
 then drop to the linearization error.}
 \label{fig:arm-reduced-zoom}
\end{figure}

\begin{figure}[H]
 \centering
 \includegraphics[width=\linewidth]{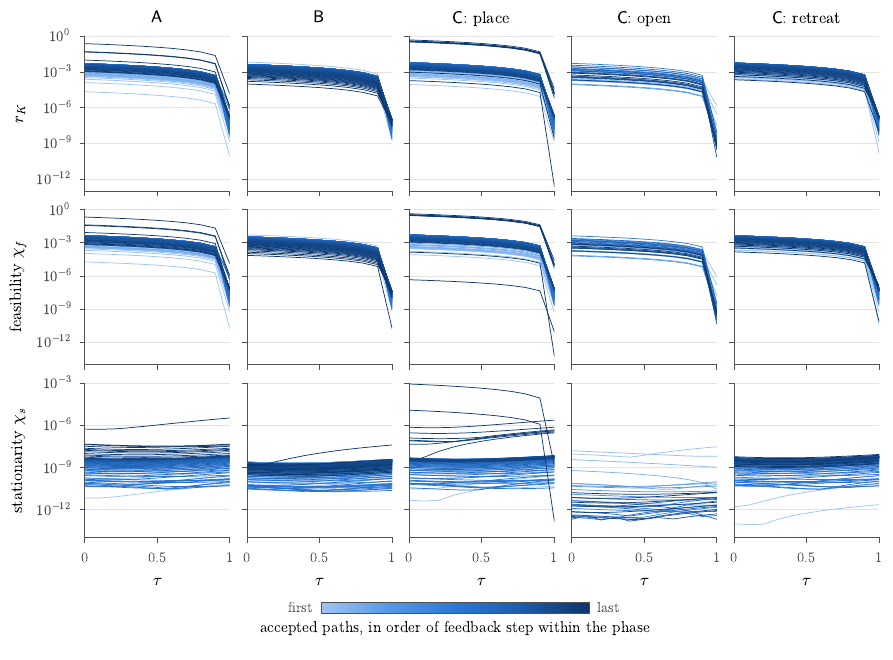}
 \caption{Arm, reduced run: residuals against the homotopy parameter for
 all $473$ accepted paths, by phase. Feasibility and the KKT residual
 follow $(1-\tau)$ times their start values until $\tau=0.9$;
 stationarity is already small at the seed and changes little.}
 \label{fig:arm-reduced-tau}
\end{figure}

\begin{figure}[H]
 \centering
 \includegraphics[width=\linewidth]{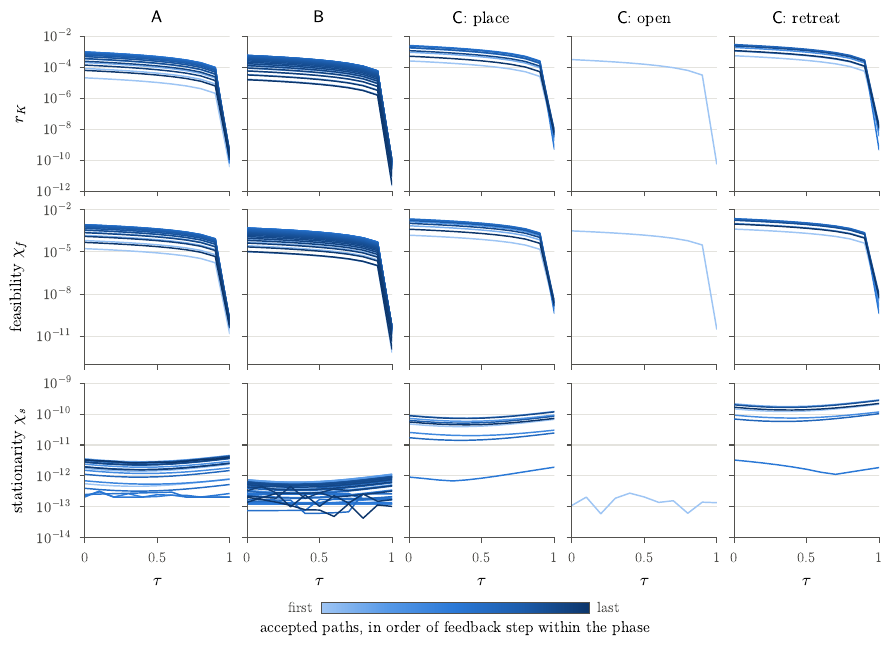}
 \caption{Arm, full-horizon run: residuals against the homotopy parameter
 for the $81$ accepted paths traced at every 25th feedback step.}
 \label{fig:arm-full-tau}
\end{figure}

\FloatBarrier
\bibliographystyle{plain}
\bibliography{refs}
\end{document}